\documentclass[10pt]{article}
\usepackage[affil-it]{authblk}
\usepackage{pstricks}
\usepackage[off]{auto-pst-pdf}
\usepackage{amssymb}
\usepackage{amsmath}
\usepackage{latexsym}
\usepackage{mathrsfs}
\usepackage{xcolor}
\usepackage{bm}
\usepackage{longtable}
\usepackage{multirow}
\usepackage{lscape}
\usepackage{hyperref}

\newcounter{theorem}
\newtheorem{theorem}{Theorem}[section]
\newtheorem{corollary}[theorem]{Corollary}
\newtheorem{condition}[theorem]{Condition}
\newtheorem{lemma}[theorem]{Lemma}
\newtheorem{notation}[theorem]{Notation}
\newtheorem{algorithm}[theorem]{Algorithm}
\newtheorem{definition}[theorem]{Definition}
\newtheorem{remark}[theorem]{Remark}
\newtheorem{example}[theorem]{Example}
\newcommand{\comment}[1]{}

\makeatletter
\ifPreview
\let\Hy@FirstPageHook\relax
\let\Hy@EveryPageAnchor\relax
\fi
\makeatother
\newenvironment{proof}[1][\unskip]{ \textbf{Proof{#1}}. }{$\Box$\medskip}
\newcommand{\doublebrace}[4]{\left\{\begin{array}{ll}#1 & #2 \\ #3 & #4\end{array}\right.}

\author{Todor Milev \\ \texttt{todor.milev@gmail.com}}
\title{
Computing vector partition functions 
}
\date{February 15, 2023\footnote{minor updates on \today}}

\begin{document}
\affil{Google Education}
\maketitle
{\abstract{ 
A vector partition function is the number of ways to write a vector as a non-negative integer-coefficient sum of the elements of a finite set of vectors $\Delta$. We present a new algorithm for computing closed-form formulas for vector partition functions as quasi-polynomials over a finite set of pointed polyhedral cones, implemented in the ``calculator'' computer algebra system.

We include an exposition of previously known theory of vector partition functions. While our results are not new, our exposition is elementary and self-contained. 

}}
\section{Introduction}

Given a finite set $\Delta$ of non-zero integral vectors with non-negative coordinates, and a vector $\bm{\gamma}$ with coordinates $\left( \gamma_1,\dots, \gamma_n \right)$, the vector partition function $P_\Delta(\bm \gamma)$ is by definition the number of ways we can split $\bm \gamma$ as an integral sum with non-negative coefficients of the vectors in $\Delta$. It is well known that there exist finitely many pointed (i.e. with walls passing through the origin) polyhedral cones, with walls parallel to hyperplanes spanned by subsets of $\Delta$ of rank $n-1$, such that $P_\Delta$ is a quasipolynomial over each as a function of the coordinates $(\gamma_1,\dots, \gamma_n)$.

The primary goal of this paper  is to describe an elementary algorithm (Algorithm \ref{algorithm:main}) for computing these quasipolynomials, introduce a light-weight notation for the final result and present pseudo-code with reference C++ implementation in \cite{Milev:calculator}. The central point of our algorithm is the computation of one possible realization of the partial fraction decomposition whose existence is asserted by \cite[Proposition 3.2]{DeConciniProcesi:GeometryToricArrangements}. Our treatment is elementary and is well suited for computer realizations.

We accompany our work with examples. In the appendix, we tabulate computations of the Kostant partition function for the root systems $A_2$, $B_2$, $C_2$, $G_2$, $A_3$, $B_3$, $C_3$, $A_4$. Here, we recall that Kostant partition function is the vector partition function of a root system -  a particular set of vectors that is of interest to Lie algebra representation theory. For example, root system $G_2$ consists of $(1,0), (0,1),(1,1), (1,2), (1,3), (2,3)$. The original motivation for the present work was the study of the characters of generalized Verma modules, which are sums of Kostant partition functions. As far as we are aware, the present article is the first one to publish closed form formulas for the Kostant partition functions of $G_2$ and $B_3$, although our software is not the first one to obtain them - credit for this most likely goes to \cite{Cochet:2005VectorPartitionFunctionAndRepresentationTheory}.

Our secondary goal is to present in Theorem \ref{theorem:mainTheory} an elementary proof of the theoretical results stated in the beginning - namely, that the quasipolynomials give closed formulas for the vector partition function over finitely many pointed polyhedral cones with walls passing through the origin and spanned by $n-1$-element subsets of $\Delta$. Although this seems to have been known and stated directly or indirectly in many texts  (\cite{DeConciniProcesiVergne:PartitionFunctions}, \cite{Sturmfels:OnVectorPartitionFunctions}, \cite{DahmenMicchelli:TheNumberOfSolutionsToDiophantine}), the earliest claimed end-to-end proof that the author was able to locate is in \cite[Theorem 3.26]{DeConciniProcesi:GeometryToricArrangements}. The proof here was done independently of these references and is self-contained.

The proof of Theorem \ref{theorem:mainTheory} proceeds as follows. The partial fraction decomposition given in Algorithm \ref{algorithm:main}  proves that the vector partition function is a quasipolynomial over a finite set of translated polyhedral cones whose walls are parallel to hyperplanes spanned by subset of rank $n-1$ of the starting set $\Delta$. The walls of an individual cone all meet at a single point, however this point depends on the cone and is not necessarily the origin. To resolve this issue, in Section \ref{secWallDist} we present an alternative elementary scheme for computing the vector partition function based on Bernoulli sums. This second algorithm gives formulas over pointed polyhedral cones whose walls pass through the origin. The walls of this second set of cones are not necessarily spanned by $n-1$-element subsets of $\Delta$. Finally, we ``glue'' the polyhedral chambers given by the first algorithm  by those produced by the second approach and that produces our main result Theorem \ref{theorem:mainTheory}.

We have realized all algorithms in this article in the ``calculator'' computer algebra system, \cite{Milev:calculator}. The present article is a rewrite of the non-representation theory aspects of the unpublished 2009 work \cite{Milev:PartialFractions}. Our original software implementation dates back to 2009. However, the present revisit introduces the semi-reduced partial fraction trick from Algorithm \ref{algorithm:main} - the 2009 version combined Algorithm \ref{algorithm:main} and \ref{algorithm:semiReducedToReduced} into a single run time. During our work on the present article, we fixed a number of software issues, added heavy unit testing which increased confidence in the correctness of our software and introduced human-readable printouts of the quasipolynomial formulas and all intermediates.

Our software is not the first one to compute algebraic expressions of vector partition functions and is predated by LattE (\cite{LattE}),  the ``barvinok'' program (up to our knowledge, no longer available on the internet; see \cite{Kevin&others:BarvinokAlgorithm}), as well as a series of MAPLE programs with focus on Kostant partition functions (\cite{Cochet:2005VectorPartitionFunctionAndRepresentationTheory},  \cite{BaldoniBeckCochetVergne:VolumeComputationForPolytopesAndPartitionFunctions}). Special credit goes to \cite{Cochet:2005VectorPartitionFunctionAndRepresentationTheory} and related works;  \cite{Cochet:2005VectorPartitionFunctionAndRepresentationTheory} works with the full quasipolynomial of the vector partition function, just as our implementation does.

\textbf{Acknowledgements.} The author would like to thank M. Vergne, W. Baldoni, T. Bliem, A. Petukhov and M. Beck for the valuable advice and discussions on vector partition functions.

\textbf{Competing interests.} The author declares no competing interests.

\color{black}
\subsection{Notation}
In the spirit of related articles, when denoting vectors, we shall use small bold-faced greek letters such as $\bm{\gamma}$, with the only novelty being our use of the bold face for vectors. We will denote the coordinates of $\bm \gamma$ by $\gamma_{1},\dots, \gamma_n$. Similarly, when $\bm \gamma_1, \dots, \bm \gamma_s$ denote $s$ different vectors, we denote the coordinates of $\bm \gamma_i$ by $\gamma_{i,1}, \dots, \gamma_{i,n}$. We reserve superscripts for regular arithmetics, that is, $\gamma_1^{2}$ is the square of the first coordinate of vector $\bm \gamma$.

\section{Linear optimization prerequisites}
In the present section, we fix the terminology with which we will describe $n$-dimensional pointed polyhedral cones. For completeness, we will also describe a few basic algorithms needed to manipulate pointed polyhedral cones in a computer program.

Throughout this paper, we will consider hyperplanes and their arrangements. Unless stated otherwise, we shall assume that all planes pass through the origin and are therefore determined by their normal vector. 

We recall that by projectivization, we can map any point in $n-1$-dimensional space to a line through the origin in $n$ dimensions. To do this, embed the $n-1$-dimensional space as the plane $x_n=1$ and assign to each point with $x_n=1$ the line through it and the origin. Therefore our geometric intuition for finite-volume polyhedra in dimension $n-1$ can be applied to understand pointed polyhedral cones in dimension $n$. 

A homogeneous system of inequalities with rational coefficients
\begin{equation}\label{eq:coneDefinition}
\begin{array}{ccccccl}
a_{1,1} x_1&+&\dots&+&a_{1,n}x_n&\geq& 0 \\
&&&&&\vdots \\
a_{m, 1}x_1&+&\dots&+&a_{m,n}x_n&\geq& 0
\end{array}
\end{equation}
cuts off an infinite cone passing through the origin of the coordinate system. More precisely, the $k^{th}$ inequality in \eqref{eq:coneDefinition} cuts a half-space bordered by the hyperplane with normal given by the $\left(a_{k,1}, \dots, a_{k, n}\right)$.

We will abbreviate \eqref{eq:coneDefinition} using the scalar ``dot product'' as 
\begin{equation}\label{eq:coneDefinitionScalarProduct}
\begin{array}{rcl}
\langle \bm a_1, \bm x\rangle &\geq& 0\\
&\vdots\\
\langle \bm a_m, \bm x\rangle &\geq& 0\\
\end{array},
\end{equation}
where $\bm a_k = \left(a_{k,1}, \dots, a_{k,n}\right)$. 

We note that by clearing denominators in the inequalities, we can guarantee that the system has integer coefficients. We furthermore note that some of the inequalities in \eqref{eq:coneDefinition} could be redundant. Here, an inequality given by \eqref{eq:coneDefinition} is redundant if omitting it does not change the set of solutions of the system \eqref{eq:coneDefinition}. For example if the system contains the inequalities $x_1\geq 0$, $x_2\geq 0$ and $x_1+x_2\geq 0$, the last of these inequalities is a consequence of the other two and can be dropped. In this example, the hyperplane with normal $(1,1)$ would intersect the cone at a single point - the origin. In more than two variables, this generalizes to saying that an inequality is redundant if it intersects the cone in a set of points of zero $n-1$-dimensional surface area. A set of points given by \eqref{eq:coneDefinition} may happen to be of zero $n$-dimensional volume. For an example in dimension $1$, the inequality $ x_1\geq 0$, $-x_1 \geq 0$ forces that $x_1=0$; this toy example can be extended to higher-dimensional examples by adding more variables.

\begin{definition}\label{definition:chamber}
We say that a subset $C$ of ${\mathbb{Q}}^n$ given by homogeneous non-strict ($\geq$) rational linear inequalities \eqref{eq:coneDefinition} is a chamber. We call the hyperplanes given by the defining inequalities to be the \emph{supporting hyperplanes} of $C$ on the condition that they intersect $C$ in a set of non-zero $n-1$-volume. The intersection of a supporting hyperplane with $C$ will be called a \emph{wall of $C$}.
\end{definition}
In Definition \ref{definition:pointedPolyhedralCone}, we will define a pointed polyhedral cone, which is equivalent to the notion of a chamber in Definition \ref{definition:chamber} by Lemma \ref{lemma:conesAreChambers}. We will therefore use the terms ``chamber'' and ``pointed polyhedral cone'' interchangeably.

We explicitly allow the chamber $C$ to be the entire space $\mathbb Q^n$ by allowing the empty set of defining rational inequalities. 

When representing $C$ on a computer, we request the following.

\begin{definition}[Normalized walls]\label{definition:coneWalls}
We say that the normal vectors $\bm a_1,\dots, \bm a_m$ defining a cone $C$ through \eqref{eq:coneDefinition} are normalized walls of $C$ if the following hold.
\begin{itemize}
\item There are no redundant inequalities in \eqref{eq:coneDefinition}.
\item All $\bm a_i$'s are rescaled so that the coefficients have integer coefficients and the coefficients have greatest common divisor $1$.
\item The vectors $\bm a_i$'s are sorted with respect to some total order on their coordinates $\left(a_{i,1}, \dots, a_{i,n}\right)$.
\end{itemize}
\end{definition}
For $C$ of non-zero volume, the normalized walls determine one  possible choice of normal vectors for each of the supporting hyperplanes of $C$.

We claim that if $C$ is of non-zero volume (i.e.,  contains $n$ linearly independent points), then the ordered list of normalized walls  $\bm a_1,\dots, \bm a_m$ exists and is unique. In Lemma \ref{lemma:conesHaveNormalWalls}, we prove the existence. We also prove uniqueness but only with the additional assumption that all coordinates of $C$ are non-negative. 

The computational significance of this definition is clear. Suppose we have two cones given by two sets of inequalities \eqref{eq:coneDefinition}. Suppose the two cones are of non-zero volume. Compute their normalized walls. Then the two cones are identical (contain the same points) if and only if their normalized wall lists are identical. 

The third condition requires that the normalized vectors be sorted with respect to some total order. Our computer realization sorts the walls using the graded colexicographic order on the vectors, in other words, we compare two vectors by the sum of their coordinates with smaller sums yielding colexigoraphically-smaller vectors. Ties are then broken by declaring the larger vector to be the that has a larger entry in the last coordinate where the two vectors are different. 

\begin{definition}[Vertex of chamber]\label{definition:vertex}
Let $C$ be a chamber. We say that a point in the chamber is a normalized vertex (or just vertex for short) if
\begin{itemize} 
\item it belongs to $n-1$ different normalized walls of the chamber and
\item it is rescaled by multiplication of a positive number so that all coordinates are  integers without a common divisor.
\end{itemize}
\end{definition}

To find all the vertices of chamber $C$, we can use the following algorithm.

\begin{algorithm}[Vertices from normals]\label{algorithm:findTheVertices}
\item[Step 1.] For every $n-1$-element subset $A$ of $\bm a_1, \dots, \bm a_m$: 
\item[Step 1.1.] Using row-echelon form matrix reduction, compute a basis $\bm v_1, \dots$ for the linear space of all vectors perpendicular to $A$.
\item[Step 1.2.] If the number of vectors $\bm v_1, \dots$ is larger than $1$, go back to Step 1.
\item[Step 1.3.] Else, we have a single vector $\bm v_1$. Rescale the vector by \eqref{definition:vertex}.
\item[Step 1.4.] If $\langle \bm v_1, \bm a_k\rangle <0$ for some $k$, set $\bm v=-\bm v_1$, else set $\bm v=\bm v_1$.
\item[Step 1.5.] If $\langle \bm v, \bm a_j\rangle <0$ for some $j$, go back to Step 1.
\item[Step 1.6.] Add $\bm v$ to the list of vertices, if it's not already there.
\item[Step 1.7.] Go back to Step 1 to proceed with the next subset, until all subsets are exhausted.
\end{algorithm}
Vectors found in Step 1.5 may indeed not be new if the chamber is given by redundant inequalities.

\begin{definition}\label{definition:pointedPolyhedralCone}
Let $\bm \alpha_1, \dots, \bm \alpha_k$ be vectors. We say that 
\begin{equation}\label{eq:coneSpan}
C=C\left( \bm\alpha_1, \dots, \bm \alpha_k\right) = \left\{ t_1\bm \alpha_1+\dots + t_k \bm \alpha_k |t_i\geq 0\right\}
\end{equation}
is the polyhedral cone generated by the vectors $\bm \alpha_1, \dots, \bm \alpha_k$.

Conversely, we will say that a set $C$ is a pointed polyhedral cone if it can be given by a formula \eqref{eq:coneSpan}.
\end{definition}

The following lemma is well-known from the subject of linear algebra and we omit a proof.
\begin{lemma}\label{lemma:conesAreChambers}
A set $C$ is a chamber (Definition \ref{definition:chamber}) if and only if it is a pointed polyhedral cone (Definition \ref{definition:pointedPolyhedralCone}).
\end{lemma}
Given a set of $\bm \alpha_i$'s of full rank, here's an algorithm to find the normalized walls of $C(\bm \alpha_1, \dots, \bm \alpha_s)$. 
\begin{algorithm}[Normals from vertices]~ \label{algorithm:normalsFromGenerators}
\begin{itemize}
\item[Step 1.] For every $n-1$-element subset $J$ of the $\bm \alpha_i$'s:
\begin{itemize}
\item[Step 1.1.] Using row-echelon form matrix reduction, compute a vector space basis for the vectors $\bm m_1, \dots,\bm m_k$ that are perpendicular to $J$. 

\item[Step 1.2.] If $k\geq 2$, the set $J$ is linearly dependent. Discard $J$ and go back to Step 1. Else, set $\bm m=\bm m_1$.
\item[Step 1.3] If two of the $\bm \alpha_i$'s have scalar products with $\bm m$ of different sign, go back to Step 1.
\item[Step 1.4] If there's an $\bm \alpha_i$ for which $\langle \bm\alpha_i, \bm m\rangle < 0$, multiply $\bm m$ by $-1$.
\item[Step 1.5] Rescale using a positive rational number the coordinates of $\bm m$ so they all have integer coefficients with greatest common divisor $1$. 
\item[Step 1.6] If $\bm m$ was already found for a previous value of $J$, go back to Step 1.
\item[Step 1.7.] If we have that $\langle \bm\alpha_i, \bm m\rangle = 0$ for all $\bm \alpha_i$, then the original set of $\bm\alpha_i$'s does not span the entire vector space. Abort the algorithm. Else, add $\bm m$ to the set of normals of walls of $C$. 
\item[Step 1.8] Go back to Step 1 to process the next set $J$.
\end{itemize}
\item[Step 2.] Sort the normals $\bm m$ using a preferred sorting function.
\end{itemize}
\end{algorithm}
The walls produced by Algorithm \ref{algorithm:normalsFromGenerators} are normalized in the sense of Definition \ref{definition:coneWalls}. We omit the proof.

We will use the following lemma, whose proof we omit.
\begin{lemma}\label{lemma:conesAreGeneratedByTheirVertices}
	
Let $C$ be a pointed polyhedral cone (chamber) that does not contain any point with a negative coordinate. Then $C$ is generated by its vertices, i.e., 
\[
C\left(\bm \alpha_1, \dots, \bm \alpha_t\right) = C
\]
where $\bm \alpha_1, \dots, \bm \alpha_t$ are the vertices of $C$.
\end{lemma}
We will also use the following lemma.
\begin{lemma}\label{lemma:conesHaveNormalWalls}
Let $C$ be a pointed polyhedral cone that contains $n$ linearly independent points. Suppose all points in $C$ have have non-negative coordinates. Then $C$ has a unique representation using normalized walls (Definition \ref{definition:coneWalls}).
\end{lemma}
\begin{proof}
We prove existence first. Compute the set of vertices $\bm \alpha_1, \dots, \bm \alpha_t$ of $C$ using Algorithm \ref{algorithm:findTheVertices} from the normals $\bm a_1, \dots, \bm a_m$ that define $C$. By Lemma \ref{lemma:conesAreGeneratedByTheirVertices} $C=C(\bm\alpha_1, \dots, \bm \alpha_t)$. By Algorithm \ref{algorithm:normalsFromGenerators}, find a new set of normals $\bm m_1, \dots, \bm m_l$. By the remarks after the algorithm, this is the desired normalized set of walls. 

To prove uniqueness, we note that each $\bm m_i$ contains $n-1$ linearly independent vertices of $C$. This means that every representation of $C$ - including possibly non-normalized ones such as $\bm a_1, \dots, \bm a_m$ - must contain a positive multiple of each of the $\bm m_i$'s. This means that every normalized set of walls contains all of $\bm m_1, \dots, \bm m_l$ up to a positive constant. By the second and third condition of Definition \ref{definition:coneWalls} this means the normalized walls are unique.
\end{proof}

\section{Partial fraction decomposition algorithm}\label{secPFdecomposition}
\subsection{Partial fraction decomposition definition}
Given an integer-coefficient vector $\bm \gamma=\left(\gamma_1, \dots, \gamma_n \right)$, we use the abbreviated notation
\[
x^{\bm \gamma} = x_1^{\gamma_1}, \dots, x_n^{\gamma_n},
\]
where $x_1, \dots, x_n$ are variables. Throughout this paper, we will be computing over the ring
\begin{equation}\label{eq:ringOfImportance}
\mathbb V= \displaystyle {\mathbb{Q}}\left[x_1, \dots,x_k, \frac{1}{x_1}, \dots,\frac{1}{x_n} \right] \left[ \frac{1}{\left(1 - x^{\bm \gamma} \right)^p }\right]_{\bm \gamma \in \mathbb{Z}^k, p\in \mathbb Z}, 
\end{equation}
i.e., the ring generated by multiplications and additions of all $x_i$'s, all $x_i^{-1}$'s and all fractions of the form 
\[
\frac{1}{\left(1-x^{\bm \gamma}\right)^p} = \frac{1}{\left(1-x_1^{\gamma_1}\dots x_n^{\gamma_n }\right)^p}
\]
In more strict terms, $\mathbb V$ is the localization of the polynomials ${\mathbb{Q}}\left[x_1, \dots,x_n \right] $ relative to the multiplicative set generated by $x_1, \dots, x_n, (1-x^{\bm \gamma})$, $\bm \gamma \in \mathbb Z^n$.
We remark that $\mathbb V$ is a Weyl algebra module, a fact we will later use to expand the elements of $\mathbb V$ into power series.

The reason we are interested in $\mathbb V$ is that, given a set $\Delta$ whose vector partition function $P_I$  we want to compute, then 

\begin{equation}\label{eq:VectorPartitionFunctionGeneratorFunction}
\begin{array}{rcll|l}
\displaystyle
\prod_{\bm \alpha \in \Delta}\frac{1}{1-x^{\bm \alpha}} &=&\displaystyle \prod_{\bm \alpha \in \Delta} \left(1+x^{\bm \alpha}+ x^{2\bm \alpha} + x^{3\bm \alpha}+\dots \right) &&\text{Geometric series}\\
&=&\displaystyle \sum_{\gamma \in \mathbb Z^n }P_\Delta(\bm \gamma) x^{\bm \gamma}
\end{array}.
\end{equation}
The last equality follows from the fact that, when we distribute the finite product of the infinite geometric series, for each way to break down $\bm \gamma = a_1\bm \alpha_1+\dots+a_s\bm \alpha_s$ we get one monomial $x^{\bm \gamma}$ by distributing $x^{a_i\bm \alpha_i}$ from the $i^{th}$ geometric series. The ability to expand elements of $\mathbb V$ in power series - in particular the ability to compute the expansion in \eqref{eq:VectorPartitionFunctionGeneratorFunction} - implies means to compute $P_\Delta(\bm \gamma)$. %We will compute the expansion above by first computing a partial fraction decomposition for \eqref{eq:VectorPartitionFunctionGeneratorFunction}.

 Consider an element of $\mathbb V$ given in the form:
\begin{equation}\label{eq:partialFractionDecompositionNonReduced}
	q=\sum_{\substack{I \in \mathbb I \\ \{\bm \alpha_1, \dots, \bm \alpha_r\} \subset \Delta}} \frac{p_{I}}{ \displaystyle \prod_{i=1}^{r} \prod_{j} \left(1- x^{b_{i, j} \bm \alpha_{i }} \right)^{m_{i , j}} }
\end{equation}
In the sum, the vectors $\bm \alpha_1, \dots, \bm \alpha_r$ lie in a subset of $\Delta$ that is a function of some indexing element $I$. The integers $r$, $b_{i,j}$ and $m_{i,j}$ as well as the Laurent polynomial $p_I$ are also functions of $ I$, but we omit the $I$ from their notation in the interest of readability.

\begin{definition}\label{definition:partialFractionDecomposition}
	~
	\begin{itemize}
		\item We call a summand in \eqref{eq:partialFractionDecompositionNonReduced} a \emph{partial fraction}. 
		\begin{itemize}
			\item If the $\bm \alpha_1, \dots, \bm \alpha_r$ of the summand are linearly dependent, we say that the partial fraction is \emph{non-reduced}.
			\item If the $\bm \alpha_1, \dots, \bm \alpha_r$ of the summand are linearly independent and span the whole space, we say that the partial fraction is \emph{semi-reduced}.
			\item If a partial fraction is semi-reduced and the $\prod_j $-product in the denominator of \eqref{eq:partialFractionDecompositionNonReduced} contains a single term, we say that the partial fraction is \emph{fully reduced}.
		\end{itemize}
		\item A formula \eqref{eq:partialFractionDecompositionNonReduced} will be called 
		\begin{itemize}
			\item a \emph{non-reduced partial fraction decomposition} of $q$ if it has at least one non-reduced partial fraction;
			\item a \emph{semi-reduced partial fraction decomposition} of $q$ if it has no non-reduced partial fractions but has at least one partial fraction that is not fully reduced;
			\item a \emph{fully reduced partial fraction decomposition} of $q$ if it only has fully reduced partial fractions.
		\end{itemize}
	\end{itemize}
\end{definition}

Existence of a fully reduced partial fraction decomposition for every $q\in \mathbb V$ is proven in \cite[Proposition 3.2]{DeConciniProcesi:GeometryToricArrangements}.  The present section will be dedicated to describing in great detail how to do such a decomposition on a computer. A note on terminology: \cite{DeConciniProcesi:GeometryToricArrangements} defines a unique ``partial fraction decomposition'' for an element $ q\in \mathbb V$, which is much stronger and more technical than our term ``fully reduced partial fraction decomposition''. In our terminology, fully reduced partial fraction decompositions are not unique.

\begin{example}\label{example:fullyReducedDecompositionIsNotUnique} The equality 
\begin{align*}
p=\frac{1}{(1-x_{1} ) (1-x_{2} ) (1-x_{1} x_{2} ) }=&	\displaystyle ~~~ \frac{-x_{2}^{-1}}{(1-x_{1} )^2 (1-x_{1} x_{2} ) }\\&
	+\displaystyle \frac{x_{2}^{-1}}{(1-x_{1} )^2 (1-x_{2} ) }\end{align*}
gives a fully reduced partial fraction decomposition of $p$. As $p$ is symmetric with respect to the two variables, a second fully reduced partial fraction decomposition of $p$ is obtained by swapping $x_1$ and $x_2$ in the right hand side of the equality.
\end{example}

For additional precision, in the following notation we describe the approximate data structure we used to represent a summand of \eqref{eq:partialFractionDecompositionNonReduced} in the computer's memory.

\begin{notation} \label{notation:partialFractionData} Every partial fraction will be represented in the computer's memory by storing the following data.
	\begin{itemize}
		\item[(a)] The numerator polynomial $p_I$.
		\item[(b)] For every vector $\bm\alpha \in \Delta$:
		\begin{itemize}
			\item If $\bm \alpha$ appears in the denominator of the fraction (see \eqref{eq:partialFractionDecompositionNonReduced}) in the form 
			\[
			\left(1- x^{b_1\bm \alpha}\right)^{m_1} \dots \left(1-x^{b_s\bm \alpha} \right)^{m_s} 
			\]
			store a list of pairs of integers $(b,m)$:
			\begin{equation}\label{eq:denominatorKeyValuePair}
				\bm	\alpha\mapsto  \left[\left(b_1, m_1\right), \dots,\left(b_s, m_s\right) \right]
			\end{equation}
			\item If $\bm \alpha$ does not appear in the denominator of the partial fraction, store the empty list.
		\end{itemize} 
	\end{itemize}
\end{notation}

We note that the symbol $\mapsto$ used in \eqref{eq:denominatorKeyValuePair} suggests that \eqref{eq:denominatorKeyValuePair} can naturally be represented by a map in the computer science sense (in our C++ realization we used a hash-map). In turn, the full sum \eqref{eq:partialFractionDecompositionNonReduced} can also be represented as a map (again, we used a hash-map) from the set of the data from Notation \ref{notation:partialFractionData}(b) to the set of polynomials $p_I$ from Notation \ref{notation:partialFractionData}(a); the data (b) serves as a key and the data (a) serves as a value of the map.

\subsection{A universal formula}

In our upcoming main Algorithm \ref{algorithm:main} we will construct a fully reduced partial fraction decomposition for an element $q\in\mathbb V$. To do so, we will apply carefully a formula due to Szenes-Vergne (\cite{SzenesVergne:ResidueFlas}) and the geometric series sum formula. To make our exposition even simpler, we have united the two formulas in a single universal formula in Lemma \ref{lemma:formulas}(c).

\begin{lemma}~\label{lemma:formulas} 
\begin{itemize}
\item[(a)](Szenes-Vergne formula) Let $\bm \alpha_1, \dots, \bm \alpha_k$ be vectors with integer coefficients such that $\displaystyle \sum_{i} \bm \alpha_i\neq 0$ 

\begin{equation}\begin{array}{rcl}
\displaystyle \frac{1}{1-x^{\bm \alpha_1}}\dots \frac{1}{ 1-x^{\bm \alpha_k }} &=&\displaystyle \frac{1}{1-x^{\sum \bm \alpha_i }} \cdot \\
&&\displaystyle \sum_{j=1 }^{j=k} \frac{ x^{\bm \alpha_1}\dots x^{\bm  \alpha_{j-1}}}{ (1 -x^{\bm \alpha_1 }) \dots (1-x^{\bm \alpha_{j - 1}})}\\&&\displaystyle \phantom{\sum} \frac{1}{ (1-x^{ \bm \alpha_{j + 1}})\dots(1-x^{\bm \alpha_{k}})} 
\end{array}
\end{equation} 
\item[(b)] (Geometric series sum formula)
\begin{equation}\label{eq:geometricSeries1}
\frac{1}{1-x^{\bm \alpha}} =\frac{1+x^{ \bm \alpha} +\dots+x^{(k-1) \bm \alpha}}{1-x^{k \bm \alpha}} =\frac{g_k( x^{\bm\alpha}) }{1 - x^{k\bm \alpha }},
\end{equation} 
for $k\in{\mathbb{Z}}_{>0}$, and  
\begin{equation}\label{eq:geometricSeries2}
\frac{1}{1-x^{\bm \alpha}}= \frac{-x^{-\bm\alpha} -\dots-x^{k\bm \alpha}}{1-x^{k \bm\alpha}}= \frac{g_k\left( x^{\bm\alpha} \right) }{1 - x^{k\bm \alpha }},
\end{equation} 
for $k\in{\mathbb{Z}}_{<0}$.
\item[(c)](generalized Szenes-Vergne formula) Let $k \geq 1$ and $\bm \alpha_1,\dots ,\bm \alpha_k$ be as above. Then
\begin{equation}\label{eq:SzenesVergneFormula}
\begin{array} {r@{~}c@{~}l}
\displaystyle 
\frac{1}{1-x^{\bm \alpha_1}}\dots \frac{1}{1-x^{\bm \alpha_k}}& =& \displaystyle \frac{ 1 }{ 1-x^{\sum a_i \bm \alpha_i}}\cdot \\
&&\displaystyle \left( \sum_{j=1}^{j=k} \frac{x^{ a_1 \bm \alpha_1} \dots x^{a_{j-1} \bm \alpha_{j-1}}}{(1-x^{\bm \alpha_1}) \dots (1-x^{\bm \alpha_{j-1}})}\cdot \right.\\
&&\displaystyle ~~~~~~\left. \frac{g_{a_j}(x^{\bm \alpha_j})} {(1-x^{\bm \alpha_{j+1}}) \dots(1- x^{\bm \alpha_{k}})} \right),
\end{array}
\end{equation} 
where $g_k\left(x^{\bm\alpha} \right)$ are the numerators in geometric series sum formula (b). 
\end{itemize}
\end{lemma}
\textbf{Remark.} When $k=1$ (c) collapses down to (b). When all of the $a_i$'s are 1, (c) collapses down to (a).

\begin{proof}[ of Lemma \ref{lemma:formulas}]
(a) is obtained from \cite[Lemma 1.8]{SzenesVergne:ResidueFlas} by setting all complex coefficients in the lemma to be 1 and rearranging some of the fractions. (c) is obtained by combining (a) and (b). 
\end{proof}

\begin{example}
A simple example of \eqref{eq:SzenesVergneFormula} is the formula
\[
\begin{array}{@{}r@{}c@{}l}
\displaystyle \frac{1}{(1-x_{2})(1-x_{1}x_{2})}&=&\displaystyle \frac{-x_{2}^{-1}}{(1-x_{1})(1-x_{1}x_{2})}+\frac{x_{2}^{ -1}}{(1 - x_{1})(1-x_{2})} 
\end{array}
\]
generated by the linear combination $(1,0)=-(0,1)+(1,1)$ with $a_1=-1$ and $a_2=1$.
\end{example} We used this formula to generate our previous Example \ref{example:fullyReducedDecompositionIsNotUnique}.

\subsection{An algorithm for computing partial fraction decompositions} \label{secTheAlg} \label{secTheTheoreticalAlgorithm} The goal of this section is to present an algorithm for computing a fully reduced partial fraction decomposition with respect to a set of non-zero vector with non-negative integer coordinates $\Delta$ (see Definition \ref{definition:partialFractionDecomposition}). The algorithm is similar in spirit to the proof of \cite[Proposition 3.2]{DeConciniProcesi:GeometryToricArrangements}, and will be carried out using consecutive applications of formula \eqref{eq:SzenesVergneFormula}.

From formula \eqref{eq:partialFractionDecompositionNonReduced} and Definition \ref{definition:partialFractionDecomposition}, we see that a semi-reduced partial fraction decomposition is a formula of the form:
\begin{equation}\label{eq:partialFractionDecompositionSemiReduced}
q=\sum_{\substack{I \in \mathbb I \\ \{\bm \alpha_1, \dots, \bm \alpha_n\} \subset \Delta} } \frac{p_{I}}{ \displaystyle \prod_{i=1}^{n} \prod_{j} \left(1- x^{b_{i, j} \bm \alpha_{i }} \right)^{m_{i , j}} }
\end{equation}
where the $\bm \alpha_i$'s are linearly independent and a fully reduced partial fraction decomposition is a formula of the form:
\begin{equation}\label{eq:partialFractionDecomposition}
	\begin{array}{rcl}
		q&=&\displaystyle \sum_{\substack{I \in \mathbb I \\ \{\bm \alpha_1, \dots, \bm \alpha_n\} \subset \Delta}}  \frac{p_{J}} {\left(1-x^{b_1\bm \alpha_{1}} \right)^{ m_1}\dots \left(1-x^{b_n\bm \alpha_{ n}} \right)^{m_n}}\\
		&=&\displaystyle \sum_{\substack{I \in \mathbb I \\ \{\bm \alpha_1, \dots, \bm \alpha_n\} \subset \Delta}}  \frac{p_{J}} {\displaystyle\prod_{i} \left(1-x^{b_i\bm \alpha_{i }} \right)^{m_i} },
	\end{array}
\end{equation}

Given $k,m >0$, any two fractions $\displaystyle \frac{1}{\left(1-x^{k \bm \alpha}\right)}$, $\displaystyle \frac{1}{\left(1-x^{m \bm \alpha}\right)}$ can be transformed to have a common denominator via the geometric series formula \eqref{eq:geometricSeries1}:
\[
\begin{array}{rcl}
\displaystyle 
\frac{1}{\left(1-x^{k \bm \alpha}\right)} &=& \displaystyle \frac{g_{\frac{r }{k}}\left( x^{k\bm \alpha } \right)} {\left(1-x^{r \bm \alpha}\right)} \\ \displaystyle  \frac{1}{\left(1-x^{m \bm \alpha}\right)}&=&\displaystyle \frac{g_{\frac{r }{m}}\left( x^{m\bm \alpha } \right)} {\left(1-x^{r \bm \alpha}\right)}, 
\end{array} 
\]
where $r=\text{lcm}(k,m) $. This shows how to convert a semi-reduced partial fraction into a fully reduced one. 

We note that this also shows how to reduce the $\prod_{j}$-product in \eqref{eq:partialFractionDecompositionNonReduced} to a single term even when the $\bm \alpha_i$'s are linearly dependent. However, these transformations increase the number of monomials in the numerator of each fraction significantly, which is bad for the memory consumption of our upcoming algorithm. Therefore it makes sense to postpone the application of these formulas to as late in our algorithm as possible. This explains why we did not get rid of the $\prod_j$ term in formula \eqref{eq:partialFractionDecompositionNonReduced}.

We now proceed with an algorithm that transforms any non-reduced partial fraction decomposition to a semi-reduced one. We are going to assume that $\Delta$ is a set of vectors with positive integer coefficients, equipped with a total order $>$. In a computer realization, the total order on the vector can of course be replaced with the order in which the vectors appear in a list.

\begin{algorithm}[Semi-reduced partial fractions] \label{algorithm:main}~
\begin{itemize}
\item Input: a non-reduced partial fraction relative to $\Delta$.
\item Output: a reduced partial fraction relative to $\Delta$.
\end{itemize}
\begin{itemize}
\item[Step 1.] Set $R=0$ and set $N$ to be the starting partial fraction, represented in the computer's memory as described in Notation \ref{notation:partialFractionData}.
\item[Step 2.] If $N$ is empty, halt; else remove an arbitrary fraction $u$ from \eqref{eq:partialFractionDecompositionNonReduced} - on a computer, that can be realized by popping out a [key, value] pair from a map. Let $u$ be given by the formula below.
\begin{equation}\label{eq:generalFormPartialFraction}
u= \frac{p}{\displaystyle \prod_{i=1}^r \prod_{j=1}^{l_i} \left(1-x^{b_{i,j} \bm \alpha_i } \right)^{ m_i,j}}
\end{equation}
and let $w$ be the subproduct in $u=v\cdot w$ that selects the highest powers $b_{i,l_i}$ for each of the inner-most multiplicands in the denominator:

\[
\begin{array}{rcl}
u&=& v\cdot w \\
w&=&\displaystyle \frac{1}{\displaystyle \prod_{i=1}^r \left(1-x^{b_{i,l_i} \bm\alpha_i } \right)^{ m_{i, l_i}}}\\
&=&\displaystyle \frac{1}{\displaystyle \left(1-x^{b_{1} \bm\alpha_1} \right)^{m_1} \dots \left( 1 -x^{b_{r} \bm\alpha_r }\right)^{ m_r }}\\
&=&\displaystyle \frac{1}{\displaystyle \left(1-x^{ \bm\beta_1} \right)^{ m_1}\dots \left( 1- x^{ \bm\beta_r}\right)^{m_r}},
\end{array}
\]
where we have set $b_i=b_{i,l_i}$, $m_i =m_{i,l_i}$ and finally 
\[\bm\beta_i = b_i\bm \alpha_i \] to simplify the notation.

In a computer realization, the construction of $v$ amounts to removing the entry with largest value of $b$ in the pair $(b,m)$ from the data of Notation \ref{notation:partialFractionData}(b), and the construction of $w$ amounts to assembling the so-removed pairs $(b,m)$ in a new partial fraction.
\item[Step 3.] If the fraction $u$ is reduced add it to $R$ and go back to Step 2, else proceed with Step 4.
\item[Step 4.] Let $\bm \alpha_1< \dots< \bm \alpha_r$ be the vectors appearing as exponents of $x$ in the denominator of $u$. Select the smallest index $s$ for which $\bm \alpha_1, \dots,\bm \alpha_s$ are linearly dependent and take a linear combination
\begin{equation}\label{eq:algorithmMain}
\begin{array}{rcl}
\displaystyle a_{t-1}\cdot b_{t-1}\bm \alpha_{t-1}&=&\displaystyle  a_{t} \cdot b_t \bm \alpha_t +\dots+a_s \cdot b_s\bm \alpha_s\\
&=&\displaystyle  a_{t} \cdot \bm \beta_t+\dots+a_s \cdot \bm \beta_s
\end{array}
\end{equation}
where $t-1$ is the smallest index so that $a_{t-1}\neq 0$ and we also have that $a_s\neq 0$ by the choice of $s$. Rescale the linear combination so that all $a_i$'s are integers, the integers have no greatest common divisor, and $a_{t-1}>0$. Finally, drop all summands of the linear combination whose coefficient is zero. To simplify our notation, let us relabel the resulting linear combination as:
\[
a_{t-1}\cdot b_{t-1}\bm \alpha_{t-1}= a_{1} \cdot \bm \delta_1+\dots+a_l \cdot \bm \delta_l,
\]
where each summand $ a_i \bm \delta_i$ equals one of the non-zero summands $a_j\cdot \bm \beta_j$ from the right hand side of \eqref{eq:algorithmMain}.

\item[Step 5.] Transform $u$ by applying the generalized Szenes-Vergne formula \eqref{eq:SzenesVergneFormula}:

\begin{equation}\label{eq:applySzenesVergneFormula}
\begin{array}{rcl}
\displaystyle u&=&v\cdot w\\
&=&
\displaystyle v\cdot \left(\frac{1}{\displaystyle \prod_{\substack{ i=1 \\ \varepsilon_i }}^{r} \left(1-x^{b_i \bm \alpha_{ i}}\right)^{m_i-\varepsilon_i} } \right) \cdot \frac{1}{\left(1-x^{a_{t-1}\cdot b_{t-1}\bm \alpha_{t-1}}\right)} \\
&&\displaystyle \cdot \sum_{j=1}^{j=l} \frac{x^{ a_t \bm \delta_t} \dots x^{\bm \delta_{j- 1}}}{(1-x^{\bm \delta_1}) \dots (1- x^{\bm \delta_{j-1}})} \frac{g_{a_j}\left(x^{\bm \delta_j}\right)} {(1-x^{\bm \delta_{j+1}}) \dots(1- x^{\bm \delta_{l}})},
\end{array} 
\end{equation}
where $\varepsilon_i =\begin{cases} 1&\text{if }b_i\bm\alpha_i\text{ participates in the right-hand side of \eqref{eq:algorithmMain}}\\0&\text{otherwise} \end{cases}$.

Put the $l$ new summands generated by \eqref{eq:applySzenesVergneFormula} back into $N$.
\item[Step 6.] Go back to Step 2.
\end{itemize}
\end{algorithm}

We note that while we did not find it very easy to code Algorithm \ref{algorithm:main} on a computer, we made use of a very simple trick for finding errors in our implementation. At the start of each computation, we choose arbitrary values for $x_1, \dots, x_n$ and substitute them in the starting non-reduced partial fraction decomposition. At the end of the algorithm, we carry out the same substitution in each summand and add up all the terms. If the two quantities do not come out equal, we have an error in our implementation. This simple error check test helped us catch every single error (to the best of our knowledge) during development. This check is computationally inexpensive and continues to run in the ``production'' version of our code. 
\begin{lemma}
Algorithm \ref{algorithm:main} will come to a halt and produce a semi-reduced partial fraction decomposition recorded in the variable $N$.
\end{lemma}
\begin{proof}
From Step 2 and Step 3 we see that, unless the algorithm has reached a semi-reduced partial fraction, it will proceed with the application of the Szenes-Vergne formula in Step 5. Therefore the algorithm can only halt when it has reached a semi-reduced partial fraction decomposition, and all we need to do is prove that it halts.

Let $u$ be a partial fraction. Its denominator has multiplicands of the form 

\[
\prod_{i=1}^{|\Delta|} \prod_{j=1}^{l_i} \left(1-x^{b_{i,j}\bm \alpha_{i}}\right)^{m_{i,j}},
\]
where, if a vector $\bm \alpha_i\in \Delta$ does not appear in the denominator, we set $l_i=1$ and $m_{i,1}=0$. For the given $u$, form the multiplicity tuple
\[
\bm m(u) = \left( \underbrace{\sum _{j=1}^{l_1} m_{1,j} ~,~ \dots ~,~ \sum_{j=1}^{l_{|\Delta|}}m_{s,j}}_{\text{one coordinate for each }\bm \alpha_i\in \Delta } \right)
\]
The tuple $\bm m(u)$ measures the total multiplicity with which each vector $\bm \alpha\in \Delta$ - or its scalar multiples - appear in the denominator of $u$. Suppose that Algorithm \ref{algorithm:main} removed the partial fraction $u$ in Step 2, and replaced it with partial fractions $ u_1, \dots ,  u_l$ in Step 5. For each of the $ u_i$'s, $\bm m(u_i)$ has its  $t-1$'st coordinate be $1$ larger than the $t-1$'st coordinate of $\bm m (u)$, with all previous coordinates equal. This means that $\bm m(u_i)$ is strictly larger in the lexicographic order than $\bm m(u)$. Since the Szenes-Vergne formula preserves the sum of the coordinates of $\bm m$, there are only finitely many possible values for the $\bm m$-tuple and so our algorithm must come to a halt. 
\end{proof}

To transition from a semi-reduced partial fraction decomposition to a fully reduced partial fraction decomposition, all we need to do is apply the geometric series sum formula. We state this final step as a separate algorithm, but omit any further details.

\begin{algorithm}[Semi-reduced to reduced partial fractions] ~\label{algorithm:semiReducedToReduced}
\begin{itemize}
\item Input: a semi-reduced partial fraction decomposition.
\item Output: a fully reduced partial fraction decomposition.
\end{itemize}
\begin{itemize}
\item[Step 1] Set $R=0$ and set $N$ to be the starting partial fraction, represented as described in Notation \ref{notation:partialFractionData}. 
\item[Step 2] If $N$ is empty, halt, else extract an element $u$ from $N$. If every vector $\bm \alpha$ appearing in the denominator of \eqref{eq:generalFormPartialFraction} appears as in the form $\left(1-x^{b\bm \alpha}\right)^m$ for one fixed multiple $b>0$, add $u$ to $R$ and repeat Step 2 with the next summand; else proceed with Step 3.
\item[Step 3] Let $\bm \alpha$ be vector that appears with different rescaling in the denominator, i.e., a vector for which the denominator has the multiplicands:

\[
u=\frac{p}{ \left(1-x^{b_1 \bm \alpha}\right)^{m_1}\dots \left(1-x^{b_l \bm \alpha}\right)^{m_l} \cdot (\text{terms not involving }\bm \alpha) }
\]
\item[Step 4] Compute the positive least common multiple $L= lcm (b_1, \dots, b_l)$. Here, we recall that negative $b_i$'s allowed. Set 
\[
c_1= \frac{L}{b_1} \qquad \dots\qquad  c_l=\frac{L}{b_l}
\]
\item[Step 5] Transform $u$ using the geometric series sum formulas for \eqref{eq:geometricSeries1}, \eqref{eq:geometricSeries2}:
\[
\begin{array} {rcl}
u&=&\displaystyle \frac{p}{\displaystyle
\left( 1-x^{b_1 \bm \alpha}\right)^{m_1}\dots \left(1-x^{b_l \bm \alpha}\right)^{m_l} \cdot (\text{terms not involving }\alpha) } \\
&=&\displaystyle \frac{\displaystyle p\cdot g_{c_1}^{m_1} \left(x^{b_1\bm\alpha}\right) \cdots g_{c_l} ^{m_l}\left(x^{b_l\bm\alpha}\right) } { \displaystyle \left(1 - x^{L \bm \alpha}\right)^{m_1}\dots \left(1 -x^{L \bm \alpha}\right)^{m_l} \cdot (\text{terms not involving }\alpha) } \\
&=&\displaystyle \frac{\displaystyle p \prod_{i} g_{c_i}^{m_i}\left(x^{b_i\bm\alpha}\right) }{ \displaystyle \left(1- x^{L \bm \alpha} \right)^{ \sum_i m_i} \cdot (\text{terms not involving }\alpha) } 
\end{array}
\]
and put the element back into $N$ for further reduction.
\end{itemize}

\end{algorithm}

\subsubsection{Examples}
In all of the examples above, we use the graded colexicographic order on all sets of vectors that need to be ordered. All examples in this section, including the latex formulas, were generated by our software.
\begin{example}
The vector partition function of $(1,0),(0,1),(1,1)$. This is also the root system $A_2$.
\end{example}
The semi-reduced partial fraction decomposition produced by Algorithm \ref{algorithm:main} is 
\[
\begin{array}{rcl}
\displaystyle 
\frac{1}{(1-x_{1})(1-x_{2})(1-x_{1}x_{2})}&=& \displaystyle 
\frac{-x_{2}^{-1}}{(1-x_{1})^{2}(1-x_{1}x_{2})}  +\frac{x_{2}^{-1}}{(1-x_{1})^{2}(1-x_{2})}\end{array}
\]
which is also reduced. The formulas above were already used in Example \eqref{example:fullyReducedDecompositionIsNotUnique}.

\begin{example}\label{example:partialFractionDecompositionB2}
The vector partition function of $(1,0),(0,1),(1,1), (1,2)$. This is the root system $B_2$ (given in simple coordinates).
\end{example}
\[
\begin{array}{rcl}  \underbrace{\frac{1}{(1-x_{1})(1-x_{2})(1-x_{1}x_{2})(1-x_{1}x_{2}^{2})}}_{(1,1)-(0,1)-(1,0)}&=& \displaystyle  \frac{x_{2}^{-1}}{(1-x_{1})^{2}(1-x_{2})(1-x_{1}x_{2}^{2})}\\&&+\underbrace{\frac{- x_{2}^{-1}}{(1-x_{1})^{2}(1-x_{1}x_{2}^{2})(1-x_{1}x_{2})}}_{(1,2)-2(1,1)+(1,0)}\\&=& \displaystyle  \frac{x_{1} x_{2}^{-1}}{(1-x_{1})^{3}(1-x_{1}x_{2})}\\
&&\displaystyle  +\frac{-x_{1}-x_{2}^{-1}}{(1-x_{1})^{3}(1-x_{1}x_{2}^{2})}\\
&&\displaystyle \displaystyle +\underbrace{\frac{x_{2}^{-1}}{(1-x_{1})^{2}(1-x_{2})(1- x_{1} x_{2}^{2} )}}_{(1,2)-2(0,1)-(1,0)}\\
&=&  \displaystyle \frac{x_{1}x_{2}^{-1}}{(1-x_{1})^{3}(1-x_{1}x_{2})}\\
&&\displaystyle  +\frac{x_{2}^{-3}}{(1-x_{1})^{3}(1-x_{2})}\\
&&\displaystyle +\frac{-x_{1}-x_{2}^{-1}-x_{2}^{-2}- x_{2}^{-3}}{(1-x_{1})^{3}(1 -x_{1}x_{2}^{2} )} \end{array}\]
We have indicated each linear combination used in our universal formula \eqref{eq:SzenesVergneFormula}.

\begin{example}\label{example:partialFractionDecompositionSmallWithElongation}
	The vector partition function of $(1,0,0), (0,1,0), (0,0,1), (2,2,2)$. 
\end{example}

\[
\begin{array}{@{}r@{~}c@{~}l}  \underbrace{\frac{1}{(1-x_{1})(1-x_{2})(1-x_{3})(1-x_{1}^{2}x_{2}^{2}x_{3}^{2})}}_{(1,1,1)-(0,0,1)-(0,1,0)-(1,0,0)}&=&  \frac{-x_{2}^{-1}-x_{2}^{-2}}{(1-x_{1})(1-x_{1}^{2})(1-x_{1}^{2}x_{2}^{2}x_{3}^{2})(1-x_{3})}\\&&+\frac{x_{2}^{-2}x_{3}^{-2}}{(1-x_{1})(1-x_{1}^{2})(1-x_{2})(1-x_{3})}\\&&+\frac{-x_{2}^{-2}x_{3}^{-1}-x_{2}^{-2}x_{3}^{-2}}{(1-x_{1})(1-x_{1}^{2})(1-x_{2})(1-x_{1}^{2}x_{2}^{2}x_{3}^{2})}\\&=&  \frac{-x_{1}x_{2}^{-1}-x_{1}x_{2}^{-2}-x_{2}^{-1}-x_{2}^{-2}}{(1-x_{1}^{2})^{2}(1-x_{1}^{2}x_{2}^{2}x_{3}^{2})(1-x_{3})}\\&&+\frac{-x_{1}x_{2}^{-2}x_{3}^{-1}-x_{1}x_{2}^{-2}x_{3}^{-2}-x_{2}^{-2}x_{3}^{-1}-x_{2}^{-2}x_{3}^{-2}}{(1-x_{1}^{2})^{2}(1-x_{2})(1-x_{1}^{2}x_{2}^{2}x_{3}^{2})}\\&&+\frac{x_{1}x_{2}^{-2}x_{3}^{-2}+x_{2}^{-2}x_{3}^{-2}}{(1-x_{1}^{2})^{2}(1-x_{2})(1-x_{3})}\end{array}
\]
The second equality illustrates Algorithm \ref{algorithm:semiReducedToReduced} in action.

\section{Brion-Vergne decomposition for fully reduced partial fractions} \label{secBVdecomposition}
In the present section, we explain how we can exploit our fully reduced partial fraction decomposition to quickly expand elements of $\mathbb V$ into power series. Let $\partial_i$ be the differential operator on the variables $x_1,\dots, x_k$ acting by:

\[
\partial_i \left(x_j\right) =\begin{cases} 1& \text{if } i=j\\0 &\text{otherwise}.\end{cases}
\]
From the chain rule, $\partial_i$ satisfies 

\[
\partial_i \frac{1}{\left(1-x^{\bm\alpha}\right)^p} = \frac{-(-p)\partial_i( x^{\bm\alpha}) }{\left(1- x^{\bm\alpha}\right)^{p+1}} = \frac{p \alpha_i}{\left(1-x^{\bm\alpha}\right)^{p+1} }
\] 
where $\alpha_i$ is the $i^{th}$ coordinate of $\bm \alpha$, and so $\partial_i$ acts on the vector space $\mathbb V$ given in \eqref{eq:ringOfImportance}.

For any vector $\bm \gamma =\left(\gamma_1, \dots, \gamma_k\right)$ define the differential operator $\partial_{\bm \gamma}$ by 
\[
\partial_{\bm\gamma}:=\sum_{k=1}^n \gamma_k x_k\partial_k
\] 
Given two vectors $\bm\alpha,\bm \beta $, define their scalar (``dot'')
product $\langle\bm\alpha, \bm\beta  \rangle$ by \[\langle\bm\alpha, \bm\beta \rangle= \sum_{k=1}^n \alpha_i\beta_i, \]
where $\alpha_i$, $\beta_i$ are the coordinates of the two vectors. It follows that 
\begin{equation}\label{eq:differentialOperatorToScalarProduct}
\partial_{\bm \alpha} \left(x^{\bm\beta}\right) = \langle\bm\alpha, \bm\beta \rangle x^{ \bm\beta}
\end{equation}
Let $\bm \alpha_1,\dots,\bm \alpha_n$ be linearly independent vectors. Let $\bm \beta_i$ be defined by:
\[
\langle \bm \beta_i, \bm \alpha_k \rangle = \begin{cases}
1&\text{if } i=k\\
0&\text{else}
\end{cases}
\]
We recall from linear algebra that if 
\[
A=\begin{pmatrix}
\alpha_{1,1} & \dots &\alpha_{1, n} \\
\vdots \\
\alpha_{n,1} & \dots &\alpha_{n, n} 
\end{pmatrix}
\]
is the matrix whose $i^{th}$ row contains the coordinates of $\bm\alpha_ i$, then the coordinates of vector $\bm \beta_j$ are given as the $j^{th}$ column of the matrix $A^{-1}$.

We will now apply this to compute the vector partition function coming from a fully reduced partial fraction decomposition. Suppose we want to compute the power series expansion of an element that appears in a partial fraction decomposition such as

\[
u= \frac{x^{\bm \delta}}{ \left(1-x^{\bm \alpha_1}\right)^{m_1} \dots \left(1-x^{\bm \alpha_n} \right)^{ m_n}},
\]
where $\bm \alpha_1, \dots, \bm \alpha_n$ are linearly independent.

From \eqref{eq:differentialOperatorToScalarProduct} it follows that $\partial_{\bm \beta_i} (x^{\bm \alpha_k })= \begin{cases} x^{\bm \alpha_k}& {\mathrm{if~} i=k}\\ {0}&{ \mathrm{ otherwise} }\end{cases}$. Then 
\begin{equation}
\left(\frac{\partial_{\bm \beta_i}}{m_i}+ 1 \right) \left(\frac{1}{(1-x^{\bm \alpha_i})^{m_i}}\right)= \frac{1}{(1-x^{\bm \alpha_i })^{ m_i+1}}
\end{equation}
Therefore we can set $\xi$ as given below and we can compute $u$ as its action:
\begin{equation}\label{eqBVdecomposition}
\begin{array}{rcl}
\xi& =&\displaystyle\prod_{i=1}^n\left(\frac{\partial_{\bm \beta_i}}{1}+1\right)\dots \left( \frac{ \partial_{\bm \beta_i }}{m_i - 1} + 1 \right)\\
&=&\displaystyle \prod_{i=1}^n\prod_{j=1}^{m_i-1}\left(\frac {\partial_{\bm \beta_i}}{j} +1 \right) \\
\displaystyle \xi \left(\frac{1}{1 -x^{\bm \alpha_1}}\dots\frac{1}{1-x^{\bm \alpha_n }}\right) &=& \displaystyle \frac{1}{(1-x^{\bm \alpha_1})^{m_1}} \dots\frac{1}{(1-x^{\bm \alpha_n})^{m_n}}\quad.
\end{array}
\end{equation}

This equality can be found in \cite{BrionVergne}. If we apply the right- hand-side to left-hand-side transformation of \eqref{eqBVdecomposition}, we can transform a fully reduced partial fraction decomposition of an element $ v\in \mathbb V$ to a sum of differential operators acting on fractions of the form $\frac{1}{ 1-x^{\bm \alpha_1}}\dots\frac{1}{1-x^{\bm \alpha_n}}$ with the $\bm \alpha_i$'s in each summand linearly independent. We call such a decomposition a Brion-Vergne decomposition of an element $v\in \mathbb V$.

In order to apply $\xi$ to power series, we also need to compute its action on a monomial $x^{\bm \alpha}$: 
\begin{equation}\label{eq:xiAppliedToMonomial}
\begin{array}{rcll|l}
\displaystyle\xi \left(x^{\bm\alpha}\right) &=&\displaystyle \prod_{i=1}^n\prod_{j=1}^{m_i -1} \left(\frac { \partial_{\bm \beta_i }}{j} + 1\right) \left(x^{ \bm \alpha} \right) \\
&=&\displaystyle \prod_{i=1}^n\prod_{j=1}^{m_i-1}\left(\frac {\langle\bm \beta_i, \bm \alpha \rangle }{j}+1\right) x^{\bm\alpha} && \eqref{eq:differentialOperatorToScalarProduct}
\end{array}
\end{equation}
The coefficient of $x^{\bm \alpha}$ is a product of linear terms of total degree $\displaystyle -n+\sum_{i=1}^n m_{i}$. 

Let $\Lambda=\Lambda(\bm \alpha_1, \dots, \bm \alpha_n)$ be the lattice generated by $\bm \alpha_1, \dots, \bm \alpha_n$, and let $\Lambda_+$ be the set of points generated by non-negative integral linear combinations, i.e.,  
\begin{equation} \label{eq:latticeDefinition}
\begin{array}{rclcl}
\displaystyle \Lambda&=&\Lambda\left(\bm \alpha_1, \dots, \bm \alpha_n\right) &=&\displaystyle \left\{\sum_{i=1}^n a_i \bm \alpha_i| a_i\in \mathbb Z \right\} \\
\displaystyle \Lambda_+&=&\Lambda_+\left(\bm \alpha_1, \dots, \bm \alpha_n\right) &=&\displaystyle \left\{\sum_{i=1}^n a_i \bm \alpha_i | a_i\in \mathbb Z_{\geq 0} \right\} 
\end{array} 
\end{equation}
Here, $\Lambda_+$ is simply the intersection of the lattice $\Lambda$ with the cone generated by the $\bm \alpha_i$'s:
\begin{equation}\label{eq:latticePlusIsLatticeIntersectedWithCone}
\Lambda_+ = \Lambda \cap C\left(\bm \alpha_1, \dots, \bm \alpha_n\right).
\end{equation}

Suppose $\bm \alpha_1, \dots, \bm \alpha_n$ are linearly independent. We now have enough notation to write down the power series expansion of the element: 
\begin{equation}\label{eq:expansionOfW}
\begin{array}{rcl}
w&=&\displaystyle \frac{1}{\left(1-x^{\bm \alpha_1}\right)\dots  \left(1-x^{\bm \alpha_n} \right)} \\
&=&\displaystyle \left(1+x^{\bm \alpha_1} + x^{2\bm \alpha_1}+\dots \right)\dots \left(1+x^{\bm \alpha_n} + x^{2\bm \alpha_n}+\dots \right)\\
&=&\displaystyle \sum_{\bm \gamma\in \Lambda_+} \iota_{\Lambda_+} (\bm\gamma) x^{\bm \gamma},
\end{array}
\end{equation}
where $\iota_{\Lambda_+} (\bm \gamma)$ is the number of ways one can write $\bm\gamma $ as a non-negative integral linear combination of the $\bm \alpha_i $'s. Since the $\bm \alpha_i$'s are linearly independent, it follows that 

\begin{equation}\label{eqQuasiNumberNE}
\iota_{\Lambda_+}(\bm \gamma) = \begin{cases}1 &\text{if~} \bm\gamma\mathrm{~lies~in~} \Lambda_{+} \\ 
0& \text{otherwise.}\end{cases}
\end{equation}
Let us extend the notation \eqref{eqQuasiNumberNE} to an arbitrary set $X$ by setting $\iota_X$ to be the indicator function of $X$, that is:
\begin{equation}\label{eq:indicatorFunction}
\iota_{X}(\bm \gamma) = \begin{cases}1 &\text{if~} \bm\gamma\mathrm{~lies~in~} X \\ 
0& \text{otherwise.}\end{cases}
\end{equation}
Abbreviate by $C$ the cone $C(\bm \alpha_1, \dots, \bm \alpha_n)$. Then \eqref{eq:latticePlusIsLatticeIntersectedWithCone} can be rewritten as:
\begin{equation}\label{eq:indicatorFunctionLambdaPlus}
\iota_{\Lambda_+} = \iota_{\Lambda}\cdot \iota _{C}
\end{equation}
Let 
\[
u=\displaystyle \sum_{\bm \gamma\in\Lambda} P(\bm\gamma)x^{\bm \gamma}.
\]
We are now ready to find a closed form formula for $P(\bm\gamma)$ and in this way expand $u$ into formal power series. 
\[
\begin{array}{rcll|l}
\\
u&=&\displaystyle \frac{x^{\bm\delta}}{\left(1-x^{\bm \alpha_1}\right)^{m_1} \dots \left(1 - x^{\bm \alpha_n}\right)^{m_n}}\\
&=&\displaystyle x^{\bm \delta } \xi \left(\frac{1}{\left(1-x^{\bm \alpha_1}\right) \dots \left( 1-x^{\bm \alpha_n }\right)} \right) & &\eqref{eqBVdecomposition} \\
&=&\displaystyle x^{\bm \delta } \xi \left(\sum_{\bm \lambda\in \Lambda_+ } \iota_{\Lambda_+}( \bm \lambda) x^{\bm \lambda}\right) && \eqref{eq:expansionOfW} \\
&=&\displaystyle x^{\bm \delta } \sum_{\bm\lambda\in \Lambda_+ } \iota_{\Lambda_+}(\bm \lambda) \xi \left(x^{\bm \lambda}\right)\\
&=&\displaystyle x^{\bm \delta } \sum_{\bm\lambda\in \Lambda_+ } \iota_{\Lambda_+}(\bm \lambda) \prod_{ i=1}^n \prod_{ j=1 }^{ m_i-1} \left( \frac {\langle\bm \beta_i, \bm \lambda \rangle}{j} +1 \right) x^{\bm\lambda}   &&\eqref{eq:xiAppliedToMonomial}\\
&=&\displaystyle \sum_{\bm \gamma-\bm \delta \in \Lambda_+ } \iota_{\Lambda_+}(\bm \gamma-\bm \delta) \prod_{i=1}^n \prod_{j=1 }^{ m_i - 1}\left(\frac {\langle\bm \beta_i, \bm \gamma-\bm \delta \rangle}{j}+1\right) x^{\bm\gamma} &&\text{set }\bm\gamma=\bm\lambda + \bm \delta
\end{array}
\] 
and, comparing coefficients for $u$, we get:
\begin{equation}\label{eq:vectorPartitionFunctionOneFraction}
\begin{array}{rcl}
\displaystyle P(\bm\gamma)&=& \displaystyle\iota_{\Lambda_+}(\bm \gamma-\bm \delta) \prod_{i=1 }^n \prod_{j=1}^{m_i-1}\left(\frac { \langle \bm \beta_i, \bm \gamma-\bm \delta \rangle }{j} +1 \right) \\
&=&\displaystyle\iota_{\Lambda_+}(\bm\gamma-\bm\delta)\prod_{i=1}^n\binom{\langle\bm\beta_i,\bm\gamma-\bm\delta\rangle+m_i-1}{m_i-1}
\end{array}
\end{equation}
where we've converted the last product into a binomial coefficient.

Suppose $\bm\gamma -\bm \delta$ is an integral but not positive linear combination of the $\bm\alpha_i$'s, that is, $\bm\gamma -\bm \delta \in \Lambda\setminus \Lambda_+$. Then $P(\bm \gamma) $ is zero. If we drop the multiplicand $\iota_{\Lambda_+}(\bm\gamma-\bm \delta )$ from formula \eqref{eq:vectorPartitionFunctionOneFraction} and $\langle\bm \beta_i, \bm\gamma -\bm\delta \rangle$ is negative enough for all $i$, then the so-modified formula \eqref{eq:vectorPartitionFunctionOneFraction} will no longer hold. However, if $\langle\bm \beta_i, \bm\gamma -\bm\delta \rangle$ is negative but small in absolute value, formula \eqref{eq:vectorPartitionFunctionOneFraction} with $\iota_{\Lambda_+}(\bm\gamma-\bm \delta)$ omitted will continue to produce correct results because  $\binom{n}{m}=0$ if $m>n$ and $n>0$. More precisely, we have the following.

\begin{remark}\label{remark:noTau}
Suppose we omit the multiplicand $\iota_{\Lambda_+}(\bm\gamma-\bm\delta)$ from \eqref{eq:vectorPartitionFunctionOneFraction}.
\begin{itemize}
\item[(a)] The so modified formula continues to produce correct values for $P(\bm \gamma)$ for all $\bm \gamma$ for which 
\begin{itemize}
\item[(1)] $\bm\gamma-\bm\delta\in \Lambda$ 
\item[(2)] $\langle \bm\beta_i, \bm\gamma-\bm \delta \rangle+m_i-1\geq 0$ for some $m_i$.
\end{itemize}
		
\item[(b)] The so modified formula ceases to produce correct values for $P(\bm \gamma)$ for $\bm \gamma$ for which 
\begin{itemize}
\item[(1)] $\bm\gamma-\bm\delta\in \Lambda$ 
\item[(2)] $\langle \bm\beta_i, \bm\gamma-\bm \delta \rangle+m_i-1< 0$ for all $m_i$.
\end{itemize}
\end{itemize}
\end{remark}

\subsection{Example}
All formulas in this section were generated by our software.
\begin{example}
The Brion-Vergne decomposition of the partial fraction of $B_2$ given in Example \ref{example:partialFractionDecompositionB2} is the following.
\end{example}
\[
\begin{array}{@{}r@{}c@{}l}  
\displaystyle \frac{1}{(1-x_{1})(1-x_{2})(1-x_{1}x_{2})(1-x_{1}x_{2}^{2})}&=&
\displaystyle x_{1}x_{2}^{-1}\cdot \frac{1}{2}\left(x_{1}\partial _{1}-x_{2}\partial _{2}\right)^{2}\cdot\\
&&\displaystyle \frac{1}{(1-x_{1})(1-x_{1}x_{2})}
\\&+&\displaystyle \left(-x_{1}-x_{2}^{-1}-x_{2}^{-2}-x_{2}^{-3}\right) \cdot 
\\&&\displaystyle \frac{1}{2}\left (x_{1}\partial _{1}-\frac{x_{2}\partial _{2}}{2}\right)^{2}\cdot
\\&&\displaystyle \frac{1}{(1-x_{1})(1-x_{1}x_{2}^{2})}
\\&+&\displaystyle x_{2}^{-3}\cdot \frac{1}{2}\left(x_{1}\partial _{1}\right)^{2}\cdot \frac{1}{(1-x_{1})(1-x_{2})}\end{array}
\]

\begin{example}
The Brion-Vergne decomposition of the partial fraction of \\
$(1,0,0),(0,1,0), (0,0,1), (2,2,2)$ given in Example \ref{example:partialFractionDecompositionSmallWithElongation} is the following.
\end{example}
\[
\renewcommand{\arraystretch}{2.3}
\begin{array}{@{}r@{}c@{}l}
\displaystyle \frac{1}{(1-x_{1})(1-x_{2})(1-x_{3})(1-x_{1}^{2}x_{2}^{2}x_{3}^{2})}&=&\displaystyle (-x_{1}x_{2}^{-1}-x_{1}x_{2}^{-2}-x_{2}^{-1}-x_{2}^{-2})\cdot 
\\&&\displaystyle \left(\frac{x_{1}\partial _{1}}{2}-\frac{x_{2}\partial _{2}}{2} \right) \cdot 
\\&& \displaystyle  \frac{1}{(1-x_{1}^{2})(1-x_{1}^{2}x_{2}^{2}x_{3}^{2})(1-x_{3})}
\\&+& \displaystyle (-x_{1}x_{2}^{-2}x_{3}^{-1}-x_{1}x_{2}^{-2}x_{3}^{-2}
\\&& \displaystyle  ~-x_{2}^{-2}x_{3}^{-1}-x_{2}^{-2}x_{3}^{-2})\cdot 
\\&& \displaystyle \left(\frac{x_{1}\partial _{1}}{2}-\frac{x_{3}\partial _{3}}{2}\right)\cdot
\\&& \displaystyle \frac{1}{(1-x_{1}^{2})(1-x_{2})(1-x_{1}^{2}x_{2}^{2}x_{3}^{2})}\\
&+& \displaystyle (x_{1}x_{2}^{-2}x_{3}^{-2}+x_{2}^{-2}x_{3}^{-2})\cdot \left(\frac{x_{1}\partial _{1}}{2}\right)\cdot
\\&& \displaystyle \frac{1}{(1-x_{1}^{2})(1-x_{2})(1-x_{3})}
\end{array}
\]
\section{The algebraic expressions for the vector partition function do not depend on the distance from the walls}\label{secWallDist}

In the previous sections, we showed how to express the generating function of a vector partition function as a sum of fully reduced partial fractions by applying in Algorithm \ref{algorithm:main} a single universal formula \eqref{eq:SzenesVergneFormula} - the generalized Szenes-Vergne formula. Then, using formula \eqref{eq:vectorPartitionFunctionOneFraction}, we showed how to express the reduced partial fractions as quasipolynomials. This resulted in a formula for the vector partition function of $\bm \gamma$ that is a sum of polynomials in $\bm \gamma$ with coefficients of the form 
\begin{equation}\label{eq:iotaLambdaPlusForm}
\iota_{\Lambda_+}(\bm \gamma-\bm \delta) = \iota_{\Lambda}(\bm\gamma-\bm \delta) \iota _{C}(\bm \gamma-\bm \delta)
\end{equation}
where $\bm\delta$ runs over a set of finitely many vectors. However, in our introduction, we were speaking of quasipolynomial chambers with walls passing through the origin, which implies formulas with coefficients of the form
\begin{equation}\label{eq:iotaLambdaForm}
\iota_{\Lambda} (\bm\gamma-\bm \delta)\iota_{C}(\bm\gamma)
\end{equation}

In this section, we will resolve this ``discrepancy'' by proving that our final formula for the vector partition function does not change values when we simultaneously replace all coefficients of the form \eqref{eq:iotaLambdaPlusForm} with coefficients of the form \eqref{eq:iotaLambdaForm}. We note that by Remark \ref{remark:noTau}(b), replacing \eqref{eq:iotaLambdaPlusForm} by \eqref{eq:iotaLambdaForm} does change the values of  individual summands, however by the claim of the present section, the simultaneous replacement of all summands keeps the total sum unchanged. More precisely, we aim to show the following theorem.

\begin{theorem}\label{theorem:mainTheory}~ 
Let $\bm \alpha_1, \dots, \bm \alpha_s$ be vectors with non-negative coordinates that span the entire space.
\begin{itemize}
\item[(a)]
There exist:

\begin{itemize}
\item a lattice $\Lambda$ (see \eqref{eq:latticeDefinition}) 
\item a collection of finitely many pointed polyhedral cones $\mathcal I=\{ D_1, \dots, D_N\}$ (see \eqref{eq:coneDefinition})
\item finitely many integer-coordinate vectors $\bm\delta_1, \dots,\bm \delta_m$ 
\item finitely many polynomials $ Q_{i,j}$ indexed by the $\bm\delta_i$'s and the $D_i$'s
\end{itemize}  
such that 
\[
\begin{array}{rcl}	
P_\Delta (\bm\gamma) &= &P_{\Delta}(\gamma_1,\dots, \gamma_n) \\
&=&\displaystyle \begin{cases}
Q_{i,j}(\gamma_1,\dots, \gamma_n) & \text{when }\bm \gamma\in D_i\cap \left(\bm\delta_j+\Lambda \right)\\
0& \text{in all other cases}
\end{cases}\\
&=&\displaystyle \sum_{i,j} \iota_{D_i}(\bm \gamma) \iota _{\delta_j + \Lambda}(\bm \gamma) Q_{i,j}(\bm \gamma)\\
\end{array}
\]
\end{itemize}
\item[(b)] The polyhedral cones $D_i$ from (a) can be chosen so that the planes giving each wall of each cone is spanned by $n-1$ vectors in $\Delta$.
\end{theorem}
The earliest proof known to the author is in \cite[Theorem 3.26]{DeConciniProcesi:GeometryToricArrangements}. In the remainder of this section, we include our own proof as we find it to be purely elementary, self-contained, and independent of the proof of \cite[Theorem 3.26]{DeConciniProcesi:GeometryToricArrangements}. A first version of our proof was discovered independently by the author and posted online in the unpublished article \cite{Milev:PartialFractions}; the present proof is a rewrite of this original work.

So far, our use of the term ``quasipolynomial'' was informal; let us now state a formal definition.

\begin{definition}
A function $f:\mathbb R^n\to \mathbb C$ is a quasipolynomial if there exist finitely many vectors $\bm \delta_1, \dots, \bm\delta_N \in \mathbb Z^n$, finitely many polynomials $P_1, \dots, P_N$ and a full-rank lattice $\bm \Lambda \subset \mathbb Z^n$  such that $f(\bm \gamma)$ is given by:

\[
f(\bm \gamma)  = \sum_{i=1}^N \iota_{\bm \delta_i+\Lambda} (\bm \gamma) P_i(\bm \gamma)
\]
\end{definition}
Indicator functions of full-rank lattices can be expressed using the indicator functions of their full-rank sub-lattices. Suppose $\Omega\subset \Lambda$ is a sub-lattice of $\Lambda$ of full rank. Then 
\begin{equation}\label{eq:latticeToSublattice}
\iota_{\bm \delta + \Lambda}(\bm \gamma) = \sum_{ \substack{ \bm\delta - \bm \mu_j\in \Lambda \\ \bm \mu_j\text{ run over representatives of } \mathbb Z^n/\Omega}} \iota_{\bm \mu_j +\Omega}(\bm \gamma)
\end{equation}
gives the desired transition from indicator function of a lattice to an indicator function of a sub-lattice. Given a finite sum of lattice translates over different lattices $\Lambda_1, \dots$, using the trick \eqref{eq:latticeToSublattice}, we can rewrite the sum as a sum that uses lattice translates of $\displaystyle \Omega = \bigcap_i \Lambda_i$, so we can assume sums of lattice indicator functions use translates of a single lattice. Furthermore, given lattices $\Lambda, \Psi\subset\mathbb Z^n$ of full rank and vectors $\bm \delta, \bm \varepsilon\in \mathbb Z^n$, we note that 
\begin{equation}\label{eq:productOfIndicatorsOfLatticeShifts}
\iota _{\bm\delta +\Psi}(\bm \gamma)\iota _{\bm\varepsilon+\Lambda}(\bm \gamma) = \sum_{\substack { \bm \mu_i \text{ runs over the representatives of } \mathbb Z^n /\Psi\cap \Lambda \\ \bm \mu_i-\bm \delta \in \Psi\\ \bm \mu_i -\bm \varepsilon \in \Lambda} } \iota _{\bm \mu_i+ \Psi\cap \Lambda }(\bm \gamma)
\end{equation}
and so products of indicator functions of translates of two lattices can be reduced to linear combination of indicator functions of translates of lattices. To summarize this discussion, the linear combinations of indicator functions of lattice translates form an algebra, and therefore quasipolynomials form an algebra.

With this terminology, Theorem \ref{theorem:mainTheory} can be rephrased as the following.

\begin{corollary}\label{corollary:mainTheory}
There exist finitely many chambers $D_i$ with disjoint interiors and such that their walls spanned by vectors in $\Delta$, such that $P_{\Delta}(\bm \gamma)$ is a quasipolynomial over each $D_i$.
\end{corollary}

\subsection{Computations with lattices}
We briefly review the algorithms needed to compute with lattices for the reader's convenience and as reference for computer implementations. Recall from \eqref{eq:latticeDefinition} that a full-rank lattice $\Lambda=\Lambda(\bm\alpha_1, \dots, \bm \alpha_n)$ is the set of all integer-coefficient linear combinations of the linearly independent vectors with rational coordinates $\bm \alpha_i$. 

\subsubsection{Integral Gaussian Elimination}
To represent a $\Lambda$ on a computer, form the vector-row matrix $L$ of the $\bm \alpha_i$'s:
\[
L=
\begin{pmatrix}
	\alpha_{1,1}&\dots& \alpha_{1,n}\\
	&\vdots&\\
	\alpha_{n,1}&\dots& \alpha_{n,n}\\
\end{pmatrix}.
\]
Now, carry out modified integral Gaussian Elimination by rows on $L$. This modified Gaussian elimination consists of the following steps.
\begin{algorithm}[Integral Gaussian Elimination]\label{algorithm:integralGaussianElimination} ~
	\begin{itemize}
		\item Input: a matrix $L$ with rational coordinates.
		\item Output: a matrix $L$ (modified in-place) with rational coordinates in upper-triangular form whose vector-rows generate the same lattice as the input.
	\end{itemize}
	\begin{itemize}
		\item[Step 1.] Find the greatest common divisor $d$ of all denominators of all entries of $L$. 
		\item[Step 2.] Multiply $L$ by $d$. 
		\item[Step 3.] Set $i=1$. The variable $i$ holds the index of the current ``pivot row'', where we have borrowed the term from regular Gaussian Elimination.
		\item[Step 4.] For each column of $L$ of index $j$:
		\begin{itemize}
			\item[Step 4.1.] Find the first index $a\geq i$ so that the $(a,j)^{th}$ entry of $L$ is non-zero. If there's no such $a$, go back to Step 2 and proceed with the next column $j$.  
			\item[Step 4.2.] If $a>i$, swap the $i^{th}$ and $a^{th}$ rows.
			\item[Step 4.3.] Set $x$ to be the $(i,j)^{th}$ entry of $L$. Here, $x$ is the analog of the pivot element from regular Gaussian Elimination.
			\item[Step 4.4.] If $x<0$, multiply the $i^{th}$ row of $L$ by $-1$ and replace $x$ by $-x$. 
			\item[Step 4.5.] For every row $k>i$:
			\begin{itemize}
				\item[Step 4.5.1.] Let $y$ be the $(k,j)^{th}$ entry of $L$. 
				\item[Step 4.5.2.] Let $y= q\cdot x+r$ be the integer division with remainder of $y$ by $x$ with $0\leq r< x$.
				\item[Step 4.5.3.] If $q=0$:
				\begin{itemize}
					\item[Step 4.5.3.1] If $r=0$ this means that $y=0$; increment $k$ and go back to Step 4.5 until done.
					\item[Step 4.5.3.2] Else, swap the $k^{th}$ row of $L$ with the $i^{th}$ row and go back to Step 4.5.1.
				\end{itemize}
				\item[Step 4.5.4.] Else $q\neq 0 $; apply the row transformation $ \bm \alpha_k \mapsto \bm \alpha_k - \bm q\cdot\bm \alpha_i$ to the $k^{th}$ row of $L$. 
				\item[Step 4.5.5.] Set $x=r$.
				\item[Step 4.5.6.] Go back to Step 4.5.1.
			\end{itemize} 
			\item[Step 4.6.] For every row $k<i$:
			\begin{itemize}
				\item[Step 4.6.1.] Let $y$ be the $(k,j)^{th}$ entry of $L$. 
				\item[Step 4.6.2.] Let $y= q\cdot x+r$ be the integer division with remainder of $y$ by $x$ with $0\leq r< x$.
				\item[Step 4.6.3.] Apply the row transformation $ \bm \alpha_k \mapsto \bm \alpha_k - \bm q\cdot\bm \alpha_i$ to the $k^{th}$ row of $L$. 
				\item[Step 4.6.4.] Increment $k$ and go back to Step 4.6 until done. 
			\end{itemize}
			\item[Step 4.7.] Increment $i$.
			\item[Step 4.8.] Increment $j$ and go back to Step 4 until done.
		\end{itemize} 
		\item[Step 5.] Divide $L$ by $d$.
	\end{itemize}
\end{algorithm}
Given two lattices $L_1$ and $L_2$, they are equal if and only if their reductions with Integral Gaussian Elimination coincide. In our computer implementation, we record lattices by storing their reduced matrices computed using the algorithm above. 

The Integral Gaussian Elimination algorithm is similar to regular Gaussian elimination. The differences are:
\begin{itemize}
	\item Division/multiplication by scalar is not allowed, unless the scalar is $-1$ (Step 4.4).
	\item Row subtractions below the main diagonal (Step 4.5.4) do not eliminate the leading entry, but rather decrease its absolute value, after which we swap rows. The elimination sequence below the main diagonal - Step 4.5.1 to Step 4.5.6 - reproduces the effect of the Euclidean algorithm carried out with starting inputs the pivot entry $x$ and the leading entry $y$ of the other row.
	\item Unlike regular Gaussian elimination, row subtractions above the pivot row (Step 4.6) may fail to eliminate the $(k,j)^{th}$ entry - they just make sure it is smaller than the $(i,j)^{th} $ entry of the pivot row. 
	\item When the starting matrix $L$ is invertible, Algorithm \ref{algorithm:integralGaussianElimination} will produce an upper triangular matrix with positive numbers on the main diagonal, and the possibly non-zero entries above the main diagonal will be non-negative and smaller than the diagonal entries in the same column.
\end{itemize} 

\subsubsection{Dual lattice computation}\label{section:dualLatticeComputation}
Given a lattice $\Lambda$ of full rank, denote by $\Lambda^*$ the dual lattice:
\[
\Lambda^*= \left\{\bm \gamma |\langle\bm \gamma,\bm \lambda\rangle\in \mathbb Z \text{ for all } \bm \lambda \in \Lambda\right\}
\]
To compute $\Lambda^*$, take the inverse transpose of the matrix of $\Lambda$ and reduce it using Integral Gaussian Elimination. 

\subsubsection{Common refinement of lattices}\label{section:commonRefinement}
Given lattices $\Psi$ and $\Lambda$ of full rank, denote by
\[
\langle \Psi, \Lambda \rangle = \text{ the lattice generated by the elements of } \Lambda\cup\Psi
\]
To compute the lattice $\langle\Psi , \Lambda\rangle$, augment the matrix of $\Psi$ by the matrix of $\Lambda$ by combining the matrix rows and run Integral Gaussian Elimination Algorithm \ref{algorithm:integralGaussianElimination}. The resulting non-zero rows of the reduced matrix hold the reduced matrix of $\langle \Psi, \Lambda \rangle$.
\subsubsection{Intersection of lattices}
To compute the intersection of two lattices $\Lambda$ and $\Psi$, we note that 
\[
\Lambda\cap \Psi = \langle\Lambda^*,\Psi^*  \rangle^*,
\]
which we can compute using Sections \ref{section:dualLatticeComputation} and \ref{section:commonRefinement}.
\subsubsection{Lattice extras}

In the following sections, we need to compute lattices given by:  
\[
\Omega(\bm a, \Lambda) = \{\bm \gamma \in \Lambda|\langle\bm a, \bm \gamma\rangle\in \mathbb Z\},
\] 
where $\Lambda$ is a full-rank lattice and $\bm a$ is a fixed rational coordinate vector. This is done with the following algorithm.
\begin{algorithm}\label{algorithm:computeSubLatticeIntegralScalarProducts}~
	\begin{itemize}
		\item Compute the matrix $M$ of $\Lambda^*$.
		\item Append $\bm a$ as the $n+1^{st}$ row of $M$. 
		\item Reduce $M$ using the Integral Gaussian Elimination (Algorithm \ref{algorithm:integralGaussianElimination}).
		\item The last row of the so reduced $M$ must be zero, drop it. Let $N$ be the resulting matrix.
		\item The inverse transpose matrix $(N^{-1})^t$ now holds the desired reduced matrix of $\Omega(\bm a, \Lambda)$.
	\end{itemize} 	
\end{algorithm}
Given lattices $\Lambda,\Theta $ of full rank and vectors $\bm b, \bm \gamma$, in the following sections, we need algorithms to compute the lattice $\Psi$ given by

\[
\Psi=\Psi(\bm b, \bm \alpha, \Lambda, \Theta)=\left\{\bm \gamma \in \Theta | \langle\bm b, \bm \gamma\rangle\bm \alpha\in \Lambda \right\}
\]
This is done with the following computation. Let $L$ be the matrix of $\Lambda$ (with coordinates of the generating vectors arranged in rows). Compute 
\[
\bm \nu=\left(L^{t}\right)^{-1}\cdot \bm \alpha, 
\]
where $\bm \alpha$ is regarded as a vector-row. Let $q$ be the least common multiple of the denominators  of the coordinates of $\bm \nu$ and let $p$ be the greatest common divisor of the numerators. Then  $\Psi$ is computed as:
\begin{equation}\label{eq:computeSubLatticeScalarProductTimesDirectionInLattice}
\Psi=\Psi(\bm b, \bm \alpha, \Lambda, \Theta)=\Omega(a,\Theta)= \frac{q}{p}\Omega(\bm a, \Theta)	
\end{equation}
where, in the last equality, multiplication by $\frac{q}{p}$ means simultaneous rescaling of all vectors in the lattice.
\subsubsection{Summing quasipolynomials}\label{section:summingQuasipolynomials}
Let $\Lambda$ and $\Omega$ be two lattices and let $\bm \delta_1, \dots, \bm \delta_N$ be representatives of $\mathbb Z^n/\Lambda$ and $\bm \nu_1, \dots, \bm \nu_M$ be representatives of $\mathbb Z^n/\Omega$. Let $Q(\bm \gamma)$ and $P(\bm \gamma)$ be quasipolynomials given by
\[
Q(\bm \gamma) = \sum_{i=1}^N \iota_{\bm \delta_i+\Lambda}(\bm \gamma) Q_i(\bm \gamma)\qquad
P(\bm \gamma) = \sum_{i=1}^M \iota_{\bm \delta_j+\Omega}(\bm \gamma) P_i(\bm \gamma)
\]
for some polynomials $Q_i, P_i$. Using formula \eqref{eq:latticeToSublattice}, we can rewrite both $Q(\bm \gamma)$ and $P(\bm \gamma)$ as a sum of polynomials with coefficients in $\iota_{\bm \mu_i+\Omega\cap \Lambda }(\bm \gamma)$, where the $\bm \mu_i$'s run over the representatives of $\mathbb Z^n/\left(\Omega\cap \Lambda\right)$. This gives us a procedure to compute the sum $Q(\bm \gamma)+P(\bm \gamma)$.
\subsection{An elementary algorithm for computing the vector partition function}

The core of our argument for proving Theorem \ref{theorem:mainTheory} is the following. 
For a while, let us set aside all considerations from the previous sections, and 
come up with an elementary algorithm for computing the vector partition function that does not use generating functions. The resulting quasipolynomial formulas will be valid over pointed polyhedral cones whose walls are not cut off by planes spanned by $n-1$-element subsets of the original vectors. 
At the same time, Section \ref{secBVdecomposition} provided us with formulas valid over translates of polyhedral cones with walls parallel to planes spanned by $n-1$-element subsets of $\Delta$. Thus, the chambers produced by the two approaches are ``glued'' by one another, which proves Theorem \ref{theorem:mainTheory}. The remainder of Section \ref{secWallDist} is dedicated to making this argument precise.

Suppose we want to compute the vector partition function $P_{\Delta}$ with respect to the vectors 
\[
\Delta=\left\{\bm\alpha_1,\dots, \bm\alpha_{s-1},\bm \alpha_s\right\}.
\] 
Suppose, by induction, we have already computed the partition function $P_{\Gamma}$ with respect to the first $s-1$ vectors 
\[
\Gamma=\left\{\bm\alpha_1,\dots, \bm\alpha_{s-1} \right\} = \Delta\setminus \{\bm \alpha_s\}
\] 
In a vector partition decomposition of $\bm \gamma$ with respect to $\Delta$, if one chooses exactly $t$ times the vector $\bm \alpha_s$, one can decompose $\bm \gamma$ in exactly $P_{\Gamma}\left(\bm \gamma-t\bm\alpha_s\right)$ ways. Therefore 
\begin{equation}\label{eq:recursiveRaySum}
P_{\Delta}(\bm\gamma)= \sum_{t=0}^{\infty}P_{\Gamma}\left(\bm \gamma- t\bm\alpha_s\right) .
\end{equation}
This is a finite sum as the partition function is zero for $t$ sufficiently large. \textbf{If} we had a polynomial formula for $P_{\Gamma}$, then the sum \eqref{eq:recursiveRaySum} would be linear combination of the expressions  \begin{equation}\label{eq:bernoulliSum}
\displaystyle B_k(X)= \sum_{t=0}^{t=X} t^{k}
\end{equation}
for some upper bound $X$. Formula \eqref{eq:bernoulliSum} is a Bernoulli sum (also known as Faulhaber's formula) that sums to a polynomial of degree $k+1$ in $X$. \textbf{If} the final point of summation $X$ were a linear polynomial in the coordinates $\gamma_1, \dots, \gamma_n$ of $\bm \gamma$, then formula \eqref{eq:bernoulliSum} would be a Bernoulli sum in which we have substituted-in a linear polynomial, so our vector partition function would be a polynomial in $\gamma_1, \dots, \gamma_n$. 

Our argument will have to be considerably more complicated as both of the highlighted \textbf{if}'s above do not hold in general. To make our argument work, we need to resolve several issues which we present in the following subsections.

\subsection{Normal chambers with respect to a vector}\label{sec:normalChambers}
For the proof of our main Theorem \ref{theorem:mainTheory}, we proceed by induction on the number of vectors $s$.

Let $\mathbb R^n_+ = \{\bm\gamma\in \mathbb R^n | \bm \gamma \text{ has non-negative coordinates}\}$. Suppose by induction hypothesis that we have managed to subdivide $\mathbb R^n_+$ into finitely many chambers such that, over each chamber $C$, $P_{\Gamma}$ is either given by a quasipolynomial or is zero over the interior of $C$. 

The sum \eqref{eq:recursiveRaySum} needs to be carried over a ray in the direction of $\bm \alpha_s$, which will intersect multiple pointed polyhedral cones given by our induction hypothesis. Therefore we seek to split the sum \eqref{eq:bernoulliSum} into multiple sums: 

\begin{equation}\label{eq:partialFractionRefinement1}
P_{\Delta}(\bm\gamma)= \sum_{j=0}^{r-1}\left( \sum_{t=t_j+1}^{t=t_{j+1}} P_{\Gamma}\left(\bm \gamma- t\bm\alpha_s\right)\right),
\end{equation}
where $t_0, t_1,t_2,\dots, t_{r}$ are integers chosen so that $t_0=-1$ and the vector $ \bm\gamma-t\bm \alpha_s$ lies entirely within one pointed polyhedral cone among the cones that give all formulas for $P_{\Gamma}$. 

For our induction basis, we start with $n$ vectors and the chamber \\ $C\left(\bm\alpha_1, \dots, \bm \alpha_n\right)$, extended to a subdivision of $\mathbb R^n_+$ in an arbitrary fashion. Over $C\left(\bm\alpha_1, \dots, \bm \alpha_n\right)$, the vector partition function equals $P_{\bm\alpha_1, \dots , \bm \alpha_n}(\bm \gamma)=$  $ \iota_{\Lambda \left(\bm\alpha_1, \dots, \bm \alpha_n \right)} (\bm \gamma)$, and the vector partition function is zero over the interiors of all other chambers in the starting subdivision.

We now proceed with the induction step. Suppose the starting point $\bm\gamma$ lies in any chamber $C$ given by normal-vector inequalities as in \eqref{eq:coneDefinitionScalarProduct} such that all points in $C$ have non-negative coordinates. Since $\Delta$ consists of vectors with non-negative integer coordinates, at least one point on the ray $\bm \gamma -t\bm \alpha_s$ has a negative coordinate, so the ray must exit our chamber through at least one wall. The ray may exit through $k$ walls simultaneously if it happens to hit an $n-k$-dimensional ``edge'' of the chamber.

Suppose a normal of the exit wall of chamber $C$ is $\bm b$. It follows that $\langle \bm b,\bm\alpha_s \rangle > 0$. 
\begin{definition}\label{def:normalChamber}
If a chamber (pointed polyhedral cone) $C$ has a unique wall with normal $\bm b$ so that 
\[
\langle \bm b,\bm\alpha \rangle > 0,
\] 
then we say that $C$ is normal with respect to $\bm\alpha$.
\end{definition}

In other words, $C$ is normal with respect to $\bm\alpha_s$ if the ray $\bm \gamma - t\bm \alpha _s$ is guaranteed to exit the chamber through a single distinguished wall  or one of its edges. 

\begin{definition}[Normally separated chambers]\label{def:normallySeparated}
We say that two chambers $C, D$ are normally separated if:
\begin{itemize}
\item[(a)] Their interiors are disjoint.
\item[(b)] If the chambers $C,D$ have $n-1$ common points that span a plane $P$, then \[C\cap P=D\cap P.\] 
\end{itemize}
\end{definition}
Condition (b) states that if two chambers are neighbors that touch along a patch of dimension $n-1$, then all points along the touching wall are common. If the ray $\bm \gamma-t\bm\alpha_s$ must exit through a wall with normal $\bm a$, then it will pass through a neighbor on the other side with normal $-\bm a$, if such a neighbor exists. Condition (b) then guarantees that if such a neighbor exists, then it is unique.
\begin{definition}[Normally separated collection]\label{def:normallySeparatedCollection}
We say that a collection of chambers is normal if every two chambers in it are normally separated.
\end{definition}

\begin{definition}[Subdivision]\label{definition:subdivision} Let $\mathcal I$ be a set of closed sets with pairwise disjoint interiors. We way that the set of closed sets $\mathcal I'$ subdivides $\mathcal I$ if:
\begin{itemize}
\item The interiors of sets in $\mathcal I' $ are disjoint.
\item 
\[
\bigcup_{D\in \mathcal I'} D = \bigcup_{D\in \mathcal I} D
\]
\item For every $D\in \mathcal I$ there exist finitely many $D_1, \dots, D_p\in \mathcal I'$ such that $D=D_1\cup\dots\cup D_p$.	
\end{itemize}

\end{definition}

Every collection of chambers with disjoint interiors can be subdivided into a collection of chambers with normal neighbors. Indeed, take the set of all hyperplanes of all walls of all the chambers. They subdivide space into finitely many chambers; then the chambers whose interiors intersect the interiors of the starting chambers give us the desired collection with normal neighbors. We will present a different algorithm for ensuring normal neighbors by subdivision in Algorithm \ref{algorithm:extendSubdivisionToNormal} below.

\begin{definition}[Normal with respect to direction]\label{def:normalChamberCollection}
We say that a collection of chambers is normal with respect to $\bm \alpha$ if the following hold.
\begin{itemize}
\item[(a)] Each chamber is normal with respect to $\bm\alpha$.
\item[(b)] The chambers are normally separated, i.e., satisfy Definition \ref{def:normalChamberCollection}.
\end{itemize}
\end{definition}

Now fix $\bm \gamma$ and assume that there exist chambers $C_1, \dots, C_r$ that satisfy the following induction setup conditions.
\begin{condition}[Induction setup] ~\label{condition:summationProperties}
\begin{itemize}
\item[(a)] $\bm \gamma\in C_1$.
\item[(b)] The ray $\bm \gamma -t\bm \alpha_s$ passes through the chambers $C_1, \dots, C_r$, in this order.
\item[(c)] The collection $\left\{ C_1, \dots, C_r\right\}$ is normal with respect to $\bm\alpha_s$.
\item[(d)] The ray $\bm \gamma -t\bm \alpha_s$ has no points outside of $\displaystyle\bigcup_{i=1}^r C_i$ but inside $C(\bm \alpha_1, \dots, \bm \alpha_s)$ (the cone generated by the first $s$ vectors).
\end{itemize}
\end{condition}
We will ensure that Condition \ref{condition:summationProperties} can be satisfied in Section \ref{sec:subdivision}.

Since the collection is normal, each chamber $C_i$ has a unique exit wall with normal $\bm b_i$ for which $\langle\bm b_i,\bm \alpha_s \rangle > 0$. The ray through $C_i$ must enter a chamber with wall $-\bm b_i $, and so we have that $-\bm b_i$ is a wall of $C_{i+1}$ for $i<r$. 

The largest value of $t$ for which the point $\bm \gamma - t\bm \alpha_s$ is in $C_i$ is computed as follows:
\begin{equation}
\begin{array}{rcl}
\displaystyle \langle \bm \gamma -t\bm \alpha_s , \bm b_i\rangle &=&0\\
t&=&\displaystyle \frac{\langle \bm \gamma , \bm b_i\rangle}{\langle \bm \alpha_s, \bm b_i \rangle}
\end{array}
\end{equation}
and so the first nested sum in \eqref{eq:partialFractionRefinement1} runs from $0$ to 

\begin{equation}\label{eq:pointTj}
t_1 =\left\lfloor \frac{\langle \bm \gamma , \bm b_1\rangle}{\langle \bm \alpha_s, \bm b_1 \rangle} \right\rfloor .
\end{equation}
We can reorder the sum \eqref{eq:partialFractionRefinement1} to:
\begin{equation}\label{eq:partialFractionRefinement2}
P_{\Delta }(\bm \gamma) = P_{\Delta}\left(\bm \gamma - (t_1+1)\bm \alpha_s \right)+ \sum_{t=0}^{t_1} P_{\Gamma}(\bm \gamma -t\bm \alpha_s )\quad .
\end{equation}
We need to show how to convert both of these summands to a linear combination of expressions of the form 
\[
\iota_{ C}(\bm \gamma) \iota_{\bm\delta_j+\Lambda}(\bm \gamma) Q_{i,j}(\bm \gamma) .
\]

To achieve this, we need to eliminate the floor function from \eqref{eq:pointTj} in favor of lattice summation operations. We will do this using two lattice summation tricks which we state in a more general form.

\subsubsection{Elimination of the floor function \texorpdfstring{from  $\iota_{\bm \delta+\Lambda}(\bm \gamma)f(\left\lfloor\left\langle \bm a, \bm\gamma\right\rangle+c\right\rfloor)$}{}}\label{section:eliminateFloorFromLinearFunction}
Let $\bm a$ be an non-zero rational-coordinate vector and let $\Lambda\subset \mathbb Z^n$ be a lattice of full rank. Let $\Omega$ be the set:

\[
\Omega = \Omega(\bm a, \Lambda) = \left\{\bm \gamma \in \Lambda | \langle \bm a,\bm \gamma \rangle \in \mathbb Z \right\}
\]
$\Omega$ is computed by Algorithm \ref{algorithm:computeSubLatticeIntegralScalarProducts}. The quotient-group $\mathbb Z^n/\Omega$ is finite and let $\bm \mu_1, \dots, \bm \mu_N$ be representative elements of this quotient. Fix a constant $c\in \mathbb Q$ and set 
\[
\begin{array}{rcl}
d_i&=& \langle\bm a, \bm \mu_i\rangle - \left\lfloor\langle\bm a, \bm \mu_i\rangle+c\right\rfloor 
\end{array}
\] 
This quantity does not change if we replace $\bm \mu_i$ with $\bm \gamma$ for $\bm \gamma \in \bm \mu_i+\Omega$ and so we have that:
\[
\left\lfloor\langle \bm a, \bm \gamma \rangle+c \right\rfloor =\langle \bm a, \bm \gamma \rangle-d_i\qquad \text{ for all } \bm \gamma \in \bm\mu_i +\Omega
\]
Let $\bm \delta\in \mathbb Z^n$ be an arbitrary vector.
For an arbitrary function $f:\mathbb C\to\mathbb C$ and the fixed constant $c\in \mathbb Q$, we can compute:
\begin{equation}\label{eq:eliminateFloor}
\renewcommand{\arraystretch}{3}
\begin{array}{rcl}
\displaystyle \iota_{\bm \delta+\Lambda}(\bm\gamma) f(\lfloor \langle \bm a, \bm \gamma\rangle + c \rfloor) &\stackrel{\text{use }\eqref{eq:latticeToSublattice}}{=}& \displaystyle \sum_{\substack{i=1\\ \bm \mu_i-\bm \delta\in \Lambda}}^N \iota_{\bm \mu_i + \Omega}(\bm \gamma ) f\left(\left\lfloor \left\langle \bm a, \bm \gamma \right \rangle + c \right\rfloor\right)\\
&=&\displaystyle \sum_{\substack{i=1\\ \bm \mu_i-\bm \delta\in \Lambda}}^N \iota_{\bm \mu_i + \Omega }(\bm\gamma) f\left( \langle \bm a, \bm \gamma \rangle -d_i\right) \\
\end{array}
\end{equation}

\subsubsection{Elimination of the floor function \texorpdfstring{from $Q(\bm \gamma -\left\lfloor \left\langle \bm b, \bm\gamma \right \rangle +c\right\rfloor \bm \alpha )$}{}}\label{section:eliminateFloorFromQuasipolynomial}
In the present section, we explain how to remove the floor function coming from \eqref{eq:pointTj}. More precisely, we explain how to remove the floor function from an arbitrary expression of the form $Q\left( \bm \gamma -\left\lfloor \left\langle \bm b, \bm\gamma \right\rangle+c\right\rfloor \bm \alpha \right)$, where $Q$ is a quasipolynomial, $\bm b\in \mathbb Q^n$, and $\bm \alpha\in \mathbb Z^n$ are constant vectors and $c\in \mathbb Q$ is a constant. Intuitively this is clear: by Section \ref{section:eliminateFloorFromLinearFunction}, the coordinate functions of the vector-valued argument $\bm \gamma - \left \lfloor \left\langle \bm b, \bm\gamma \right\rangle +c\right\rfloor \bm \alpha$ of $Q$ are quasipolynomial with respect to $\bm \gamma$, and since quasipolynomials form an algebra, $Q\left( \bm \gamma - \left \lfloor \left\langle \bm b, \bm\gamma \right\rangle+c\right\rfloor \bm \alpha \right )$ must be a quasipolynomial with respect to $\bm \gamma$ as well.

However, let us work out the computational details in order to provide a reference for computer implementations.  Let $Q$ be given by 
\begin{equation}\label{eq:toBeSubstitutedWithFloors}
Q(\bm \gamma)=\sum_{j=1}^M\iota_{\bm \delta_j +\Lambda}(\bm \gamma) R_{j}(\bm \gamma)
\end{equation}
for some polynomials $R_{j}$, some vectors $\bm \delta_j\in \mathbb Z^n$ and a full rank lattice $\Lambda$.
As in the previous section, let $\Omega $ be the lattice 
\[
\Omega = \Omega(\bm b, \Lambda)= \left\{\bm \gamma\in \Lambda| \left\langle \bm b, \bm \gamma\right\rangle \in \mathbb Z \right\}.
\]
As in the previous section, let $\bm \mu_1, \dots, \bm \mu_N$ be representatives of the quotient-lattice $\mathbb Z^n/\Omega$. Let $d_i\in \mathbb Q$ be the numbers given by
\[
\begin{array}{rcl}
d_i&=& \left\langle\bm b, \bm \mu_i \right\rangle-\left\lfloor \left\langle\bm b, \bm \mu_i \right\rangle+c\right\rfloor \\
\left \lfloor \left\langle\bm b, \bm \mu_i \right\rangle+c\right\rfloor &=&\left\langle\bm b, \bm \mu_i \right\rangle-d_i
\end{array}
\]
We can compute:
\begin{equation}\label{eq:substitutionWithFloorfunction1}
\begin{array}{rcl}
\displaystyle Q(\bm \gamma -\left\lfloor \left\langle \bm b, \bm \gamma \right\rangle+c \right\rfloor \bm \alpha)&=&\displaystyle  \sum_{ i=1}^N\iota_{\bm\mu_i + \Omega}(\bm \gamma)Q(\bm \gamma - \left( \left\langle \bm b, \bm \gamma \right\rangle - d_i\right) \bm \alpha)
\\
&=&\displaystyle \sum_{i=1}^N \iota_{\bm\mu_i+\Omega}(\bm \gamma) \sum_{j=1}^M \iota_{\bm \delta_j +\Lambda}\left( \bm \gamma - \left( \left\langle \bm b, \bm \gamma \right\rangle - d_i\right) \bm \alpha \right) \\
&&\displaystyle  ~~~~~~~~~~~~~~~~~~~~~~ R_{j}\left( \bm \gamma - \left( \left\langle \bm b, \bm \gamma \right\rangle - d_i\right) \bm \alpha \right)
\end{array}
\end{equation}
Let $\Psi$ be the lattice
\[ 
\Psi = \Psi(\bm b, \bm \alpha, \Lambda) =\left\{\bm \gamma\in ~ \Lambda \left| \left \langle \bm b, \bm \gamma \right \rangle \bm \alpha \in \Lambda \right. \right\} 
\]
Let $\bm\nu_1,\dots, \bm \nu_P$ be representatives of the quotient-lattice $\mathbb Z^n/\Psi $. Then 
\[
\begin{array}{rcl}
\displaystyle 
\iota_{\bm \delta_j +\Lambda}\left( \bm \gamma+d_i\bm \alpha- \langle\bm b, \bm \gamma\rangle \bm \alpha\right)
&=& \displaystyle \sum_{k=1}^P \iota_{\bm \nu_k+\Psi}(\bm \gamma) \iota_{\bm \delta_j + \Lambda } \left( \bm \gamma - \left( \left \langle \bm b, \bm \gamma \right\rangle - d_i\right) \bm \alpha \right)\\
&=&\displaystyle \sum_{k=1}^P\iota_{\bm \nu_k+\Psi}(\bm \gamma) \iota_{\bm \delta_j+\left( \langle\bm b, \bm \nu_k \rangle-d_i\right) \bm \alpha+\Lambda}\left( \bm \gamma\right)
\end{array}
\]
and, combining with \eqref{eq:substitutionWithFloorfunction1}, we get that
\begin{equation}\label{eq:substitutionWithFloorfunctionFinal}
\begin{array}{@{}r@{~}c@{~}l}
Q(\bm \gamma -\left\lfloor \left\langle \bm b, \bm \gamma \right\rangle\bm \right\rfloor \bm \alpha) & =& \displaystyle \sum_{ i = 1}^N \sum_{j=1}^M \sum_{k=1 }^P \iota_{\bm\mu_i + \Omega}(\bm \gamma) \iota_{ \bm \nu_k+\Psi}(\bm \gamma)\iota_{\bm \delta_j+\left( \langle\bm b, \bm \nu_k \rangle-d_i\right) \bm \alpha+\Lambda}\left( \bm \gamma\right)\\
&&\displaystyle ~~~~~~~~~~~~~~~ R_{j}\left( \bm \gamma+\left(d_i- \left\langle \bm b, \bm \gamma \right \rangle \right)\bm \alpha\right)\\
\end{array}
\end{equation}
The product of indicator functions can be reduced to a single indicator function using \eqref{eq:productOfIndicatorsOfLatticeShifts}. This shows how to transform \eqref{eq:substitutionWithFloorfunction1} to a quasi-polynomial.

\subsubsection{Induction step}
\begin{figure}[h!]
\psset{xunit = 0.03cm, yunit = 0.03cm}
\begin{pspicture}(-100,-100)(100,100)
\psline[linecolor=black, linestyle=dashed](0,0)(-40.82480000000001,70.7107)
\rput[rb](-40.82480000000001,70.7107){(1, 0, 0)}
\pscircle*(-40.82480000000001,70.7107){0.12}
\psline[linecolor=black, linestyle=dashed](0,0)(81.6497,0)
\rput[l](81.6497,0){(0, 1, 0)}
\pscircle*(81.6497,0){0.12}
\psline[linecolor=black, linestyle=dashed](0,0)(-40.82480000000001,-70.71069999999997)
\rput[rt](-40.82480000000001,-70.71069999999997){(0, 0, 1)}
\pscircle*(-40.82480000000001,-70.71069999999997){0.12}
\psline(20.412450000000007,-35.35534999999999)(81.6497,0)
\psline(20.412450000000007,-35.35534999999999)(20.412450000000007,0)
\psline(81.6497,0)(20.412450000000007,0)
\rput(40.82489524406,-11.785119023690015){\color{red} $2$}
\psline(20.412450000000007,35.35534999999999)(81.6497,0)
\psline(20.412450000000007,35.35534999999999)(20.412450000000007,0)
\psline(81.6497,0)(20.412450000000007,0)
\rput(40.82489524406,11.785119023690015){\color{red} $1$}
\psline(20.412450000000007,35.35534999999999)(0.00003333329999577472,0)
\psline(20.412450000000007,35.35534999999999)(20.412450000000007,0)
\psline(0.00003333329999577472,0)(20.412450000000007,0)
\rput(13.608297502800013,11.785140236899991){\color{red} $3$}
\psline(0.00003333329999577472,0)(20.412450000000007,35.35534999999999)
\psline(0.00003333329999577472,0)(-10.206175000000002,17.677674999999994)
\psline(20.412450000000007,35.35534999999999)(-10.206175000000002,17.677674999999994)
\rput(3.4021299942999974,17.677674999999994){\color{red} $5$}
\psline(20.412450000000007,-35.35534999999999)(0.00003333329999577472,0)
\psline(20.412450000000007,-35.35534999999999)(20.412450000000007,0)
\psline(0.00003333329999577472,0)(20.412450000000007,0)
\rput(13.608297502799985,-11.785140236899991){\color{red} $4$}
\psline(-40.82480000000001,70.7107)(20.412450000000007,35.35534999999999)
\psline(-40.82480000000001,70.7107)(-10.206175000000002,17.677674999999994)
\psline(20.412450000000007,35.35534999999999)(-10.206175000000002,17.677674999999994)
\rput[rb](-10.206187247439999,41.24793426059){\color{red} $7$}
\psline[linecolor=black](-40.82480000000001,70.7107)(-40.82480000000001,0)
\psline[linecolor=black](-40.82480000000001,70.7107)(-10.206175000000002,17.677674999999994)
\psline[linecolor=black](-40.82480000000001,0)(-10.206175000000002,17.677674999999994)
\rput[rb](-30.618607996589986,29.462815236899985){$12$}
\psline[linecolor=black](0.00003333329999577472,0)(-10.206175000000002,-17.677674999999994)
\psline[linecolor=black](0.00003333329999577472,0)(20.412450000000007,-35.35534999999999)
\psline[linecolor=black](20.412450000000007,-35.35534999999999)(-10.206175000000002,-17.677674999999994)
\rput(3.4021299942999974,-17.677674999999994){\color{red} $6$}
\psline[linecolor=black](-40.82480000000001,-70.71069999999997)(20.412450000000007,-35.35534999999999)
\psline[linecolor=black](-40.82480000000001,-70.71069999999997)(-10.206175000000002,-17.677674999999994)
\psline[linecolor=black](20.412450000000007,-35.35534999999999)(-10.206175000000002,-17.677674999999994)
\rput[rt](-10.206187247439999,-41.24793426059){\color{red} $8$}
\psline[linecolor=black](-40.82480000000001,-70.71069999999997)(-40.82480000000001,0)
\psline[linecolor=black](-40.82480000000001,-70.71069999999997)(-10.206175000000002,-17.677674999999994)
\psline[linecolor=black](-40.82480000000001,0)(-10.206175000000002,-17.677674999999994)
\rput[rt](-30.618607996589986,-29.462815236899985){$10$}
\psline[linecolor=black](0.00003333329999577472,0)(-10.206175000000002,-17.677674999999994)
\psline[linecolor=black](0.00003333329999577472,0)(-10.206175000000002,17.677674999999994)
\psline[linecolor=black](-10.206175000000002,17.677674999999994)(-40.82480000000001,0)
\psline[linecolor=black](-40.82480000000001,0)(-10.206175000000002,-17.677674999999994)
\rput[r](-15.309279166700009,0){\color{blue} $11$}
\pscircle*(-40.82480000000001,70.7107){0.09}
\pscircle*(81.6497,0){0.09}
\pscircle*(-40.82480000000001,-70.71069999999997){0.09}
\pscircle*(20.412450000000007,35.35534999999999){0.09}
\pscircle*(20.412450000000007,-35.35534999999999){0.09}
\pscircle*(-40.82480000000001,0){0.09}
\end{pspicture}
\psset{xunit = 0.03cm, yunit = 0.03cm}
\begin{pspicture}(-100,-100)(100,100)
\psline[linecolor=black, linestyle=dashed](0,0)(-40.82480000000001,70.7107)
\rput[rb](-40.82480000000001,70.7107){(1, 0, 0)}
\pscircle*(-40.82480000000001,70.7107){0.12}
\psline[linecolor=black, linestyle=dashed](0,0)(81.6497,0)
\rput[l](81.6497,0){(0, 1, 0)}
\pscircle*(81.6497,0){0.12}
\psline[linecolor=black, linestyle=dashed](0,0)(-40.82480000000001,-70.71069999999997)
\rput[rt](-40.82480000000001,-70.71069999999997){(0, 0, 1)}
\pscircle*(-40.82480000000001,-70.71069999999997){0.12}
\psline(20.412450000000007,-35.35534999999999)(81.6497,0)
\psline(20.412450000000007,-35.35534999999999)(20.412450000000007,0)
\psline(81.6497,0)(20.412450000000007,0)
\rput(40.82489524406,-11.785119023690015){\color{red} $2$}
\psline(20.412450000000007,35.35534999999999)(81.6497,0)
\psline(20.412450000000007,35.35534999999999)(20.412450000000007,0)
\psline(81.6497,0)(20.412450000000007,0)
\rput(40.82489524406,11.785119023690015){\color{red} $1$}
\psline(20.412450000000007,35.35534999999999)(0.00003333329999577472,0)
\psline(20.412450000000007,35.35534999999999)(20.412450000000007,0)
\psline(0.00003333329999577472,0)(20.412450000000007,0)
\rput(13.608297502800013,11.785140236899991){\color{red} $3$}
\psline(0.00003333329999577472,0)(20.412450000000007,35.35534999999999)
\psline(0.00003333329999577472,0)(-10.206175000000002,17.677674999999994)
\psline(20.412450000000007,35.35534999999999)(-10.206175000000002,17.677674999999994)
\rput(3.4021299942999974,17.677674999999994){\color{red} $5$}
\psline(20.412450000000007,-35.35534999999999)(0.00003333329999577472,0)
\psline(20.412450000000007,-35.35534999999999)(20.412450000000007,0)
\psline(0.00003333329999577472,0)(20.412450000000007,0)
\rput(13.608297502799985,-11.785140236899991){\color{red} $4$}
\psline(-40.82480000000001,70.7107)(20.412450000000007,35.35534999999999)
\psline(-40.82480000000001,70.7107)(-10.206175000000002,17.677674999999994)
\psline(20.412450000000007,35.35534999999999)(-10.206175000000002,17.677674999999994)
\rput[rb](-10.206187247439999,41.24793426059){\color{red} $7$}
\psline(0.00003333329999577472,0)(-10.206175000000002,-17.677674999999994)
\psline(0.00003333329999577472,0)(20.412450000000007,-35.35534999999999)
\psline(20.412450000000007,-35.35534999999999)(-10.206175000000002,-17.677674999999994)
\rput(3.4021299942999974,-17.677674999999994){\color{red} $6$}
\psline(-40.82480000000001,70.7107)(-40.82480000000001,0)
\psline(-40.82480000000001,70.7107)(-10.206175000000002,17.677674999999994)
\psline(-40.82480000000001,0)(-10.206175000000002,17.677674999999994)
\rput[rb](-30.618607996589986,29.462815236899985){$12$}
\psline(-40.82480000000001,-70.71069999999997)(20.412450000000007,-35.35534999999999)
\psline(-40.82480000000001,-70.71069999999997)(-10.206175000000002,-17.677674999999994)
\psline(20.412450000000007,-35.35534999999999)(-10.206175000000002,-17.677674999999994)
\rput[rt](-10.206187247439999,-41.24793426059){\color{red} $8$}
\psline[linecolor=black](-40.82480000000001,-70.71069999999997)(-40.82480000000001,0)
\psline[linecolor=black](-40.82480000000001,-70.71069999999997)(-10.206175000000002,-17.677674999999994)
\psline[linecolor=black](-40.82480000000001,0)(-10.206175000000002,-17.677674999999994)
\rput[rt](-30.618607996589986,-29.462815236899985){$10$}
\psline[linecolor=black](-10.206175000000002,-17.677674999999994)(0.00003333329999577472,0)
\psline[linecolor=black](-10.206175000000002,-17.677674999999994)(-40.82480000000001,0)
\psline[linecolor=black](0.00003333329999577472,0)(-40.82480000000001,0)
\rput(-17.010327497199995,-5.892534763100002){\color{blue} $11''$}
\psline[linecolor=black](-10.206175000000002,17.677674999999994)(0.00003333329999577472,0)
\psline[linecolor=black](-10.206175000000002,17.677674999999994)(-40.82480000000001,0)
\psline[linewidth=3pt, linecolor=orange, arrows=>](0.00003333329999577472,0)(-40.82480000000001,0)
\rput(-17.010327497199995,5.892534763100002){\color{blue} $11'$}
\pscircle*(-40.82480000000001,70.7107){0.09}
\pscircle*(81.6497,0){0.09}
\pscircle*(-40.82480000000001,-70.71069999999997){0.09}
\pscircle*(20.412450000000007,35.35534999999999){0.09}
\pscircle*(20.412450000000007,-35.35534999999999){0.09}
\pscircle*(-40.82480000000001,0){0.09}
\end{pspicture}
\caption{Computing the combinatorial chambers for $(1,0,0)$, $(0,1,0)$, $(0,0,1)$, $(1,1,0)$, $(0,1,1)$, $(1,0,1)$. The graph is in barycentric coordinates. $P_{\Delta}(\bm \gamma)$ is computed over the chambers with red label, and we're computing it over the chamber with blue label. The arrow shows the vector $-(1,0,1)$. It exits the blue chamber through a vertex that belongs to Chambers 3, 4, 5 and 6.}\label{figure:showMultipleChamberExit}	
\end{figure}
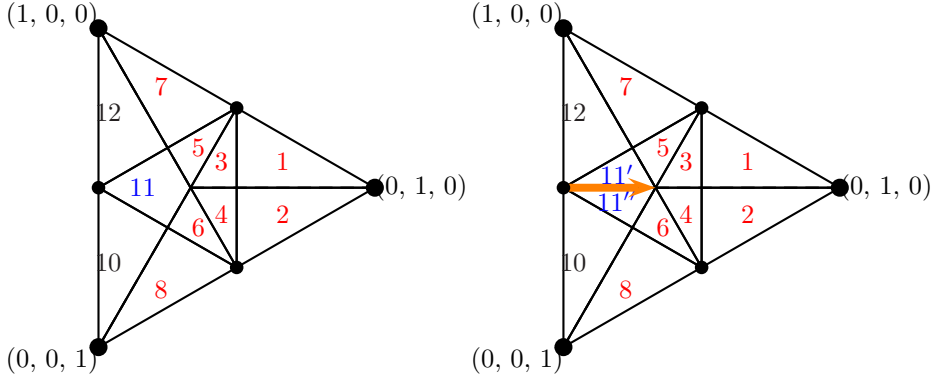

Suppose by induction we have computed $P_{\Delta}(\bm \gamma)$ for $\bm\alpha\in C_2, \dots, C_r$, i.e., suppose we have computed a lattice $\Lambda$, a set of vectors $\bm\delta_1, \dots, \bm \delta _m$ and polynomials $Q_{i, j}$ such that

\begin{equation}\label{eq:partialFractionRefinement3}
\begin{array}{rcll}
Q_i(\bm \gamma)&=&\displaystyle \sum_{j=1}^m \iota_{\bm \delta_j + \Lambda}(\bm \gamma) Q_{i, j}(\bm \gamma) \\
P_{\Delta}(\bm \gamma) &=&\displaystyle Q_i(\bm \gamma)
& \text{ for }i\geq 2 \text{ and all } \bm \gamma \in C_i,
\end{array}
\end{equation}
If we could guarantee that $\bm \gamma+ (t_1+1)\bm \alpha_s \in C_2$, this would compute the first summand in \eqref{eq:partialFractionRefinement2} as $Q_2(\bm \gamma)$. However, if the ray $\bm \gamma + t\bm \alpha_s$ intersects $C_2$ closely to an $n-2$-dimensional edge, it can happen that $t_2 = t_1$ and so $\bm \gamma+ (t_1+1)\bm \alpha_s\notin C_2$ (for example, see Figure \ref{figure:showMultipleChamberExit}).

This case can be taken care of by expressing $P_\Delta$ as an alternating sum: 

\begin{equation}\label{eq:vectorPartitionFunctionAlternatingSumComponent}
\begin{array}{rcl}
P_\Delta \left(\bm \gamma -(t_1+1)\bm \alpha_s\right)&=&\phantom{+} Q_2(\bm \gamma -(t_{1 \phantom{-1}}+1)\bm \alpha_s) - Q_2(\bm \gamma -(t_2+1)\bm \alpha_s)  \\
&&+Q_3(\bm \gamma -(t_{2\phantom{-1}}+1)\bm \alpha_s) - Q_3(\bm \gamma -(t_3+1)\bm \alpha_s) \\
&&+\dots \\
&&+Q_r(\bm \gamma -(t_{r-1}+1)\bm \alpha_s)- Q_r(\bm \gamma -(t_{r}+1)\bm \alpha_s).
\end{array}
\end{equation}
The Bernoulli sum $B_k(t)$ vanishes at $t=-1$ and at $t=0$, so the last summand $-Q_r\left(\bm \gamma -(t_{r}+1)\bm \alpha_s\right)$ - which originates from a Bernoulli sum - is zero, but we added it for symmetry.
When we substitute the values of $t_i$ from \eqref{eq:pointTj}, the above summands are of the form
\[
Q_i\left( \bm \gamma -\left\lfloor\left\langle\bm \gamma ,\frac{\bm b_i}{\langle\bm \alpha_s,\bm b_i \rangle}\right\rangle +1 \right \rfloor \bm \alpha_s \right)
\] 
which can be reduced to quasipolynomial using \eqref{eq:substitutionWithFloorfunctionFinal}.

Suppose we have computed the vector partition function $P_{\Delta}(\bm \gamma)$ over the chamber $D$.
In this equality, the first summand counts all vector partitions that have the vector $\bm \alpha_s$ participate at least $t_1+1$ times.

In our induction hypothesis, we supposed that either $P_{\Gamma}$ is zero over the interior of $C_1$, or we have a quasipolynomial formula for $P_\Gamma$. 

\subsubsection{The case of zero $P_\Gamma$ over the interior of $C_1$}\label{section:zeroPGamma}
Suppose first that $P_\Gamma$ is zero over the interior of $C_1$. Then all summands in $\displaystyle \sum_{t=0}^{t_1} P_{\Gamma}(\bm \gamma+t\bm \alpha_s)$ except possibly the last are zero. So, the second summand of \eqref{eq:partialFractionRefinement2} becomes 
\begin{equation}\label{eq:zeroCaseForPGamma}
P_{\Gamma}\left( \bm \gamma -\left\lfloor \frac{\left\langle\bm \gamma, \bm b_1\right\rangle}{\left\langle\bm \alpha_s, \bm b_1 \right \rangle} \right\rfloor\bm \alpha_s \right)
\end{equation}
and can be non-zero only when the argument belongs to the wall separating $C_1$ and $C_2$, in other words when $\frac{\left\langle\bm \gamma, \bm b_1\right\rangle}{ \left\langle\bm \alpha_s, \bm b_1 \right \rangle} \in \mathbb Z$. Let $\Theta$ be the lattice 
\[ 
\Theta =\left\{\bm \gamma\in \mathbb Z^n \left| \frac{\left\langle\bm \gamma, \bm b_1\right\rangle}{\left\langle\bm \alpha_s, \bm b_1 \right \rangle}\in \mathbb Z\right. \right\}.
\]
Then \eqref{eq:zeroCaseForPGamma} is computed using \eqref{eq:substitutionWithFloorfunctionFinal} in the formula for $P_{\Gamma}$ over $C_2$, and multiplying that by $\iota_{\Theta}(\bm \gamma)$.

\subsubsection{The case of non-zero $P_\Gamma$ over the interior of $C_1$}\label{section:nonZeroPGamma}

Suppose alternatively that over $C_1$, we have a quasi-polynomial formula for $P_{\Gamma}$ given by:

\begin{equation}\label{eq:sumPGammaStartingForm}
P_{\Gamma}(\bm \gamma) = \sum_{i=1}^N\iota_{\bm \delta_i+\Lambda} (\bm \gamma ) q_i\left(\bm \gamma \right)
\end{equation}
where the $q_i$'s are polynomials in $n$ variables, $\Lambda\subset \mathbb Z^n$ is a lattice and $\bm \delta_1, \dots,\bm \delta_N$ are vectors representative of elements of the finite quotient group $\mathbb Z^n / \Lambda$. It follows that

\[
P_\Gamma\left(\bm\gamma-t\bm\alpha_s\right) = \sum_{i=1}^N \iota_{\bm \delta_i+\Lambda} (\bm \gamma -t\bm \alpha_s ) q_i\left(\bm \gamma -t\bm \alpha_s\right).
\]
Let $a$ be the smallest positive integer for which 
\begin{equation}\label{eq:aScale}
a\bm \alpha_s\in \Lambda.
\end{equation}
It follows that  
\[
\iota_{\bm\delta_i+\Lambda} \left( \bm \gamma +a\bm\alpha_s \right) = \iota_{\bm\delta_i + \Lambda} \left( \bm \gamma\right), 
\]
for every vector $\bm \gamma$ and so we can sum:
\begin{equation}\label{eq:sumPGammaGeneric}
\begin{array}{@{}r@{~~}c@{~~}l}
\displaystyle \sum_{t=0}^{t=t_1} P_\Gamma\left(\bm\gamma-t\bm\alpha_s\right)& =& \displaystyle \sum_{ t=0}^{t=t_{1}} \sum_{i=1}^N \iota_{ \bm \delta_i+\Lambda} (\bm \gamma -t\bm \alpha_s ) q_i\left(\bm \gamma -t\bm \alpha_s\right)\\
&&\text{substitute } t=au+w
\\
&=&\displaystyle \sum_{w=0}^{a-1} \sum_{u=0 }^{\left\lfloor \frac{t_1-w}{a}\right\rfloor}  \sum_{i=1}^N\iota_{\bm \delta_i+\Lambda} ( \bm \gamma -(au+w)\bm \alpha_s ) q_i\left(\bm \gamma -(au+w)\bm \alpha_s\right)\\~\\
&=&\displaystyle \sum_{w=0}^{a-1} \sum_{u=0 }^{\left\lfloor \frac{t_{1}-w}{a}\right\rfloor}  \sum_{i=1}^N\iota_{\bm \delta_i+w\bm \alpha_s +\Lambda} (\bm \gamma )q_i\left(\bm \gamma -(au+w)\bm \alpha_s\right)\\
&=&\displaystyle \sum_{w=0}^{a-1}\sum_{i=1}^{N}\iota_{\bm \delta_i+w\bm \alpha_s+\Lambda}(\bm \gamma ) \sum_{u=0 }^{\left\lfloor \frac{t_1 -w}{a}\right\rfloor} q_i \left(\bm \gamma -(au+w)\bm \alpha_s\right)
\end{array}
\end{equation}
The inner-most sum above 
\[
\displaystyle \sum_{u=0 }^{\left\lfloor \frac{t_1 -w}{a}\right\rfloor} q_i \left(\bm \gamma -(au+w)\bm \alpha_s\right)
\] is polynomial in $u$ and therefore a linear combination of the Bernoulli sums of the form:
\begin{equation}
B_k\left( \left\lfloor \frac{t_1-w}{a}\right\rfloor\right)
\end{equation}
with coefficients that are products of polynomials in  $\bm \gamma$. Here, $k$ is the degree of the Bernoulli sum equal to $s-n$, i.e., the number of partitioning vectors less the dimension of the space. Substitute \eqref{eq:pointTj} in the above to get a linear combination of elements of the form: 
\begin{equation}\label{eq:bernoulliSumSubstituted}
\renewcommand{\arraystretch}{3}
\begin{array}{r@{~}c@{~}l}
\displaystyle  B_k\left( \left\lfloor \frac{ \left \lfloor \frac{\langle \bm \gamma , \bm b_1\rangle}{\langle \bm \alpha_s,  \bm b_1 \rangle} \right\rfloor -w}{a}\right\rfloor\right) &=& \displaystyle  B_k\left( \left\lfloor \frac{ \left \lfloor \frac{\langle \bm \gamma , \bm b_1\rangle}{\langle \bm \alpha_s,  \bm b_1 \rangle} -w \right \rfloor }{ a} \right \rfloor \right)\\
&=&\displaystyle  B_k\left( \left\lfloor  \frac{\langle \bm \gamma , \bm b_1\rangle}{a\langle \bm \alpha_s,  \bm b_1 \rangle} -\frac{ w }{a}\right\rfloor\right)
\end{array}
\end{equation}
with coefficients of the form $\iota_{\bm \delta_i + w\bm \alpha_s +\Lambda}(\bm \gamma)$  where 
\[
a>0,0 \leq w< a
\]
and the simplification follow from the fact that $a,w$ are integers. The Bernoulli sum $B_k$ is a polynomial. We can now apply \eqref{eq:eliminateFloor} with $f=B_k$,  $\displaystyle \bm a=\frac{\bm b_1}{a \left \langle \bm \alpha_s, \bm b_1 \right\rangle}$, $c=-\frac{w}{a}$ and $\bm \delta = \bm \delta_i +w \bm \alpha_s$ to eliminate the floor functions in 
\[
\displaystyle \iota_{\bm \delta_i + w\bm \alpha_s +\Lambda}(\bm \gamma)B_k\left( \left\lfloor  \frac{\langle \bm \gamma , \bm b_1 \rangle }{a\langle \bm \alpha_s,  \bm b_1 \rangle} -\frac{w }{a}\right\rfloor\right),
\] and thus obtain the quasipolynomials promised in Theorem \ref{theorem:mainTheory}(a).

\subsubsection{Example}\label{section:whyWeNeedCompression}
Consider the vector partition function for $(2,2), (1,0),(0,1)$.
	
	\psset{xunit = 0.01cm, yunit = 0.01cm}
	\begin{pspicture}(-100,-100)(100,100)
		
		\psline[linecolor=orange, linewidth=2pt, arrows=->](52, 104)(52,52)
		\pscircle*[linecolor=blue](52,104){0.07}
		\psline[linecolor=red, linewidth=2pt, arrows=->](52,52)(52,0)
		
		\psline[linecolor=orange, linewidth=2pt, arrows=->](104, 208)(104,104)
		\pscircle*[linecolor=blue](104,208){0.07}
		\psline[linecolor=red, linewidth=2pt, arrows=->](104,104)(104,52)

		\psline[linecolor=black, linestyle=dashed](0,0)(500,0)
		\psline[linecolor=black, linestyle=dashed](0,0)(0,200)
		\psline[linecolor=black, linestyle=dashed](0,0)(52,0)
		\psline[linecolor=black, linestyle=dashed](0,0)(0,52)
		\pscircle*(0,0){0.09}
		\pscircle*(52.000000000000014,0){0.09}
		\pscircle*(52.000000000000014,0){0.09}
		\pscircle*(104.00000000000001,0){0.09}
		\pscircle*(104.00000000000001,0){0.09}
		\pscircle*(104.00000000000001,104){0.09}
		\pscircle*(156,0){0.09}
		\pscircle*(156,0){0.09}
		\pscircle*(156,104){0.09}
		\pscircle*(208,0){0.09}
		\pscircle*(208,0){0.09}
		\pscircle*(208,104){0.09}
		\pscircle*(260,0){0.09}
		\pscircle*(260,0){0.09}
		\pscircle*(208,208){0.09}
		\pscircle*(260,104){0.09}
		\pscircle*(312,0){0.09}
		\pscircle*(312,0){0.09}
		\pscircle*(260,208){0.09}
		\pscircle*(312,104){0.09}
		\pscircle*(364,0){0.09}
		\pscircle*(364,0){0.09}
		\pscircle*(312,208){0.09}
		\pscircle*(364,104){0.09}
		\pscircle*(416,0){0.09}
		\pscircle*(416,0){0.09}
		\pscircle*(156,0){0.09}
		\pscircle*(52.000000000000014,0){0.09}
		\pscircle*(52,52){0.04}
		\pscircle*(104,52){0.04}
		\pscircle*(104,52){0.04}
		\pscircle*(156,52){0.04}
		\pscircle*(156,156){0.04}
		\pscircle*(208,52){0.04}
		\pscircle*(208,52){0.04}
		\pscircle*(208,156){0.04}
		\pscircle*(260,52){0.04}
		\pscircle*(260,52){0.04}
		\pscircle*(260,156){0.04}
		\pscircle*(312,52){0.04}
		\pscircle*(312,52){0.04}
		\pscircle*(260,260){0.04}
		\pscircle*(312,156){0.04}
		\pscircle*(364,52){0.04}
		\pscircle*(364,52){0.04}
		\pscircle*(312,260){0.04}
		\pscircle*(364,156){0.04}
		\pscircle*(416,52){0.04}
		\pscircle*(416,52){0.04}
		\psline[linecolor=blue, linestyle=dashed](0,0)(500,0)
		\psline[linecolor=blue, linestyle=dashed](0,0)(200,200)
		\pscircle*(52.000000000000014,0){0.09}
	\end{pspicture}

We now illustrate the summation techniques of the preceding subsections. In this example, we have deliberately chosen to account for vector $(2,2)$ before vectors $(1,0)$ and $(0,1)$, so that the last vector $(0,1)$ is outside of the convex hull of the previous ones. In dimension $2$, cones are generated by at most two vectors so this convex hull issue could be avoided by putting the two generator vectors first - in the case of the present example, those would be vectors $(1,0)$ and $(0,1)$. However, in dimension $3$, the issue can no longer be avoided: pick up any infinite cone with four sides and any of the three vectors generating the cone will not contain the fourth. In view of the fact that the issue is unavoidable in dimension $3$ and above, we keep our unusual vector order in dimension $2$ in the name of a simpler visualization.

Let us now show how to compute $P_\Delta$ according to Sections \ref{section:zeroPGamma} and \ref{section:nonZeroPGamma}. Let $C_1 =\left\{(x_1, x_2)|x_1\geq x_2, \geq 0\right\}$ and $C_2 =\left\{(x_1, x_2)|x_2\geq x_1, \geq 0\right\}$ be the two chambers drawn in the figure. We immediately see that
\begin{equation}\label{eq:example1Quasipolynomial1}
P_{\{(2,2),(1,0)\} }(x_1, x_2) = 
\begin{cases}
1 & \left (x_1,x_2 \right) \in \Lambda((1,0),(0,2)) \text{ and } x_2\leq x_1\\
0&\text{otherwise}
\end{cases}
\end{equation}
We are computing $P_{\Delta}$ over chamber $C_1$. The  $a$ from \eqref{eq:aScale} equals $2$, the vector $\bm \alpha_s$ is $(0,1)$ and $\Lambda \mapsto \Lambda((1,0),(0,2))$. The representatives $\bm \delta_1, \bm \delta_2$ of $\mathbb Z^2 /\Lambda $ are $\bm \delta_1\mapsto(0,0), \bm \delta_2\mapsto (0,1)$. The quasipolynomials $q_1, q_2$ from \eqref{eq:sumPGammaStartingForm} are $q_1\mapsto 1$ and $q_2\mapsto 0$. The normal of the exit wall with respect to $(0,1)$ is $\bm b_1$ is $(0,1)$. The argument of \eqref{eq:bernoulliSumSubstituted} becomes: $\left\lfloor  \frac{\langle \bm (x_1,x_2) , \bm b_1 \rangle }{a\langle \bm \alpha_s,  \bm b_1 \rangle} -\frac{w }{a}\right\rfloor\mapsto \left\lfloor \frac{x_2}{2}-\frac{w}{2}\right\rfloor$.

Therefore for $(x_1,x_2)\in C_1$, \eqref{eq:sumPGammaGeneric} becomes:  
\[
\begin{array}{@{}r@{~}c@{~}l}
P_\Delta(x_1, x_2)_{|C_1} &=&\displaystyle \sum_{w=0}^{w=1}\sum_{i=1}^{i=2}\iota_{\bm \delta_i+w\bm (0,1)+\Lambda}(x_1, x_2 ) \sum_{u=0 }^{\left\lfloor \frac{x_2}{2}-\frac{w}{2}\right\rfloor} q_i \left((x_1,x_2)-(au+w)\bm (0,1)\right)\\
&=& \displaystyle \sum_{w=0}^{w=1} \iota_{(0,w)+\Lambda}(x_1, x_2 ) \sum_{u=0 }^{\left\lfloor \frac{x_2}{2}-\frac{w}{2}\right\rfloor} 1\\
&=&\displaystyle \sum_{w=0}^{w=1} \iota_{(0,w)+\Lambda}(x_1, x_2 ) \cdot (\underbrace{ t+1}_{\text{Bernoulli sum } B_1(t)})_{|t=\left\lfloor \frac{x_2}{2}-\frac{w}{2}\right\rfloor}\\
&=&\displaystyle \sum_{w=0}^{w=1} \iota_{(0,w)+\Lambda}(x_1, x_2 ) \left(\left\lfloor \frac{x_2}{2}-\frac{w}{2}\right\rfloor+1\right)\\
&=&\begin{cases}
\frac{x_2}{2}+1 & \text{ for } (x_1, x_2) \in \Lambda \\
\frac{x_2}{2}+\frac{1}{2} & \text{ for } (x_1, x_2) \in (0,1)+\Lambda
\end{cases}
\end{array}
\]

Now let us compute $P_\Delta$ over chamber $C_2$. We fall in the case of Section \ref{section:zeroPGamma}. By Section \ref{section:zeroPGamma}, the first summand in equation \eqref{eq:partialFractionRefinement2} is obtained by computing $Q_1((x_1,x_2)-\lfloor\langle(-1, 1), (x_1,x_2)\rangle +0\rfloor(0, 1))$, where $Q_1$ is the quasipolynomial of \eqref{eq:example1Quasipolynomial1}. By Section \ref{section:eliminateFloorFromQuasipolynomial}, this summand can be computed to be:

$S_1 = \begin{cases}1 & \text{over }(0,0), (0,1), (0,2), (0,3)+\Lambda((2,0), (4,0))\\ 0&\text{otherwise}\end{cases}$

Similarly, the non-zero summand $S_2$ from the alternating sum in \eqref{eq:vectorPartitionFunctionAlternatingSumComponent} can be computed to be $P_{\Delta} ((x_1, x_2)-\lfloor\langle(-1, 1), (x_1, x_2) \rangle +1\rfloor(0, 1))$ which we computed above. Again using the machinery for eliminating the floor function from \ref{section:eliminateFloorFromQuasipolynomial} we get that 

\[
S_2(x_1, x_2) =\begin{cases}
\frac{x_1}{2}+\frac{1}{2} &\text{for }(x_1,x_2)\in (1,0), (1,1), (1,2), (1,3)+\Lambda((2,0), (4,0))\\
\frac{x_1}{2} &\text{for }(x_1,x_2)\in (0,0), (0,1), (0,2), (0,3)+\Lambda((2,0), (4,0))
\end{cases} 
\]
Summing the quantities we get that for $(x_1, x_2)\in C_2$:
\[
P_{\Delta}(x_1,x_2)_{C_2} =\begin{cases}
	\frac{x_1}{2}+1 &\text{for }(x_1,x_2)\in (0,0), (0,1), (0,2), (0,3)+\Lambda((2,0), (4,0))\\
	\frac{x_1}{2}+\frac{1}{2} &\text{for }(x_1,x_2)\in (1,0), (1,1), (1,2), (1,3)+\Lambda((2,0), (4,0))
\end{cases} 
\]
As the careful reader would immediately see, the formula above can be compressed to:
\[
P_{\Delta}(x_1,x_2)_{|C_2} =\begin{cases}
	\frac{x_1}{2}+1 &\text{for }(x_1,x_2)\in \Lambda((2,0), (0,1))\\
	\frac{x_1}{2}+\frac{1}{2} &\text{for }(x_1,x_2)\in (1,0)+ \Lambda((2,0), (0,1)),
\end{cases} 
\]
which can be also seen immediately by the formula over chamber $C_1$ and by the symmetry of the variables $x_1, x_2$.

Now, let us compute $P_\Delta(x_1, x_2)$ using our main algorithm. A Brion-Vergne decomposition is:
\begin{align*}&~~~
\frac{1}{(1-x_{1} ) (1-x_{2} ) (1-x_{1}^2x_{2}^2) }\\=&
\left(-x_{1} x_{2}^{-1}-x_{1} x_{2}^{-2}-x_{2}^{-1}-x_{2}^{-2}\right)\cdot \left(\frac{x_{1} \partial_{1} }{2}-\frac{x_{2} \partial_{2} }{2}\right)\cdot\frac{1}{(1-x_{1}^2) (1-x_{1}^2x_{2}^2) }\\&
\displaystyle +\left(x_{1} x_{2}^{-2}+x_{2}^{-2}\right)\cdot \left(\frac{x_{1} \partial_{1} }{2}\right)\cdot\frac{1}{(1-x_{1}^2) (1-x_{2} ) }
\end{align*}

We omit the computations over $C_2$. Over the interior of $C_1$, by \eqref{eq:vectorPartitionFunctionOneFraction}, the first fraction contributes to the vector partition function:
\[
\begin{array}{rcl}
	P_1 &=&
\begin{cases}
\underbrace{- \frac{1}{2} \left( (x_1-1)-(x_2+1)\right)}_{ \text{contrib. by. } -x_1x_2^{-1}}& \text{for }(x_1,x_2)\in (1,1)+ \Lambda((2,0), (0,2))\\
\underbrace{- \frac{1}{2} \left( (x_1-1)-(x_2+2)\right)}_{ \text{contrib. by. } -x_1x_2^{-2}}& \text{for }(x_1,x_2)\in (1,0)+ \Lambda((2,0), (0,2))\\
\underbrace{- \frac{1}{2} \left( x_1-(x_2+1)\right)}_{ \text{contrib. by. } -x_2^{-1}}& \text{for }(x_1,x_2)\in (0,1)+ \Lambda((2,0), (0,2))\\
\underbrace{- \frac{1}{2} \left( x_1-(x_2+2)\right)}_{ \text{contrib. by. } -x_2^{-2}}& \text{for }(x_1,x_2)\in  \Lambda((2,0), (0,2))\\
\end{cases}
\end{array}
\]
Similarly, the second fraction contributes to the vector partition function:

\[
\begin{array}{rcl}
	P_2 &=&
	\begin{cases}
		\underbrace{ \frac{1}{2} \left( x_1-1\right)}_{ \text{contrib. by. } x_1x_2^{-2}}& \text{for }(x_1,x_2)\in (1,0)+ \Lambda((2,0), (0,1))\\
		\underbrace{ \frac{1}{2} x_1}_{ \text{contrib. by. } x_2^{-2}}& \text{for }(x_1,x_2)\in   \Lambda((2,0), (0,1))
	\end{cases}
\end{array}
\]
and summing the two we get the formula 
\[
P_{\Delta}(x_1, x_2)_{|C_1}=\begin{cases}
\frac{x_2}{2}+1& \text{for } (x_1,x_2)\in (1,0),(0,0)+ \Lambda((2,0), (0,2))\\
\frac{x_2}{2}+\frac{1}{2}& \text{for } (x_1, x_2)\in (0,1),(1,1)+ \Lambda((2,0), (0,2))
\end{cases}
\]
Again, the final formula for $P_{\Delta}(x_1,x_2) _{|C_1}$ is not optimal, but can be compressed. 

\subsection{Compressing the quasipolynomial formulas}
In the example of Section \ref{section:whyWeNeedCompression}, we saw that both the elementary algorithm and the standard algorithm for computing vector partition functions yield sub-optimal formulas using lattices that are too rough. 

These formulas can be compressed in a straightforward fashion. Indeed, suppose that the polynomial formulas over two different lattice shifts are the same. If the difference of the lattice shifts is $\bm \delta$, we need to check whether the polynomial formulas are the same for every two sets of polynomials that differ by the same shift. If so, then we can compress our formulas by refining our lattice to include $\bm \delta$ and dropping the redundant formulas (see Section \ref{section:commonRefinement}).

\subsection{Subdividing into normal chambers with respect to a vector}\label{sec:subdivision}
In the previous sections, we demonstrated how to sum the vector partition function $P_{\Delta}(\bm \gamma)$ to the form Theorem \ref{theorem:mainTheory} (a), under the assumption that we can select a collection of chambers for which we have computed $P_{ \Delta \setminus \{\bm \alpha_s\}}$ and such that each chamber is normal with respect to $\bm\alpha_s$.

We are free to subdivide our chambers to smaller pieces, without changing the formula for the vector partition function $P_{\Delta \setminus \bm \alpha_s}$. It is intuitively clear that such subdivision should be sufficient to ensure the normality of the chambers with respect to $\bm \alpha_s$. In the present section, we aim to introduce an algorithm to do this.

We recall from Definition \ref{def:normalChamber} that a chamber $C$ is normal with respect to a vector $\bm \nu$ if it has a unique wall $\bm b$ with $\langle \bm b, \bm \nu\rangle > 0$. We also recall from Definition \ref{def:normallySeparatedCollection} that a collection of chambers is normally separated if every chamber has at most $1$ neighbor with $n-1$ linearly independent common points along each of its walls. Finally, Definition \ref{def:normalChamberCollection} states that a collection of chambers is normal with respect to $\bm v$ if it is normally separated and each of the chambers is normal with respect to $\bm v$.

\begin{theorem}\label{theorem:subdivision}
Let $\bm \nu$ be a vector with positive coordinates and $C_1, \dots, C_k$ be a set of normally separated chambers (Definition \ref{def:normallySeparatedCollection}), such that all chamber points have positive coordinates. Then the collection can be subdivided into chambers so that the final collection is normal with respect to $\bm \nu$.
\end{theorem}

Our proof of Theorem \ref{theorem:subdivision} is constructive. We furthermore report that the resulting algorithm is practical and we have implemented it in our software.

For every chamber $C$, we store the following information. First, we store the normalized walls of each chamber (see Definition \ref{definition:coneWalls}). Next, for the purposes of our computer implementation, we also store the number of chambers created since the start of the computation to serve as a unique identifier of each chamber. The vertices of a chamber can be computed from its walls by Algorithm \ref{algorithm:findTheVertices}. However, we compute and store all vertices of each chamber immediately upon creation. Our first use of chamber vertices is to determine quickly when a given plane splits a given chamber in two pieces of non-zero volume. Namely, we will need the following observations, whose proofs we omit.
\begin{lemma}[Subdivision by plane]\label{lemma:subdivisionByPlane}
Let $C$ be a chamber that is generated by its vertices $\bm \alpha_1, \dots, \bm \alpha_t$ and defined by the normals $\bm a_1, \dots, \bm a_m$. Let $\bm a$ be a non-zero vector. The following hold.
\begin{itemize}
\item[(a)] The plane with normal $\bm a$ splits $C$ into two chambers of non-zero volume if and only if there exist $\bm \alpha_i, \bm \alpha_j$ for which $\langle \bm a, \bm \alpha_i\rangle <0$ and  $\langle \bm a,\bm \alpha_j\rangle >0$.
\item[(b)] The two chambers $C_+\cup C_-$ from (a) are given by:
\[
\begin{array}{rcl}
C_+&=&\left\{ \bm \gamma| \langle \bm \gamma ,\bm a_i \rangle\geq 0 \text{ for all }i \text{ and } \langle  \bm \gamma ,\bm a\rangle \geq 0\right\}\\
C_-&=&\left\{ \bm \gamma| \langle \bm\gamma ,\bm a_i \rangle\geq 0 \text{ for all }i \text{ and } \langle  \bm \gamma ,-\bm a\rangle \geq 0\right\}\\
\end{array}
\]
\end{itemize}
\end{lemma}
The inequalities in (b) are redundant and, when represented by our computer implementation, need to be normalized to satisfy Definition \ref{definition:coneWalls}.

For the following lemma, we recall that vertices were defined in Definition \ref{definition:vertex} and normalized walls were defined in Definition \ref{definition:coneWalls}. The lemma gives a criterion for when two chambers are normally separated.
\begin{lemma} \label{lemma:separatedNormally}
Let $C$ and $D$ be two chambers that are generated by their vertices. Then they are normally separated if one of the two mutually exclusive possibilities hold.

\begin{itemize}
\item There  exist two linearly independent vectors $\bm a, \bm b$ such that $\langle \bm a, \bm \alpha\rangle\geq 0$, $\langle \bm b, \bm \alpha\rangle\geq 0$ for all $\bm\alpha\in C$ and $\langle \bm a, \bm \beta\rangle\leq 0$, $\langle \bm b, \bm \beta \rangle \leq 0$ for all $\bm\beta\in D$.
\item There exists a vector $\bm a$ such that:
\begin{itemize}
\item $\bm a$ is a normalized wall of $C$
\item $-\bm a$ is a normalized wall of $D$
\item all vertices $ \bm \alpha$ of $C$ for which $\langle \bm a,\bm \alpha\rangle = 0$ are also vertices of $D$.
\end{itemize}  
\end{itemize}

\end{lemma}
We omit the proof of this lemma. The next lemma gives a criterion for when two chambers are not normally separated.
\begin{lemma}\label{lemma:notSeparatedNormally}~
\begin{itemize}
\item[(a)] Suppose two chambers $C$, $D$ have non-empty intersection of their interiors but do not coincide. Then there exists a wall of $C$ whose plane splits $D$ into two non-zero volume chambers, or there exists a wall of $D$ whose plane splits $C$ into two non-zero volume chambers.
\item[(b)] Suppose two chambers $C,D$ are generated by their vertices, have disjoint interiors, but are not normally separated. Then the following hold.
\begin{itemize}
\item[(1)] There exists a normalized wall with normal vector $\bm a$ of $C$ such that $-\bm a$ is a the normal vector of a normalized wall of $D$.
\item[(2)] At least one of the following two holds.
\begin{itemize}
\item There exists a wall of $C$ with normal vector $\bm v\neq \bm a$ so that the plane with normal vector $\bm v$ splits $D$ into two non-zero volume chambers. There exists a vertex $\bm \alpha$ of $C$ that lies in the walls with normals $ \bm a$ and $\bm v$.
\item There exists a wall of $D$ with normal vector $\bm v\neq-\bm a$ so that the plane with normal vector $\bm v$ splits $C$ into two non-zero volume chambers. There exists a vertex $\bm \alpha$ of $D$ that lies in the walls with normals $-\bm a$ and $\bm v$.
\end{itemize}
\end{itemize}
\end{itemize}
\end{lemma}
\begin{proof} We omit the proof of (a) and only prove (b). The fact that two chambers are not normally separated but have disjoint interiors implies that they have a common supporting plane, whose appropriately rescaled normal is the vector $\bm a$ from (b.1).
	
Let $P$ be the plane with normal vector $\bm a$. To prove (b.2), we note that $C\cap P$ and $D\cap P$ are two $n-1$-dimensional cones whose $n-1$-dimensional interiors intersect. By (a) either an $n-2$-dimensional wall of $C\cap P$ splits $D\cap P$ or an $n-2$-dimensional wall of $D\cap P$ splits $C\cap P$. The $n-2$-dimensional walls are the intersection of $P$ and an $n-1$-dimensional wall of either $C$, which establishes the first alternative of (b.2) or $D$ which establishes the second alternative of (b.2). The vertex $\bm \alpha$ exists because the wall with normal $\bm v$ intersects the wall with normal $\pm \bm a$ in a point that belongs to one of the two chambers.
\end{proof}

We now explain how to compute when two chambers are neighbors.

\begin{definition}[Chamber Neighbors]\label{definition:chamberNeighbors}
Given a wall $\bm a$ of $C$ and a set $\mathcal I$ of chambers with $C\in \mathcal I$, let  
\[
N(C, \bm a, \mathcal I) = \left\{D\in \mathcal I \left|
\begin{array}{l} D,C \text{ have }n-1 \text{ common points of rank }n-1 \\
\text{ belonging to the plane with normal } \bm a 
\end{array}\right.
\right\}
\]
We call $N(C, \bm a, \mathcal I)$ the set of neighbors of $C$ in the collection $\mathcal I$ along the wall with normal $\bm a$. 

We say that  $M(C, \bm a, \mathcal I)$ is an upper bound for the neighbors if $N(C, \bm a, \mathcal I)\subset M(C, \bm a, \mathcal I)$ for every $C\in \mathcal I$ and every $\bm a$.
\end{definition}

For each wall $\bm a_i$ of $C$, we store the set of neighbors $N(C, \bm a_i, \mathcal I)$ by recording the identifiers of all neighbors of $C$.

Fix a collection of chambers $\mathcal I$ with points with non-negative coordinates, a distinguished chamber $D\in \mathcal I$, and a subdivision of $D=D_1\cup\dots\cup D_p$ into chambers with disjoint interiors. 

\begin{definition}[Chamber replacement]\label{definition:chamberReplacement}
We say that \\ $\mathcal I'= \mathcal I \setminus \left\{D\right\}\cup \left\{ D_1,\dots, D_p \right\}$ is the replacement of $D$ by $D_1, \dots, D_p$ in $I$. 
\end{definition}

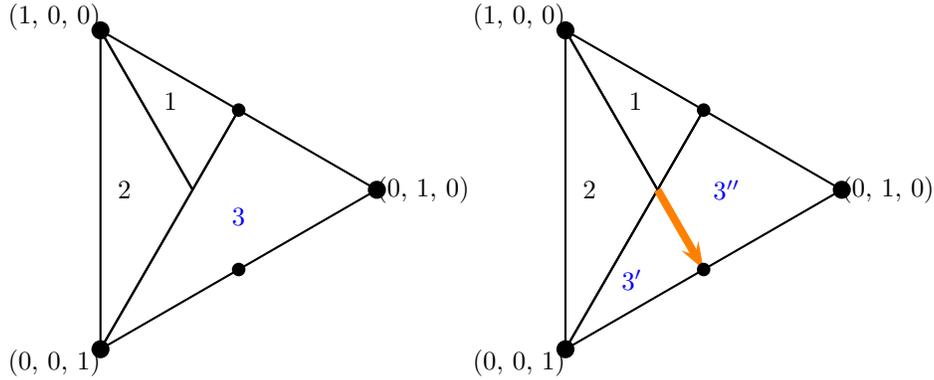
\begin{figure}[h!]

	\psset{xunit = 0.03cm, yunit = 0.03cm}
	\begin{pspicture}(-100,-100)(100,100)
		\rput[rb](-40.82480000000001,70.7107){(1, 0, 0)}
		\pscircle*(-40.82480000000001,70.7107){0.12}
		\rput[l](81.6497,0){(0, 1, 0)}
		\pscircle*(81.6497,0){0.12}
		\rput[rt](-40.82480000000001,-70.71069999999997){(0, 0, 1)}
		\pscircle*(-40.82480000000001,-70.71069999999997){0.12}
		\psline[linecolor=black](20.412450000000007,35.35534999999999)(-40.82480000000001,70.7107)
		\psline[linecolor=black](20.412450000000007,35.35534999999999)(0.00003333329999577472,0)
		\psline[linecolor=black](-40.82480000000001,70.7107)(0.00003333329999577472,0)
		\rput[rb](-6.804078339,35.35534999999999){$1$}
		\psline[linecolor=black](-40.82480000000001,-70.71069999999997)(-40.82480000000001,70.7107)
		\psline[linecolor=black](-40.82480000000001,-70.71069999999997)(0.00003333329999577472,0)
		\psline[linecolor=black](-40.82480000000001,70.7107)(0.00003333329999577472,0)
		\rput[r](-27.21649500570001,0){$2$}
		\psline[linecolor=black](81.6497,0)(-40.82480000000001,-70.71069999999997)
		\psline[linecolor=black](81.6497,0)(20.412450000000007,35.35534999999999)
		\psline[linecolor=black](-40.82480000000001,-70.71069999999997)(20.412450000000007,35.35534999999999)
		\rput(20.412450000000007,-11.785069526200004){\color{blue} $3$}
		\pscircle*(-40.82480000000001,70.7107){0.09}
		\pscircle*(81.6497,0){0.09}
		\pscircle*(-40.82480000000001,-70.71069999999997){0.09}
		\pscircle*(20.412450000000007,35.35534999999999){0.09}
		\pscircle*(20.412450000000007,-35.35534999999999){0.09}
	\end{pspicture}
	\psset{xunit = 0.03cm, yunit = 0.03cm}
	\begin{pspicture}(-100,-100)(100,100)
		\rput[rb](-40.82480000000001,70.7107){(1, 0, 0)}
		\pscircle*(-40.82480000000001,70.7107){0.12}
		\rput[l](81.6497,0){(0, 1, 0)}
		\pscircle*(81.6497,0){0.12}
		\rput[rt](-40.82480000000001,-70.71069999999997){(0, 0, 1)}
		\pscircle*(-40.82480000000001,-70.71069999999997){0.12}
		\psline[linecolor=black](20.412450000000007,35.35534999999999)(-40.82480000000001,70.7107)
		\psline[linecolor=black](20.412450000000007,35.35534999999999)(0.00003333329999577472,0)
		\psline[linecolor=black](-40.82480000000001,70.7107)(0.00003333329999577472,0)
		\rput[rb](-6.804078339,35.35534999999999){$1$}
		\psline[linecolor=black](-40.82480000000001,-70.71069999999997)(-40.82480000000001,70.7107)
		\psline[linecolor=black](-40.82480000000001,-70.71069999999997)(0.00003333329999577472,0)
		\psline[linecolor=black](-40.82480000000001,70.7107)(0.00003333329999577472,0)
		\rput[r](-27.21649500570001,0){$2$}
		\psline[linecolor=black](-40.82480000000001,-70.71069999999997)(20.412450000000007,-35.35534999999999)
		\psline[linecolor=black](-40.82480000000001,-70.71069999999997)(0.00003333329999577472,0)
		\psline[linecolor=black](20.412450000000007,-35.35534999999999)(0.00003333329999577472,0)
		\rput[rt](-6.804078339,-35.35534999999999){\color{blue} $3'$}
		\psline[linecolor=black](81.6497,0)(20.412450000000007,-35.35534999999999)
		\psline[linecolor=black](81.6497,0)(20.412450000000007,35.35534999999999)
		\psline[linecolor=black](20.412450000000007,35.35534999999999)(0.00003333329999577472,0)
		\psline[linewidth=3pt, linecolor=black, linecolor=orange, arrows=>](20.412450000000007,-35.35534999999999)(0.00003333329999577472,0)
		\rput(30.618658333300004,0){\color{blue} $3''$}
		\pscircle*(-40.82480000000001,70.7107){0.09}
		\pscircle*(81.6497,0){0.09}
		\pscircle*(-40.82480000000001,-70.71069999999997){0.09}
		\pscircle*(20.412450000000007,35.35534999999999){0.09}
		\pscircle*(20.412450000000007,-35.35534999999999){0.09}
	\end{pspicture}
	\caption{	\label{figure:extraneousNeighbors} Algorithm \ref{algorithm:extendSubdivisionToNormal}, illustrated on the vector partition function of $(1,0,0)$, $(0,1,0)$, $(0,0,1)$, $(1,1,0)$, $(0,1,1)$ in barycentric coordinates. Step 5.1.1.4 splits chamber $3$ into chambers $3'$ and $3''$. Steps 5.1.1.5 - 5.1.1.7 will mark chamber $3', 3''$ as neighbors of both Chambers $1,2$, however the neighborhood relations $3'\rightleftarrows 1$ and  $3''\rightleftarrows 2$ are extraneous. These are extranesous relations are removed in Step 5.1.1.2.}
\end{figure}

Suppose we have already computed all neighbors in a collection of chambers, and we seek to update this neighborhood information for the replacement of a single chamber.
More precisely, suppose we have an upper bound for the map of neighbors $M\left(C, \bm a, \mathcal I\right)$ for all walls $\bm a$ and all chambers $C\in \mathcal I$. Suppose we have also computed $N\left(D_i, \bm a, \left\{D_1, \dots, D_p\right\}\right)$. Then we can compute an upper bound for $M\left(C, \bm a, \mathcal I'\right)$ by the following.
\begin{itemize}
\item For $C\notin \left\{D_1,\dots, D_p\right\}$, an upper bound $M\left(C, \bm a, \mathcal I'\right)$  can be computed by replacing all occurrences of $D$ in $M\left(C, -\bm a, \mathcal I\right)$ by the $D_i$'s for which $\bm a $ is a normalized wall.
\item For $D_i$, an upper bound $M\left(D_i, \bm a, \mathcal I'\right)$ is given by the union \\ $M\left(D_i, \bm a, \left\{D_1, \dots, D_p\right\}\right)\cup M\left(D, \bm a, \mathcal I\right)$.
\end{itemize}
The upper bound $M\left(C, \bm a, \mathcal I'\right)$ computed above is not exact. This can be seen in Figure \ref{figure:extraneousNeighbors}.

In order for our chambers $\mathcal I$ to be normally separated (Definitions \ref{def:normallySeparated}), we need the following.
\begin{itemize}
\item For all chambers to have disjoint interiors (Definition \ref{def:normallySeparated}(a)), all neighbors of every $C$ along each wall $\bm a_i$ must have $-\bm a_i$ as a wall normal.
\item For our chambers to be normally separated by Definition \ref{def:normallySeparated}(b), every $C$ must have one or zero neighbors along each of its walls, i.e., we need to have that $\left|N(C, \bm a, \mathcal I)\right|\leq 1$ for all $\bm a$.
\end{itemize}
Previously, we noted that any given collection of chambers with disjoint interiors can be subdivided into a normally separated collection by computing all chambers in the complement of all hyperplanes of all walls of the chambers. This process has the drawback that if the starting chambers were all normal with respect to a given vector $\bm \nu$, the resulting chambers need not be. 

We now describe an algorithm for subdividing a collection with normal neighbors to a new collection with normal neighbors that respects a fixed subdivision of one chamber. 
\begin{algorithm}[Extend subdivision to normal neighbors] \label{algorithm:extendSubdivisionToNormal}~
\begin{itemize}
\item Input: 
\begin{itemize}
\item A collection of chambers $\mathcal R$ with normal neighbors. The chambers $\mathcal R$ are required to have no points with negative coordinates.
\item A distinguished chamber $D\in I$ and a subdivision $\left\{D_1, \dots, D_p \right\}$ of $D$ into chambers with pairwise disjoint interiors such that the collection $\left\{D_1, \dots, D_p\right\}$ has normal neighbors. 
\item The map $N(C, \bm a, \mathcal R)$ for all $\bm a$ and all $C\in \mathcal R$.
\item The map $N(D_i, \bm a, \mathcal \{D_1, \dots, D_p\})$ for all $\bm a$ and $D_i\in \mathcal \{D_1, \dots, D_p\}$.
\end{itemize}  
\item Output: a subdivision $\mathcal O$ of $\mathcal R\setminus\{D\}\cup \left\{ D_1, \dots, D_p \right\}$ with normally separated chambers and a neighbor map $N(C, \bm a, \mathcal O)$ computed for all $C, \bm a$.
\end{itemize}
\begin{itemize}
\item[Step 1.] Initialize a collection $\mathcal O$ with the elements of $\mathcal R \setminus \{ D \} \cup \left\{D_1, \dots, D_p \right\}$.
\item[Step 2.] Set $M(C,\bm a)$ to be the upper bound for the neighbor map obtained by replacing $D$ with the $D_i$'s in the neighbor map $N(C, \bm a, \mathcal R)$ as described after Definition \ref{definition:chamberReplacement}.
\item[Step 3.] Form a queue of all chambers  $C$ for which $\left|M(C, \bm a)\right| >1$ for some $\bm a$. 
\item[Step 4.] If the queue is empty, skip to Step 6, else remove the first element $C$ of the queue.
\item[Step 5.] For each chamber $E$ in $M(C, \bm a)$:
\begin{itemize}
\item[Step 5.1.] For every vertex $\bm \beta$ of $E$ for which $\langle \bm \beta, \bm a\rangle =0$:
\begin{itemize}
\item[Step 5.1.1.] For every wall with normal vector $\bm b$ of $E$ that passes through $\bm \beta$:
\begin{itemize}
\item[Step 5.1.1.1.] Compute the cone $C_-$ with walls obtained by adding the wall $-\bm b$ to the walls of $C$ (see Lemma \ref{lemma:subdivisionByPlane}).
\item[Step 5.1.1.2.] If $C_-$ equals $C$, remove $E$ from $M(C, \bm a)$ and remove $C$ from $M(E, -\bm a)$. Then, go back to Step 3.
\item[Step 5.1.1.3.] Compute the cone $C_+$ with walls obtained by adding the wall $\bm b$ to the walls of $C$. If plane with normal $\bm b$ does not split $C$ using the criterion of Lemma \ref{lemma:subdivisionByPlane}, go back to Step 5.1.1.
\item[Step 5.1.1.4.] Remove $C$ form $\mathcal O$ and add $C_+, C_-$ to it.
\item[Step 5.1.1.5.] Set $M( C_+, \bm b)= C_-$ and $M(C_-, -\bm b)= C_+$. 
\item[Step 5.1.1.6.] For every normal $\bm a$ of a wall of $C$ that is a wall of $C_+$, set $M(C_+, \bm a) = M(C, \bm a)$. 
\item[Step 5.1.1.7.] For every normal $\bm a$ of a wall of $C$ that is a wall of $C_-$, set $M(C_-, \bm a) = M(C, \bm a)$. 
\item[Step 5.1.1.8.] Go back to Step 3.
\end{itemize}
\item[Step 5.1.2.] If we exhausted all normals $\bm b$ without going back to Step 3, bo back to Step 5.1.
\end{itemize}
\item[Step 5.2] If we exhausted all vertices $\bm\beta  $ without going back to Step 3, go back to Step 5.
\end{itemize}
\item[Step 5.] If we exhausted chambers $E$ without going back to Step 3, remove $C$ from the queue of chambers and go back to Step 4.
\item[Step 6.] Set $N(C,\bm a, \mathcal O) = M(C, \bm a)$ for all $C\in \mathcal O$.
\end{itemize}
\end{algorithm} 
Let us justify the correctness of the algorithm. The starting chambers have disjoint interiors. The only place where new chambers are generated is the subdivision in Step 5.1.1.4, and since the chamber from which the subdivision originated is removed, this guarantees that the collection continues to have chambers with disjoint interiors.
Suppose now at some point the algorithm above generates two chambers $C, D$ that are not normally separated. By Lemma \ref{lemma:notSeparatedNormally}, a wall of $C$ must split $D$ or a wall of $D$ must split $C$, which guarantees that one of the two chambers will be removed by Step 5.1.1.4. Since there are only finitely many chambers, and therefore finitely many supporting planes, this guarantees that the algorithm terminates.

The map $M(c, \bm a)$ is guaranteed to provide an upper bound for the neighbors of $\mathcal O$ as neighbors are conservatively accounted for by Steps 5.1.1.5, 5.1.1.6 and 5.1.1.7. This process may introduce extraneous neighbors in $M(C, \bm a)$, see Figure \ref{figure:extraneousNeighbors} for example. However, by Lemma \ref{lemma:separatedNormally}, these extraneous neighbor elements will be detected and removed in Step 5.1.1.2.

Before we proceed with the proof of Theorem \ref{theorem:subdivision}, we need a few remarks on what it means for one convex body to ``occlude'' another convex body in the direction of $\nu$. We first define the $t_{\bm \gamma}$ to be the distance between a point $\bm \gamma$ to a chamber $C$ in the direction of $-\bm\nu$. More precisely, let $\bm \gamma$ be any point that does not lie on the plane of the wall of $C$. Define $t_{\bm \gamma}: \text{set of open sets in }\mathbb R^n\to \mathbb R\cup \{\textbf{undefined}\}$ as
\[
t_{\bm \gamma}(C) = 
\begin{cases}
\inf \{t \geq 0|\bm\gamma- t \bm \nu\in C\}&\text{ if } \{t \geq 0|\bm\gamma- t \bm \nu\in C\} \neq \emptyset \\
\textbf{undefined} & \text{ if } \{t \geq 0|\bm\gamma- t \bm \nu\in C\} =\emptyset
\end{cases}.
\]
Here, $\textbf{undefined}$ stands for a formal quantity that is non-comparable to a number.

\begin{definition} We say that the open set $D$ occludes the open set $C$ in the direction of $\bm \nu$ if $t_{\bm\alpha}(C) \leq  t_{\bm \alpha}(D)$ for all $\bm \alpha$ for which both quantities are defined, and there's at least one such $\bm \alpha$. 
\end{definition}
\begin{lemma}[Occlusion] \label{lemma:occlusionIsPartialOrder}
Let $C$ and $D$ be two open convex sets that do not intersect. Let $\bm \nu$ be a non-zero vector. Suppose there's a point $\bm\alpha$ for which both $t_{\bm \alpha}(C)$ and $t_{\bm \alpha}(D)$ are defined. Then either $C$ occludes $D$ or $D$ occludes $C$ in the direction of $\bm\nu$.
\end{lemma}
\begin{proof}
Let $C, D$ be interiors of two convex sets. If $t_{\bm\alpha}(C)= t_{\bm \alpha}(D)$ for all $\bm \alpha$ for which the quantities are defined, there's nothing to prove. So, without loss, suppose that 
\begin{equation} \label{eq:occlusion1}
t_{\bm\alpha}(C) < t_{\bm \alpha}(D)
\end{equation} 
for some $\bm \alpha$. Suppose, contrary to the lemma, there exists $\bm \beta$ such that 
\begin{equation} \label{eq:occlusion2}
t_{\bm\beta}(C) > t_{\bm \alpha}(D).
\end{equation} Let $\varepsilon_1, \varepsilon_2, \varepsilon_3, \varepsilon_4$ be small quantities (to be chosen next) and set
\[
\begin{array}{rcl}
u_1 = t_{\bm \alpha}(C)+\varepsilon_1,& ~~~~&  u_2 = t_{\bm \alpha}(D)+\varepsilon_2\\
v_1 = t_{\bm \beta} (C)+\varepsilon_3,& ~~~~&  v_2 = t_{\bm \beta }(D)+\varepsilon_4
\end{array}
\]
Now choose $\varepsilon_1, \varepsilon_2, \varepsilon_3, \varepsilon_4$ so that
\[
\begin{array}{rcl}
\bm \alpha+u_1 \bm \nu\in C,  &~~~& \bm \alpha +u_2\bm \nu \in D\\
\bm \beta+ v_1 \bm \nu \in C, && \bm \beta+v_2\bm \nu \in D
\end{array}
\] 
and so that 
\begin{equation}
u_1<u_2,\qquad v_1 > v_2.
\end{equation} 
This choice is possible because of the infimum in the definition of $t_{\bm \alpha}$ and the fact that both inequalities \eqref{eq:occlusion1} and \eqref{eq:occlusion2} are strict. Set 
\[
a=u_2-u_1 > 0, \qquad b=v_1-v_2 > 0.
\]

The four points above form a 2d-trapezoid (embedded in $n$-dimensional space). Let $\bm \gamma$ be the intersection of the diagonals of the trapezoid:

\[
\begin{array}{rcl}
\bm \gamma&=& \displaystyle\frac{b}{a+b} \left(\bm \alpha + u_2\bm \nu\right) + \frac{a}{a+b} \left(\bm \beta + v_2 \bm \nu\right) \\
&=& \displaystyle \frac{b}{a+b} \left( \bm \alpha + (a+u_1)\bm \nu\right) + \frac{a}{a+b} \left(\bm \beta + (v_1-b) \bm \nu\right) \\
&=& \displaystyle\frac{b}{a+b}\left (\bm \alpha + u_1\bm \nu\right) + \frac{a}{a+b} \left(\bm \beta + v_1 \bm \nu \right) \\
\end{array}
\]
Since both $a,b$ are positive, in the first line, we have a positive linear combination with coefficients that sum to $1$ of two points in a convex set. In other words, the first line shows that $\bm\gamma\in D$. Likewise, the last line shows that $\bm\gamma\in C$. However, $C$ and $D$ have disjoint interiors which is a contradiction. This proves the lemma.
\end{proof}

The following lemma states that a chamber $D$ with points with positive coordinates can be subdivided into finitely many chambers that are normal with respect to $\bm \nu$. The statement of this lemma is intuitively clear, however we present a proof as we will use the explicit subdivision formulas in an upcoming algorithm statement.
\begin{lemma} [Directional Subdivision]\label{lemma:directionalSubdivision}
Let $\nu$ be a non-zero vector with positive coordinates and let $D$ be a chamber consisting of points with positive coordinates that has $p$ walls with normals $\bm b_1, \dots, \bm b_p$ for which $\langle \bm\nu,\bm b_i \rangle > 0 $. 

Then there exist $p$ chambers $D_1, \dots, D_p$ with disjoint interiors such that $\bigcup_{i} D_i=D$ and such that each $D_i$ is normal with respect to $\bm \nu$.
\end{lemma}
\begin{proof}
Denote by $B_i$ the plane given by the normal vector $\bm b_i$. For each point $\bm \gamma$ in $D$, $\langle\bm \gamma, \bm b_i\rangle \geq 0$, and so $\bm \gamma -t\bm \nu$ must intersect $B_i$  at $t=\frac{\langle\bm \gamma, \bm b_i \rangle}{\langle \bm \nu ,\bm b_i  \rangle } $. The points for which $\bm \gamma -t\bm \nu$ intersects $B_i $ before it intersects $B_j$ are given by the inequality 
\begin{equation}\label{eq:directionalSubdivisionWalls}
\left\langle \bm \gamma , \frac{\bm b_j}{\left\langle\bm\nu, \bm b_j \right\rangle }- \frac{\bm b_i}{\left\langle\bm\nu, \bm b_i \right \rangle }\right \rangle\geq 0
\end{equation}
Set \begin{equation}\label{eq:directionalSubdivisionNewWalls}
\bm c_{i,j}=\frac{\bm b_j}{\left\langle\bm\nu, \bm b_j \right\rangle }- \frac{\bm b_i}{\left\langle\bm\nu, \bm b_i \right\rangle }	
\end{equation}
 from the inequality above. Let $\bm a_1, \dots, \bm a_q$ be the remaining walls of $D$. Let $D_{i}$ be defined by the inequalities: 
\begin{equation}\label{eq:directionalSubdivisionChambers}
\begin{array}{rcl}
\langle\bm \gamma, \bm b_i\rangle&\geq& 0 \\
\langle\bm \gamma, \bm a_1\rangle&\geq& 0 \\
&\vdots \\
\langle\bm \gamma, \bm a_q\rangle&\geq& 0 \\
\langle\bm \gamma, \bm c_{i,1}\rangle&\geq& 0 \\
& \vdots \\
\langle\bm \gamma, \bm c_{i,p}\rangle&\geq& 0 
\end{array}.
\end{equation}
Every two chambers $D_i, D_j$ have pairwise disjoint interiors as $D_i$ contains the inequality $\langle\bm \gamma, \bm c_{i,j} \rangle \geq 0$ and $D_j$ - its opposite. 
Since the ray $\bm \gamma+t\bm \nu$ can only exit through the walls with normals $\bm b_1, \dots, \bm b_p$ it follows that $D_1\cup \dots \cup D_p=D$. 
\end{proof}

\begin{proof} [~of Theorem \ref{theorem:subdivision}]
Our proof is algorithmic - see \ref{algorithm:subdivisionByArbitrarySlicesOnce} for a detailed pseudocode formulation. Introduce the collections $\mathcal N_0= \left\{C_1, \dots, C_k\right\}$, $\mathcal \mathcal R_0 =\emptyset$. Here, $\mathcal N$ stands for ``non-refined'', and $\mathcal R$ stands for ``refined''.
The sets $\mathcal N_0, \mathcal R_0$ serve as the starting $0^{th}$ step of our algorithm. 
In the $j+1^{st}$ step, the sets $\mathcal N_{j+1}, \mathcal I_{j+1}$ will be constructed from $ \mathcal N_{j}, \mathcal R_{j}$. At any point of the algorithm, 
\begin{itemize}
\item the chambers in the union of the two collections $\mathcal N_j\cup \mathcal R_j$ have pairwise disjoint interiors
\item the union of the chambers in the two collections is equal to the union of the starting chambers , i.e., 
\[
\bigcup_{D\in \mathcal N_j \cup \mathcal R_j }D=\bigcup_i C_i 
\]
\item the collection of refined chamber $\mathcal R_j$ is normal with respect to $\bm \nu$. 
\end{itemize}

By the Occlusion Lemma \ref{lemma:occlusionIsPartialOrder}, there exists a chamber $D\in \mathcal N_{j}$ whose interior is not occluded by the interior of any chamber in $\mathcal N_j $. By the Directional Subdivision Lemma \ref{lemma:directionalSubdivision} subdivide $D$ into $p$ chambers $D_1,\dots, D_k$ with disjoint interiors so each $D_i$ is normal with respect to $\bm\nu$. Set 
\[
\mathcal N'_{j+1} = \mathcal N_j \setminus D
\]
and \[\mathcal R'_{j+1} = \mathcal R_j\cup \left\{D_1, \dots, D_k\right\}.\] By Algorithm \ref{algorithm:extendSubdivisionToNormal}, subdivide the collection $\mathcal N'_{j+1}\cup \mathcal R'_{j+1}$ so it is a normally separated collection. Put all resulting chambers $C$ that were subdivided from elements of $\mathcal N'_{j+1}$ into $\mathcal N_{j+1}$. Put the remaining chambers in $\mathcal R_{j+1}$. 

Algorithm \ref{algorithm:extendSubdivisionToNormal} splits cones only using the new supporting hyperplanes introduced by the Directional Subdivision Lemma \ref{lemma:directionalSubdivision}. The new supporting hyperplanes have normal vectors given by \eqref{eq:directionalSubdivisionWalls} and so all of them contain $\bm \nu$. In this way both the lemma and the algorithm only introduce walls whose scalar products with $\bm\nu$ is zero. This means that none of these walls participate in generating new walls through Lemma \ref{lemma:directionalSubdivision}. In this way, the only possible normal vectors of walls that can appear in this algorithm are the vectors obtained by \eqref{eq:directionalSubdivisionWalls} as the pairs $\bm b_i, \bm b_j$ run over all walls of all of the original chambers $\left\{C_1, \dots, C_k\right\}$ with negative scalar product with $\bm \nu$. There are finitely many such vectors, and therefore only finitely many new supporting hyperplanes introduced by the algorithm. The hyperplane arrangement of all such hyperplanes, together with the supporting hyperplanes of the original chambers splits space into finitely many indivisible chambers. The union of the original set $\displaystyle \bigcup_{C\in \mathcal N_0}C$ contains finitely many of those. Every time we compute $\mathcal N'_{j+1}$ from $\mathcal N_j$, the union $\displaystyle \bigcup _{C\in\mathcal N'_{j+1}} C $ contains fewer indivisible chambers. Since $\mathcal N_{j+1}$ is a subdivision of $\mathcal N'_{j+1}$, the two collections have the same union. This shows that the number of indivisible chambers in $\mathcal N_{j}$ decreases strictly as $j$ increases. Therefore the algorithm must terminate. The final $\mathcal R_j$ in this algorithm is the normal collection whose existence we asserted.
\end{proof}

We implemented the proof above on computer. We state the pseudo-code of our implementation for reference.
\begin{algorithm}[Subdivision with respect to direction] ~\label{algorithm:subdivisionByArbitrarySlicesOnce}
\begin{itemize}
\item Input.
\begin{itemize}
\item A vector $\bm \nu$.
\item A collection of chambers $\mathcal N$ with disjoint interiors.
\item The map of neighbors $N(C, \bm a, \mathcal N)$, computed for all chambers $C$ and all vectors $\bm a$. We recall from Definition \ref{definition:chamberNeighbors} that $N(C, \bm a, \mathcal N)$ is the set of all neighbors of $C$ along $\bm a$ in the collection $\mathcal N$.
\end{itemize} 
\item Output.
\begin{itemize}
\item A collection of chambers $\mathcal R$ with disjoint interiors that is a subdivision of $\mathcal N$, such that $\mathcal R$ is normal in the direction of $\bm \nu$ (Definition \ref{def:normalChamberCollection}).
\item The map of neighbors $N(C, \bm a, \mathcal R)$, computed for all chambers $C\in \mathcal R$ and all vectors $\bm a$. 
\end{itemize} 
\end{itemize}
\begin{itemize}
\item[Step 1.] Store all elements of $\mathcal N $ in a queue (first in, first out). We continue to refer to that queue by the letter $\mathcal N$.
\item[Step 2.] Initialize the set of chambers $\mathcal R$ to the empty set.
\item[Step 3.] Initialize the map of neighbors $N(C, \bm a)$ using $N(C, \bm a, \mathcal N)$.
\item[Step 4.] If $\mathcal N$ is empty, stop the program. The output chamber collection is contained in $\mathcal R$, and the output neighbor map is $N(C, \bm a)=N(C, \bm a, \mathcal R)$.
\item[Step 5.] Pop the first element $C$ in the queue $\mathcal N$. 
\item[Step 6.] If $C$ has a wall $\bm b$ with $\langle \bm b, \bm \nu\rangle>0$ for which one of the neighbors $N(C, \bm b)$ is not in $\mathcal R$, push $C$ to the back of the queue $\mathcal N$ and go back to Step 4.
\item[Step 7.] Else, if $C$ has a unique wall $\bm b$ with $\langle \bm b, \bm \nu\rangle>0$, add $C$ to the collection $\mathcal R$ and go back to Step 4.
\item[Step 8.] Else, if $C$ has multiple walls $\bm b$ with $\langle \bm b, \bm \nu\rangle>0$: \begin{itemize}
\item[Step 8.1.] Subdivide $C$ using the Directional Subdivision Lemma \ref{lemma:directionalSubdivision} as follows.
\begin{itemize}
	\item[Step 8.1.1.] Let the walls of $C$ that have positive scalar product with $\bm \nu$ be $\bm b_1, \dots, \bm b_p$.
	
	\item[Step 8.1.2.] For each index $i$:
	\begin{itemize}
\item[Step 8.1.2.1.] Let $D_i$ be the cone whose walls are given by \eqref{eq:directionalSubdivisionChambers}.
		\item[Step 8.1.2.2.] Compute the vertices of $D_i$ using Algorithm \ref{algorithm:findTheVertices}.
		\item[Step 8.1.2.3.] Compute the normalized walls of $D_i$ using Algorithm \ref{algorithm:normalsFromGenerators}.
		\item[Step 8.1.2.4.] For every $j$, if the wall with normal $\bm c_{i,j}$ defined in \eqref{eq:directionalSubdivisionChambers} remains a wall after the normalization of $D_i$, set $D_i$ and $D_j$ as neighbors, i.e., set $N(D_i, \bm c_{i,j})=\{D_j\}$, $N(D_j, \bm c_{i,j})=\{D_i\}$.
	\end{itemize}
	\item[Step 8.1.3.] Remove $D$ from $\mathcal N$ and push $\{D_1, \dots, D_p\}$ to the back of the queue $\mathcal N$.
	\item[Step 8.1.4.] Use Algorithm \ref{algorithm:extendSubdivisionToNormal} to subdivide the collection $\mathcal R\cup \mathcal N$ into a collection with normal neighbors.
	\begin{itemize}
		\item Place all chambers that were in $\mathcal R$ and are not subdivided by Algorithm \ref{algorithm:extendSubdivisionToNormal} back in $\mathcal R$. 
		\item Place all chambers that were in $\mathcal N$ and are not subdivided by Algorithm \ref{algorithm:extendSubdivisionToNormal} back in $\mathcal N$. 
		\item Update the neighbors map $N(C, \bm a)$ to the output neighbor map of Algorithm \ref{algorithm:extendSubdivisionToNormal}.
	\end{itemize}
\end{itemize}
\end{itemize}
\item[Step 9.] Go back to Step 4.
\end{itemize}
\end{algorithm}

\subsection{Proof of Theorem \ref{theorem:mainTheory}}
We are now ready to prove our main theorem. As described in Section \ref{sec:normalChambers},  we proceed by induction on the number of direction vectors $\bm\alpha_k$, $k\leq s$. The vector partition function does not change when we shuffle the $\bm\alpha_i$'s, so reorder the vectors so that the first $n$ vectors $\bm\alpha_1, \dots, \bm\alpha_n$ are linearly independent.

We start our induction at $k=n$. For our starting setup, take the cone generated by all vectors 
\[
X=C(\bm \alpha_1, \dots, \bm \alpha_s).
\] 
Subdivide it with the planes spanned by the $n-1$ element subsets of the first $n$ vectors $\bm \alpha_1, \dots, \bm \alpha_n$ using the Subdivision by Plane Lemma \ref{lemma:subdivisionByPlane}. The cone generated by the first $n$ vectors $C(\bm \alpha_1, \dots, \bm\alpha_n)$ is in the subdivision. Set 
\[
\begin{array}{rcl}
P_{\left\{\bm \alpha_1, \dots, \bm \alpha_n\right\}}(\bm \gamma) &=& \begin{cases}
1&\text{ when }\bm \gamma \in \Lambda(\bm\alpha_1, \dots, \bm \alpha_n)\cap C(\bm\alpha_1, \dots, \bm \alpha_n)\\
0&\text{ otherwise}
\end{cases}\\
&=&\iota_{C(\bm \alpha_1,\dots, \bm \alpha_n)}(\bm \gamma)\iota_{\Lambda(\bm \alpha_1,\dots, \bm \alpha_n)}(\bm \gamma)
\end{array}
\]
Set $P_{\left\{\bm \alpha_1, \dots, \bm \alpha_n\right\} }$ to be zero outside of $C(\bm \alpha_1, \dots, \bm \alpha_n)$. In case $n=s$, it is immediate that the formula we produced in this fashion gives quasipolynomial formula promised by Theorem \ref{theorem:mainTheory}.

Suppose now we have proven Theorem \ref{theorem:mainTheory} for $k-1$ vectors. Suppose the induction hypothesis gives us quasipolynomial formulas over $r$ chambers $C_1, \dots, C_r$. As additional induction hypotheses, we suppose that:

\begin{itemize}
\item $\displaystyle \bigcup_{i=1}^r C_i= X$, i.e., the cones cover the full cone.
\item The $C_i$'s are a set with normally separated neighbors (Definition \ref{def:normallySeparatedCollection}), that is, their interiors are disjoint. 
\item For each $i$, we have a quasipolynomial $P_i$ such that one of the alternatives hold:
\begin{itemize}
\item $P_i(\bm \gamma)\neq 0$ and $P_i= P_{\bm\alpha_1, \dots, \bm\alpha_{k-1}}(\bm \gamma)$ over $C_i$.
\item $P_i(\bm \gamma)=0$.
\end{itemize}  
\end{itemize}
In other words, when $P_i =0$, we allow it to not reproduce the values of $P_{\bm\alpha_1, \dots, \bm\alpha_{k-1}}$ over the relevant chamber. That is to be expected. Indeed, let $D$ be a cone from the base of our induction hypothesis that is a neighbor of $C(\bm \alpha_1, \dots, \bm \alpha_n)$. Then the vector partition function takes on the zero value over the interior of $D$, but takes on the value of $1$ over the intersection of one of the walls of $D$ with a lattice. 

We are now ready to prove the induction step. By Theorem \ref{theorem:subdivision}, we can subdivide the collection $\left\{C_1, \dots, C_r\right\}$ into a collection $\mathcal I = \left\{D_1, \dots, D_R\right\}$ that is normal with respect to $\bm\alpha_k$. 

Let $\mathcal K\subset \mathcal I$ be the set of cones over which we have computed the quasipolynomial formulas promised by our induction step, and let $\mathcal J=\mathcal I\setminus \mathcal K$ be the set over which we have not computed the formulas. By the Occlusion Lemma \ref{lemma:occlusionIsPartialOrder}, there exists a cone $D$ in $\mathcal J$ that is not occluded in the direction of $\bm\alpha_k$ inside $\mathcal J$. Since $D$ is not occluded, for every $\bm \alpha$ in the interior of $D$, the ray $\{\bm \gamma -t\bm \alpha_k|t\geq 0\}$ intersects interiors of cones lying in $\mathcal K$, i.e., cones for which we have computed the vector partition function as a quasipolynomial. Therefore we can compute $P_{ \bm\alpha_1, \dots, \bm\alpha_{k-1} }$ as a quasipolynomial using Bernoulli sums and the lattice computations described in Section \ref{section:summingQuasipolynomials}. Then we can move $D$ from the collection $\mathcal J$ into the collection $\mathcal K$. In this way, we can compute the quasipolynomial vector partition function for every chamber in $I$. Since Theorem \ref{theorem:subdivision} produces normally separated neighbors, the conditions of our induction are satisfied by this construction. This proves Theorem \ref{theorem:mainTheory}(a).

We are now proving Theorem \ref{theorem:mainTheory}(b). By Algorithm \ref{algorithm:main}, there exist a reduced partial fraction relative to $\Delta$ that gives us quasipolynomial formulas that are valid over translates of pointed polyhedral cones such that each of their walls are parallel to hyperplanes spanned by $n-1$ vectors from the original set $\Delta$. On the other hand, Section \ref{section:summingQuasipolynomials} gives us pointed quasipolynomial formulas that are valid over cones whose walls pass through the origin. 

If two polynomials have values that coincide over a translate of a pointed polyhedral cone of non-zero volume intersected with a lattice, then the two polynomials coincide over the intersection of a lattice with a box of arbitrary size. This can only happen if the polynomials are identical. 

Therefore the quasipolynomial formulas from Algorithm \ref{algorithm:main} and Section \ref{section:summingQuasipolynomials} are equal whenever the regions over which they are defined have intersection of non-zero volume. In other words, the chambers produced by one algorithm ``glue'' the chambers produced by the other and the other way round, which completes the proof of our theorem.

\subsection{Computing chambers of quasipolynomiality}
For reference, we formulate the algorithm for computing the chambers of quasipolynomiality used in the proof of Theorem \ref{theorem:mainTheory}(a). 

\begin{figure}[t]
	{	\tiny
		\psset{xunit = 0.01cm, yunit = 0.01cm}
		\begin{pspicture}(-100,-100)(100,100)
			\rput[rb](-40.82480000000001,70.7107){$(1, 0, 0)$}
			\pscircle*(-40.82480000000001,70.7107){0.12}
			\rput[l](81.6497,0){$(0, 1, 0)$}
			\pscircle*(81.6497,0){0.12}
			\rput[rt](-40.82480000000001,-70.71069999999997){$(0, 0, 1)$}
			\pscircle*(-40.82480000000001,-70.71069999999997){0.12}
			\psline[linecolor=black](81.6497,0)(-40.82480000000001,-70.71069999999997)
			\psline[linecolor=black](81.6497,0)(20.412450000000007,35.35534999999999)
			\psline[linecolor=black](-40.82480000000001,70.7107)(-40.82480000000001,-70.71069999999997)
			\psline[linecolor=black](-40.82480000000001,70.7107)(20.412450000000007,35.35534999999999)
			\pscircle*(-40.82480000000001,70.7107){0.09}
			\pscircle*(81.6497,0){0.09}
			\pscircle*(-40.82480000000001,-70.71069999999997){0.09}
			\pscircle*(20.412450000000007,35.35534999999999){0.09}
			\pscircle*(20.412450000000007,-35.35534999999999){0.09}
			\pscircle*(0.00003333329999577472,0){0.09}
		\end{pspicture}
		\begin{pspicture}(-100,-100)(100,100)
			\pscircle*(-40.82480000000001,70.7107){0.12}
			\pscircle*(81.6497,0){0.12}
			\pscircle*(-40.82480000000001,-70.71069999999997){0.12}
			\psline[linecolor=black](81.6497,0)(-40.82480000000001,-70.71069999999997)
			\psline[linecolor=black](81.6497,0)(20.412450000000007,35.35534999999999)
			\psline[linecolor=black](-40.82480000000001,-70.71069999999997)(20.412450000000007,35.35534999999999)
			\rput(20.412450000000007,-11.785069526200004){$1$}
			\psline[linecolor=black](-40.82480000000001,70.7107)(-40.82480000000001,-70.71069999999997)
			\psline[linecolor=black](-40.82480000000001,70.7107)(20.412450000000007,35.35534999999999)
			\psline[linecolor=black](-40.82480000000001,-70.71069999999997)(20.412450000000007,35.35534999999999)
			\rput[r](-20.41234250849999,11.785140236899991){$2$}
			\pscircle*(-40.82480000000001,70.7107){0.09}
			\pscircle*(81.6497,0){0.09}
			\pscircle*(-40.82480000000001,-70.71069999999997){0.09}
			\pscircle*(20.412450000000007,35.35534999999999){0.09}
			\pscircle*(20.412450000000007,-35.35534999999999){0.09}
			\pscircle*(0.00003333329999577472,0){0.09}
		\end{pspicture}
		\begin{pspicture}(-100,-100)(100,100)
			\pscircle*(-40.82480000000001,70.7107){0.12}
			\pscircle*(81.6497,0){0.12}
			\pscircle*(-40.82480000000001,-70.71069999999997){0.12}
			\psline[linecolor=black](20.412450000000007,35.35534999999999)(-40.82480000000001,70.7107)
			\psline[linecolor=black](20.412450000000007,35.35534999999999)(0.00003333329999577472,0)
			\psline[linecolor=black](-40.82480000000001,70.7107)(0.00003333329999577472,0)
			\rput[rb](-6.804078339,35.35534999999999){$3$}
			\psline[linecolor=black](-40.82480000000001,-70.71069999999997)(-40.82480000000001,70.7107)
			\psline[linecolor=black](-40.82480000000001,-70.71069999999997)(0.00003333329999577472,0)
			\psline[linecolor=black](-40.82480000000001,70.7107)(0.00003333329999577472,0)
			\rput[r](-27.21649500570001,0){$4$}
			\psline[linecolor=black](81.6497,0)(-40.82480000000001,-70.71069999999997)
			\psline[linecolor=black](81.6497,0)(20.412450000000007,35.35534999999999)
			\psline[linecolor=black](-40.82480000000001,-70.71069999999997)(20.412450000000007,35.35534999999999)
			\rput(20.412450000000007,-11.785069526200004){$1$}
			\pscircle*(-40.82480000000001,70.7107){0.09}
			\pscircle*(81.6497,0){0.09}
			\pscircle*(-40.82480000000001,-70.71069999999997){0.09}
			\pscircle*(20.412450000000007,35.35534999999999){0.09}
			\pscircle*(20.412450000000007,-35.35534999999999){0.09}
			\pscircle*(0.00003333329999577472,0){0.09}
		\end{pspicture}
		\begin{pspicture}(-100,-100)(100,100)
			\pscircle*(-40.82480000000001,70.7107){0.12}
			\pscircle*(81.6497,0){0.12}
			\pscircle*(-40.82480000000001,-70.71069999999997){0.12}
			\psline[linecolor=black](20.412450000000007,35.35534999999999)(-40.82480000000001,70.7107)
			\psline[linecolor=black](20.412450000000007,35.35534999999999)(0.00003333329999577472,0)
			\psline[linecolor=black](-40.82480000000001,70.7107)(0.00003333329999577472,0)
			\rput[rb](-6.804078339,35.35534999999999){$3$}
			\psline[linecolor=black](-40.82480000000001,-70.71069999999997)(-40.82480000000001,70.7107)
			\psline[linecolor=black](-40.82480000000001,-70.71069999999997)(0.00003333329999577472,0)
			\psline[linecolor=black](-40.82480000000001,70.7107)(0.00003333329999577472,0)
			\rput[r](-27.21649500570001,0){$4$}
			\psline[linecolor=black](-40.82480000000001,-70.71069999999997)(20.412450000000007,-35.35534999999999)
			\psline[linecolor=black](-40.82480000000001,-70.71069999999997)(0.00003333329999577472,0)
			\psline[linecolor=black](20.412450000000007,-35.35534999999999)(0.00003333329999577472,0)
			\rput[rt](-6.804078339,-35.35534999999999){$5$}
			\psline[linecolor=black](81.6497,0)(20.412450000000007,-35.35534999999999)
			\psline[linecolor=black](81.6497,0)(20.412450000000007,35.35534999999999)
			\psline[linecolor=black](20.412450000000007,35.35534999999999)(0.00003333329999577472,0)
			\psline[linecolor=black](20.412450000000007,-35.35534999999999)(0.00003333329999577472,0)
			\rput(30.618658333300004,0){$6$}
			\pscircle*(-40.82480000000001,70.7107){0.09}
			\pscircle*(81.6497,0){0.09}
			\pscircle*(-40.82480000000001,-70.71069999999997){0.09}
			\pscircle*(20.412450000000007,35.35534999999999){0.09}
			\pscircle*(20.412450000000007,-35.35534999999999){0.09}
			\pscircle*(0.00003333329999577472,0){0.09}
		\end{pspicture}
		\begin{pspicture}(-100,-100)(100,100)
			\pscircle*(-40.82480000000001,70.7107){0.12}
			\pscircle*(81.6497,0){0.12}
			\pscircle*(-40.82480000000001,-70.71069999999997){0.12}
			\psline(20.412450000000007,35.35534999999999)(-40.82480000000001,70.7107)
			\psline(20.412450000000007,35.35534999999999)(0.00003333329999577472,0)
			\psline(-40.82480000000001,70.7107)(0.00003333329999577472,0)
			\rput[rb](-6.804078339,35.35534999999999){$3$}
			\psline(-40.82480000000001,-70.71069999999997)(-40.82480000000001,70.7107)
			\psline(-40.82480000000001,-70.71069999999997)(0.00003333329999577472,0)
			\psline(-40.82480000000001,70.7107)(0.00003333329999577472,0)
			\rput[r](-27.21649500570001,0){$4$}
			\psline[linecolor=black](-40.82480000000001,-70.71069999999997)(20.412450000000007,-35.35534999999999)
			\psline[linecolor=black](-40.82480000000001,-70.71069999999997)(0.00003333329999577472,0)
			\psline[linecolor=black](20.412450000000007,-35.35534999999999)(0.00003333329999577472,0)
			\rput[rt](-6.804078339,-35.35534999999999){$5$}
			\psline[linecolor=black](81.6497,0)(20.412450000000007,-35.35534999999999)
			\psline[linecolor=black](81.6497,0)(20.412450000000007,35.35534999999999)
			\psline[linecolor=black](20.412450000000007,35.35534999999999)(20.412450000000007,-35.35534999999999)
			\rput(40.82486666670002,0){$8$}
			\psline[linecolor=black](0.00003333329999577472,0)(20.412450000000007,35.35534999999999)
			\psline[linecolor=black](0.00003333329999577472,0)(20.412450000000007,-35.35534999999999)
			\psline[linecolor=black](20.412450000000007,35.35534999999999)(20.412450000000007,-35.35534999999999)
			\rput(13.608256678000004,0){$7$}
			\pscircle*(-40.82480000000001,70.7107){0.09}
			\pscircle*(81.6497,0){0.09}
			\pscircle*(-40.82480000000001,-70.71069999999997){0.09}
			\pscircle*(20.412450000000007,35.35534999999999){0.09}
			\pscircle*(20.412450000000007,-35.35534999999999){0.09}
			\pscircle*(0.00003333329999577472,0){0.09}
		\end{pspicture}
		\begin{pspicture}(-100,-100)(100,100)
			\pscircle*(-40.82480000000001,70.7107){0.12}
			\pscircle*(81.6497,0){0.12}
			\pscircle*(-40.82480000000001,-70.71069999999997){0.12}
			\psline(20.412450000000007,35.35534999999999)(-40.82480000000001,70.7107)
			\psline(20.412450000000007,35.35534999999999)(0.00003333329999577472,0)
			\psline(-40.82480000000001,70.7107)(0.00003333329999577472,0)
			\rput[rb](-6.804078339,35.35534999999999){$3$}
			\psline(-40.82480000000001,-70.71069999999997)(-40.82480000000001,70.7107)
			\psline(-40.82480000000001,-70.71069999999997)(0.00003333329999577472,0)
			\psline(-40.82480000000001,70.7107)(0.00003333329999577472,0)
			\rput[r](-27.21649500570001,0){$4$}
			\psline[linecolor=black](-40.82480000000001,-70.71069999999997)(20.412450000000007,-35.35534999999999)
			\psline[linecolor=black](-40.82480000000001,-70.71069999999997)(0.00003333329999577472,0)
			\psline[linecolor=black](20.412450000000007,-35.35534999999999)(0.00003333329999577472,0)
			\rput[rt](-6.804078339,-35.35534999999999){$5$}
			\psline[linecolor=black](20.412450000000007,35.35534999999999)(81.6497,0)
			\psline[linecolor=black](20.412450000000007,35.35534999999999)(20.412450000000007,0)
			\psline[linecolor=black](81.6497,0)(20.412450000000007,0)
			\rput(40.82489524406,11.785119023690015){$9$}
			\psline[linecolor=black](20.412450000000007,-35.35534999999999)(81.6497,0)
			\psline[linecolor=black](20.412450000000007,-35.35534999999999)(20.412450000000007,0)
			\psline[linecolor=black](81.6497,0)(20.412450000000007,0)
			\rput(40.82489524406,-11.785119023690015){$10$}
			\psline[linecolor=black](0.00003333329999577472,0)(20.412450000000007,35.35534999999999)
			\psline[linecolor=black](0.00003333329999577472,0)(20.412450000000007,-35.35534999999999)
			\psline[linecolor=black](20.412450000000007,35.35534999999999)(20.412450000000007,-35.35534999999999)
			\rput(13.608256678000004,0){$7$}
			\pscircle*(-40.82480000000001,70.7107){0.09}
			\pscircle*(81.6497,0){0.09}
			\pscircle*(-40.82480000000001,-70.71069999999997){0.09}
			\pscircle*(20.412450000000007,35.35534999999999){0.09}
			\pscircle*(20.412450000000007,-35.35534999999999){0.09}
			\pscircle*(0.00003333329999577472,0){0.09}
		\end{pspicture}
		\begin{pspicture}(-100,-100)(100,100)
			\pscircle*(-40.82480000000001,70.7107){0.12}
			\pscircle*(81.6497,0){0.12}
			\pscircle*(-40.82480000000001,-70.71069999999997){0.12}
			\psline(20.412450000000007,35.35534999999999)(-40.82480000000001,70.7107)
			\psline(20.412450000000007,35.35534999999999)(0.00003333329999577472,0)
			\psline(-40.82480000000001,70.7107)(0.00003333329999577472,0)
			\rput[rb](-6.804078339,35.35534999999999){$3$}
			\psline(-40.82480000000001,-70.71069999999997)(-40.82480000000001,70.7107)
			\psline(-40.82480000000001,-70.71069999999997)(0.00003333329999577472,0)
			\psline(-40.82480000000001,70.7107)(0.00003333329999577472,0)
			\rput[r](-27.21649500570001,0){$4$}
			\psline[linecolor=black](-40.82480000000001,-70.71069999999997)(20.412450000000007,-35.35534999999999)
			\psline[linecolor=black](-40.82480000000001,-70.71069999999997)(0.00003333329999577472,0)
			\psline[linecolor=black](20.412450000000007,-35.35534999999999)(0.00003333329999577472,0)
			\rput[rt](-6.804078339,-35.35534999999999){$5$}
			\psline[linecolor=black](20.412450000000007,35.35534999999999)(81.6497,0)
			\psline[linecolor=black](20.412450000000007,35.35534999999999)(20.412450000000007,0)
			\psline[linecolor=black](81.6497,0)(20.412450000000007,0)
			\rput(40.82489524406,11.785119023690015){$9$}
			\psline[linecolor=black](20.412450000000007,-35.35534999999999)(81.6497,0)
			\psline[linecolor=black](20.412450000000007,-35.35534999999999)(20.412450000000007,0)
			\psline[linecolor=black](81.6497,0)(20.412450000000007,0)
			\rput(40.82489524406,-11.785119023690015){$10$}
			\psline[linecolor=black](20.412450000000007,-35.35534999999999)(0.00003333329999577472,0)
			\psline[linecolor=black](20.412450000000007,-35.35534999999999)(20.412450000000007,0)
			\psline[linecolor=black](0.00003333329999577472,0)(20.412450000000007,0)
			\rput(13.608297502799985,-11.785140236899991){$11$}
			\psline[linecolor=black](20.412450000000007,35.35534999999999)(0.00003333329999577472,0)
			\psline[linecolor=black](20.412450000000007,35.35534999999999)(20.412450000000007,0)
			\psline[linecolor=black](0.00003333329999577472,0)(20.412450000000007,0)
			\rput(13.608297502800013,11.785140236899991){$12$}
			\pscircle*(-40.82480000000001,70.7107){0.09}
			\pscircle*(81.6497,0){0.09}
			\pscircle*(-40.82480000000001,-70.71069999999997){0.09}
			\pscircle*(20.412450000000007,35.35534999999999){0.09}
			\pscircle*(20.412450000000007,-35.35534999999999){0.09}
			\pscircle*(0.00003333329999577472,0){0.09}
		\end{pspicture}
		
	}
	\caption{Algorithm \ref{algorithm:subdivisionByArbitrarySlices} shown on the vector partition function of $(1,0,0)$, $(0,1,0)$, $(0,0,1)$, $(1,1,0)$, $(0,1,1)$, $(1,1,1)$ in barycentric coordinates. The second triangle shows the steps of Algorithm \ref{algorithm:subdivisionByArbitrarySlicesOnce} with respect to vector $(1,1,0)$, the third, fourth and fifth triangle show the steps of Algorithm \ref{algorithm:subdivisionByArbitrarySlicesOnce} with respect to $(0,1,1)$, and the final two triangles show the steps of Algorithm \ref{algorithm:subdivisionByArbitrarySlicesOnce} with respect to $(1,1,1)$. \label{figure:algorithmArbitrarySlicedChambers}}
\end{figure}
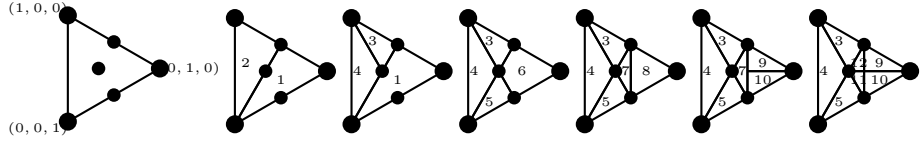

\begin{algorithm}[Arbitrarily sliced chambers of quasipolynomiality]~\label{algorithm:subdivisionByArbitrarySlices}
\begin{itemize}
\item Input: a set of vectors $\Delta=\bm \alpha_1, \dots, \bm \alpha_s$ of full rank with positive integer coordinates.
\item Output: a collection of chambers with normal neighbors $C_1, \dots, C_N$ so that the vector partition function $P_{\Delta}(\bm \gamma)$ is given by a quasipolynomial in $\bm \gamma$ for each $ C_i$ and is equal to $0$ whenever $\bm \alpha\notin C_1\cup \dots C_N$
\end{itemize}
\begin{itemize}
\item[Step 1.] Reorder $\bm \alpha_1, \dots, \bm \alpha_s$ so that the first $n$ vectors are linearly independent.
\item[Step 2.] Compute the cone $C= C(\bm \alpha_1, \dots, \bm \alpha_n)$ using Algorithm \ref{algorithm:normalsFromGenerators}. 
\item[Step 3.] Initialize $\mathcal N=\{C\}$.
\item[Step 4.] For $i$ from $1$ to $s$ do:
\begin{itemize}
\item[Step 4.1.] Subdivide $\mathcal N$ to a collection normal with respect to $\bm\alpha_i$ using Algorithm \ref{algorithm:subdivisionByArbitrarySlicesOnce} and write the resulting subdivision back into $\mathcal N$.
\end{itemize}
\end{itemize}
\end{algorithm}
The output of this algorithm (from Step 2 on) depends on the order of vectors $\bm \alpha_1, \dots, \bm \alpha_s$. The algorithm is illustrated in Figure \ref{figure:algorithmArbitrarySlicedChambers}.

\subsection{Computing the largest possible chambers}
In Theorem \ref{theorem:mainTheory}(a) we saw that a single quasipolynomial formula is guaranteed to hold over a chamber computed by induction as described in the proof of Theorem \ref{theorem:subdivision}. To achieve the subdivision of Theorem \ref{theorem:subdivision}, we used the Directional Subdivision Lemma \ref{lemma:directionalSubdivision} in Step 8.1 of Algorithm \ref{algorithm:subdivisionByArbitrarySlicesOnce}.

The Directional Subdivision Lemma is realized by refining the starting chamber into smaller chambers with new walls with normal vectors $\bm c_{i,j}$ given by \eqref{eq:directionalSubdivisionNewWalls}. These are not necessarily spanned by $n-1$-element subsets of $\Delta$, which was promised by Theorem \ref{theorem:mainTheory}(b). So, we can optimize our chamber computations by replacing the Directional Subdivision Lemma in Step 8.1 of Algorithm \ref{algorithm:subdivisionByArbitrarySlices} with the following procedure.
\begin{algorithm}[Drop-in replacement for Directional Subdivision.]~\label{algorithm:directionalSubdivisionReplacement}
\begin{itemize}
\item[Step 1.] Compute all normals $\bm c_{i,j}$ needed by the Directional Subdivision Lemma. 
\item[Step 2.] Find the first vector $\bm c_{i,j}$ that gives a normal to a plane spanned by $n-1$-element subsets of $\Delta$. 
\begin{itemize}
\item[Step 2.1.] If there's no such vector, do not subdivide the chamber. Add the chamber to the set of chambers $\mathcal R$ from Algorithm \ref{algorithm:subdivisionByArbitrarySlicesOnce}.
\item[Step 2.2.] If there's such a vector, use it to split the chamber in two using the Subdivision by Plane Lemma \ref{lemma:subdivisionByPlane}	
\end{itemize}
\end{itemize}	
\end{algorithm}
There will be chambers that have multiple walls whose normals have positive scalar products with the direction vector $\bm \alpha_s$ but which will be left non-divided by Step 2.1 above. Future subdivisions of such chambers may still arise through slices that come from neighboring chambers through Algorithm \ref{algorithm:extendSubdivisionToNormal} through a subdivision of a neighboring chamber to normal.

This leads to the following algorithm for computing a smaller set of (larger) chambers of quasipolynomiality. The algorithm can be regarded as al alternative to Algorithm \ref{algorithm:subdivisionByArbitrarySlices}.
\begin{algorithm}[Properly sliced chambers of quasipolynomiality]\label{algorithm:subdivisionProper}~
\begin{itemize}
\item The algorithm has same input and output as Algorithm \ref{algorithm:subdivisionByArbitrarySlices}.
\item The algorithm has the same implementation as Algorithm \ref{algorithm:subdivisionByArbitrarySlices} except that it uses \ref{algorithm:directionalSubdivisionReplacement} in place of Step 8.2 of the sub-routine of Algorithm \ref{algorithm:subdivisionByArbitrarySlicesOnce}.
\end{itemize}
\end{algorithm}

The chambers computed in Algorithm \ref{algorithm:subdivisionProper} are still finer than they could be made, due to some excessive subdivision with hyperplanes that are linearly spanned by $n-1$ vectors of $\Delta$ but the hyperplane slice lies outside of any cone spanned by $n$ vectors. This is explained in \cite[Example 4.15]{BaldoniBeckCochetVergne:VolumeComputationForPolytopesAndPartitionFunctions} and \cite[Section 3]{BaldoniDeLoeraVergne:CountingIntegerFlowsInNetworks}. To elaborate on this matter, we need a couple of definitions.

\begin{definition}[Sufficiently refined by $\Delta$]
We say that a set of chambers $C_1, \dots, C_k$ is sufficiently refined by $\Delta$ if, for every cone $C$ spanned by $n$ vectors in $\Delta$ and every chamber $C_i$, the interior of $C_i$ lies entirely in $C$ or lies entirely outside of $C$. 
\end{definition}
\begin{definition}[Combinatorially separated]
We say that two chambers $E$ and $D$ are combinatorially separated by $\Delta$ if the following hold.
\begin{itemize}
\item The set $\{E,D\}$ is sufficiently refined by $\Delta$.
\item There exists a cone $C$ spanned by $n$ vectors in $\Delta$ such that $C$ contains the interior of one of the cones but does not contain the interior of the other.
\end{itemize} 
\end{definition}

Inspection of the formulas produced by Algorithm \ref{algorithm:main} suggest that two chambers may yield different quasipolynomial formulas only if they are combinatorially separated by $\Delta$. Therefore, to obtain larger combinatorial chambers, we can amalgamate chambers that are not combinatorially separated with the following algorithm.

\begin{algorithm}[Amalgamation]~\label{algorithm:amalgamate}
\begin{itemize}
\item Input.
\begin{itemize}
\item A set of chambers $\mathcal I$ that are combinatorially separated by $\Delta$
\end{itemize}
\item Output: collection of chambers $\mathcal D$ with the following properties.
\begin{itemize}
	\item $\mathcal D$ is combinatorially separated by $\Delta$.
	\item  $|\mathcal D|<|\mathcal C|$
	\item Each $D\in \mathcal D$ is the union of finitely many of the elements of $\mathcal I$.
\end{itemize} 
\end{itemize}
\begin{itemize}
\item[Step 1.] Initialize a map $\mathcal M$ from $\mathcal I$ to the set of finite sequences of $0$'s and $1$'s.
\item[Step 2.] For each $C\in \mathcal I$:
\begin{itemize}
	\item[Step 2.1.] Initialize $L$ to be the empty list.
	\item[Step 2.2.] For the $j^{th}$ cone $E$ spanned by $n$ vectors in $\Delta$, enumerated in the same order for each $C$:
	\begin{itemize}
	\item[Step 2.2.1.] If $C\subset E$, append $1$ to $L$. Else, append $0$ to the list $L$.
	\item[Step 2.2.2.](Optional) If $L$ is too large to comfortably hold in the computer's memory, compress it using a hash function that sends different sequences to different compressions with practical probability $1$. 
	\item[Step 2.2.3.] Go back to Step 2.2 until done.
	\end{itemize}
	\item[Step 2.3.] Set $\mathcal M(C)=L$.
	\item[Step 2.4.] Go back to step 2 until done.
\end{itemize}
\item[Step 3.] Initialize $\mathcal D$ with the empty list.
\item[Step 4.] For every possible value $L$ of the map $\mathcal M$:
\begin{itemize}
	\item[Step 4.1.] Let $C_{1}, \dots, C_{k}$ be all cones for which $\mathcal M(C_{i})= L$.
	\item[Step 4.2.] If the union of the $C_{i}$'s is convex, combine them into a single chamber and put it into $\mathcal D$.
	\item[Step 4.3.] Else, put all the $C_{i}$'s unchanged into $\mathcal D$.
\end{itemize}
\end{itemize}
	
\end{algorithm}

We conclude with a final algorithm for computing chambers of quasipolynomiality. This can be regarded as an alternative to Algorithms \ref{algorithm:subdivisionByArbitrarySlices} and \ref{algorithm:subdivisionProper}.
\begin{algorithm}[Amalgamated chambers of quasipolynomiality]\label{algorithm:subdivisionProperAmalgamated} ~
\begin{itemize}
	\item The input and output are the same as in Algorithms \ref{algorithm:subdivisionByArbitrarySlices} and \ref{algorithm:subdivisionProper}.
\end{itemize}
\begin{itemize}
	\item[Step 1.] Run Algorithm \ref{algorithm:subdivisionProper} to generate a starting set of chambers of quasipolynomiality.
	\item[Step 2.] Run Algorithm \ref{algorithm:amalgamate} to amalgamate the chambers into larger ones.
\end{itemize}
\end{algorithm}

Our initial observations suggest that the bulk of the run time is taken by the Amalgamation Algorithm. So, if a particular application of the vector partition function computations does not require the most optimal set of chambers, it makes sense to skip the Amalgamation Algorithm. 

The work \cite{BaldoniDeLoeraVergne:CountingIntegerFlowsInNetworks} describes a scheme for computing the chambers of quasipolynomiality that yields the same result as Algorithm \ref{algorithm:subdivisionProperAmalgamated}. The scheme is quite different from what we describe here: \cite{BaldoniDeLoeraVergne:CountingIntegerFlowsInNetworks} starts from one fully refined chamber, and work their way up to a decomposition of the entire space by finding a fully refined neighbor along a given wall. Our approach goes in the opposite direction: we start with the convex hull of the vectors in $\Delta$, and subdivide it down to get our chambers. In this way, our algorithm takes advantage of the information of which chambers are neighbors of one another, which helps speed up computations.

The chambers counts for the particular case of root systems are reported in  \cite{BaldoniBeckCochetVergne:VolumeComputationForPolytopesAndPartitionFunctions} (based on \cite{BaldoniDeLoeraVergne:CountingIntegerFlowsInNetworks}). A glance at the reported computational times may suggest that our final Algorithm \ref{algorithm:subdivisionProperAmalgamated} is more efficient, however it is hard to establish that with certainty based on the stopwatch alone, given that the algorithms were coded in different programming environments and ran on different hardware. Certainly, the two algorithms are very different: ours is top-down (and, arguably, technically more compute), whereas the algorithm in \cite{BaldoniDeLoeraVergne:CountingIntegerFlowsInNetworks} is bottom-up. The algorithm of \cite{BaldoniDeLoeraVergne:CountingIntegerFlowsInNetworks} is still very important, as it allows to compute a combinatorial chamber containing a given interior point, without computing the rest of the chamber complex. 

At any rate, we can report that root systems $B_5$ and $C_5$ have $138061$ combinatorial chambers each, a number that does not appear in \cite{BaldoniBeckCochetVergne:VolumeComputationForPolytopesAndPartitionFunctions}. The two numbers must be equal as each of the sets can be transformed to the other by rescaling the longer (respectively, shorter) vectors. The fact that this is confirmed by the software provides evidence for its correctness.

Our computations are summarized in the following Table. We recall that the output of Algorithms \ref{algorithm:subdivisionByArbitrarySlices} and \ref{algorithm:subdivisionProper} depend on the order of our vectors; our results use the graded lexicographic order on the simple-coordinate representations of these root systems. For readers not familiar with root systems, we recall that the subscript of the root system gives the dimension of the vector space, for example the root system of $F_4$ is in $4$-dimensional space.

\begin{tabular}{|c|c|c|c|c|}\hline
&&\multicolumn{3}{c|}{Number of chambers + run time}	\\
R. Sys. & $|\Delta|$& Alg. \ref{algorithm:subdivisionByArbitrarySlices}& Alg. \ref{algorithm:subdivisionProper}& Alg. \ref{algorithm:subdivisionProperAmalgamated}\\\hline
$A_2$&$3$&  \begin{tabular}{c}$2$ \\$5\text{ms}$ \end{tabular} &\begin{tabular}{c}$2$ \\$5\text{ms}$ \end{tabular}  &\begin{tabular}{c}$2$ \\$5\text{ms}$ \end{tabular}\\\hline
$B_2$&$4$&  \begin{tabular}{c}$3$ \\$10\text{ms}$ \end{tabular} &\begin{tabular}{c}$3$ \\$10\text{ms}$ \end{tabular}  &\begin{tabular}{c}$3$ \\$10\text{ms}$ \end{tabular}\\\hline
$C_2$&$4$&  \begin{tabular}{c}$3$ \\$18\text{ms}$ \end{tabular} &\begin{tabular}{c}$3$ \\$10\text{ms}$ \end{tabular}  &\begin{tabular}{c}3 \\$10\text{ms}$ \end{tabular}\\\hline
$G_2$&$6$&  \begin{tabular}{c} $5$\\$18\text{ms} $ \end{tabular} &\begin{tabular}{c}$5$ \\$18\text{ms} $ \end{tabular}  &\begin{tabular}{c} $5 $\\$18\text{ms} $ \end{tabular}\\\hline 
$A_3$&$6$&  \begin{tabular}{c}$ 7$ \\$45\text{ms} $ \end{tabular} &\begin{tabular}{c}$7$ \\$45\text{ms}$ \end{tabular}  &\begin{tabular}{c}$7 $ \\$50\text{ms} $ \end{tabular}\\\hline
$B_3$&$9$&  \begin{tabular}{c}$ 45$ \\$ 1080\text{ms}$ \end{tabular} &\begin{tabular}{c}$23 $ \\$316 \text{ms}$ \end{tabular}  &\begin{tabular}{c}$ 23$ \\$320\text{ms} $ \end{tabular}\\\hline
$C_3$&$9$&  \begin{tabular}{c} $31 $\\$581\text{ms} $ \end{tabular} &\begin{tabular}{c}$23 $ \\$308ms $ \end{tabular}  &\begin{tabular}{c}$ 23$ \\$ 309\text{ms}$ \end{tabular}\\\hline
$A_4$&$10$&  \begin{tabular}{c}$56 $\\$3.9\text{s} $ \end{tabular} &\begin{tabular}{c}$48 $ \\$2.8\text{s} $ \end{tabular}  &\begin{tabular}{c}$48 $ \\$2.9\text{s} $ \end{tabular}\\\hline
$B_4$&$16$&  \begin{tabular}{c}$ >1000000$ \\
first 13 vect. only\\
$2\text{h} $ \end{tabular} &\begin{tabular}{c}$1111 $ \\$15.5s $ \end{tabular}  &\begin{tabular}{c}$695 $ \\$20.8\text{s} $ \end{tabular}\\\hline
$C_4$&$16$&  \begin{tabular}{c}$ N/A$ \\$N/A $ \end{tabular} &\begin{tabular}{c}$1035$ \\$16.5s $ \end{tabular}  &\begin{tabular}{c}$695 $ \\$ 19.8\text{s}$ \end{tabular}\\\hline
$D_4$&$12$&  \begin{tabular}{c}$12721$ \\$50s $ \end{tabular} &\begin{tabular}{c}$ 145$ \\$14.3s $ \end{tabular}  &\begin{tabular}{c}$133 $ \\$14.8\text{s} $ \end{tabular}\\\hline
$F_4$&24&N/A &\begin{tabular}{c}$39058 $ \\$3 \text{min }43 \text{s} $ \end{tabular}  &\begin{tabular}{c}$12946 $ \\$ 12\text{min } 19 \text{s}$ \end{tabular}\\\hline
$A_5$&15&  \begin{tabular}{c}$>1800000 $ \\
first 14 vect. only\\
$2\text{h} $ \end{tabular} &\begin{tabular}{c}$1006 $ \\$ 37\text{s}$ \end{tabular}  &\begin{tabular}{c}$820 $ \\$ 47\text{s}$ \end{tabular}\\\hline
$B_5$ &$25$&N/A&\begin{tabular}{c}
$472236$\\
$1\text{h}$
\end{tabular} &\begin{tabular}{c} $138061$\\$5\text{h } 40\text{min}$  \end{tabular}\\\hline
$C_5$ &$25$&N/A&\begin{tabular}{c} $479674$\\$56\text{min} $ \end{tabular}  &\begin{tabular}{c}$ 138061$ \\$10 \text{h }27\text{min} $ \end{tabular}\\\hline
$A_6$&$21$ &N/A &\begin{tabular}{c} $142908$\\$11\text{min }35\text{s} $ \end{tabular}  &\begin{tabular}{c}$44288$ \\$1\text{h } 26\text{min} $ \end{tabular}\\\hline
\end{tabular}

\nocite{LatE} \nocite{SphericalExplorer}
\bibliographystyle{plain}
%\bibliography{../bibliography}
% Generated using \bibliography{../bibliography}
% Bibliography start

% Bibliography end.

\appendix
\section{Tables with the Kostant Partition Function }
In the appendix, we include computer printouts of the vector partition functions of the positive root systems $A_2, B_2, C_2,G_2 A_3, B_3,C_3, A_4$. We recall that the root systems are sets of vectors that are of special interest to Lie Representation theory. Their vector partition functions are called Kostant Partition Functions. Their study was the original inspiration for the present article, with a number of desired future applications. In view of the length of the present article, it appears wise to postpone such applications to future work.

To the best of our knowledge, this is the first work to publish formulas for the vector partition functions of $G_2, B_3, C_3$ and $A_4$, although we seem to not be the first one to obtain such a printout - the credit for this should be given, to the best of our knowledge, to \cite{Cochet:2005VectorPartitionFunctionAndRepresentationTheory}. Formulas for vector partition functions have been available in the software ``barvinok'' \cite{SvenVerdoolaege:Barvinok}, especially in the unimodular case of $A_3, A_4$ - i.e., the case where the ambient lattice is always $\mathbb Z^n$; we are not aware if those have ever been printed in an easy to read format.

The tables are organized as follows. First, we plot a picture of the combinatorial chambers of the vector partition function. The chambers are labeled with an identifier. 

Next, we present the Kostant partition function over each chamber. In the first column, we present the chamber identifier. In the second column, we show the polynomial formula over each shift of the lattice $\Lambda$. If the lattice of our quasipolynomial is $\mathbb Z^n$, we omit it, else we indicated the generators of the lattice in the first row of the given chamber. The third column gives the lattice shifts for which each of the polynomials is defined. When one polynomial formula is valid over multiple lattice shifts, we list multiple lattice shifts in the last column. For a fixed chamber, the polynomials with different lattice shifts have higher-order terms equal. We print those equal parts in black, and print all terms that vary (``error terms'') in red. The error terms in $G_2$ appear to be linear, in $B_3$ - cubic and in $C_3$ - quadratic.

Third, we list a detailed definition of each of the chambers. The first column holds the identifier; the second holds the defining inequalities. The next column holds the vertices of each chamber. Next, we present an internal point, obtained by summing all vertices and rescaling down to relatively prime coordinates. The final column holds the identifiers of all neighbors that are neighbors to the given one. Next to each vector, we indicate the vector partition quasipolynomial formula (from the previous table), evaluated on the vector.

In order to check our work, for small enough vectors, we compute the vector partition function by enumerating all vector partition with separate program we wrote for the purpose. A checkbox following the number of vector partitions indicates that we checked this independent program to verify our computation. Lack of checkbox means our software skipped the check to speed up computations. We would like to report that this mechanism helped us find errors with our processing of lattice shifts that were only showing up on fairly large examples (namely, in type $B_3$). This error-check gives us a reasonable confidence in the implementation of the relatively long algorithms of the present article. 

\pagebreak
\begin{table}[h!]
	\begin{center}
		{\tiny
			\psset{xunit=4cm, yunit=4cm}
			\begin{pspicture}(0,0)(1,1)
				\psline[linecolor=black, linestyle=dashed](0,0)(1,0)
				\rput[lt](1,0){(1, 0)}
				\pscircle*[linecolor=black](1,0){0.05}
				\psline[linecolor=black, linestyle=dashed](0,0)(0,1)
				\rput[rb](0,1){(0, 1)}
				\pscircle*[linecolor=black](0,1){0.05}
				\psline[linecolor=blue, linestyle=dashed](0,0)(0,1)
				\psline[linecolor=blue, linestyle=dashed](0,0)(0.5,0.5)
				\psline[linecolor=blue, linestyle=dashed](0,0)(1,0)
				\rput[rb](0.25,0.75){1}
				\rput[lt](0.75,0.25){2}
				\pscircle*[linecolor=blue](1,0){0.0375}
				\pscircle*[linecolor=blue](0,1){0.0375}
				\pscircle*[linecolor=blue](0.5,0.5){0.0375}
		\end{pspicture}}
		\caption{$2$ combinatorial chambers of \(A_2\):(1, 0), (0, 1), (1, 1)}
	\end{center}
\end{table}
\begin{longtable}{|c|c|c|}\caption{\footnotesize V.p.f. of \(A_2\):(1, 0), (0, 1), (1, 1)}\\\hline N & Polynomial/Lattice & Shift(s)\\ \hline
	\endfirsthead\multicolumn{3}{c}{{\bfseries \tablename\ \thetable{} -- continued from previous page}} \\
	\hline  N &  Polynomial/Lattice & Shift(s) \\ \hline
	\endhead\hline \multicolumn{3}{|c|}{{Continued on next page}} \\ \hline\endfoot\endlastfoot1&\(x_{1} +1\)&-\\
	\hline\hline 
	2&\(x_{2} +1\)&-\\
	\hline\end{longtable}

\begin{longtable}{|ccccc|}\caption{\footnotesize V.p.f. of \(A_2\):(1, 0), (0, 1), (1, 1)}\\\hline N & Defining inequalities & Vertices& Int. Pt.& Neighbors \\ \hline
	\endfirsthead\multicolumn{5}{c}{{\bfseries \tablename\ \thetable{} -- continued from previous page}} \\
	\hline N & Defining inequalities & Vertices&Int. Pt. & Neighbors\\ \hline
	\endhead\hline \multicolumn{5}{|c|}{{Continued on next page}} \\ \hline\endfoot\endlastfoot1&
	\(\begin{array}{rcl}x_{1}&\geq& 0\\
		-x_{1}+x_{2}&\geq& 0\\
	\end{array}\)&\(\begin{array}{l}(0, 1): 1\checkmark\\(1, 1): 2\checkmark\end{array}\)&(1, 2): 2\checkmark& [2]\\\hline
	2&
	\(\begin{array}{rcl}x_{2}&\geq& 0\\
		x_{1}-x_{2}&\geq& 0\\
	\end{array}\)&\(\begin{array}{l}(1, 0): 1\checkmark\\(1, 1): 2\checkmark\end{array}\)&(2, 1): 2\checkmark& [1]\\\hline
\end{longtable}

\allowdisplaybreaks\begin{align*}&~~~
	\frac{1}{(1-x_{1} ) (1-x_{2} ) (1-x_{1} x_{2} ) }\\=&~~~
	\displaystyle \frac{-x_{2}^{-1}}{(1-x_{1} )^2 (1-x_{1} x_{2} ) }\\&
	+\displaystyle \frac{x_{2}^{-1}}{(1-x_{1} )^2 (1-x_{2} ) }\\=&
	-x_{2}^{-1}\cdot \left(x_{1} \partial_{1} -x_{2} \partial_{2} \right)\cdot\frac{1}{(1-x_{1} ) (1-x_{1} x_{2} ) }\\&
	\displaystyle +x_{2}^{-1}\cdot \left(x_{1} \partial_{1} \right)\cdot\frac{1}{(1-x_{1} ) (1-x_{2} ) }\end{align*}\begin{table}[h!]
	\begin{center}
		{\tiny
			\psset{xunit=4cm, yunit=4cm}
			\begin{pspicture}(0,0)(1,1)
				\psline[linecolor=black, linestyle=dashed](0,0)(1,0)
				\rput[lt](1,0){(1, 0)}
				\pscircle*[linecolor=black](1,0){0.05}
				\psline[linecolor=black, linestyle=dashed](0,0)(0,1)
				\rput[rb](0,1){(0, 1)}
				\pscircle*[linecolor=black](0,1){0.05}
				\psline[linecolor=blue, linestyle=dashed](0,0)(0,1)
				\psline[linecolor=blue, linestyle=dashed](0,0)(0.33333333,0.66666667)
				\psline[linecolor=blue, linestyle=dashed](0,0)(1,0)
				\psline[linecolor=blue, linestyle=dashed](0,0)(0.5,0.5)
				\rput[rb](0.41666667,0.58333333){1}
				\rput[rb](0.16666667,0.83333333){2}
				\rput[lt](0.75,0.25){3}
				\pscircle*[linecolor=blue](1,0){0.0375}
				\pscircle*[linecolor=blue](0,1){0.0375}
				\pscircle*[linecolor=blue](0.5,0.5){0.0375}
				\pscircle*[linecolor=blue](0.33333333,0.66666667){0.0375}
		\end{pspicture}}
		\caption{$3$ combinatorial chambers of \(B_2\):(1, 0), (0, 1), (1, 1), (1, 2)}
	\end{center}
\end{table}

\pagebreak

\begin{longtable}{|c|c|c|}\caption{\footnotesize V.p.f. of \(B_2\):(1, 0), (0, 1), (1, 1), (1, 2)}\\\hline N & Polynomial/Lattice & Shift(s)\\ \hline
	\endfirsthead\multicolumn{3}{c}{{\bfseries \tablename\ \thetable{} -- continued from previous page}} \\
	\hline  N &  Polynomial/Lattice & Shift(s) \\ \hline
	\endhead\hline \multicolumn{3}{|c|}{{Continued on next page}} \\ \hline\endfoot\endlastfoot\multirow{3}{*}{1}&\(\Lambda=\langle((1, 0),(0, 2))\rangle\)&\\\hline&\(\begin{array}{l}-\frac{x_{1}^2}{2}+x_{1} x_{2} -\frac{x_{2}^2}{4}+\frac{x_{1} }{2}+\frac{x_{2} }{2}\\
		\color{red}
		+1\color{black}
	\end{array}
	\)&\((0, 0)\)\\
	\cline{2-3}&\(\begin{array}{l}-\frac{x_{1}^2}{2}+x_{1} x_{2} -\frac{x_{2}^2}{4}+\frac{x_{1} }{2}+\frac{x_{2} }{2}\\
		\color{red}
		+\frac{3}{4}\color{black}
	\end{array}
	\)&\((0, 1)\)\\
	\hline\hline 
	2&\(\frac{x_{1}^2}{2}+\frac{3}{2}x_{1} +1\)&-\\
	\hline\hline 
	\multirow{3}{*}{3}&\(\Lambda=\langle((1, 0),(0, 2))\rangle\)&\\\hline&\(\begin{array}{l}\frac{x_{2}^2}{4}+x_{2} \\
		\color{red}
		+1\color{black}
	\end{array}
	\)&\((0, 0)\)\\
	\cline{2-3}&\(\begin{array}{l}\frac{x_{2}^2}{4}+x_{2} \\
		\color{red}
		+\frac{3}{4}\color{black}
	\end{array}
	\)&\((0, 1)\)\\
	\hline\end{longtable}

\pagebreak

\begin{longtable}{|ccccc|}\caption{\footnotesize V.p.f. of \(B_2\):(1, 0), (0, 1), (1, 1), (1, 2)}\\\hline N & Defining inequalities & Vertices& Int. Pt.& Neighbors \\ \hline
	\endfirsthead\multicolumn{5}{c}{{\bfseries \tablename\ \thetable{} -- continued from previous page}} \\
	\hline N & Defining inequalities & Vertices&Int. Pt. & Neighbors\\ \hline
	\endhead\hline \multicolumn{5}{|c|}{{Continued on next page}} \\ \hline\endfoot\endlastfoot1&
	\(\begin{array}{rcl}-x_{1}+x_{2}&\geq& 0\\
		2x_{1}-x_{2}&\geq& 0\\
	\end{array}\)&\(\begin{array}{l}(1, 1): 2\checkmark\\(1, 2): 3\checkmark\end{array}\)&(2, 3): 5\checkmark& [2, 3]\\\hline
	2&
	\(\begin{array}{rcl}x_{1}&\geq& 0\\
		-2x_{1}+x_{2}&\geq& 0\\
	\end{array}\)&\(\begin{array}{l}(0, 1): 1\checkmark\\(1, 2): 3\checkmark\end{array}\)&(1, 3): 3\checkmark& [4]\\\hline
	3&
	\(\begin{array}{rcl}x_{2}&\geq& 0\\
		x_{1}-x_{2}&\geq& 0\\
	\end{array}\)&\(\begin{array}{l}(1, 0): 1\checkmark\\(1, 1): 2\checkmark\end{array}\)&(2, 1): 2\checkmark& [4]\\\hline
\end{longtable}

\allowdisplaybreaks\begin{align*}&~~~
	\frac{1}{(1-x_{1} ) (1-x_{2} ) (1-x_{1} x_{2} ) (1-x_{1} x_{2}^2) }\\=&~~~
	\displaystyle \frac{x_{1} x_{2}^{-1}}{(1-x_{1} )^3 (1-x_{1} x_{2} ) }\\&
	+\displaystyle \frac{-x_{1} -x_{2}^{-1}-x_{2}^{-2}-x_{2}^{-3}}{(1-x_{1} )^3 (1-x_{1} x_{2}^2) }\\&
	+\displaystyle \frac{x_{2}^{-3}}{(1-x_{1} )^3 (1-x_{2} ) }\\=&
	x_{1} x_{2}^{-1}\cdot \frac{1}{2}\left(x_{1} \partial_{1} -x_{2} \partial_{2} \right)^2\cdot\frac{1}{(1-x_{1} ) (1-x_{1} x_{2} ) }\\&
	\displaystyle +\left(-x_{1} -x_{2}^{-1}-x_{2}^{-2}-x_{2}^{-3}\right)\cdot \frac{1}{2}\left(x_{1} \partial_{1} -\frac{x_{2} \partial_{2} }{2}\right)^2\cdot\frac{1}{(1-x_{1} ) (1-x_{1} x_{2}^2) }\\&
	\displaystyle +x_{2}^{-3}\cdot \frac{1}{2}\left(x_{1} \partial_{1} \right)^2\cdot\frac{1}{(1-x_{1} ) (1-x_{2} ) }\end{align*}\begin{table}[h!]
	\begin{center}
		{\tiny
			\psset{xunit=4cm, yunit=4cm}
			\begin{pspicture}(0,0)(1,1)
				\psline[linecolor=black, linestyle=dashed](0,0)(1,0)
				\rput[lt](1,0){(1, 0)}
				\pscircle*[linecolor=black](1,0){0.05}
				\psline[linecolor=black, linestyle=dashed](0,0)(0,1)
				\rput[rb](0,1){(0, 1)}
				\pscircle*[linecolor=black](0,1){0.05}
				\psline[linecolor=blue, linestyle=dashed](0,0)(0,1)
				\psline[linecolor=blue, linestyle=dashed](0,0)(0.5,0.5)
				\psline[linecolor=blue, linestyle=dashed](0,0)(1,0)
				\psline[linecolor=blue, linestyle=dashed](0,0)(0.66666667,0.33333333)
				\rput[lt](0.58333333,0.41666667){1}
				\rput[rb](0.25,0.75){2}
				\rput[lt](0.83333333,0.16666667){3}
				\pscircle*[linecolor=blue](1,0){0.0375}
				\pscircle*[linecolor=blue](0,1){0.0375}
				\pscircle*[linecolor=blue](0.5,0.5){0.0375}
				\pscircle*[linecolor=blue](0.66666667,0.33333333){0.0375}
		\end{pspicture}}
		\caption{$3$ combinatorial chambers of \(C_2\):(1, 0), (0, 1), (1, 1), (2, 1)}
	\end{center}
\end{table}
\begin{longtable}{|c|c|c|}\caption{\footnotesize V.p.f. of \(C_2\):(1, 0), (0, 1), (1, 1), (2, 1)}\\\hline N & Polynomial/Lattice & Shift(s)\\ \hline
	\endfirsthead\multicolumn{3}{c}{{\bfseries \tablename\ \thetable{} -- continued from previous page}} \\
	\hline  N &  Polynomial/Lattice & Shift(s) \\ \hline
	\endhead\hline \multicolumn{3}{|c|}{{Continued on next page}} \\ \hline\endfoot\endlastfoot\multirow{3}{*}{1}&\(\Lambda=\langle((2, 0),(0, 1))\rangle\)&\\\hline&\(\begin{array}{l}-\frac{x_{1}^2}{4}+x_{1} x_{2} -\frac{x_{2}^2}{2}+\frac{x_{1} }{2}+\frac{x_{2} }{2}\\
		\color{red}
		+1\color{black}
	\end{array}
	\)&\((0, 0)\)\\
	\cline{2-3}&\(\begin{array}{l}-\frac{x_{1}^2}{4}+x_{1} x_{2} -\frac{x_{2}^2}{2}+\frac{x_{1} }{2}+\frac{x_{2} }{2}\\
		\color{red}
		+\frac{3}{4}\color{black}
	\end{array}
	\)&\((1, 0)\)\\
	\hline\hline 
	\multirow{3}{*}{2}&\(\Lambda=\langle((2, 0),(0, 1))\rangle\)&\\\hline&\(\begin{array}{l}\frac{x_{1}^2}{4}+x_{1} \\
		\color{red}
		+1\color{black}
	\end{array}
	\)&\((0, 0)\)\\
	\cline{2-3}&\(\begin{array}{l}\frac{x_{1}^2}{4}+x_{1} \\
		\color{red}
		+\frac{3}{4}\color{black}
	\end{array}
	\)&\((1, 0)\)\\
	\hline\hline 
	3&\(\frac{x_{2}^2}{2}+\frac{3}{2}x_{2} +1\)&-\\
	\hline\end{longtable}

\begin{longtable}{|ccccc|}\caption{\footnotesize V.p.f. of \(C_2\):(1, 0), (0, 1), (1, 1), (2, 1)}\\\hline N & Defining inequalities & Vertices& Int. Pt.& Neighbors \\ \hline
	\endfirsthead\multicolumn{5}{c}{{\bfseries \tablename\ \thetable{} -- continued from previous page}} \\
	\hline N & Defining inequalities & Vertices&Int. Pt. & Neighbors\\ \hline
	\endhead\hline \multicolumn{5}{|c|}{{Continued on next page}} \\ \hline\endfoot\endlastfoot1&
	\(\begin{array}{rcl}x_{1}-x_{2}&\geq& 0\\
		-x_{1}+2x_{2}&\geq& 0\\
	\end{array}\)&\(\begin{array}{l}(1, 1): 2\checkmark\\(2, 1): 3\checkmark\end{array}\)&(3, 2): 5\checkmark& [1, 3]\\\hline
	2&
	\(\begin{array}{rcl}x_{1}&\geq& 0\\
		-x_{1}+x_{2}&\geq& 0\\
	\end{array}\)&\(\begin{array}{l}(0, 1): 1\checkmark\\(1, 1): 2\checkmark\end{array}\)&(1, 2): 2\checkmark& [4]\\\hline
	3&
	\(\begin{array}{rcl}x_{2}&\geq& 0\\
		x_{1}-2x_{2}&\geq& 0\\
	\end{array}\)&\(\begin{array}{l}(1, 0): 1\checkmark\\(2, 1): 3\checkmark\end{array}\)&(3, 1): 3\checkmark& [4]\\\hline
\end{longtable}

\allowdisplaybreaks\begin{align*}&~~~
	\frac{1}{(1-x_{1} ) (1-x_{2} ) (1-x_{1} x_{2} ) (1-x_{1}^2x_{2} ) }\\=&~~~
	\displaystyle \frac{-x_{1}^{-1}x_{2}^{-2}}{(1-x_{1} )^3 (1-x_{1} x_{2} ) }\\&
	+\displaystyle \frac{x_{2}^{-2}}{(1-x_{1} )^2(1-x_{1}^2) (1-x_{2} ) }\\&
	+\displaystyle \frac{-x_{2}^{-2}}{(1-x_{1} )^2(1-x_{1}^2) (1-x_{1}^2x_{2} ) }\\&
	+\displaystyle \frac{x_{1}^{-1}x_{2}^{-2}}{(1-x_{1} )^3 (1-x_{1}^2x_{2} ) }\\=&~~~
	\displaystyle \frac{-x_{1}^{-1}x_{2}^{-2}}{(1-x_{1} )^3 (1-x_{1} x_{2} ) }\\&
	+\displaystyle \frac{x_{1}^{-1}x_{2}^{-2}}{(1-x_{1} )^3 (1-x_{1}^2x_{2} ) }\\&
	+\displaystyle \frac{-x_{1}^2x_{2}^{-2}-2x_{1} x_{2}^{-2}-x_{2}^{-2}}{(1-x_{1}^2)^3 (1-x_{1}^2x_{2} ) }\\&
	+\displaystyle \frac{x_{1}^2x_{2}^{-2}+2x_{1} x_{2}^{-2}+x_{2}^{-2}}{(1-x_{1}^2)^3 (1-x_{2} ) }\\=&
	-x_{1}^{-1}x_{2}^{-2}\cdot \frac{1}{2}\left(x_{1} \partial_{1} -x_{2} \partial_{2} \right)^2\cdot\frac{1}{(1-x_{1} ) (1-x_{1} x_{2} ) }\\&
	\displaystyle +x_{1}^{-1}x_{2}^{-2}\cdot \frac{1}{2}\left(x_{1} \partial_{1} -2x_{2} \partial_{2} \right)^2\cdot\frac{1}{(1-x_{1} ) (1-x_{1}^2x_{2} ) }\\&
	\displaystyle +\left(-x_{1}^2x_{2}^{-2}-2x_{1} x_{2}^{-2}-x_{2}^{-2}\right)\cdot \frac{1}{2}\left(\frac{x_{1} \partial_{1} }{2}-x_{2} \partial_{2} \right)^2\cdot\frac{1}{(1-x_{1}^2) (1-x_{1}^2x_{2} ) }\\&
	\displaystyle +\left(x_{1}^2x_{2}^{-2}+2x_{1} x_{2}^{-2}+x_{2}^{-2}\right)\cdot \frac{1}{2}\left(\frac{x_{1} \partial_{1} }{2}\right)^2\cdot\frac{1}{(1-x_{1}^2) (1-x_{2} ) }\end{align*}\begin{table}[h!]
	\begin{center}
		{\tiny
			\psset{xunit=4cm, yunit=4cm}
			\begin{pspicture}(0,0)(1,1)
				\psline[linecolor=black, linestyle=dashed](0,0)(1,0)
				\rput[lt](1,0){(1, 0)}
				\pscircle*[linecolor=black](1,0){0.05}
				\psline[linecolor=black, linestyle=dashed](0,0)(0,1)
				\rput[rb](0,1){(0, 1)}
				\pscircle*[linecolor=black](0,1){0.05}
				\psline[linecolor=blue, linestyle=dashed](0,0)(0.66666667,0.33333333)
				\psline[linecolor=blue, linestyle=dashed](0,0)(0.6,0.4)
				\psline[linecolor=blue, linestyle=dashed](0,0)(0,1)
				\psline[linecolor=blue, linestyle=dashed](0,0)(0.5,0.5)
				\psline[linecolor=blue, linestyle=dashed](0,0)(1,0)
				\psline[linecolor=blue, linestyle=dashed](0,0)(0.75,0.25)
				\rput[lt](0.70833333,0.29166667){1}
				\rput[lt](0.55,0.45){2}
				\rput[rb](0.25,0.75){3}
				\rput[lt](0.875,0.125){4}
				\rput[lt](0.63333333,0.36666667){5}
				\pscircle*[linecolor=blue](1,0){0.0375}
				\pscircle*[linecolor=blue](0,1){0.0375}
				\pscircle*[linecolor=blue](0.5,0.5){0.0375}
				\pscircle*[linecolor=blue](0.66666667,0.33333333){0.0375}
				\pscircle*[linecolor=blue](0.75,0.25){0.0375}
				\pscircle*[linecolor=blue](0.6,0.4){0.0375}
		\end{pspicture}}
		\caption{$5$ combinatorial chambers of \(G_2\):(1, 0), (0, 1), (1, 1), (2, 1), (3, 1), (3, 2)}
	\end{center}
\end{table}
\begin{longtable}{|c|c|c|}\caption{\footnotesize V.p.f. of \(G_2\):(1, 0), (0, 1), (1, 1), (2, 1), (3, 1), (3, 2)}\\\hline N & Polynomial/Lattice & Shift(s)\\ \hline
	\endfirsthead\multicolumn{3}{c}{{\bfseries \tablename\ \thetable{} -- continued from previous page}} \\
	\hline  N &  Polynomial/Lattice & Shift(s) \\ \hline
	\endhead\hline \multicolumn{3}{|c|}{{Continued on next page}} \\ \hline\endfoot\endlastfoot\multirow{13}{*}{1}&\(\Lambda=\langle((6, 0),(0, 2))\rangle\)&\\\hline&\(\begin{array}{l}\begin{array}{l}-\frac{x_{1}^4}{432}+\frac{x_{1}^3x_{2} }{36}-\frac{x_{1}^2x_{2}^2}{8}+\frac{x_{1} x_{2}^3}{4}-\frac{x_{2}^4}{6}+\frac{x_{1}^3}{27}-\frac{x_{1}^2x_{2} }{3} \\ +x_{1} x_{2}^2-\frac{3}{4}x_{2}^3-\frac{7}{36}x_{1}^2+\frac{7}{6}x_{1} x_{2} -\frac{17}{24}x_{2}^2\end{array}\\
		\color{red}
		+\frac{x_{1} }{3}+\frac{3}{4}x_{2} +1\color{black}
	\end{array}
	\)&\((0, 0)\)\\
	\cline{2-3}&\(\begin{array}{l}\begin{array}{l}-\frac{x_{1}^4}{432}+\frac{x_{1}^3x_{2} }{36}-\frac{x_{1}^2x_{2}^2}{8}+\frac{x_{1} x_{2}^3}{4}-\frac{x_{2}^4}{6}+\frac{x_{1}^3}{27}-\frac{x_{1}^2x_{2} }{3} \\ +x_{1} x_{2}^2-\frac{3}{4}x_{2}^3-\frac{7}{36}x_{1}^2+\frac{7}{6}x_{1} x_{2} -\frac{17}{24}x_{2}^2\end{array}\\
		\color{red}
		+\frac{10}{27}x_{1} +\frac{23}{36}x_{2} +\frac{341}{432}\color{black}
	\end{array}
	\)&\(\begin{array}{l}
		(1, 0)\\
		(1, 1)\end{array}
	\)\\
	\cline{2-3}&\(\begin{array}{l}\begin{array}{l}-\frac{x_{1}^4}{432}+\frac{x_{1}^3x_{2} }{36}-\frac{x_{1}^2x_{2}^2}{8}+\frac{x_{1} x_{2}^3}{4}-\frac{x_{2}^4}{6}+\frac{x_{1}^3}{27}-\frac{x_{1}^2x_{2} }{3} \\ +x_{1} x_{2}^2-\frac{3}{4}x_{2}^3-\frac{7}{36}x_{1}^2+\frac{7}{6}x_{1} x_{2} -\frac{17}{24}x_{2}^2\end{array}\\
		\color{red}
		+\frac{x_{1} }{3}+\frac{3}{4}x_{2} +\frac{7}{8}\color{black}
	\end{array}
	\)&\((0, 1)\)\\
	\cline{2-3}&\(\begin{array}{l}\begin{array}{l}-\frac{x_{1}^4}{432}+\frac{x_{1}^3x_{2} }{36}-\frac{x_{1}^2x_{2}^2}{8}+\frac{x_{1} x_{2}^3}{4}-\frac{x_{2}^4}{6}+\frac{x_{1}^3}{27}-\frac{x_{1}^2x_{2} }{3} \\ +x_{1} x_{2}^2-\frac{3}{4}x_{2}^3-\frac{7}{36}x_{1}^2+\frac{7}{6}x_{1} x_{2} -\frac{17}{24}x_{2}^2\end{array}\\
		\color{red}
		+\frac{11}{27}x_{1} +\frac{19}{36}x_{2} +\frac{19}{27}\color{black}
	\end{array}
	\)&\((2, 0)\)\\
	\cline{2-3}&\(\begin{array}{l}\begin{array}{l}-\frac{x_{1}^4}{432}+\frac{x_{1}^3x_{2} }{36}-\frac{x_{1}^2x_{2}^2}{8}+\frac{x_{1} x_{2}^3}{4}-\frac{x_{2}^4}{6}+\frac{x_{1}^3}{27}-\frac{x_{1}^2x_{2} }{3} \\ +x_{1} x_{2}^2-\frac{3}{4}x_{2}^3-\frac{7}{36}x_{1}^2+\frac{7}{6}x_{1} x_{2} -\frac{17}{24}x_{2}^2\end{array}\\
		\color{red}
		+\frac{x_{1} }{3}+\frac{3}{4}x_{2} +\frac{15}{16}\color{black}
	\end{array}
	\)&\(\begin{array}{l}
		(3, 0)\\
		(3, 1)\end{array}
	\)\\
	\cline{2-3}&\(\begin{array}{l}\begin{array}{l}-\frac{x_{1}^4}{432}+\frac{x_{1}^3x_{2} }{36}-\frac{x_{1}^2x_{2}^2}{8}+\frac{x_{1} x_{2}^3}{4}-\frac{x_{2}^4}{6}+\frac{x_{1}^3}{27}-\frac{x_{1}^2x_{2} }{3} \\ +x_{1} x_{2}^2-\frac{3}{4}x_{2}^3-\frac{7}{36}x_{1}^2+\frac{7}{6}x_{1} x_{2} -\frac{17}{24}x_{2}^2\end{array}\\
		\color{red}
		+\frac{11}{27}x_{1} +\frac{19}{36}x_{2} +\frac{125}{216}\color{black}
	\end{array}
	\)&\((2, 1)\)\\
	\cline{2-3}&\(\begin{array}{l}\begin{array}{l}-\frac{x_{1}^4}{432}+\frac{x_{1}^3x_{2} }{36}-\frac{x_{1}^2x_{2}^2}{8}+\frac{x_{1} x_{2}^3}{4}-\frac{x_{2}^4}{6}+\frac{x_{1}^3}{27}-\frac{x_{1}^2x_{2} }{3} \\ +x_{1} x_{2}^2-\frac{3}{4}x_{2}^3-\frac{7}{36}x_{1}^2+\frac{7}{6}x_{1} x_{2} -\frac{17}{24}x_{2}^2\end{array}\\
		\color{red}
		+\frac{10}{27}x_{1} +\frac{23}{36}x_{2} +\frac{23}{27}\color{black}
	\end{array}
	\)&\((4, 0)\)\\
	\cline{2-3}&\(\begin{array}{l}\begin{array}{l}-\frac{x_{1}^4}{432}+\frac{x_{1}^3x_{2} }{36}-\frac{x_{1}^2x_{2}^2}{8}+\frac{x_{1} x_{2}^3}{4}-\frac{x_{2}^4}{6}+\frac{x_{1}^3}{27}-\frac{x_{1}^2x_{2} }{3} \\ +x_{1} x_{2}^2-\frac{3}{4}x_{2}^3-\frac{7}{36}x_{1}^2+\frac{7}{6}x_{1} x_{2} -\frac{17}{24}x_{2}^2\end{array}\\
		\color{red}
		+\frac{11}{27}x_{1} +\frac{19}{36}x_{2} +\frac{277}{432}\color{black}
	\end{array}
	\)&\(\begin{array}{l}
		(5, 0)\\
		(5, 1)\end{array}
	\)\\
	\cline{2-3}&\(\begin{array}{l}\begin{array}{l}-\frac{x_{1}^4}{432}+\frac{x_{1}^3x_{2} }{36}-\frac{x_{1}^2x_{2}^2}{8}+\frac{x_{1} x_{2}^3}{4}-\frac{x_{2}^4}{6}+\frac{x_{1}^3}{27}-\frac{x_{1}^2x_{2} }{3} \\ +x_{1} x_{2}^2-\frac{3}{4}x_{2}^3-\frac{7}{36}x_{1}^2+\frac{7}{6}x_{1} x_{2} -\frac{17}{24}x_{2}^2\end{array}\\
		\color{red}
		+\frac{10}{27}x_{1} +\frac{23}{36}x_{2} +\frac{157}{216}\color{black}
	\end{array}
	\)&\((4, 1)\)\\
	\cline{2-3}\hline 
	\multirow{13}{*}{2}&\(\Lambda=\langle((6, 0),(0, 2))\rangle\)&\\\hline&\(\begin{array}{l}\begin{array}{l}-\frac{x_{1}^4}{54}+\frac{x_{1}^3x_{2} }{12}-\frac{x_{1}^2x_{2}^2}{8}+\frac{x_{1} x_{2}^3}{12}-\frac{x_{2}^4}{48}-\frac{13}{108}x_{1}^3+\frac{x_{1}^2x_{2} }{2} \\ -\frac{x_{1} x_{2}^2}{2}+\frac{x_{2}^3}{6}-\frac{7}{72}x_{1}^2+\frac{5}{6}x_{1} x_{2} -\frac{5}{12}x_{2}^2+\frac{x_{2} }{3}\end{array}\\
		\color{red}
		+\frac{7}{12}x_{1} +1\color{black}
	\end{array}
	\)&\((0, 0)\)\\
	\cline{2-3}&\(\begin{array}{l}\begin{array}{l}-\frac{x_{1}^4}{54}+\frac{x_{1}^3x_{2} }{12}-\frac{x_{1}^2x_{2}^2}{8}+\frac{x_{1} x_{2}^3}{12}-\frac{x_{2}^4}{48}-\frac{13}{108}x_{1}^3+\frac{x_{1}^2x_{2} }{2} \\ -\frac{x_{1} x_{2}^2}{2}+\frac{x_{2}^3}{6}-\frac{7}{72}x_{1}^2+\frac{5}{6}x_{1} x_{2} -\frac{5}{12}x_{2}^2+\frac{x_{2} }{3}\end{array}\\
		\color{red}
		+\frac{59}{108}x_{1} +\frac{149}{216}\color{black}
	\end{array}
	\)&\((1, 0)\)\\
	\cline{2-3}&\(\begin{array}{l}\begin{array}{l}-\frac{x_{1}^4}{54}+\frac{x_{1}^3x_{2} }{12}-\frac{x_{1}^2x_{2}^2}{8}+\frac{x_{1} x_{2}^3}{12}-\frac{x_{2}^4}{48}-\frac{13}{108}x_{1}^3+\frac{x_{1}^2x_{2} }{2} \\ -\frac{x_{1} x_{2}^2}{2}+\frac{x_{2}^3}{6}-\frac{7}{72}x_{1}^2+\frac{5}{6}x_{1} x_{2} -\frac{5}{12}x_{2}^2+\frac{x_{2} }{3}\end{array}\\
		\color{red}
		+\frac{7}{12}x_{1} +\frac{15}{16}\color{black}
	\end{array}
	\)&\(\begin{array}{l}
		(0, 1)\\
		(3, 1)\end{array}
	\)\\
	\cline{2-3}&\(\begin{array}{l}\begin{array}{l}-\frac{x_{1}^4}{54}+\frac{x_{1}^3x_{2} }{12}-\frac{x_{1}^2x_{2}^2}{8}+\frac{x_{1} x_{2}^3}{12}-\frac{x_{2}^4}{48}-\frac{13}{108}x_{1}^3+\frac{x_{1}^2x_{2} }{2} \\ -\frac{x_{1} x_{2}^2}{2}+\frac{x_{2}^3}{6}-\frac{7}{72}x_{1}^2+\frac{5}{6}x_{1} x_{2} -\frac{5}{12}x_{2}^2+\frac{x_{2} }{3}\end{array}\\
		\color{red}
		+\frac{55}{108}x_{1} +\frac{17}{27}\color{black}
	\end{array}
	\)&\((2, 0)\)\\
	\cline{2-3}&\(\begin{array}{l}\begin{array}{l}-\frac{x_{1}^4}{54}+\frac{x_{1}^3x_{2} }{12}-\frac{x_{1}^2x_{2}^2}{8}+\frac{x_{1} x_{2}^3}{12}-\frac{x_{2}^4}{48}-\frac{13}{108}x_{1}^3+\frac{x_{1}^2x_{2} }{2} \\ -\frac{x_{1} x_{2}^2}{2}+\frac{x_{2}^3}{6}-\frac{7}{72}x_{1}^2+\frac{5}{6}x_{1} x_{2} -\frac{5}{12}x_{2}^2+\frac{x_{2} }{3}\end{array}\\
		\color{red}
		+\frac{59}{108}x_{1} +\frac{325}{432}\color{black}
	\end{array}
	\)&\(\begin{array}{l}
		(1, 1)\\
		(4, 1)\end{array}
	\)\\
	\cline{2-3}&\(\begin{array}{l}\begin{array}{l}-\frac{x_{1}^4}{54}+\frac{x_{1}^3x_{2} }{12}-\frac{x_{1}^2x_{2}^2}{8}+\frac{x_{1} x_{2}^3}{12}-\frac{x_{2}^4}{48}-\frac{13}{108}x_{1}^3+\frac{x_{1}^2x_{2} }{2} \\ -\frac{x_{1} x_{2}^2}{2}+\frac{x_{2}^3}{6}-\frac{7}{72}x_{1}^2+\frac{5}{6}x_{1} x_{2} -\frac{5}{12}x_{2}^2+\frac{x_{2} }{3}\end{array}\\
		\color{red}
		+\frac{7}{12}x_{1} +\frac{7}{8}\color{black}
	\end{array}
	\)&\((3, 0)\)\\
	\cline{2-3}&\(\begin{array}{l}\begin{array}{l}-\frac{x_{1}^4}{54}+\frac{x_{1}^3x_{2} }{12}-\frac{x_{1}^2x_{2}^2}{8}+\frac{x_{1} x_{2}^3}{12}-\frac{x_{2}^4}{48}-\frac{13}{108}x_{1}^3+\frac{x_{1}^2x_{2} }{2} \\ -\frac{x_{1} x_{2}^2}{2}+\frac{x_{2}^3}{6}-\frac{7}{72}x_{1}^2+\frac{5}{6}x_{1} x_{2} -\frac{5}{12}x_{2}^2+\frac{x_{2} }{3}\end{array}\\
		\color{red}
		+\frac{55}{108}x_{1} +\frac{245}{432}\color{black}
	\end{array}
	\)&\(\begin{array}{l}
		(2, 1)\\
		(5, 1)\end{array}
	\)\\
	\cline{2-3}&\(\begin{array}{l}\begin{array}{l}-\frac{x_{1}^4}{54}+\frac{x_{1}^3x_{2} }{12}-\frac{x_{1}^2x_{2}^2}{8}+\frac{x_{1} x_{2}^3}{12}-\frac{x_{2}^4}{48}-\frac{13}{108}x_{1}^3+\frac{x_{1}^2x_{2} }{2} \\ -\frac{x_{1} x_{2}^2}{2}+\frac{x_{2}^3}{6}-\frac{7}{72}x_{1}^2+\frac{5}{6}x_{1} x_{2} -\frac{5}{12}x_{2}^2+\frac{x_{2} }{3}\end{array}\\
		\color{red}
		+\frac{59}{108}x_{1} +\frac{22}{27}\color{black}
	\end{array}
	\)&\((4, 0)\)\\
	\cline{2-3}&\(\begin{array}{l}\begin{array}{l}-\frac{x_{1}^4}{54}+\frac{x_{1}^3x_{2} }{12}-\frac{x_{1}^2x_{2}^2}{8}+\frac{x_{1} x_{2}^3}{12}-\frac{x_{2}^4}{48}-\frac{13}{108}x_{1}^3+\frac{x_{1}^2x_{2} }{2} \\ -\frac{x_{1} x_{2}^2}{2}+\frac{x_{2}^3}{6}-\frac{7}{72}x_{1}^2+\frac{5}{6}x_{1} x_{2} -\frac{5}{12}x_{2}^2+\frac{x_{2} }{3}\end{array}\\
		\color{red}
		+\frac{55}{108}x_{1} +\frac{109}{216}\color{black}
	\end{array}
	\)&\((5, 0)\)\\
	\cline{2-3}\hline 
	\multirow{7}{*}{3}&\(\Lambda=\langle((6, 0),(0, 1))\rangle\)&\\\hline&\(\begin{array}{l}\frac{x_{1}^4}{432}+\frac{5}{108}x_{1}^3+\frac{23}{72}x_{1}^2\\
		\color{red}
		+\frac{11}{12}x_{1} +1\color{black}
	\end{array}
	\)&\((0, 0)\)\\
	\cline{2-3}&\(\begin{array}{l}\frac{x_{1}^4}{432}+\frac{5}{108}x_{1}^3+\frac{23}{72}x_{1}^2\\
		\color{red}
		+\frac{95}{108}x_{1} +\frac{325}{432}\color{black}
	\end{array}
	\)&\((1, 0)\)\\
	\cline{2-3}&\(\begin{array}{l}\frac{x_{1}^4}{432}+\frac{5}{108}x_{1}^3+\frac{23}{72}x_{1}^2\\
		\color{red}
		+\frac{91}{108}x_{1} +\frac{17}{27}\color{black}
	\end{array}
	\)&\((2, 0)\)\\
	\cline{2-3}&\(\begin{array}{l}\frac{x_{1}^4}{432}+\frac{5}{108}x_{1}^3+\frac{23}{72}x_{1}^2\\
		\color{red}
		+\frac{11}{12}x_{1} +\frac{15}{16}\color{black}
	\end{array}
	\)&\((3, 0)\)\\
	\cline{2-3}&\(\begin{array}{l}\frac{x_{1}^4}{432}+\frac{5}{108}x_{1}^3+\frac{23}{72}x_{1}^2\\
		\color{red}
		+\frac{95}{108}x_{1} +\frac{22}{27}\color{black}
	\end{array}
	\)&\((4, 0)\)\\
	\cline{2-3}&\(\begin{array}{l}\frac{x_{1}^4}{432}+\frac{5}{108}x_{1}^3+\frac{23}{72}x_{1}^2\\
		\color{red}
		+\frac{91}{108}x_{1} +\frac{245}{432}\color{black}
	\end{array}
	\)&\((5, 0)\)\\
	\hline\hline 
	\multirow{3}{*}{4}&\(\Lambda=\langle((1, 0),(0, 2))\rangle\)&\\\hline&\(\begin{array}{l}\frac{x_{2}^4}{48}+\frac{x_{2}^3}{4}+\frac{25}{24}x_{2}^2+\frac{7}{4}x_{2} \\
		\color{red}
		+1\color{black}
	\end{array}
	\)&\((0, 0)\)\\
	\cline{2-3}&\(\begin{array}{l}\frac{x_{2}^4}{48}+\frac{x_{2}^3}{4}+\frac{25}{24}x_{2}^2+\frac{7}{4}x_{2} \\
		\color{red}
		+\frac{15}{16}\color{black}
	\end{array}
	\)&\((0, 1)\)\\
	\hline\hline 
	\multirow{13}{*}{5}&\(\Lambda=\langle((6, 0),(0, 2))\rangle\)&\\\hline&\(\begin{array}{l}\begin{array}{l}\frac{x_{1}^4}{54}-\frac{5}{36}x_{1}^3x_{2} +\frac{3}{8}x_{1}^2x_{2}^2-\frac{5}{12}x_{1} x_{2}^3+\frac{x_{2}^4}{6}-\frac{5}{108}x_{1}^3+\frac{x_{1}^2x_{2} }{6} \\ -\frac{x_{2}^3}{12}-\frac{11}{72}x_{1}^2+x_{1} x_{2} -\frac{13}{24}x_{2}^2\end{array}\\
		\color{red}
		+\frac{5}{12}x_{1} +\frac{7}{12}x_{2} +1\color{black}
	\end{array}
	\)&\((0, 0)\)\\
	\cline{2-3}&\(\begin{array}{l}\begin{array}{l}\frac{x_{1}^4}{54}-\frac{5}{36}x_{1}^3x_{2} +\frac{3}{8}x_{1}^2x_{2}^2-\frac{5}{12}x_{1} x_{2}^3+\frac{x_{2}^4}{6}-\frac{5}{108}x_{1}^3+\frac{x_{1}^2x_{2} }{6} \\ -\frac{x_{2}^3}{12}-\frac{11}{72}x_{1}^2+x_{1} x_{2} -\frac{13}{24}x_{2}^2\end{array}\\
		\color{red}
		+\frac{49}{108}x_{1} +\frac{17}{36}x_{2} +\frac{157}{216}\color{black}
	\end{array}
	\)&\(\begin{array}{l}
		(1, 0)\\
		(1, 1)\\
		(4, 1)\end{array}
	\)\\
	\cline{2-3}&\(\begin{array}{l}\begin{array}{l}\frac{x_{1}^4}{54}-\frac{5}{36}x_{1}^3x_{2} +\frac{3}{8}x_{1}^2x_{2}^2-\frac{5}{12}x_{1} x_{2}^3+\frac{x_{2}^4}{6}-\frac{5}{108}x_{1}^3+\frac{x_{1}^2x_{2} }{6} \\ -\frac{x_{2}^3}{12}-\frac{11}{72}x_{1}^2+x_{1} x_{2} -\frac{13}{24}x_{2}^2\end{array}\\
		\color{red}
		+\frac{5}{12}x_{1} +\frac{7}{12}x_{2} +\frac{7}{8}\color{black}
	\end{array}
	\)&\(\begin{array}{l}
		(0, 1)\\
		(3, 0)\\
		(3, 1)\end{array}
	\)\\
	\cline{2-3}&\(\begin{array}{l}\begin{array}{l}\frac{x_{1}^4}{54}-\frac{5}{36}x_{1}^3x_{2} +\frac{3}{8}x_{1}^2x_{2}^2-\frac{5}{12}x_{1} x_{2}^3+\frac{x_{2}^4}{6}-\frac{5}{108}x_{1}^3+\frac{x_{1}^2x_{2} }{6} \\ -\frac{x_{2}^3}{12}-\frac{11}{72}x_{1}^2+x_{1} x_{2} -\frac{13}{24}x_{2}^2\end{array}\\
		\color{red}
		+\frac{53}{108}x_{1} +\frac{13}{36}x_{2} +\frac{19}{27}\color{black}
	\end{array}
	\)&\((2, 0)\)\\
	\cline{2-3}&\(\begin{array}{l}\begin{array}{l}\frac{x_{1}^4}{54}-\frac{5}{36}x_{1}^3x_{2} +\frac{3}{8}x_{1}^2x_{2}^2-\frac{5}{12}x_{1} x_{2}^3+\frac{x_{2}^4}{6}-\frac{5}{108}x_{1}^3+\frac{x_{1}^2x_{2} }{6} \\ -\frac{x_{2}^3}{12}-\frac{11}{72}x_{1}^2+x_{1} x_{2} -\frac{13}{24}x_{2}^2\end{array}\\
		\color{red}
		+\frac{53}{108}x_{1} +\frac{13}{36}x_{2} +\frac{125}{216}\color{black}
	\end{array}
	\)&\(\begin{array}{l}
		(2, 1)\\
		(5, 0)\\
		(5, 1)\end{array}
	\)\\
	\cline{2-3}&\(\begin{array}{l}\begin{array}{l}\frac{x_{1}^4}{54}-\frac{5}{36}x_{1}^3x_{2} +\frac{3}{8}x_{1}^2x_{2}^2-\frac{5}{12}x_{1} x_{2}^3+\frac{x_{2}^4}{6}-\frac{5}{108}x_{1}^3+\frac{x_{1}^2x_{2} }{6} \\ -\frac{x_{2}^3}{12}-\frac{11}{72}x_{1}^2+x_{1} x_{2} -\frac{13}{24}x_{2}^2\end{array}\\
		\color{red}
		+\frac{49}{108}x_{1} +\frac{17}{36}x_{2} +\frac{23}{27}\color{black}
	\end{array}
	\)&\((4, 0)\)\\
	\cline{2-3}\end{longtable}

\begin{longtable}{|ccccc|}\caption{\footnotesize V.p.f. of \(G_2\):(1, 0), (0, 1), (1, 1), (2, 1), (3, 1), (3, 2)}\\\hline N & Defining inequalities & Vertices& Int. Pt.& Neighbors \\ \hline
	\endfirsthead\multicolumn{5}{c}{{\bfseries \tablename\ \thetable{} -- continued from previous page}} \\
	\hline N & Defining inequalities & Vertices&Int. Pt. & Neighbors\\ \hline
	\endhead\hline \multicolumn{5}{|c|}{{Continued on next page}} \\ \hline\endfoot\endlastfoot1&
	\(\begin{array}{rcl}x_{1}-2x_{2}&\geq& 0\\
		-x_{1}+3x_{2}&\geq& 0\\
	\end{array}\)&\(\begin{array}{l}(2, 1): 3\checkmark\\(3, 1): 4\checkmark\end{array}\)&(5, 2): 10\checkmark& [8, 5]\\\hline
	2&
	\(\begin{array}{rcl}x_{1}-x_{2}&\geq& 0\\
		-2x_{1}+3x_{2}&\geq& 0\\
	\end{array}\)&\(\begin{array}{l}(1, 1): 2\checkmark\\(3, 2): 7\checkmark\end{array}\)&(4, 3): 12\checkmark& [1, 8]\\\hline
	3&
	\(\begin{array}{rcl}x_{1}&\geq& 0\\
		-x_{1}+x_{2}&\geq& 0\\
	\end{array}\)&\(\begin{array}{l}(0, 1): 1\checkmark\\(1, 1): 2\checkmark\end{array}\)&(1, 2): 2\checkmark& [7]\\\hline
	4&
	\(\begin{array}{rcl}x_{2}&\geq& 0\\
		x_{1}-3x_{2}&\geq& 0\\
	\end{array}\)&\(\begin{array}{l}(1, 0): 1\checkmark\\(3, 1): 4\checkmark\end{array}\)&(4, 1): 4\checkmark& [6]\\\hline
	5&
	\(\begin{array}{rcl}-x_{1}+2x_{2}&\geq& 0\\
		2x_{1}-3x_{2}&\geq& 0\\
	\end{array}\)&\(\begin{array}{l}(2, 1): 3\checkmark\\(3, 2): 7\checkmark\end{array}\)&(5, 3): 16\checkmark& [6, 7]\\\hline
\end{longtable}

\allowdisplaybreaks\begin{align*}&~~~
	\frac{1}{(1-x_{1} ) (1-x_{2} ) (1-x_{1} x_{2} ) (1-x_{1}^2x_{2} ) (1-x_{1}^3x_{2} ) (1-x_{1}^3x_{2}^2) }\\=&~~~
	\displaystyle \frac{-x_{1}^{-4}x_{2}^{-5}}{(1-x_{1} )^4(1-x_{1}^2) (1-x_{1} x_{2} ) }\\&
	+\displaystyle \frac{x_{2}^{-5}}{(1-x_{1} )^2(1-x_{1}^2)(1-x_{1}^3)^2 (1-x_{2} ) }\\&
	+\displaystyle \frac{-x_{1}^3x_{2}^{-2}-x_{2}^{-3}-x_{2}^{-4}-x_{2}^{-5}}{(1-x_{1} )^2(1-x_{1}^2)(1-x_{1}^3)^2 (1-x_{1}^3x_{2}^2) }\\&
	+\displaystyle \frac{-x_{2}^{-2}-x_{1}^{-2}x_{2}^{-3}+x_{1}^{-3}x_{2}^{-4}+x_{1}^{-4}x_{2}^{-5}}{(1-x_{1} )^4(1-x_{1}^2) (1-x_{1}^3x_{2}^2) }\\&
	+\displaystyle \frac{x_{1}^3x_{2}^{-3}}{(1-x_{1} )^2(1-x_{1}^2)(1-x_{1}^3)^2 (1-x_{1}^3x_{2} ) }\\&
	+\displaystyle \frac{-2x_{1} x_{2}^{-3}}{(1-x_{1} )^3(1-x_{1}^2)(1-x_{1}^3) (1-x_{1}^3x_{2} ) }\\&
	+\displaystyle \frac{2x_{1} x_{2}^{-2}+2x_{1}^{-2}x_{2}^{-3}}{(1-x_{1} )^3(1-x_{1}^2)(1-x_{1}^3) (1-x_{1}^3x_{2}^2) }\\&
	+\displaystyle \frac{x_{1}^{-1}x_{2}^{-3}}{(1-x_{1} )^4(1-x_{1}^2) (1-x_{1}^2x_{2} ) }\\&
	+\displaystyle \frac{-x_{1}^{-2}x_{2}^{-3}}{(1-x_{1} )^5 (1-x_{1}^2x_{2} ) }\\&
	+\displaystyle \frac{x_{1}^{-1}x_{2}^{-2}+x_{1}^{-3}x_{2}^{-3}}{(1-x_{1} )^5 (1-x_{1}^3x_{2}^2) }\\&
	+\displaystyle \frac{x_{2}^{-3}}{(1-x_{1} )^4(1-x_{1}^3) (1-x_{1}^3x_{2} ) }\\&
	+\displaystyle \frac{-x_{2}^{-2}-x_{1}^{-3}x_{2}^{-3}}{(1-x_{1} )^4(1-x_{1}^3) (1-x_{1}^3x_{2}^2) }\\=&~~~
	\displaystyle \frac{-x_{2}^{-5}-4x_{1}^{-1}x_{2}^{-5}-6x_{1}^{-2}x_{2}^{-5}-4x_{1}^{-3}x_{2}^{-5}-x_{1}^{-4}x_{2}^{-5}}{(1-x_{1}^2)^5 (1-x_{1} x_{2} ) }\\&
	+\displaystyle \frac{~~~~\begin{array}{l}-x_{1}^8x_{2}^{-2}-4x_{1}^7x_{2}^{-2}-10x_{1}^6x_{2}^{-2}-16x_{1}^5x_{2}^{-2}-x_{1}^5x_{2}^{-3}-19x_{1}^4x_{2}^{-2}-4x_{1}^4x_{2}^{-3} \\ \hline  -16x_{1}^3x_{2}^{-2}-10x_{1}^3x_{2}^{-3}-10x_{1}^2x_{2}^{-2}-16x_{1}^2x_{2}^{-3}-4x_{1} x_{2}^{-2}-19x_{1} x_{2}^{-3}-x_{2}^{-2} \\ \hline  -16x_{2}^{-3}-10x_{1}^{-1}x_{2}^{-3}-4x_{1}^{-2}x_{2}^{-3}-x_{1}^{-3}x_{2}^{-3}\end{array}~~~~}{(1-x_{1}^3)^5 (1-x_{1}^3x_{2}^2) }\\&
	+\displaystyle \frac{x_{1}^{-1}x_{2}^{-2}+x_{1}^{-3}x_{2}^{-3}}{(1-x_{1} )^5 (1-x_{1}^3x_{2}^2) }\\&
	+\displaystyle \frac{~~~~\begin{array}{l}x_{1}^8x_{2}^{-3}+4x_{1}^7x_{2}^{-3}+10x_{1}^6x_{2}^{-3}+16x_{1}^5x_{2}^{-3}+19x_{1}^4x_{2}^{-3}+16x_{1}^3x_{2}^{-3}+10x_{1}^2x_{2}^{-3} \\ \hline  +4x_{1} x_{2}^{-3}+x_{2}^{-3}\end{array}~~~~}{(1-x_{1}^3)^5 (1-x_{1}^3x_{2} ) }\\&
	+\displaystyle \frac{-x_{1}^{-2}x_{2}^{-3}}{(1-x_{1} )^5 (1-x_{1}^2x_{2} ) }\\&
	+\displaystyle \frac{x_{1}^3x_{2}^{-3}+4x_{1}^2x_{2}^{-3}+6x_{1} x_{2}^{-3}+4x_{2}^{-3}+x_{1}^{-1}x_{2}^{-3}}{(1-x_{1}^2)^5 (1-x_{1}^2x_{2} ) }\\&
	+\displaystyle \frac{~~~~\begin{array}{l}x_{1}^{23}x_{2}^{-2}+4x_{1}^{22}x_{2}^{-2}+10x_{1}^{21}x_{2}^{-2}+20x_{1}^{20}x_{2}^{-2}+x_{1}^{20}x_{2}^{-3}+37x_{1}^{19}x_{2}^{-2} \\ \hline  -x_{1}^{20}x_{2}^{-4}+4x_{1}^{19}x_{2}^{-3}+62x_{1}^{18}x_{2}^{-2}-x_{1}^{20}x_{2}^{-5}-2x_{1}^{19}x_{2}^{-4}+10x_{1}^{18}x_{2}^{-3} \\ \hline  +92x_{1}^{17}x_{2}^{-2}-2x_{1}^{19}x_{2}^{-5}-4x_{1}^{18}x_{2}^{-4}+20x_{1}^{17}x_{2}^{-3}+126x_{1}^{16}x_{2}^{-2}-4x_{1}^{18}x_{2}^{-5} \\ \hline  -8x_{1}^{17}x_{2}^{-4}+37x_{1}^{16}x_{2}^{-3}+161x_{1}^{15}x_{2}^{-2}-8x_{1}^{17}x_{2}^{-5}-13x_{1}^{16}x_{2}^{-4}+62x_{1}^{15}x_{2}^{-3} \\ \hline  +194x_{1}^{14}x_{2}^{-2}-13x_{1}^{16}x_{2}^{-5}-20x_{1}^{15}x_{2}^{-4}+92x_{1}^{14}x_{2}^{-3}+216x_{1}^{13}x_{2}^{-2}-20x_{1}^{15}x_{2}^{-5} \\ \hline  -26x_{1}^{14}x_{2}^{-4}+126x_{1}^{13}x_{2}^{-3}+224x_{1}^{12}x_{2}^{-2}-26x_{1}^{14}x_{2}^{-5}-34x_{1}^{13}x_{2}^{-4}+161x_{1}^{12}x_{2}^{-3} \\ \hline  +221x_{1}^{11}x_{2}^{-2}-34x_{1}^{13}x_{2}^{-5}-41x_{1}^{12}x_{2}^{-4}+194x_{1}^{11}x_{2}^{-3}+204x_{1}^{10}x_{2}^{-2}-41x_{1}^{12}x_{2}^{-5} \\ \hline  -44x_{1}^{11}x_{2}^{-4}+216x_{1}^{10}x_{2}^{-3}+176x_{1}^9x_{2}^{-2}-44x_{1}^{11}x_{2}^{-5}-46x_{1}^{10}x_{2}^{-4}+224x_{1}^9x_{2}^{-3} \\ \hline  +140x_{1}^8x_{2}^{-2}-46x_{1}^{10}x_{2}^{-5}-44x_{1}^9x_{2}^{-4}+221x_{1}^8x_{2}^{-3}+105x_{1}^7x_{2}^{-2}-44x_{1}^9x_{2}^{-5} \\ \hline  -41x_{1}^8x_{2}^{-4}+204x_{1}^7x_{2}^{-3}+74x_{1}^6x_{2}^{-2}-41x_{1}^8x_{2}^{-5}-34x_{1}^7x_{2}^{-4}+176x_{1}^6x_{2}^{-3} \\ \hline  +46x_{1}^5x_{2}^{-2}-34x_{1}^7x_{2}^{-5}-26x_{1}^6x_{2}^{-4}+140x_{1}^5x_{2}^{-3}+26x_{1}^4x_{2}^{-2}-26x_{1}^6x_{2}^{-5} \\ \hline  -20x_{1}^5x_{2}^{-4}+105x_{1}^4x_{2}^{-3}+13x_{1}^3x_{2}^{-2}-20x_{1}^5x_{2}^{-5}-13x_{1}^4x_{2}^{-4}+74x_{1}^3x_{2}^{-3} \\ \hline  +6x_{1}^2x_{2}^{-2}-13x_{1}^4x_{2}^{-5}-8x_{1}^3x_{2}^{-4}+46x_{1}^2x_{2}^{-3}+2x_{1} x_{2}^{-2}-8x_{1}^3x_{2}^{-5}-4x_{1}^2x_{2}^{-4} \\ \hline  +26x_{1} x_{2}^{-3}-4x_{1}^2x_{2}^{-5}-2x_{1} x_{2}^{-4}+13x_{2}^{-3}-2x_{1} x_{2}^{-5}-x_{2}^{-4}+6x_{1}^{-1}x_{2}^{-3}-x_{2}^{-5} \\ \hline  +2x_{1}^{-2}x_{2}^{-3}\end{array}~~~~}{(1-x_{1}^6)^5 (1-x_{1}^3x_{2}^2) }\\&
	+\displaystyle \frac{~~~~\begin{array}{l}-x_{1}^{23}x_{2}^{-3}-4x_{1}^{22}x_{2}^{-3}-10x_{1}^{21}x_{2}^{-3}-20x_{1}^{20}x_{2}^{-3}-37x_{1}^{19}x_{2}^{-3}-62x_{1}^{18}x_{2}^{-3} \\ \hline  -92x_{1}^{17}x_{2}^{-3}-126x_{1}^{16}x_{2}^{-3}-161x_{1}^{15}x_{2}^{-3}-194x_{1}^{14}x_{2}^{-3}-216x_{1}^{13}x_{2}^{-3}-224x_{1}^{12}x_{2}^{-3} \\ \hline  -221x_{1}^{11}x_{2}^{-3}-204x_{1}^{10}x_{2}^{-3}-176x_{1}^9x_{2}^{-3}-140x_{1}^8x_{2}^{-3}-105x_{1}^7x_{2}^{-3}-74x_{1}^6x_{2}^{-3} \\ \hline  -46x_{1}^5x_{2}^{-3}-26x_{1}^4x_{2}^{-3}-13x_{1}^3x_{2}^{-3}-6x_{1}^2x_{2}^{-3}-2x_{1} x_{2}^{-3}\end{array}~~~~}{(1-x_{1}^6)^5 (1-x_{1}^3x_{2} ) }\\&
	+\displaystyle \frac{~~~~\begin{array}{l}-x_{1}^4x_{2}^{-2}-4x_{1}^3x_{2}^{-2}-6x_{1}^2x_{2}^{-2}-x_{1}^2x_{2}^{-3}-4x_{1} x_{2}^{-2}-4x_{1} x_{2}^{-3}-x_{2}^{-2} \\ \hline  +x_{1} x_{2}^{-4}-6x_{2}^{-3}+4x_{2}^{-4}-4x_{1}^{-1}x_{2}^{-3}+x_{2}^{-5}+6x_{1}^{-1}x_{2}^{-4}-x_{1}^{-2}x_{2}^{-3}+4x_{1}^{-1}x_{2}^{-5} \\ \hline  +4x_{1}^{-2}x_{2}^{-4}+6x_{1}^{-2}x_{2}^{-5}+x_{1}^{-3}x_{2}^{-4}+4x_{1}^{-3}x_{2}^{-5}+x_{1}^{-4}x_{2}^{-5}\end{array}~~~~}{(1-x_{1}^2)^5 (1-x_{1}^3x_{2}^2) }\\&
	+\displaystyle \frac{~~~~\begin{array}{l}x_{1}^{20}x_{2}^{-5}+2x_{1}^{19}x_{2}^{-5}+4x_{1}^{18}x_{2}^{-5}+8x_{1}^{17}x_{2}^{-5}+13x_{1}^{16}x_{2}^{-5}+20x_{1}^{15}x_{2}^{-5} \\ \hline  +26x_{1}^{14}x_{2}^{-5}+34x_{1}^{13}x_{2}^{-5}+41x_{1}^{12}x_{2}^{-5}+44x_{1}^{11}x_{2}^{-5}+46x_{1}^{10}x_{2}^{-5}+44x_{1}^9x_{2}^{-5} \\ \hline  +41x_{1}^8x_{2}^{-5}+34x_{1}^7x_{2}^{-5}+26x_{1}^6x_{2}^{-5}+20x_{1}^5x_{2}^{-5}+13x_{1}^4x_{2}^{-5}+8x_{1}^3x_{2}^{-5}+4x_{1}^2x_{2}^{-5} \\ \hline  +2x_{1} x_{2}^{-5}+x_{2}^{-5}\end{array}~~~~}{(1-x_{1}^6)^5 (1-x_{2} ) }\\=&
	\left(\begin{array}{l}-x_{2}^{-5}-4x_{1}^{-1}x_{2}^{-5}-6x_{1}^{-2}x_{2}^{-5} \\ -4x_{1}^{-3}x_{2}^{-5}-x_{1}^{-4}x_{2}^{-5}\end{array}\right)\begin{array}{l}\cdot \frac{1}{24}\left(\frac{x_{1} \partial_{1} }{2}-\frac{x_{2} \partial_{2} }{2}\right)^4\\
		~~\cdot\frac{1}{(1-x_{1}^2) (1-x_{1} x_{2} ) }\end{array}\\&
	\displaystyle +\left(\begin{array}{l}-x_{1}^8x_{2}^{-2}-4x_{1}^7x_{2}^{-2}-10x_{1}^6x_{2}^{-2} \\ -16x_{1}^5x_{2}^{-2}-x_{1}^5x_{2}^{-3}-19x_{1}^4x_{2}^{-2} \\ -4x_{1}^4x_{2}^{-3}-16x_{1}^3x_{2}^{-2}-10x_{1}^3x_{2}^{-3} \\ -10x_{1}^2x_{2}^{-2}-16x_{1}^2x_{2}^{-3}-4x_{1} x_{2}^{-2} \\ -19x_{1} x_{2}^{-3}-x_{2}^{-2}-16x_{2}^{-3}-10x_{1}^{-1}x_{2}^{-3} \\ -4x_{1}^{-2}x_{2}^{-3}-x_{1}^{-3}x_{2}^{-3}\end{array}\right)\begin{array}{l}\cdot \frac{1}{24}\left(\frac{x_{1} \partial_{1} }{3}-\frac{x_{2} \partial_{2} }{2}\right)^4\\
		~~\cdot\frac{1}{(1-x_{1}^3) (1-x_{1}^3x_{2}^2) }\end{array}\\&
	\displaystyle +\left(x_{1}^{-1}x_{2}^{-2}+x_{1}^{-3}x_{2}^{-3}\right)\cdot \frac{1}{24}\left(x_{1} \partial_{1} -\frac{3}{2}x_{2} \partial_{2} \right)^4\cdot\frac{1}{(1-x_{1} ) (1-x_{1}^3x_{2}^2) }\\&
	\displaystyle +\left(\begin{array}{l}x_{1}^8x_{2}^{-3}+4x_{1}^7x_{2}^{-3}+10x_{1}^6x_{2}^{-3} \\ +16x_{1}^5x_{2}^{-3}+19x_{1}^4x_{2}^{-3}+16x_{1}^3x_{2}^{-3} \\ +10x_{1}^2x_{2}^{-3}+4x_{1} x_{2}^{-3}+x_{2}^{-3}\end{array}\right)\begin{array}{l}\cdot \frac{1}{24}\left(\frac{x_{1} \partial_{1} }{3}-x_{2} \partial_{2} \right)^4\\
		~~\cdot\frac{1}{(1-x_{1}^3) (1-x_{1}^3x_{2} ) }\end{array}\\&
	\displaystyle -x_{1}^{-2}x_{2}^{-3}\cdot \frac{1}{24}\left(x_{1} \partial_{1} -2x_{2} \partial_{2} \right)^4\cdot\frac{1}{(1-x_{1} ) (1-x_{1}^2x_{2} ) }\\&
	\displaystyle +\left(\begin{array}{l}x_{1}^3x_{2}^{-3}+4x_{1}^2x_{2}^{-3}+6x_{1} x_{2}^{-3} \\ +4x_{2}^{-3}+x_{1}^{-1}x_{2}^{-3}\end{array}\right)\begin{array}{l}\cdot \frac{1}{24}\left(\frac{x_{1} \partial_{1} }{2}-x_{2} \partial_{2} \right)^4\\
		~~\cdot\frac{1}{(1-x_{1}^2) (1-x_{1}^2x_{2} ) }\end{array}\\&
	\displaystyle +\left(\begin{array}{l}x_{1}^{23}x_{2}^{-2}+4x_{1}^{22}x_{2}^{-2}+10x_{1}^{21}x_{2}^{-2} \\ +20x_{1}^{20}x_{2}^{-2}+x_{1}^{20}x_{2}^{-3}+37x_{1}^{19}x_{2}^{-2} \\ -x_{1}^{20}x_{2}^{-4}+4x_{1}^{19}x_{2}^{-3}+62x_{1}^{18}x_{2}^{-2} \\ -x_{1}^{20}x_{2}^{-5}-2x_{1}^{19}x_{2}^{-4}+10x_{1}^{18}x_{2}^{-3} \\ +92x_{1}^{17}x_{2}^{-2}-2x_{1}^{19}x_{2}^{-5}-4x_{1}^{18}x_{2}^{-4} \\ +20x_{1}^{17}x_{2}^{-3}+126x_{1}^{16}x_{2}^{-2}-4x_{1}^{18}x_{2}^{-5} \\ -8x_{1}^{17}x_{2}^{-4}+37x_{1}^{16}x_{2}^{-3}+161x_{1}^{15}x_{2}^{-2} \\ -8x_{1}^{17}x_{2}^{-5}-13x_{1}^{16}x_{2}^{-4}+62x_{1}^{15}x_{2}^{-3} \\ +194x_{1}^{14}x_{2}^{-2}-13x_{1}^{16}x_{2}^{-5}-20x_{1}^{15}x_{2}^{-4} \\ +92x_{1}^{14}x_{2}^{-3}+216x_{1}^{13}x_{2}^{-2}-20x_{1}^{15}x_{2}^{-5} \\ -26x_{1}^{14}x_{2}^{-4}+126x_{1}^{13}x_{2}^{-3}+224x_{1}^{12}x_{2}^{-2} \\ -26x_{1}^{14}x_{2}^{-5}-34x_{1}^{13}x_{2}^{-4}+161x_{1}^{12}x_{2}^{-3} \\ +221x_{1}^{11}x_{2}^{-2}-34x_{1}^{13}x_{2}^{-5}-41x_{1}^{12}x_{2}^{-4} \\ +194x_{1}^{11}x_{2}^{-3}+204x_{1}^{10}x_{2}^{-2}-41x_{1}^{12}x_{2}^{-5} \\ -44x_{1}^{11}x_{2}^{-4}+216x_{1}^{10}x_{2}^{-3}+176x_{1}^9x_{2}^{-2} \\ -44x_{1}^{11}x_{2}^{-5}-46x_{1}^{10}x_{2}^{-4}+224x_{1}^9x_{2}^{-3} \\ +140x_{1}^8x_{2}^{-2}-46x_{1}^{10}x_{2}^{-5}-44x_{1}^9x_{2}^{-4} \\ +221x_{1}^8x_{2}^{-3}+105x_{1}^7x_{2}^{-2}-44x_{1}^9x_{2}^{-5} \\ -41x_{1}^8x_{2}^{-4}+204x_{1}^7x_{2}^{-3}+74x_{1}^6x_{2}^{-2} \\ -41x_{1}^8x_{2}^{-5}-34x_{1}^7x_{2}^{-4}+176x_{1}^6x_{2}^{-3} \\ +46x_{1}^5x_{2}^{-2}-34x_{1}^7x_{2}^{-5}-26x_{1}^6x_{2}^{-4} \\ +140x_{1}^5x_{2}^{-3}+26x_{1}^4x_{2}^{-2}-26x_{1}^6x_{2}^{-5} \\ -20x_{1}^5x_{2}^{-4}+105x_{1}^4x_{2}^{-3}+13x_{1}^3x_{2}^{-2} \\ -20x_{1}^5x_{2}^{-5}-13x_{1}^4x_{2}^{-4}+74x_{1}^3x_{2}^{-3} \\ +6x_{1}^2x_{2}^{-2}-13x_{1}^4x_{2}^{-5}-8x_{1}^3x_{2}^{-4} \\ +46x_{1}^2x_{2}^{-3}+2x_{1} x_{2}^{-2}-8x_{1}^3x_{2}^{-5} \\ -4x_{1}^2x_{2}^{-4}+26x_{1} x_{2}^{-3}-4x_{1}^2x_{2}^{-5} \\ -2x_{1} x_{2}^{-4}+13x_{2}^{-3}-2x_{1} x_{2}^{-5}-x_{2}^{-4} \\ +6x_{1}^{-1}x_{2}^{-3}-x_{2}^{-5}+2x_{1}^{-2}x_{2}^{-3}\end{array}\right)\begin{array}{l}\cdot \frac{1}{24}\left(\frac{x_{1} \partial_{1} }{6}-\frac{x_{2} \partial_{2} }{4}\right)^4\\
		~~\cdot\frac{1}{(1-x_{1}^6) (1-x_{1}^3x_{2}^2) }\end{array}\\&
	\displaystyle +\left(\begin{array}{l}-x_{1}^{23}x_{2}^{-3}-4x_{1}^{22}x_{2}^{-3}-10x_{1}^{21}x_{2}^{-3} \\ -20x_{1}^{20}x_{2}^{-3}-37x_{1}^{19}x_{2}^{-3}-62x_{1}^{18}x_{2}^{-3} \\ -92x_{1}^{17}x_{2}^{-3}-126x_{1}^{16}x_{2}^{-3}-161x_{1}^{15}x_{2}^{-3} \\ -194x_{1}^{14}x_{2}^{-3}-216x_{1}^{13}x_{2}^{-3}-224x_{1}^{12}x_{2}^{-3} \\ -221x_{1}^{11}x_{2}^{-3}-204x_{1}^{10}x_{2}^{-3}-176x_{1}^9x_{2}^{-3} \\ -140x_{1}^8x_{2}^{-3}-105x_{1}^7x_{2}^{-3}-74x_{1}^6x_{2}^{-3} \\ -46x_{1}^5x_{2}^{-3}-26x_{1}^4x_{2}^{-3}-13x_{1}^3x_{2}^{-3} \\ -6x_{1}^2x_{2}^{-3}-2x_{1} x_{2}^{-3}\end{array}\right)\begin{array}{l}\cdot \frac{1}{24}\left(\frac{x_{1} \partial_{1} }{6}-\frac{x_{2} \partial_{2} }{2}\right)^4\\
		~~\cdot\frac{1}{(1-x_{1}^6) (1-x_{1}^3x_{2} ) }\end{array}\\&
	\displaystyle +\left(\begin{array}{l}-x_{1}^4x_{2}^{-2}-4x_{1}^3x_{2}^{-2}-6x_{1}^2x_{2}^{-2} \\ -x_{1}^2x_{2}^{-3}-4x_{1} x_{2}^{-2}-4x_{1} x_{2}^{-3} \\ -x_{2}^{-2}+x_{1} x_{2}^{-4}-6x_{2}^{-3}+4x_{2}^{-4} \\ -4x_{1}^{-1}x_{2}^{-3}+x_{2}^{-5}+6x_{1}^{-1}x_{2}^{-4} \\ -x_{1}^{-2}x_{2}^{-3}+4x_{1}^{-1}x_{2}^{-5}+4x_{1}^{-2}x_{2}^{-4} \\ +6x_{1}^{-2}x_{2}^{-5}+x_{1}^{-3}x_{2}^{-4}+4x_{1}^{-3}x_{2}^{-5} \\ +x_{1}^{-4}x_{2}^{-5}\end{array}\right)\begin{array}{l}\cdot \frac{1}{24}\left(\frac{x_{1} \partial_{1} }{2}-\frac{3}{4}x_{2} \partial_{2} \right)^4\\
		~~\cdot\frac{1}{(1-x_{1}^2) (1-x_{1}^3x_{2}^2) }\end{array}\\&
	\displaystyle +\left(\begin{array}{l}x_{1}^{20}x_{2}^{-5}+2x_{1}^{19}x_{2}^{-5}+4x_{1}^{18}x_{2}^{-5} \\ +8x_{1}^{17}x_{2}^{-5}+13x_{1}^{16}x_{2}^{-5}+20x_{1}^{15}x_{2}^{-5} \\ +26x_{1}^{14}x_{2}^{-5}+34x_{1}^{13}x_{2}^{-5}+41x_{1}^{12}x_{2}^{-5} \\ +44x_{1}^{11}x_{2}^{-5}+46x_{1}^{10}x_{2}^{-5}+44x_{1}^9x_{2}^{-5} \\ +41x_{1}^8x_{2}^{-5}+34x_{1}^7x_{2}^{-5}+26x_{1}^6x_{2}^{-5} \\ +20x_{1}^5x_{2}^{-5}+13x_{1}^4x_{2}^{-5}+8x_{1}^3x_{2}^{-5} \\ +4x_{1}^2x_{2}^{-5}+2x_{1} x_{2}^{-5}+x_{2}^{-5}\end{array}\right)\begin{array}{l}\cdot \frac{1}{24}\left(\frac{x_{1} \partial_{1} }{6}\right)^4\\
		~~\cdot\frac{1}{(1-x_{1}^6) (1-x_{2} ) }\end{array}\end{align*}\begin{table}[h!]
	\begin{center}
		{\tiny
			\psset{xunit=2.82842712cm, yunit=2.82842712cm}
			\begin{pspicture}(-0.40824829,-0.70710678)(0.81649658,0.70710678)
				\psline[linecolor=black, linestyle=dashed](0,0)(-0.40824829,0.70710678)
				\rput[rb](-0.40824829,0.70710678){(1, 0, 0)}
				\pscircle*[linecolor=black](-0.40824829,0.70710678){0.07071068}
				\psline[linecolor=black, linestyle=dashed](0,0)(0.81649658,0)
				\rput[lt](0.81649658,0){(0, 1, 0)}
				\pscircle*[linecolor=black](0.81649658,0){0.07071068}
				\psline[linecolor=black, linestyle=dashed](0,0)(-0.40824829,-0.70710678)
				\rput[rt](-0.40824829,-0.70710678){(0, 0, 1)}
				\pscircle*[linecolor=black](-0.40824829,-0.70710678){0.07071068}
				\psline[linecolor=brown](0,0)(0.20412415,0)
				\rput[rt](0.13608276,-0.11785113){1}
				\rput[rt](0.13608276,0.11785113){2}
				\psline[linecolor=brown](0.20412415,-0.35355339)(0.81649658,0)
				\psline[linecolor=brown](0.20412415,-0.35355339)(0.20412415,0)
				\psline[linecolor=brown](0.81649658,0)(0.20412415,0)
				\rput[rt](0.40824829,-0.11785113){3}
				\psline[linecolor=brown](-0.40824829,-0.70710678)(0.20412415,-0.35355339)
				\psline[linecolor=brown](0.20412415,-0.35355339)(0,0)
				\rput[rt](-0.06804138,-0.35355339){4}
				\psline[linecolor=brown](-0.40824829,-0.70710678)(-0.40824829,0.70710678)
				\psline[linecolor=brown](-0.40824829,-0.70710678)(0,0)
				\psline[linecolor=brown](-0.40824829,0.70710678)(0,0)
				\rput[rt](-0.27216553,0){5}
				\psline[linecolor=brown](0.20412415,0.35355339)(0.81649658,0)
				\psline[linecolor=brown](0.20412415,0.35355339)(0.20412415,0)
				\rput[rt](0.40824829,0.11785113){6}
				\psline[linecolor=brown](0.20412415,0.35355339)(-0.40824829,0.70710678)
				\psline[linecolor=brown](0.20412415,0.35355339)(0,0)
				\rput[rt](-0.06804138,0.35355339){7}
				\pscircle*[linecolor=blue](-0.40824829,0.70710678){0.05303301}
				\pscircle*[linecolor=blue](0.81649658,0){0.05303301}
				\pscircle*[linecolor=blue](-0.40824829,-0.70710678){0.05303301}
				\pscircle*[linecolor=blue](0.20412415,0.35355339){0.05303301}
				\pscircle*[linecolor=blue](0.20412415,-0.35355339){0.05303301}
				\pscircle*[linecolor=blue](0,0){0.05303301}
		\end{pspicture}}
		\caption{$7$ combinatorial chambers of \(A_3\):(1, 0, 0), (0, 1, 0), (0, 0, 1), (1, 1, 0), (0, 1, 1), (1, 1, 1)}
	\end{center}
\end{table}
\begin{longtable}{|c|c|c|}\caption{\footnotesize V.p.f. of \(A_3\):(1, 0, 0), (0, 1, 0), (0, 0, 1), (1, 1, 0), (0, 1, 1), (1, 1, 1)}\\\hline N & Polynomial/Lattice & Shift(s)\\ \hline
	\endfirsthead\multicolumn{3}{c}{{\bfseries \tablename\ \thetable{} -- continued from previous page}} \\
	\hline  N &  Polynomial/Lattice & Shift(s) \\ \hline
	\endhead\hline \multicolumn{3}{|c|}{{Continued on next page}} \\ \hline\endfoot\endlastfoot1&\(\begin{array}{l}-\frac{x_{1}^3}{3}+\frac{x_{1}^2x_{2} }{2}-\frac{x_{1} x_{2}^2}{2}+x_{1} x_{2} x_{3} -\frac{x_{1} x_{3}^2}{2}+\frac{x_{2}^3}{6}-\frac{x_{2}^2x_{3} }{2} \\ +\frac{x_{2} x_{3}^2}{2}-\frac{x_{3}^3}{6}-\frac{x_{1}^2}{2}+x_{1} x_{2} +\frac{x_{1} x_{3} }{2}-\frac{x_{2}^2}{2}+x_{2} x_{3} -\frac{x_{3}^2}{2}+\frac{5}{6}x_{1}  \\ +\frac{x_{2} }{3}+\frac{2}{3}x_{3} +1\end{array}\)&-\\
	\hline\hline 
	2&\(\begin{array}{l}-\frac{x_{1}^3}{6}+\frac{x_{1}^2x_{2} }{2}-\frac{x_{1}^2x_{3} }{2}-\frac{x_{1} x_{2}^2}{2}+x_{1} x_{2} x_{3} +\frac{x_{2}^3}{6}-\frac{x_{2}^2x_{3} }{2} \\ +\frac{x_{2} x_{3}^2}{2}-\frac{x_{3}^3}{3}-\frac{x_{1}^2}{2}+x_{1} x_{2} +\frac{x_{1} x_{3} }{2}-\frac{x_{2}^2}{2}+x_{2} x_{3} -\frac{x_{3}^2}{2}+\frac{2}{3}x_{1}  \\ +\frac{x_{2} }{3}+\frac{5}{6}x_{3} +1\end{array}\)&-\\
	\hline\hline 
	3&\(-\frac{x_{1}^3}{6}+\frac{x_{1}^2x_{3} }{2}+\frac{3}{2}x_{1} x_{3} +\frac{7}{6}x_{1} +x_{3} +1\)&-\\
	\hline\hline 
	4&\(-\frac{x_{1}^3}{3}+\frac{x_{1}^2x_{2} }{2}-\frac{x_{1}^2}{2}+\frac{3}{2}x_{1} x_{2} +\frac{5}{6}x_{1} +x_{2} +1\)&-\\
	\hline\hline 
	5&\(\frac{x_{2}^3}{6}+x_{2}^2+\frac{11}{6}x_{2} +1\)&-\\
	\hline\hline 
	6&\(\frac{x_{1} x_{3}^2}{2}-\frac{x_{3}^3}{6}+\frac{3}{2}x_{1} x_{3} +x_{1} +\frac{7}{6}x_{3} +1\)&-\\
	\hline\hline 
	7&\(\frac{x_{2} x_{3}^2}{2}-\frac{x_{3}^3}{3}+\frac{3}{2}x_{2} x_{3} -\frac{x_{3}^2}{2}+x_{2} +\frac{5}{6}x_{3} +1\)&-\\
	\hline\end{longtable}

\begin{longtable}{|ccccc|}\caption{\footnotesize V.p.f. of \(A_3\):(1, 0, 0), (0, 1, 0), (0, 0, 1), (1, 1, 0), (0, 1, 1), (1, 1, 1)}\\\hline N & Defining inequalities & Vertices& Int. Pt.& Neighbors \\ \hline
	\endfirsthead\multicolumn{5}{c}{{\bfseries \tablename\ \thetable{} -- continued from previous page}} \\
	\hline N & Defining inequalities & Vertices&Int. Pt. & Neighbors\\ \hline
	\endhead\hline \multicolumn{5}{|c|}{{Continued on next page}} \\ \hline\endfoot\endlastfoot1&
	\(\begin{array}{rcl}x_{2}-x_{3}&\geq& 0\\
		x_{1}-x_{2}+x_{3}&\geq& 0\\
		-x_{1}+x_{3}&\geq& 0\\
	\end{array}\)&\(\begin{array}{l}(0, 1, 1): 2\checkmark\\(1, 1, 1): 4\checkmark\\(1, 2, 1): 5\checkmark\end{array}\)&(2, 4, 3): 19\checkmark& [7, 11, 16]\\\hline
	2&
	\(\begin{array}{rcl}-x_{1}+x_{2}&\geq& 0\\
		x_{1}-x_{2}+x_{3}&\geq& 0\\
		x_{1}-x_{3}&\geq& 0\\
	\end{array}\)&\(\begin{array}{l}(1, 1, 0): 2\checkmark\\(1, 1, 1): 4\checkmark\\(1, 2, 1): 5\checkmark\end{array}\)&(3, 4, 2): 19\checkmark& [4, 12, 15]\\\hline
	3&
	\(\begin{array}{rcl}x_{1}&\geq& 0\\
		-x_{1}+x_{2}-x_{3}&\geq& 0\\
		-x_{1}+x_{3}&\geq& 0\\
	\end{array}\)&\(\begin{array}{l}(0, 1, 1): 2\checkmark\\(0, 1, 0): 1\checkmark\\(1, 2, 1): 5\checkmark\end{array}\)&(1, 4, 2): 8\checkmark& [15, 12]\\\hline
	4&
	\(\begin{array}{rcl}x_{1}&\geq& 0\\
		-x_{1}+x_{2}&\geq& 0\\
		-x_{2}+x_{3}&\geq& 0\\
	\end{array}\)&\(\begin{array}{l}(0, 0, 1): 1\checkmark\\(0, 1, 1): 2\checkmark\\(1, 1, 1): 4\checkmark\end{array}\)&(1, 2, 3): 7\checkmark& [3, 15]\\\hline
	5&
	\(\begin{array}{rcl}x_{2}&\geq& 0\\
		x_{1}-x_{2}&\geq& 0\\
		-x_{2}+x_{3}&\geq& 0\\
	\end{array}\)&\(\begin{array}{l}(0, 0, 1): 1\checkmark\\(1, 0, 0): 1\checkmark\\(1, 1, 1): 4\checkmark\end{array}\)&(2, 1, 2): 4\checkmark& [7, 4]\\\hline
	6&
	\(\begin{array}{rcl}x_{3}&\geq& 0\\
		-x_{1}+x_{2}-x_{3}&\geq& 0\\
		x_{1}-x_{3}&\geq& 0\\
	\end{array}\)&\(\begin{array}{l}(1, 1, 0): 2\checkmark\\(0, 1, 0): 1\checkmark\\(1, 2, 1): 5\checkmark\end{array}\)&(2, 4, 1): 8\checkmark& [16, 11]\\\hline
	7&
	\(\begin{array}{rcl}x_{3}&\geq& 0\\
		x_{1}-x_{2}&\geq& 0\\
		x_{2}-x_{3}&\geq& 0\\
	\end{array}\)&\(\begin{array}{l}(1, 1, 0): 2\checkmark\\(1, 0, 0): 1\checkmark\\(1, 1, 1): 4\checkmark\end{array}\)&(3, 2, 1): 7\checkmark& [16, 3]\\\hline
\end{longtable}

\allowdisplaybreaks\begin{align*}&~~~
	\frac{1}{(1-x_{1} ) (1-x_{2} ) (1-x_{3} ) (1-x_{1} x_{2} ) (1-x_{2} x_{3} ) (1-x_{1} x_{2} x_{3} ) }\\=&~~~
	\displaystyle \frac{x_{2}^{-2}}{(1-x_{1} )^3 (1-x_{1} x_{2} ) (1-x_{3} )^2 }\\&
	+\displaystyle \frac{-x_{2}^{-2}x_{3}^{-1}}{(1-x_{1} )^4 (1-x_{1} x_{2} ) (1-x_{2} x_{3} ) }\\&
	+\displaystyle \frac{x_{2}^{-2}x_{3}^{-1}}{(1-x_{1} )^4 (1-x_{1} x_{2} x_{3} ) (1-x_{1} x_{2} ) }\\&
	+\displaystyle \frac{-x_{2}^{-2}x_{3}^{-2}}{(1-x_{1} )^3 (1-x_{2} )^2 (1-x_{2} x_{3} ) }\\&
	+\displaystyle \frac{x_{2}^{-2}+x_{2}^{-3}x_{3}^{-1}}{(1-x_{1} )^4 (1-x_{1} x_{2} x_{3} ) (1-x_{3} ) }\\&
	+\displaystyle \frac{x_{2}^{-3}x_{3}^{-2}}{(1-x_{1} )^4 (1-x_{2} ) (1-x_{1} x_{2} x_{3} ) }\\&
	+\displaystyle \frac{-x_{2}^{-3}x_{3}^{-2}}{(1-x_{1} )^4 (1-x_{2} ) (1-x_{3} ) }\\&
	+\displaystyle \frac{-x_{2}^{-2}}{(1-x_{1} )^3 (1-x_{1} x_{2} x_{3} ) (1-x_{3} )^2 }\\&
	+\displaystyle \frac{-x_{2}^{-2}}{(1-x_{1} )^4 (1-x_{2} x_{3} ) (1-x_{3} ) }\\&
	+\displaystyle \frac{x_{2}^{-2}x_{3}^{-2}}{(1-x_{1} )^3 (1-x_{2} )^2 (1-x_{3} ) }\\=&
	x_{2}^{-2}\cdot \frac{1}{2}\left(x_{1} \partial_{1} -x_{2} \partial_{2} \right)^2\left(x_{3} \partial_{3} \right)\cdot\frac{1}{(1-x_{1} ) (1-x_{1} x_{2} ) (1-x_{3} ) }\\&
	\displaystyle -x_{2}^{-2}x_{3}^{-1}\cdot \frac{1}{6}\left(x_{1} \partial_{1} -x_{2} \partial_{2} +x_{3} \partial_{3} \right)^3\cdot\frac{1}{(1-x_{1} ) (1-x_{1} x_{2} ) (1-x_{2} x_{3} ) }\\&
	\displaystyle +x_{2}^{-2}x_{3}^{-1}\cdot \frac{1}{6}\left(x_{1} \partial_{1} -x_{2} \partial_{2} \right)^3\cdot\frac{1}{(1-x_{1} ) (1-x_{1} x_{2} x_{3} ) (1-x_{1} x_{2} ) }\\&
	\displaystyle -x_{2}^{-2}x_{3}^{-2}\cdot \frac{1}{2}\left(x_{1} \partial_{1} \right)^2\left(x_{2} \partial_{2} -x_{3} \partial_{3} \right)\cdot\frac{1}{(1-x_{1} ) (1-x_{2} ) (1-x_{2} x_{3} ) }\\&
	\displaystyle +\left(x_{2}^{-2}+x_{2}^{-3}x_{3}^{-1}\right)\cdot \frac{1}{6}\left(x_{1} \partial_{1} -x_{2} \partial_{2} \right)^3\cdot\frac{1}{(1-x_{1} ) (1-x_{1} x_{2} x_{3} ) (1-x_{3} ) }\\&
	\displaystyle +x_{2}^{-3}x_{3}^{-2}\cdot \frac{1}{6}\left(x_{1} \partial_{1} -x_{3} \partial_{3} \right)^3\cdot\frac{1}{(1-x_{1} ) (1-x_{2} ) (1-x_{1} x_{2} x_{3} ) }\\&
	\displaystyle -x_{2}^{-3}x_{3}^{-2}\cdot \frac{1}{6}\left(x_{1} \partial_{1} \right)^3\cdot\frac{1}{(1-x_{1} ) (1-x_{2} ) (1-x_{3} ) }\\&
	\displaystyle -x_{2}^{-2}\cdot \frac{1}{2}\left(x_{1} \partial_{1} -x_{2} \partial_{2} \right)^2\left(-x_{2} \partial_{2} +x_{3} \partial_{3} \right)\cdot\frac{1}{(1-x_{1} ) (1-x_{1} x_{2} x_{3} ) (1-x_{3} ) }\\&
	\displaystyle -x_{2}^{-2}\cdot \frac{1}{6}\left(x_{1} \partial_{1} \right)^3\cdot\frac{1}{(1-x_{1} ) (1-x_{2} x_{3} ) (1-x_{3} ) }\\&
	\displaystyle +x_{2}^{-2}x_{3}^{-2}\cdot \frac{1}{2}\left(x_{1} \partial_{1} \right)^2\left(x_{2} \partial_{2} \right)\cdot\frac{1}{(1-x_{1} ) (1-x_{2} ) (1-x_{3} ) }\end{align*}\begin{table}[h!]
	\begin{center}
		{\tiny
			\psset{xunit=7cm, yunit=7cm}
}
		\caption{$23$ combinatorial chambers of \(B_3\):(1, 0, 0), (0, 1, 0), (0, 0, 1), (1, 1, 0), (0, 1, 1), (1, 1, 1), (0, 1, 2), (1, 1, 2), (1, 2, 2)}
	\end{center}
\end{table}
\begin{longtable}{|c|c|c|}\caption{\footnotesize V.p.f. of \(B_3\):(1, 0, 0), (0, 1, 0), (0, 0, 1), (1, 1, 0), (0, 1, 1), (1, 1, 1), (0, 1, 2), (1, 1, 2), (1, 2, 2)}\\\hline N & Polynomial/Lattice & Shift(s)\\ \hline
	\endfirsthead\multicolumn{3}{c}{{\bfseries \tablename\ \thetable{} -- continued from previous page}} \\
	\hline  N &  Polynomial/Lattice & Shift(s) \\ \hline
	\endhead\hline \multicolumn{3}{|c|}{{Continued on next page}} \\ \hline\endfoot\endlastfoot\multirow{5}{*}{1}&\(\Lambda=\langle((2, 0, 0),(0, 1, 0),(0, 0, 2))\rangle\)&\\\hline&\(\begin{array}{l}\begin{array}{l}\frac{x_{1}^6}{1440}-\frac{x_{1}^5x_{3} }{240}+\frac{x_{1}^4x_{3}^2}{96}-\frac{x_{1}^3x_{3}^3}{72}+\frac{x_{1}^2x_{3}^4}{96}-\frac{x_{1} x_{3}^5}{320}+\frac{x_{3}^6}{2880} \\ -\frac{x_{1}^5}{120}+\frac{x_{1}^4x_{3} }{24}-\frac{x_{1}^3x_{3}^2}{12}+\frac{x_{1}^2x_{3}^3}{12}-\frac{7}{384}x_{1} x_{3}^4+\frac{x_{3}^5}{640}+\frac{19}{576}x_{1}^4 \\ -\frac{19}{144}x_{1}^3x_{3} +\frac{19}{96}x_{1}^2x_{3}^2+\frac{19}{288}x_{1} x_{3}^3-\frac{11}{1152}x_{3}^4-\frac{x_{1}^3}{24}+\frac{x_{1}^2x_{3} }{8}-\frac{13}{720}x_{1}^2\end{array}\\
		\color{red}
		+\frac{21}{32}x_{1} x_{3}^2-\frac{x_{3}^3}{48}+\frac{133}{90}x_{1} x_{3} +\frac{407}{1440}x_{3}^2+\frac{21}{20}x_{1} +\frac{127}{120}x_{3} +1\color{black}
	\end{array}
	\)&\((0, 0, 0)\)\\
	\cline{2-3}&\(\begin{array}{l}\begin{array}{l}\frac{x_{1}^6}{1440}-\frac{x_{1}^5x_{3} }{240}+\frac{x_{1}^4x_{3}^2}{96}-\frac{x_{1}^3x_{3}^3}{72}+\frac{x_{1}^2x_{3}^4}{96}-\frac{x_{1} x_{3}^5}{320}+\frac{x_{3}^6}{2880} \\ -\frac{x_{1}^5}{120}+\frac{x_{1}^4x_{3} }{24}-\frac{x_{1}^3x_{3}^2}{12}+\frac{x_{1}^2x_{3}^3}{12}-\frac{7}{384}x_{1} x_{3}^4+\frac{x_{3}^5}{640}+\frac{19}{576}x_{1}^4 \\ -\frac{19}{144}x_{1}^3x_{3} +\frac{19}{96}x_{1}^2x_{3}^2+\frac{19}{288}x_{1} x_{3}^3-\frac{11}{1152}x_{3}^4-\frac{x_{1}^3}{24}+\frac{x_{1}^2x_{3} }{8}-\frac{13}{720}x_{1}^2\end{array}\\
		\color{red}
		+\frac{21}{32}x_{1} x_{3}^2-\frac{x_{3}^3}{48}+\frac{133}{90}x_{1} x_{3} +\frac{407}{1440}x_{3}^2+\frac{21}{20}x_{1} +\frac{127}{120}x_{3} +\frac{63}{64}\color{black}
	\end{array}
	\)&\((1, 0, 0)\)\\
\end{longtable}

\begin{longtable}{|ccccc|}\caption{\footnotesize V.p.f. of \(B_3\):(1, 0, 0), (0, 1, 0), (0, 0, 1), (1, 1, 0), (0, 1, 1), (1, 1, 1), (0, 1, 2), (1, 1, 2), (1, 2, 2)}\\\hline N & Defining inequalities & Vertices& Int. Pt.& Neighbors \\ \hline
	\endfirsthead\multicolumn{5}{c}{{\bfseries \tablename\ \thetable{} -- continued from previous page}} \\
	\hline N & Defining inequalities & Vertices&Int. Pt. & Neighbors\\ \hline
	\endhead\hline \multicolumn{5}{|c|}{{Continued on next page}} \\ \hline\endfoot\endlastfoot1&
	\(\begin{array}{rcl}-x_{1}+x_{2}-x_{3}&\geq& 0\\
		-x_{1}+x_{3}&\geq& 0\\
		2x_{1}-x_{3}&\geq& 0\\
	\end{array}\)&\(\begin{array}{l}(1, 2, 1): 5\checkmark\\(0, 1, 0): 1\checkmark\\(1, 3, 2): 12\checkmark\end{array}\)&(2, 6, 3): 38\checkmark& [40, 12, 35]\\\hline
	2&
	\(\begin{array}{rcl}x_{2}-x_{3}&\geq& 0\\
		-x_{1}+x_{3}&\geq& 0\\
		2x_{1}-2x_{2}+x_{3}&\geq& 0\\
	\end{array}\)&\(\begin{array}{l}(1, 1, 1): 4\checkmark\\(1, 2, 2): 11\checkmark\\(2, 3, 2): 20\checkmark\end{array}\)&(4, 6, 5): 208& [51, 23, 40]\\\hline
	3&
	\(\begin{array}{rcl}x_{2}-x_{3}&\geq& 0\\
		x_{1}-x_{2}+x_{3}&\geq& 0\\
		-2x_{1}+x_{3}&\geq& 0\\
	\end{array}\)&\(\begin{array}{l}(0, 1, 1): 2\checkmark\\(1, 2, 2): 11\checkmark\\(1, 3, 2): 12\checkmark\end{array}\)&(2, 6, 5): 97\checkmark& [42, 35, 40]\\\hline
	4&
	\(\begin{array}{rcl}-x_{1}+2x_{2}-x_{3}&\geq& 0\\
		-2x_{1}+x_{3}&\geq& 0\\
		2x_{1}-x_{2}&\geq& 0\\
	\end{array}\)&\(\begin{array}{l}(2, 3, 4): 44\checkmark\\(1, 2, 3): 14\checkmark\\(1, 2, 2): 11\checkmark\end{array}\)&(4, 7, 9): 642& [56, 51, 64]\\\hline
	5&
	\(\begin{array}{rcl}x_{1}-x_{2}&\geq& 0\\
		-x_{2}+x_{3}&\geq& 0\\
		2x_{2}-x_{3}&\geq& 0\\
	\end{array}\)&\(\begin{array}{l}(1, 1, 1): 4\checkmark\\(1, 0, 0): 1\checkmark\\(1, 1, 2): 6\checkmark\end{array}\)&(3, 2, 3): 19\checkmark& [55, 4, 17]\\\hline
	6&
	\(\begin{array}{rcl}-x_{1}+x_{2}&\geq& 0\\
		-2x_{2}+x_{3}&\geq& 0\\
		2x_{1}-x_{2}&\geq& 0\\
	\end{array}\)&\(\begin{array}{l}(1, 1, 2): 6\checkmark\\(0, 0, 1): 1\checkmark\\(1, 2, 4): 16\checkmark\end{array}\)&(2, 3, 7): 52\checkmark& [17, 62, 57]\\\hline
	7&
	\(\begin{array}{rcl}-x_{1}+x_{2}&\geq& 0\\
		x_{1}-x_{3}&\geq& 0\\
		2x_{1}-2x_{2}+x_{3}&\geq& 0\\
	\end{array}\)&\(\begin{array}{l}(1, 1, 1): 4\checkmark\\(1, 1, 0): 2\checkmark\\(2, 3, 2): 20\checkmark\end{array}\)&(4, 5, 3): 70\checkmark& [4, 27, 24]\\
\hline
\end{longtable}

\allowdisplaybreaks\begin{align*}&~~~
	\frac{1}{(1-x_{1} ) (1-x_{2} ) (1-x_{3} ) (1-x_{1} x_{2} ) (1-x_{2} x_{3} ) (1-x_{1} x_{2} x_{3} ) (1-x_{2} x_{3}^2) (1-x_{1} x_{2} x_{3}^2) (1-x_{1} x_{2}^2x_{3}^2) }\\=&~~~
	\displaystyle \frac{x_{1} x_{2}^{-3}}{(1-x_{1} )^5 (1-x_{1} x_{2} ) (1-x_{3} )^2(1-x_{3}^2) }\\&
	+\displaystyle \frac{x_{2}^{-5}x_{3}^{-8}}{(1-x_{1} )^5 (1-x_{2} )^3 (1-x_{3} ) }\\&
	+\displaystyle \frac{~~~~\begin{array}{l}-x_{2}^{-3}x_{3}^{-5}+x_{1}^{-1}x_{2}^{-2}x_{3}^{-5}+x_{2}^{-4}x_{3}^{-5}-x_{2}^{-4}x_{3}^{-6}+x_{1}^{-1}x_{2}^{-3}x_{3}^{-6} \\ \hline  +x_{2}^{-5}x_{3}^{-6}-x_{1}^{-1}x_{2}^{-5}x_{3}^{-7}-x_{1}^{-1}x_{2}^{-5}x_{3}^{-8}\end{array}~~~~}{(1-x_{1} )^5 (1-x_{2} )^3 (1-x_{1} x_{2}^2x_{3}^2) }\\&
	+\displaystyle \frac{-x_{1} x_{2}^{-2}x_{3}^2+x_{1} x_{2}^{-2}-x_{2}^{-3}x_{3}^2+x_{2}^{-3}x_{3}^{-2}}{(1-x_{1} )^5 (1-x_{1} x_{2}^2x_{3}^2) (1-x_{3} )^2(1-x_{3}^2) }\\&
	+\displaystyle \frac{-3x_{2}^{-6}x_{3}^{-8}-x_{2}^{-7}x_{3}^{-8}}{(1-x_{1} )^6 (1-x_{2} )^2 (1-x_{3} ) }\\&
	+\displaystyle \frac{~~~~\begin{array}{l}-x_{2}^{-3}x_{3}^{-5}+x_{1}^{-1}x_{2}^{-2}x_{3}^{-5}-x_{2}^{-4}x_{3}^{-6}-x_{2}^{-5}x_{3}^{-5}+x_{1}^{-1}x_{2}^{-3}x_{3}^{-6} \\ \hline  -x_{2}^{-6}x_{3}^{-6}+x_{2}^{-6}x_{3}^{-7}+x_{2}^{-6}x_{3}^{-8}+x_{2}^{-7}x_{3}^{-7}+2x_{1}^{-1}x_{2}^{-6}x_{3}^{-7}+x_{2}^{-7}x_{3}^{-8} \\ \hline  +2x_{1}^{-1}x_{2}^{-6}x_{3}^{-8}\end{array}~~~~}{(1-x_{1} )^6 (1-x_{2} )^2 (1-x_{1} x_{2}^2x_{3}^2) }\\&
	+\displaystyle \frac{5x_{2}^{-7}x_{3}^{-8}+4x_{2}^{-8}x_{3}^{-8}+x_{2}^{-9}x_{3}^{-8}}{(1-x_{1} )^7 (1-x_{2} ) (1-x_{3} ) }\\&
	+\displaystyle \frac{x_{1} x_{2}^{-3}}{(1-x_{1} )^6 (1-x_{1} x_{2} ) (1-x_{3} )^2 }\\&
	+\displaystyle \frac{~~~~\begin{array}{l}-x_{1} x_{2}^{-2}x_{3}^2-x_{1} x_{2}^{-2}-x_{2}^{-3}x_{3}^2+x_{1} x_{2}^{-2}x_{3}^{-1}-x_{1} x_{2}^{-3}x_{3}^{-1}+x_{1} x_{2}^{-2}x_{3}^{-3} \\ \hline  +x_{1} x_{2}^{-3}x_{3}^{-2}-x_{2}^{-4}-x_{2}^{-3}x_{3}^{-2}+x_{2}^{-4}x_{3}^{-1}+x_{2}^{-3}x_{3}^{-3}-x_{2}^{-4}x_{3}^{-2} \\ \hline  -x_{2}^{-5}x_{3}^{-2}+x_{2}^{-3}x_{3}^{-5}+x_{2}^{-4}x_{3}^{-4}+x_{2}^{-5}x_{3}^{-4}\end{array}~~~~}{(1-x_{1} )^6 (1-x_{1} x_{2}^2x_{3}^2) (1-x_{3} )^2 }\\&
	+\displaystyle \frac{~~~~\begin{array}{l}-x_{2}^{-3}x_{3}^{-5}-x_{2}^{-4}x_{3}^{-6}+x_{2}^{-6}x_{3}^{-5}+x_{2}^{-7}x_{3}^{-6}-3x_{2}^{-7}x_{3}^{-7}-3x_{2}^{-7}x_{3}^{-8} \\ \hline  -4x_{2}^{-8}x_{3}^{-7}-4x_{2}^{-8}x_{3}^{-8}-x_{2}^{-9}x_{3}^{-7}-x_{2}^{-9}x_{3}^{-8}\end{array}~~~~}{(1-x_{1} )^7 (1-x_{2} ) (1-x_{1} x_{2}^2x_{3}^2) }\\&
	+\displaystyle \frac{2x_{1} x_{2}^{-3}x_{3} +2x_{1} x_{2}^{-3}}{(1-x_{1} )^7 (1-x_{1} x_{2} ) (1-x_{3} ) }\\&
	+\displaystyle \frac{~~~~\begin{array}{l}-x_{1}^3-x_{1}^2x_{2}^{-1}x_{3}^{-1}-x_{1} x_{2}^{-2}-2x_{2}^{-3}x_{3}^2-2x_{1} x_{2}^{-2}x_{3}^{-1}-4x_{2}^{-3}x_{3} -2x_{1} x_{2}^{-2}x_{3}^{-2} \\ \hline  -2x_{2}^{-3}-x_{1}^2x_{2}^{-2}x_{3}^{-4}-x_{1}^2x_{2}^{-3}x_{3}^{-5}-x_{2}^{-3}x_{3}^{-3}\end{array}~~~~}{(1-x_{1} )^7 (1-x_{1} x_{2}^2x_{3}^2) (1-x_{1} x_{2} ) }\\&
	+\displaystyle \frac{~~~~\begin{array}{l}-2x_{1} x_{2}^{-2}x_{3}^3-2x_{1} x_{2}^{-2}x_{3}^2-x_{1} x_{2}^{-2}x_{3} -2x_{2}^{-3}x_{3}^3-x_{1} x_{2}^{-2}-2x_{2}^{-3}x_{3}^2 \\ \hline  -x_{1} x_{2}^{-2}x_{3}^{-2}-2x_{2}^{-4}x_{3} -x_{1} x_{2}^{-2}x_{3}^{-3}-x_{1} x_{2}^{-3}x_{3}^{-2}-x_{2}^{-3}x_{3}^{-1}-x_{2}^{-4} \\ \hline  -x_{1} x_{2}^{-2}x_{3}^{-4}-x_{1} x_{2}^{-3}x_{3}^{-3}-x_{2}^{-3}x_{3}^{-2}+x_{2}^{-4}x_{3}^{-1}-x_{1} x_{2}^{-2}x_{3}^{-5} \\ \hline  -2x_{1} x_{2}^{-3}x_{3}^{-4}-2x_{2}^{-5}x_{3}^{-1}-2x_{1} x_{2}^{-4}x_{3}^{-4}-x_{2}^{-3}x_{3}^{-4}-3x_{2}^{-5}x_{3}^{-2} \\ \hline  -x_{2}^{-3}x_{3}^{-5}-x_{2}^{-4}x_{3}^{-4}-x_{2}^{-3}x_{3}^{-6}-x_{2}^{-4}x_{3}^{-5}-x_{2}^{-3}x_{3}^{-7}-2x_{2}^{-4}x_{3}^{-6} \\ \hline  -x_{2}^{-5}x_{3}^{-5}-2x_{2}^{-6}x_{3}^{-4}-2x_{2}^{-5}x_{3}^{-6}-x_{2}^{-6}x_{3}^{-5}-6x_{2}^{-6}x_{3}^{-6}-9x_{2}^{-7}x_{3}^{-6} \\ \hline  -5x_{2}^{-8}x_{3}^{-6}-x_{2}^{-9}x_{3}^{-6}\end{array}~~~~}{(1-x_{1} )^7 (1-x_{1} x_{2}^2x_{3}^2) (1-x_{3} ) }\\&
	+\displaystyle \frac{-x_{2}^{-4}x_{3}^{-5}-x_{2}^{-5}x_{3}^{-6}-x_{2}^{-5}x_{3}^{-7}-x_{2}^{-5}x_{3}^{-8}}{(1-x_{1} )^5 (1-x_{2} )^3 (1-x_{2} x_{3}^2) }\\&
	+\displaystyle \frac{-x_{1}^{-1}x_{2}^{-2}x_{3}^{-5}-x_{1}^{-1}x_{2}^{-3}x_{3}^{-6}+x_{1}^{-1}x_{2}^{-5}x_{3}^{-7}+x_{1}^{-1}x_{2}^{-5}x_{3}^{-8}}{(1-x_{1} )^5 (1-x_{2} )^3 (1-x_{1} x_{2} x_{3}^2) }\\&
	+\displaystyle \frac{-x_{1}^{-1}x_{2}^{-2}x_{3}^{-5}-x_{1}^{-1}x_{2}^{-3}x_{3}^{-6}-2x_{1}^{-1}x_{2}^{-6}x_{3}^{-7}-2x_{1}^{-1}x_{2}^{-6}x_{3}^{-8}}{(1-x_{1} )^6 (1-x_{2} )^2 (1-x_{1} x_{2} x_{3}^2) }\\&
	+\displaystyle \frac{~~~~\begin{array}{l}x_{1} x_{2}^{-3}x_{3} +x_{1} x_{2}^{-3}+x_{1} x_{2}^{-3}x_{3}^{-2}+x_{1} x_{2}^{-3}x_{3}^{-3}+x_{1} x_{2}^{-4}x_{3}^{-2}+x_{1} x_{2}^{-3}x_{3}^{-4} \\ \hline  +x_{1} x_{2}^{-4}x_{3}^{-3}+x_{1} x_{2}^{-3}x_{3}^{-5}+2x_{1} x_{2}^{-4}x_{3}^{-4}+2x_{1} x_{2}^{-5}x_{3}^{-4}\end{array}~~~~}{(1-x_{1} )^7 (1-x_{1} x_{2} x_{3}^2) (1-x_{3} ) }\\&
	+\displaystyle \frac{~~~~\begin{array}{l}-x_{1}^2x_{2}^{-1}-x_{1} x_{2}^{-2}x_{3}^{-1}-x_{1} x_{2}^{-3}x_{3}^{-2}-x_{2}^{-2}x_{3}^{-3}-x_{1} x_{2}^{-3}x_{3}^{-4}+x_{2}^{-3}x_{3}^{-3} \\ \hline  +x_{2}^{-3}x_{3}^{-4}-x_{1} x_{2}^{-4}x_{3}^{-5}+2x_{2}^{-3}x_{3}^{-5}+2x_{2}^{-4}x_{3}^{-4}+2x_{2}^{-3}x_{3}^{-6}+3x_{2}^{-4}x_{3}^{-5} \\ \hline  +x_{2}^{-3}x_{3}^{-7}+2x_{2}^{-4}x_{3}^{-6}+2x_{2}^{-5}x_{3}^{-5}+2x_{2}^{-5}x_{3}^{-6}\end{array}~~~~}{(1-x_{1} )^7 (1-x_{1} x_{2} x_{3}^2) (1-x_{1} x_{2}^2x_{3}^2) }\\&
	+\displaystyle \frac{-x_{1} x_{2}^{-3}}{(1-x_{1} )^5 (1-x_{1} x_{2} x_{3}^2) (1-x_{3} )^2(1-x_{3}^2) }\\&
	+\displaystyle \frac{x_{2}^{-5}x_{3}^{-5}+x_{2}^{-6}x_{3}^{-6}+2x_{2}^{-6}x_{3}^{-7}+2x_{2}^{-6}x_{3}^{-8}}{(1-x_{1} )^6 (1-x_{2} )^2 (1-x_{2} x_{3}^2) }\\&
	+\displaystyle \frac{-x_{2}^{-6}x_{3}^{-5}-x_{2}^{-7}x_{3}^{-6}-2x_{2}^{-7}x_{3}^{-7}-2x_{2}^{-7}x_{3}^{-8}}{(1-x_{1} )^7 (1-x_{2} ) (1-x_{2} x_{3}^2) }\\&
	+\displaystyle \frac{x_{1} x_{2}^{-3}-x_{1} x_{2}^{-3}x_{3}^{-1}-x_{1} x_{2}^{-3}x_{3}^{-3}-x_{1} x_{2}^{-4}x_{3}^{-2}}{(1-x_{1} )^6 (1-x_{1} x_{2} x_{3}^2) (1-x_{3} )^2 }\\&
	+\displaystyle \frac{~~~~\begin{array}{l}2x_{2}^{-3}x_{3}^2+3x_{2}^{-3}x_{3} +x_{1} x_{2}^{-3}x_{3}^{-2}-x_{2}^{-3}x_{3}^{-1}-x_{2}^{-3}x_{3}^{-2}+2x_{2}^{-4}x_{3}^{-1} \\ \hline  +3x_{2}^{-4}x_{3}^{-2}+2x_{2}^{-4}x_{3}^{-3}+x_{2}^{-4}x_{3}^{-4}+x_{2}^{-5}x_{3}^{-3}-x_{1}^{-1}x_{2}^{-4}x_{3}^{-3}+x_{2}^{-5}x_{3}^{-5} \\ \hline  +x_{2}^{-6}x_{3}^{-4}+2x_{2}^{-6}x_{3}^{-5}+2x_{2}^{-6}x_{3}^{-6}\end{array}~~~~}{(1-x_{1} )^7 (1-x_{1} x_{2}^2x_{3}^2) (1-x_{2} x_{3}^2) }\\&
	+\displaystyle \frac{x_{2}^{-5}x_{3}^{-2}-x_{2}^{-5}x_{3}^{-3}+x_{2}^{-6}x_{3}^{-5}}{(1-x_{1} )^6 (1-x_{2} x_{3}^2) (1-x_{3} )^2 }\\&
	+\displaystyle \frac{x_{1} x_{2}^{-3}x_{3} +x_{1} x_{2}^{-3}}{(1-x_{1} )^6 (1-x_{1} x_{2} ) (1-x_{3} )(1-x_{3}^2) }\\&
	+\displaystyle \frac{-x_{1} x_{2}^{-3}x_{3} -x_{1} x_{2}^{-3}}{(1-x_{1} )^6 (1-x_{1} x_{2} x_{3}^2) (1-x_{3} )(1-x_{3}^2) }\\&
	+\displaystyle \frac{~~~~\begin{array}{l}-x_{1} x_{2}^{-2}x_{3}^3-x_{1} x_{2}^{-2}x_{3}^2+x_{1} x_{2}^{-2}x_{3} -x_{2}^{-3}x_{3}^3+x_{1} x_{2}^{-2}-x_{2}^{-3}x_{3}^2 \\ \hline  +x_{2}^{-3}x_{3}^{-1}+x_{2}^{-3}x_{3}^{-2}\end{array}~~~~}{(1-x_{1} )^6 (1-x_{1} x_{2}^2x_{3}^2) (1-x_{3} )(1-x_{3}^2) }\\&
	+\displaystyle \frac{2x_{2}^{-5}x_{3}^{-1}+x_{2}^{-5}x_{3}^{-2}-x_{2}^{-5}x_{3}^{-3}+2x_{2}^{-6}x_{3}^{-4}+x_{2}^{-6}x_{3}^{-5}}{(1-x_{1} )^7 (1-x_{2} x_{3}^2) (1-x_{3} ) }\\&
	+\displaystyle \frac{x_{2}^{-4}x_{3}^{-2}+2x_{2}^{-4}x_{3}^{-3}+2x_{2}^{-4}x_{3}^{-4}+x_{1} x_{2}^{-4}x_{3}^{-6}+x_{1} x_{2}^{-5}x_{3}^{-7}+x_{1}^{-1}x_{2}^{-5}x_{3}^{-5}}{(1-x_{1} )^6 (1-x_{1} x_{2} )^2 (1-x_{2} x_{3}^2) }\\&
	+\displaystyle \frac{-x_{2}^{-4}x_{3}^{-2}-2x_{2}^{-4}x_{3}^{-3}-2x_{2}^{-4}x_{3}^{-4}-x_{1} x_{2}^{-4}x_{3}^{-6}-x_{1} x_{2}^{-5}x_{3}^{-7}-x_{1}^{-1}x_{2}^{-5}x_{3}^{-5}}{(1-x_{1} )^6 (1-x_{1} x_{2}^2x_{3}^2) (1-x_{1} x_{2} )^2 }\\&
	+\displaystyle \frac{~~~~\begin{array}{l}x_{1}^3x_{2}^{-1}+x_{1}^2x_{2}^{-2}x_{3}^{-1}+x_{1} x_{2}^{-3}+2x_{1} x_{2}^{-3}x_{3}^{-1}+2x_{1} x_{2}^{-3}x_{3}^{-2}+x_{1}^2x_{2}^{-3}x_{3}^{-4} \\ \hline  +x_{1}^2x_{2}^{-4}x_{3}^{-5}+x_{2}^{-4}x_{3}^{-3}\end{array}~~~~}{(1-x_{1} )^7 (1-x_{1} x_{2} ) (1-x_{1} x_{2} x_{3}^2) }\\&
	+\displaystyle \frac{x_{2}^{-4}x_{3}^{-6}}{(1-x_{1} )^7 (1-x_{2} ) (1-x_{2} x_{3} ) }\\&
	+\displaystyle \frac{x_{1}^2x_{2}^{-4}x_{3}^{-3}+x_{1} x_{2}^{-5}x_{3}^{-4}+x_{2}^{-5}x_{3}^{-5}+x_{1}^{-1}x_{2}^{-6}x_{3}^{-6}}{(1-x_{1} )^5(1-x_{1}^2) (1-x_{2} )^2 (1-x_{1} x_{2}^2x_{3}^2) }\\&
	+\displaystyle \frac{x_{1}^{-1}x_{2}^{-2}x_{3}^{-5}-x_{1}^{-1}x_{2}^{-3}x_{3}^{-5}+x_{1}^{-1}x_{2}^{-3}x_{3}^{-6}-x_{1}^{-1}x_{2}^{-4}x_{3}^{-6}}{(1-x_{1} )^4 (1-x_{2} )^4 (1-x_{1} x_{2}^2x_{3}^2) }\\&
	+\displaystyle \frac{-x_{1}^{-1}x_{2}^{-2}x_{3}^{-5}+x_{1}^{-1}x_{2}^{-3}x_{3}^{-5}-x_{1}^{-1}x_{2}^{-3}x_{3}^{-6}+x_{1}^{-1}x_{2}^{-4}x_{3}^{-6}}{(1-x_{1} )^4 (1-x_{2} )^4 (1-x_{1} x_{2} x_{3}^2) }\\&
	+\displaystyle \frac{-x_{1}^4x_{2}^{-3}x_{3}^{-2}}{(1-x_{1} )^6(1-x_{1}^2) (1-x_{2} ) (1-x_{1} x_{2} x_{3} ) }\\&
	+\displaystyle \frac{-x_{1}^2x_{2}^{-4}x_{3}^{-3}-x_{1} x_{2}^{-5}x_{3}^{-4}-x_{2}^{-5}x_{3}^{-5}-x_{1}^{-1}x_{2}^{-6}x_{3}^{-6}}{(1-x_{1} )^5(1-x_{1}^2) (1-x_{2} )^2 (1-x_{1} x_{2} x_{3}^2) }\\&
	+\displaystyle \frac{x_{1}^4x_{2}^{-2}x_{3}^{-1}+x_{1}^3x_{2}^{-3}x_{3}^{-2}-x_{1} x_{2}^{-6}x_{3}^{-5}-x_{2}^{-7}x_{3}^{-6}}{(1-x_{1} )^6(1-x_{1}^2) (1-x_{2} ) (1-x_{1} x_{2}^2x_{3}^2) }\\&
	+\displaystyle \frac{-x_{1}^3x_{2}^{-4}x_{3}^{-2}-x_{1} x_{2}^{-5}x_{3}^{-4}}{(1-x_{1} )^6(1-x_{1}^2) (1-x_{1} x_{2} x_{3}^2) (1-x_{1} x_{2} x_{3} ) }\\&
	+\displaystyle \frac{~~~~\begin{array}{l}x_{1}^3x_{2}^{-3}x_{3}^{-1}+x_{1}^2x_{2}^{-4}x_{3}^{-2}+x_{1}^2x_{2}^{-3}x_{3}^{-4}+x_{1} x_{2}^{-4}x_{3}^{-3}+x_{1}^2x_{2}^{-4}x_{3}^{-5} \\ \hline  +x_{1} x_{2}^{-4}x_{3}^{-4}+x_{1} x_{2}^{-5}x_{3}^{-5}+x_{2}^{-5}x_{3}^{-4}\end{array}~~~~}{(1-x_{1} )^6(1-x_{1}^2) (1-x_{1} x_{2} x_{3}^2) (1-x_{1} x_{2}^2x_{3}^2) }\\&
	+\displaystyle \frac{x_{1} x_{2}^{-4}x_{3}^{-2}}{(1-x_{1} )^6 (1-x_{1} x_{2} x_{3} ) (1-x_{3} )^2 }\\&
	+\displaystyle \frac{-x_{1} x_{2}^{-4}x_{3}^{-3}+x_{2}^{-4}x_{3}^{-3}}{(1-x_{1} )^7 (1-x_{2} x_{3}^2) (1-x_{1} x_{2} x_{3} ) }\\&
	+\displaystyle \frac{x_{1} x_{2}^{-4}x_{3}^{-3}-x_{2}^{-4}x_{3}^{-3}}{(1-x_{1} )^7 (1-x_{1} x_{2} x_{3}^2) (1-x_{1} x_{2} x_{3} ) }\\&
	+\displaystyle \frac{x_{1} x_{2}^{-6}x_{3}^{-5}+x_{2}^{-7}x_{3}^{-6}}{(1-x_{1} )^6(1-x_{1}^2) (1-x_{2} ) (1-x_{2} x_{3}^2) }\\&
	+\displaystyle \frac{x_{1} x_{2}^{-5}x_{3}^{-4}}{(1-x_{1} )^6(1-x_{1}^2) (1-x_{1} x_{2} x_{3} ) (1-x_{2} x_{3}^2) }\\&
	+\displaystyle \frac{-x_{1} x_{2}^{-4}x_{3}^{-3}-x_{1} x_{2}^{-5}x_{3}^{-3}-x_{2}^{-5}x_{3}^{-4}-x_{2}^{-6}x_{3}^{-4}}{(1-x_{1} )^6(1-x_{1}^2) (1-x_{1} x_{2}^2x_{3}^2) (1-x_{2} x_{3}^2) }\\&
	+\displaystyle \frac{x_{2}^{-4}x_{3}^{-6}}{(1-x_{1} )^6 (1-x_{2} )^2 (1-x_{2} x_{3} ) }\\&
	+\displaystyle \frac{x_{2}^{-3}x_{3}^{-4}}{(1-x_{1} )^7 (1-x_{1} x_{2} x_{3}^2) (1-x_{2} x_{3} ) }\\&
	+\displaystyle \frac{x_{2}^{-4}x_{3}^{-6}}{(1-x_{1} )^5 (1-x_{2} )^3 (1-x_{2} x_{3} ) }\\&
	+\displaystyle \frac{x_{1} x_{2}^{-3}x_{3}^{-1}}{(1-x_{1} )^5 (1-x_{1} x_{2} x_{3}^2) (1-x_{3} )^3 }\\&
	+\displaystyle \frac{-x_{2}^{-6}x_{3}^{-5}}{(1-x_{1} )^6 (1-x_{2} x_{3} ) (1-x_{3} )^2 }\\&
	+\displaystyle \frac{-x_{1} x_{2}^{-5}x_{3}^{-5}}{(1-x_{1} )^6(1-x_{1}^2) (1-x_{2} x_{3} ) (1-x_{1} x_{2} x_{3}^2) }\\&
	+\displaystyle \frac{x_{1}^3x_{2}^{-1}+x_{1}^2x_{2}^{-2}x_{3}^{-1}}{(1-x_{1} )^6 (1-x_{1} x_{2} )^2 (1-x_{1} x_{2} x_{3}^2) }\\&
	+\displaystyle \frac{x_{2}^{-6}x_{3}^{-6}}{(1-x_{1} )^7 (1-x_{2} x_{3} ) (1-x_{3} ) }\\&
	+\displaystyle \frac{-x_{2}^{-5}x_{3}^{-5}}{(1-x_{1} )^7 (1-x_{2} x_{3}^2) (1-x_{2} x_{3} ) }\\&
	+\displaystyle \frac{x_{1} x_{2}^{-2}-x_{1} x_{2}^{-2}x_{3}^{-1}+x_{2}^{-3}x_{3}^{-1}-x_{2}^{-3}x_{3}^{-3}}{(1-x_{1} )^5 (1-x_{1} x_{2}^2x_{3}^2) (1-x_{3} )^3 }\\&
	+\displaystyle \frac{-x_{1} x_{2}^{-3}x_{3}^{-1}}{(1-x_{1} )^5 (1-x_{1} x_{2} x_{3} ) (1-x_{3} )^3 }\\&
	+\displaystyle \frac{-x_{1}^2x_{2}^{-6}x_{3}^{-7}}{(1-x_{1} )^6(1-x_{1}^2) (1-x_{1} x_{2} ) (1-x_{2} x_{3} ) }\\&
	+\displaystyle \frac{x_{1}^3x_{2}^{-2}x_{3}^{-4}+x_{1}^3x_{2}^{-3}x_{3}^{-5}+x_{1}^2x_{2}^{-5}x_{3}^{-6}+x_{1}^2x_{2}^{-6}x_{3}^{-7}}{(1-x_{1} )^6(1-x_{1}^2) (1-x_{1} x_{2}^2x_{3}^2) (1-x_{1} x_{2} ) }\\&
	+\displaystyle \frac{-x_{1}^3x_{2}^{-3}x_{3}^{-4}-x_{1}^3x_{2}^{-4}x_{3}^{-5}}{(1-x_{1} )^6(1-x_{1}^2) (1-x_{1} x_{2} ) (1-x_{1} x_{2} x_{3}^2) }\\&
	+\displaystyle \frac{-x_{1}^3x_{2}^{-1}x_{3}^{-1}}{(1-x_{1} )^6 (1-x_{1} x_{2} )^2 (1-x_{1} x_{2} x_{3} ) }\\=&~~~
	\displaystyle \frac{-x_{1}^3x_{2}^{-1}x_{3}^{-1}}{(1-x_{1} )^6 (1-x_{1} x_{2} )^2 (1-x_{1} x_{2} x_{3} ) }\\&
	+\displaystyle \frac{x_{1} x_{2}^{-3}x_{3}^2+2x_{1} x_{2}^{-3}x_{3} +x_{1} x_{2}^{-3}}{(1-x_{1} )^5 (1-x_{1} x_{2} ) (1-x_{3}^2)^3 }\\&
	+\displaystyle \frac{~~~~\begin{array}{l}-x_{1}^9x_{2}^{-3}x_{3}^{-4}-6x_{1}^8x_{2}^{-3}x_{3}^{-4}-x_{1}^9x_{2}^{-4}x_{3}^{-5}-15x_{1}^7x_{2}^{-3}x_{3}^{-4}-6x_{1}^8x_{2}^{-4}x_{3}^{-5} \\ \hline  -20x_{1}^6x_{2}^{-3}x_{3}^{-4}-15x_{1}^7x_{2}^{-4}x_{3}^{-5}-15x_{1}^5x_{2}^{-3}x_{3}^{-4}-20x_{1}^6x_{2}^{-4}x_{3}^{-5}-6x_{1}^4x_{2}^{-3}x_{3}^{-4} \\ \hline  -15x_{1}^5x_{2}^{-4}x_{3}^{-5}-x_{1}^3x_{2}^{-3}x_{3}^{-4}-6x_{1}^4x_{2}^{-4}x_{3}^{-5}-x_{1}^3x_{2}^{-4}x_{3}^{-5}\end{array}~~~~}{(1-x_{1}^2)^7 (1-x_{1} x_{2} ) (1-x_{1} x_{2} x_{3}^2) }\\&
	+\displaystyle \frac{~~~~\begin{array}{l}x_{1}^9x_{2}^{-2}x_{3}^{-4}+6x_{1}^8x_{2}^{-2}x_{3}^{-4}+x_{1}^9x_{2}^{-3}x_{3}^{-5}+15x_{1}^7x_{2}^{-2}x_{3}^{-4}+6x_{1}^8x_{2}^{-3}x_{3}^{-5} \\ \hline  +20x_{1}^6x_{2}^{-2}x_{3}^{-4}+15x_{1}^7x_{2}^{-3}x_{3}^{-5}+15x_{1}^5x_{2}^{-2}x_{3}^{-4}+20x_{1}^6x_{2}^{-3}x_{3}^{-5}+6x_{1}^4x_{2}^{-2}x_{3}^{-4} \\ \hline  +x_{1}^8x_{2}^{-5}x_{3}^{-6}+15x_{1}^5x_{2}^{-3}x_{3}^{-5}+x_{1}^3x_{2}^{-2}x_{3}^{-4}+6x_{1}^7x_{2}^{-5}x_{3}^{-6}+6x_{1}^4x_{2}^{-3}x_{3}^{-5} \\ \hline  +x_{1}^8x_{2}^{-6}x_{3}^{-7}+15x_{1}^6x_{2}^{-5}x_{3}^{-6}+x_{1}^3x_{2}^{-3}x_{3}^{-5}+6x_{1}^7x_{2}^{-6}x_{3}^{-7}+20x_{1}^5x_{2}^{-5}x_{3}^{-6} \\ \hline  +15x_{1}^6x_{2}^{-6}x_{3}^{-7}+15x_{1}^4x_{2}^{-5}x_{3}^{-6}+20x_{1}^5x_{2}^{-6}x_{3}^{-7}+6x_{1}^3x_{2}^{-5}x_{3}^{-6}+15x_{1}^4x_{2}^{-6}x_{3}^{-7} \\ \hline  +x_{1}^2x_{2}^{-5}x_{3}^{-6}+6x_{1}^3x_{2}^{-6}x_{3}^{-7}+x_{1}^2x_{2}^{-6}x_{3}^{-7}\end{array}~~~~}{(1-x_{1}^2)^7 (1-x_{1} x_{2}^2x_{3}^2) (1-x_{1} x_{2} ) }\\&
	+\displaystyle \frac{-x_{1} x_{2}^{-3}x_{3}^{-1}}{(1-x_{1} )^5 (1-x_{1} x_{2} x_{3} ) (1-x_{3} )^3 }\\&
	+\displaystyle \frac{~~~~\begin{array}{l}-x_{1}^8x_{2}^{-6}x_{3}^{-7}-6x_{1}^7x_{2}^{-6}x_{3}^{-7}-15x_{1}^6x_{2}^{-6}x_{3}^{-7}-20x_{1}^5x_{2}^{-6}x_{3}^{-7}-15x_{1}^4x_{2}^{-6}x_{3}^{-7} \\ \hline  -6x_{1}^3x_{2}^{-6}x_{3}^{-7}-x_{1}^2x_{2}^{-6}x_{3}^{-7}\end{array}~~~~}{(1-x_{1}^2)^7 (1-x_{1} x_{2} ) (1-x_{2} x_{3} ) }\\&
	+\displaystyle \frac{x_{1} x_{2}^{-2}-x_{1} x_{2}^{-2}x_{3}^{-1}+x_{2}^{-3}x_{3}^{-1}-x_{2}^{-3}x_{3}^{-3}}{(1-x_{1} )^5 (1-x_{1} x_{2}^2x_{3}^2) (1-x_{3} )^3 }\\&
	+\displaystyle \frac{-x_{2}^{-5}x_{3}^{-5}}{(1-x_{1} )^7 (1-x_{2} x_{3}^2) (1-x_{2} x_{3} ) }\\&
	+\displaystyle \frac{x_{2}^{-6}x_{3}^{-6}}{(1-x_{1} )^7 (1-x_{2} x_{3} ) (1-x_{3} ) }\\&
	+\displaystyle \frac{x_{1}^3x_{2}^{-1}+x_{1}^2x_{2}^{-2}x_{3}^{-1}}{(1-x_{1} )^6 (1-x_{1} x_{2} )^2 (1-x_{1} x_{2} x_{3}^2) }\\&
	+\displaystyle \frac{-x_{2}^{-6}x_{3}^{-5}}{(1-x_{1} )^6 (1-x_{2} x_{3} ) (1-x_{3} )^2 }\\&
	+\displaystyle \frac{~~~~\begin{array}{l}-x_{1}^7x_{2}^{-5}x_{3}^{-5}-6x_{1}^6x_{2}^{-5}x_{3}^{-5}-15x_{1}^5x_{2}^{-5}x_{3}^{-5}-20x_{1}^4x_{2}^{-5}x_{3}^{-5}-15x_{1}^3x_{2}^{-5}x_{3}^{-5} \\ \hline  -6x_{1}^2x_{2}^{-5}x_{3}^{-5}-x_{1} x_{2}^{-5}x_{3}^{-5}\end{array}~~~~}{(1-x_{1}^2)^7 (1-x_{2} x_{3} ) (1-x_{1} x_{2} x_{3}^2) }\\&
	+\displaystyle \frac{x_{1} x_{2}^{-3}x_{3}^{-1}}{(1-x_{1} )^5 (1-x_{1} x_{2} x_{3}^2) (1-x_{3} )^3 }\\&
	+\displaystyle \frac{x_{2}^{-4}x_{3}^{-6}}{(1-x_{1} )^5 (1-x_{2} )^3 (1-x_{2} x_{3} ) }\\&
	+\displaystyle \frac{x_{2}^{-3}x_{3}^{-4}}{(1-x_{1} )^7 (1-x_{1} x_{2} x_{3}^2) (1-x_{2} x_{3} ) }\\&
	+\displaystyle \frac{x_{2}^{-4}x_{3}^{-6}}{(1-x_{1} )^6 (1-x_{2} )^2 (1-x_{2} x_{3} ) }\\&
	+\displaystyle \frac{~~~~\begin{array}{l}-x_{1}^7x_{2}^{-4}x_{3}^{-3}-x_{1}^7x_{2}^{-5}x_{3}^{-3}-6x_{1}^6x_{2}^{-4}x_{3}^{-3}-6x_{1}^6x_{2}^{-5}x_{3}^{-3}-15x_{1}^5x_{2}^{-4}x_{3}^{-3} \\ \hline  -x_{1}^6x_{2}^{-5}x_{3}^{-4}-15x_{1}^5x_{2}^{-5}x_{3}^{-3}-20x_{1}^4x_{2}^{-4}x_{3}^{-3}-x_{1}^6x_{2}^{-6}x_{3}^{-4}-6x_{1}^5x_{2}^{-5}x_{3}^{-4} \\ \hline  -20x_{1}^4x_{2}^{-5}x_{3}^{-3}-15x_{1}^3x_{2}^{-4}x_{3}^{-3}-6x_{1}^5x_{2}^{-6}x_{3}^{-4}-15x_{1}^4x_{2}^{-5}x_{3}^{-4}-15x_{1}^3x_{2}^{-5}x_{3}^{-3} \\ \hline  -6x_{1}^2x_{2}^{-4}x_{3}^{-3}-15x_{1}^4x_{2}^{-6}x_{3}^{-4}-20x_{1}^3x_{2}^{-5}x_{3}^{-4}-6x_{1}^2x_{2}^{-5}x_{3}^{-3}-x_{1} x_{2}^{-4}x_{3}^{-3} \\ \hline  -20x_{1}^3x_{2}^{-6}x_{3}^{-4}-15x_{1}^2x_{2}^{-5}x_{3}^{-4}-x_{1} x_{2}^{-5}x_{3}^{-3}-15x_{1}^2x_{2}^{-6}x_{3}^{-4}-6x_{1} x_{2}^{-5}x_{3}^{-4} \\ \hline  -6x_{1} x_{2}^{-6}x_{3}^{-4}-x_{2}^{-5}x_{3}^{-4}-x_{2}^{-6}x_{3}^{-4}\end{array}~~~~}{(1-x_{1}^2)^7 (1-x_{1} x_{2}^2x_{3}^2) (1-x_{2} x_{3}^2) }\\&
	+\displaystyle \frac{~~~~\begin{array}{l}x_{1}^7x_{2}^{-5}x_{3}^{-4}+6x_{1}^6x_{2}^{-5}x_{3}^{-4}+15x_{1}^5x_{2}^{-5}x_{3}^{-4}+20x_{1}^4x_{2}^{-5}x_{3}^{-4}+15x_{1}^3x_{2}^{-5}x_{3}^{-4} \\ \hline  +6x_{1}^2x_{2}^{-5}x_{3}^{-4}+x_{1} x_{2}^{-5}x_{3}^{-4}\end{array}~~~~}{(1-x_{1}^2)^7 (1-x_{1} x_{2} x_{3} ) (1-x_{2} x_{3}^2) }\\&
	+\displaystyle \frac{x_{1} x_{2}^{-4}x_{3}^{-3}-x_{2}^{-4}x_{3}^{-3}}{(1-x_{1} )^7 (1-x_{1} x_{2} x_{3}^2) (1-x_{1} x_{2} x_{3} ) }\\&
	+\displaystyle \frac{~~~~\begin{array}{l}x_{1}^7x_{2}^{-6}x_{3}^{-5}+6x_{1}^6x_{2}^{-6}x_{3}^{-5}+15x_{1}^5x_{2}^{-6}x_{3}^{-5}+x_{1}^6x_{2}^{-7}x_{3}^{-6}+20x_{1}^4x_{2}^{-6}x_{3}^{-5} \\ \hline  +6x_{1}^5x_{2}^{-7}x_{3}^{-6}+15x_{1}^3x_{2}^{-6}x_{3}^{-5}+15x_{1}^4x_{2}^{-7}x_{3}^{-6}+6x_{1}^2x_{2}^{-6}x_{3}^{-5}+20x_{1}^3x_{2}^{-7}x_{3}^{-6} \\ \hline  +x_{1} x_{2}^{-6}x_{3}^{-5}+15x_{1}^2x_{2}^{-7}x_{3}^{-6}+6x_{1} x_{2}^{-7}x_{3}^{-6}+x_{2}^{-7}x_{3}^{-6}\end{array}~~~~}{(1-x_{1}^2)^7 (1-x_{2} ) (1-x_{2} x_{3}^2) }\\&
	+\displaystyle \frac{-x_{1} x_{2}^{-4}x_{3}^{-3}+x_{2}^{-4}x_{3}^{-3}}{(1-x_{1} )^7 (1-x_{2} x_{3}^2) (1-x_{1} x_{2} x_{3} ) }\\&
	+\displaystyle \frac{x_{1} x_{2}^{-4}x_{3}^{-2}}{(1-x_{1} )^6 (1-x_{1} x_{2} x_{3} ) (1-x_{3} )^2 }\\&
	+\displaystyle \frac{~~~~\begin{array}{l}x_{1}^9x_{2}^{-3}x_{3}^{-1}+6x_{1}^8x_{2}^{-3}x_{3}^{-1}+15x_{1}^7x_{2}^{-3}x_{3}^{-1}+x_{1}^8x_{2}^{-4}x_{3}^{-2}+20x_{1}^6x_{2}^{-3}x_{3}^{-1} \\ \hline  +x_{1}^8x_{2}^{-3}x_{3}^{-4}+6x_{1}^7x_{2}^{-4}x_{3}^{-2}+15x_{1}^5x_{2}^{-3}x_{3}^{-1}+6x_{1}^7x_{2}^{-3}x_{3}^{-4}+x_{1}^7x_{2}^{-4}x_{3}^{-3} \\ \hline  +15x_{1}^6x_{2}^{-4}x_{3}^{-2}+6x_{1}^4x_{2}^{-3}x_{3}^{-1}+x_{1}^8x_{2}^{-4}x_{3}^{-5}+x_{1}^7x_{2}^{-4}x_{3}^{-4}+15x_{1}^6x_{2}^{-3}x_{3}^{-4} \\ \hline  +6x_{1}^6x_{2}^{-4}x_{3}^{-3}+20x_{1}^5x_{2}^{-4}x_{3}^{-2}+x_{1}^3x_{2}^{-3}x_{3}^{-1}+6x_{1}^7x_{2}^{-4}x_{3}^{-5}+6x_{1}^6x_{2}^{-4}x_{3}^{-4} \\ \hline  +20x_{1}^5x_{2}^{-3}x_{3}^{-4}+15x_{1}^5x_{2}^{-4}x_{3}^{-3}+15x_{1}^4x_{2}^{-4}x_{3}^{-2}+x_{1}^7x_{2}^{-5}x_{3}^{-5}+15x_{1}^6x_{2}^{-4}x_{3}^{-5} \\ \hline  +x_{1}^6x_{2}^{-5}x_{3}^{-4}+15x_{1}^5x_{2}^{-4}x_{3}^{-4}+15x_{1}^4x_{2}^{-3}x_{3}^{-4}+20x_{1}^4x_{2}^{-4}x_{3}^{-3}+6x_{1}^3x_{2}^{-4}x_{3}^{-2} \\ \hline  +6x_{1}^6x_{2}^{-5}x_{3}^{-5}+20x_{1}^5x_{2}^{-4}x_{3}^{-5}+6x_{1}^5x_{2}^{-5}x_{3}^{-4}+20x_{1}^4x_{2}^{-4}x_{3}^{-4}+6x_{1}^3x_{2}^{-3}x_{3}^{-4} \\ \hline  +15x_{1}^3x_{2}^{-4}x_{3}^{-3}+x_{1}^2x_{2}^{-4}x_{3}^{-2}+15x_{1}^5x_{2}^{-5}x_{3}^{-5}+15x_{1}^4x_{2}^{-4}x_{3}^{-5}+15x_{1}^4x_{2}^{-5}x_{3}^{-4} \\ \hline  +15x_{1}^3x_{2}^{-4}x_{3}^{-4}+x_{1}^2x_{2}^{-3}x_{3}^{-4}+6x_{1}^2x_{2}^{-4}x_{3}^{-3}+20x_{1}^4x_{2}^{-5}x_{3}^{-5}+6x_{1}^3x_{2}^{-4}x_{3}^{-5} \\ \hline  +20x_{1}^3x_{2}^{-5}x_{3}^{-4}+6x_{1}^2x_{2}^{-4}x_{3}^{-4}+x_{1} x_{2}^{-4}x_{3}^{-3}+15x_{1}^3x_{2}^{-5}x_{3}^{-5}+x_{1}^2x_{2}^{-4}x_{3}^{-5} \\ \hline  +15x_{1}^2x_{2}^{-5}x_{3}^{-4}+x_{1} x_{2}^{-4}x_{3}^{-4}+6x_{1}^2x_{2}^{-5}x_{3}^{-5}+6x_{1} x_{2}^{-5}x_{3}^{-4}+x_{1} x_{2}^{-5}x_{3}^{-5} \\ \hline  +x_{2}^{-5}x_{3}^{-4}\end{array}~~~~}{(1-x_{1}^2)^7 (1-x_{1} x_{2} x_{3}^2) (1-x_{1} x_{2}^2x_{3}^2) }\\&
	+\displaystyle \frac{~~~~\begin{array}{l}-x_{1}^9x_{2}^{-4}x_{3}^{-2}-6x_{1}^8x_{2}^{-4}x_{3}^{-2}-15x_{1}^7x_{2}^{-4}x_{3}^{-2}-20x_{1}^6x_{2}^{-4}x_{3}^{-2}-15x_{1}^5x_{2}^{-4}x_{3}^{-2} \\ \hline  -x_{1}^7x_{2}^{-5}x_{3}^{-4}-6x_{1}^4x_{2}^{-4}x_{3}^{-2}-6x_{1}^6x_{2}^{-5}x_{3}^{-4}-x_{1}^3x_{2}^{-4}x_{3}^{-2}-15x_{1}^5x_{2}^{-5}x_{3}^{-4} \\ \hline  -20x_{1}^4x_{2}^{-5}x_{3}^{-4}-15x_{1}^3x_{2}^{-5}x_{3}^{-4}-6x_{1}^2x_{2}^{-5}x_{3}^{-4}-x_{1} x_{2}^{-5}x_{3}^{-4}\end{array}~~~~}{(1-x_{1}^2)^7 (1-x_{1} x_{2} x_{3}^2) (1-x_{1} x_{2} x_{3} ) }\\&
	+\displaystyle \frac{~~~~\begin{array}{l}x_{1}^{10}x_{2}^{-2}x_{3}^{-1}+6x_{1}^9x_{2}^{-2}x_{3}^{-1}+15x_{1}^8x_{2}^{-2}x_{3}^{-1}+x_{1}^9x_{2}^{-3}x_{3}^{-2}+20x_{1}^7x_{2}^{-2}x_{3}^{-1} \\ \hline  +6x_{1}^8x_{2}^{-3}x_{3}^{-2}+15x_{1}^6x_{2}^{-2}x_{3}^{-1}+15x_{1}^7x_{2}^{-3}x_{3}^{-2}+6x_{1}^5x_{2}^{-2}x_{3}^{-1}+20x_{1}^6x_{2}^{-3}x_{3}^{-2} \\ \hline  +x_{1}^4x_{2}^{-2}x_{3}^{-1}+15x_{1}^5x_{2}^{-3}x_{3}^{-2}+6x_{1}^4x_{2}^{-3}x_{3}^{-2}+x_{1}^3x_{2}^{-3}x_{3}^{-2}-x_{1}^7x_{2}^{-6}x_{3}^{-5} \\ \hline  -6x_{1}^6x_{2}^{-6}x_{3}^{-5}-15x_{1}^5x_{2}^{-6}x_{3}^{-5}-x_{1}^6x_{2}^{-7}x_{3}^{-6}-20x_{1}^4x_{2}^{-6}x_{3}^{-5}-6x_{1}^5x_{2}^{-7}x_{3}^{-6} \\ \hline  -15x_{1}^3x_{2}^{-6}x_{3}^{-5}-15x_{1}^4x_{2}^{-7}x_{3}^{-6}-6x_{1}^2x_{2}^{-6}x_{3}^{-5}-20x_{1}^3x_{2}^{-7}x_{3}^{-6}-x_{1} x_{2}^{-6}x_{3}^{-5} \\ \hline  -15x_{1}^2x_{2}^{-7}x_{3}^{-6}-6x_{1} x_{2}^{-7}x_{3}^{-6}-x_{2}^{-7}x_{3}^{-6}\end{array}~~~~}{(1-x_{1}^2)^7 (1-x_{2} ) (1-x_{1} x_{2}^2x_{3}^2) }\\&
	+\displaystyle \frac{~~~~\begin{array}{l}-x_{1}^7x_{2}^{-4}x_{3}^{-3}-5x_{1}^6x_{2}^{-4}x_{3}^{-3}-10x_{1}^5x_{2}^{-4}x_{3}^{-3}-x_{1}^6x_{2}^{-5}x_{3}^{-4}-10x_{1}^4x_{2}^{-4}x_{3}^{-3} \\ \hline  -5x_{1}^5x_{2}^{-5}x_{3}^{-4}-5x_{1}^3x_{2}^{-4}x_{3}^{-3}-x_{1}^5x_{2}^{-5}x_{3}^{-5}-10x_{1}^4x_{2}^{-5}x_{3}^{-4}-x_{1}^2x_{2}^{-4}x_{3}^{-3} \\ \hline  -5x_{1}^4x_{2}^{-5}x_{3}^{-5}-10x_{1}^3x_{2}^{-5}x_{3}^{-4}-10x_{1}^3x_{2}^{-5}x_{3}^{-5}-5x_{1}^2x_{2}^{-5}x_{3}^{-4}-x_{1}^4x_{2}^{-6}x_{3}^{-6} \\ \hline  -10x_{1}^2x_{2}^{-5}x_{3}^{-5}-x_{1} x_{2}^{-5}x_{3}^{-4}-5x_{1}^3x_{2}^{-6}x_{3}^{-6}-5x_{1} x_{2}^{-5}x_{3}^{-5}-10x_{1}^2x_{2}^{-6}x_{3}^{-6} \\ \hline  -x_{2}^{-5}x_{3}^{-5}-10x_{1} x_{2}^{-6}x_{3}^{-6}-5x_{2}^{-6}x_{3}^{-6}-x_{1}^{-1}x_{2}^{-6}x_{3}^{-6}\end{array}~~~~}{(1-x_{1}^2)^6 (1-x_{2} )^2 (1-x_{1} x_{2} x_{3}^2) }\\&
	+\displaystyle \frac{-x_{1}^{-1}x_{2}^{-2}x_{3}^{-5}+x_{1}^{-1}x_{2}^{-3}x_{3}^{-5}-x_{1}^{-1}x_{2}^{-3}x_{3}^{-6}+x_{1}^{-1}x_{2}^{-4}x_{3}^{-6}}{(1-x_{1} )^4 (1-x_{2} )^4 (1-x_{1} x_{2} x_{3}^2) }\\&
	+\displaystyle \frac{~~~~\begin{array}{l}-x_{1}^{10}x_{2}^{-3}x_{3}^{-2}-6x_{1}^9x_{2}^{-3}x_{3}^{-2}-15x_{1}^8x_{2}^{-3}x_{3}^{-2}-20x_{1}^7x_{2}^{-3}x_{3}^{-2}-15x_{1}^6x_{2}^{-3}x_{3}^{-2} \\ \hline  -6x_{1}^5x_{2}^{-3}x_{3}^{-2}-x_{1}^4x_{2}^{-3}x_{3}^{-2}\end{array}~~~~}{(1-x_{1}^2)^7 (1-x_{2} ) (1-x_{1} x_{2} x_{3} ) }\\&
	+\displaystyle \frac{x_{1}^{-1}x_{2}^{-2}x_{3}^{-5}-x_{1}^{-1}x_{2}^{-3}x_{3}^{-5}+x_{1}^{-1}x_{2}^{-3}x_{3}^{-6}-x_{1}^{-1}x_{2}^{-4}x_{3}^{-6}}{(1-x_{1} )^4 (1-x_{2} )^4 (1-x_{1} x_{2}^2x_{3}^2) }\\&
	+\displaystyle \frac{x_{2}^{-4}x_{3}^{-6}}{(1-x_{1} )^7 (1-x_{2} ) (1-x_{2} x_{3} ) }\\&
	+\displaystyle \frac{~~~~\begin{array}{l}x_{1}^7x_{2}^{-4}x_{3}^{-3}+5x_{1}^6x_{2}^{-4}x_{3}^{-3}+10x_{1}^5x_{2}^{-4}x_{3}^{-3}+x_{1}^6x_{2}^{-5}x_{3}^{-4}+10x_{1}^4x_{2}^{-4}x_{3}^{-3} \\ \hline  +5x_{1}^5x_{2}^{-5}x_{3}^{-4}+5x_{1}^3x_{2}^{-4}x_{3}^{-3}+x_{1}^5x_{2}^{-5}x_{3}^{-5}+10x_{1}^4x_{2}^{-5}x_{3}^{-4}+x_{1}^2x_{2}^{-4}x_{3}^{-3} \\ \hline  +5x_{1}^4x_{2}^{-5}x_{3}^{-5}+10x_{1}^3x_{2}^{-5}x_{3}^{-4}+10x_{1}^3x_{2}^{-5}x_{3}^{-5}+5x_{1}^2x_{2}^{-5}x_{3}^{-4}+x_{1}^4x_{2}^{-6}x_{3}^{-6} \\ \hline  +10x_{1}^2x_{2}^{-5}x_{3}^{-5}+x_{1} x_{2}^{-5}x_{3}^{-4}+5x_{1}^3x_{2}^{-6}x_{3}^{-6}+5x_{1} x_{2}^{-5}x_{3}^{-5}+10x_{1}^2x_{2}^{-6}x_{3}^{-6} \\ \hline  +x_{2}^{-5}x_{3}^{-5}+10x_{1} x_{2}^{-6}x_{3}^{-6}+5x_{2}^{-6}x_{3}^{-6}+x_{1}^{-1}x_{2}^{-6}x_{3}^{-6}\end{array}~~~~}{(1-x_{1}^2)^6 (1-x_{2} )^2 (1-x_{1} x_{2}^2x_{3}^2) }\\&
	+\displaystyle \frac{~~~~\begin{array}{l}x_{1}^3x_{2}^{-1}+x_{1}^2x_{2}^{-2}x_{3}^{-1}+x_{1} x_{2}^{-3}+2x_{1} x_{2}^{-3}x_{3}^{-1}+2x_{1} x_{2}^{-3}x_{3}^{-2}+x_{1}^2x_{2}^{-3}x_{3}^{-4} \\ \hline  +x_{1}^2x_{2}^{-4}x_{3}^{-5}+x_{2}^{-4}x_{3}^{-3}\end{array}~~~~}{(1-x_{1} )^7 (1-x_{1} x_{2} ) (1-x_{1} x_{2} x_{3}^2) }\\&
	+\displaystyle \frac{-x_{2}^{-4}x_{3}^{-2}-2x_{2}^{-4}x_{3}^{-3}-2x_{2}^{-4}x_{3}^{-4}-x_{1} x_{2}^{-4}x_{3}^{-6}-x_{1} x_{2}^{-5}x_{3}^{-7}-x_{1}^{-1}x_{2}^{-5}x_{3}^{-5}}{(1-x_{1} )^6 (1-x_{1} x_{2}^2x_{3}^2) (1-x_{1} x_{2} )^2 }\\&
	+\displaystyle \frac{x_{2}^{-4}x_{3}^{-2}+2x_{2}^{-4}x_{3}^{-3}+2x_{2}^{-4}x_{3}^{-4}+x_{1} x_{2}^{-4}x_{3}^{-6}+x_{1} x_{2}^{-5}x_{3}^{-7}+x_{1}^{-1}x_{2}^{-5}x_{3}^{-5}}{(1-x_{1} )^6 (1-x_{1} x_{2} )^2 (1-x_{2} x_{3}^2) }\\&
	+\displaystyle \frac{2x_{2}^{-5}x_{3}^{-1}+x_{2}^{-5}x_{3}^{-2}-x_{2}^{-5}x_{3}^{-3}+2x_{2}^{-6}x_{3}^{-4}+x_{2}^{-6}x_{3}^{-5}}{(1-x_{1} )^7 (1-x_{2} x_{3}^2) (1-x_{3} ) }\\&
	+\displaystyle \frac{~~~~\begin{array}{l}-x_{1} x_{2}^{-2}x_{3}^4-2x_{1} x_{2}^{-2}x_{3}^3-x_{2}^{-3}x_{3}^4+2x_{1} x_{2}^{-2}x_{3} -2x_{2}^{-3}x_{3}^3+x_{1} x_{2}^{-2} \\ \hline  -x_{2}^{-3}x_{3}^2+x_{2}^{-3}+2x_{2}^{-3}x_{3}^{-1}+x_{2}^{-3}x_{3}^{-2}\end{array}~~~~}{(1-x_{1} )^6 (1-x_{1} x_{2}^2x_{3}^2) (1-x_{3}^2)^2 }\\&
	+\displaystyle \frac{-x_{1} x_{2}^{-3}x_{3}^2-2x_{1} x_{2}^{-3}x_{3} -x_{1} x_{2}^{-3}}{(1-x_{1} )^6 (1-x_{1} x_{2} x_{3}^2) (1-x_{3}^2)^2 }\\&
	+\displaystyle \frac{x_{2}^{-5}x_{3}^{-2}-x_{2}^{-5}x_{3}^{-3}+x_{2}^{-6}x_{3}^{-5}}{(1-x_{1} )^6 (1-x_{2} x_{3}^2) (1-x_{3} )^2 }\\&
	+\displaystyle \frac{x_{1} x_{2}^{-3}x_{3}^2+2x_{1} x_{2}^{-3}x_{3} +x_{1} x_{2}^{-3}}{(1-x_{1} )^6 (1-x_{1} x_{2} ) (1-x_{3}^2)^2 }\\&
	+\displaystyle \frac{~~~~\begin{array}{l}2x_{2}^{-3}x_{3}^2+3x_{2}^{-3}x_{3} +x_{1} x_{2}^{-3}x_{3}^{-2}-x_{2}^{-3}x_{3}^{-1}-x_{2}^{-3}x_{3}^{-2}+2x_{2}^{-4}x_{3}^{-1} \\ \hline  +3x_{2}^{-4}x_{3}^{-2}+2x_{2}^{-4}x_{3}^{-3}+x_{2}^{-4}x_{3}^{-4}+x_{2}^{-5}x_{3}^{-3}-x_{1}^{-1}x_{2}^{-4}x_{3}^{-3}+x_{2}^{-5}x_{3}^{-5} \\ \hline  +x_{2}^{-6}x_{3}^{-4}+2x_{2}^{-6}x_{3}^{-5}+2x_{2}^{-6}x_{3}^{-6}\end{array}~~~~}{(1-x_{1} )^7 (1-x_{1} x_{2}^2x_{3}^2) (1-x_{2} x_{3}^2) }\\&
	+\displaystyle \frac{x_{1} x_{2}^{-3}-x_{1} x_{2}^{-3}x_{3}^{-1}-x_{1} x_{2}^{-3}x_{3}^{-3}-x_{1} x_{2}^{-4}x_{3}^{-2}}{(1-x_{1} )^6 (1-x_{1} x_{2} x_{3}^2) (1-x_{3} )^2 }\\&
	+\displaystyle \frac{-x_{2}^{-6}x_{3}^{-5}-x_{2}^{-7}x_{3}^{-6}-2x_{2}^{-7}x_{3}^{-7}-2x_{2}^{-7}x_{3}^{-8}}{(1-x_{1} )^7 (1-x_{2} ) (1-x_{2} x_{3}^2) }\\&
	+\displaystyle \frac{x_{2}^{-5}x_{3}^{-5}+x_{2}^{-6}x_{3}^{-6}+2x_{2}^{-6}x_{3}^{-7}+2x_{2}^{-6}x_{3}^{-8}}{(1-x_{1} )^6 (1-x_{2} )^2 (1-x_{2} x_{3}^2) }\\&
	+\displaystyle \frac{~~~~\begin{array}{l}-x_{1}^2x_{2}^{-1}-x_{1} x_{2}^{-2}x_{3}^{-1}-x_{1} x_{2}^{-3}x_{3}^{-2}-x_{2}^{-2}x_{3}^{-3}-x_{1} x_{2}^{-3}x_{3}^{-4}+x_{2}^{-3}x_{3}^{-3} \\ \hline  +x_{2}^{-3}x_{3}^{-4}-x_{1} x_{2}^{-4}x_{3}^{-5}+2x_{2}^{-3}x_{3}^{-5}+2x_{2}^{-4}x_{3}^{-4}+2x_{2}^{-3}x_{3}^{-6}+3x_{2}^{-4}x_{3}^{-5} \\ \hline  +x_{2}^{-3}x_{3}^{-7}+2x_{2}^{-4}x_{3}^{-6}+2x_{2}^{-5}x_{3}^{-5}+2x_{2}^{-5}x_{3}^{-6}\end{array}~~~~}{(1-x_{1} )^7 (1-x_{1} x_{2} x_{3}^2) (1-x_{1} x_{2}^2x_{3}^2) }\\&
	+\displaystyle \frac{-x_{1} x_{2}^{-3}x_{3}^2-2x_{1} x_{2}^{-3}x_{3} -x_{1} x_{2}^{-3}}{(1-x_{1} )^5 (1-x_{1} x_{2} x_{3}^2) (1-x_{3}^2)^3 }\\&
	+\displaystyle \frac{~~~~\begin{array}{l}x_{1} x_{2}^{-3}x_{3} +x_{1} x_{2}^{-3}+x_{1} x_{2}^{-3}x_{3}^{-2}+x_{1} x_{2}^{-3}x_{3}^{-3}+x_{1} x_{2}^{-4}x_{3}^{-2}+x_{1} x_{2}^{-3}x_{3}^{-4} \\ \hline  +x_{1} x_{2}^{-4}x_{3}^{-3}+x_{1} x_{2}^{-3}x_{3}^{-5}+2x_{1} x_{2}^{-4}x_{3}^{-4}+2x_{1} x_{2}^{-5}x_{3}^{-4}\end{array}~~~~}{(1-x_{1} )^7 (1-x_{1} x_{2} x_{3}^2) (1-x_{3} ) }\\&
	+\displaystyle \frac{-x_{1}^{-1}x_{2}^{-2}x_{3}^{-5}-x_{1}^{-1}x_{2}^{-3}x_{3}^{-6}-2x_{1}^{-1}x_{2}^{-6}x_{3}^{-7}-2x_{1}^{-1}x_{2}^{-6}x_{3}^{-8}}{(1-x_{1} )^6 (1-x_{2} )^2 (1-x_{1} x_{2} x_{3}^2) }\\&
	+\displaystyle \frac{-x_{1}^{-1}x_{2}^{-2}x_{3}^{-5}-x_{1}^{-1}x_{2}^{-3}x_{3}^{-6}+x_{1}^{-1}x_{2}^{-5}x_{3}^{-7}+x_{1}^{-1}x_{2}^{-5}x_{3}^{-8}}{(1-x_{1} )^5 (1-x_{2} )^3 (1-x_{1} x_{2} x_{3}^2) }\\&
	+\displaystyle \frac{-x_{2}^{-4}x_{3}^{-5}-x_{2}^{-5}x_{3}^{-6}-x_{2}^{-5}x_{3}^{-7}-x_{2}^{-5}x_{3}^{-8}}{(1-x_{1} )^5 (1-x_{2} )^3 (1-x_{2} x_{3}^2) }\\&
	+\displaystyle \frac{~~~~\begin{array}{l}-2x_{1} x_{2}^{-2}x_{3}^3-2x_{1} x_{2}^{-2}x_{3}^2-x_{1} x_{2}^{-2}x_{3} -2x_{2}^{-3}x_{3}^3-x_{1} x_{2}^{-2}-2x_{2}^{-3}x_{3}^2 \\ \hline  -x_{1} x_{2}^{-2}x_{3}^{-2}-2x_{2}^{-4}x_{3} -x_{1} x_{2}^{-2}x_{3}^{-3}-x_{1} x_{2}^{-3}x_{3}^{-2}-x_{2}^{-3}x_{3}^{-1}-x_{2}^{-4} \\ \hline  -x_{1} x_{2}^{-2}x_{3}^{-4}-x_{1} x_{2}^{-3}x_{3}^{-3}-x_{2}^{-3}x_{3}^{-2}+x_{2}^{-4}x_{3}^{-1}-x_{1} x_{2}^{-2}x_{3}^{-5} \\ \hline  -2x_{1} x_{2}^{-3}x_{3}^{-4}-2x_{2}^{-5}x_{3}^{-1}-2x_{1} x_{2}^{-4}x_{3}^{-4}-x_{2}^{-3}x_{3}^{-4}-3x_{2}^{-5}x_{3}^{-2} \\ \hline  -x_{2}^{-3}x_{3}^{-5}-x_{2}^{-4}x_{3}^{-4}-x_{2}^{-3}x_{3}^{-6}-x_{2}^{-4}x_{3}^{-5}-x_{2}^{-3}x_{3}^{-7}-2x_{2}^{-4}x_{3}^{-6} \\ \hline  -x_{2}^{-5}x_{3}^{-5}-2x_{2}^{-6}x_{3}^{-4}-2x_{2}^{-5}x_{3}^{-6}-x_{2}^{-6}x_{3}^{-5}-6x_{2}^{-6}x_{3}^{-6}-9x_{2}^{-7}x_{3}^{-6} \\ \hline  -5x_{2}^{-8}x_{3}^{-6}-x_{2}^{-9}x_{3}^{-6}\end{array}~~~~}{(1-x_{1} )^7 (1-x_{1} x_{2}^2x_{3}^2) (1-x_{3} ) }\\&
	+\displaystyle \frac{~~~~\begin{array}{l}-x_{1}^3-x_{1}^2x_{2}^{-1}x_{3}^{-1}-x_{1} x_{2}^{-2}-2x_{2}^{-3}x_{3}^2-2x_{1} x_{2}^{-2}x_{3}^{-1}-4x_{2}^{-3}x_{3} -2x_{1} x_{2}^{-2}x_{3}^{-2} \\ \hline  -2x_{2}^{-3}-x_{1}^2x_{2}^{-2}x_{3}^{-4}-x_{1}^2x_{2}^{-3}x_{3}^{-5}-x_{2}^{-3}x_{3}^{-3}\end{array}~~~~}{(1-x_{1} )^7 (1-x_{1} x_{2}^2x_{3}^2) (1-x_{1} x_{2} ) }\\&
	+\displaystyle \frac{2x_{1} x_{2}^{-3}x_{3} +2x_{1} x_{2}^{-3}}{(1-x_{1} )^7 (1-x_{1} x_{2} ) (1-x_{3} ) }\\&
	+\displaystyle \frac{~~~~\begin{array}{l}-x_{2}^{-3}x_{3}^{-5}-x_{2}^{-4}x_{3}^{-6}+x_{2}^{-6}x_{3}^{-5}+x_{2}^{-7}x_{3}^{-6}-3x_{2}^{-7}x_{3}^{-7}-3x_{2}^{-7}x_{3}^{-8} \\ \hline  -4x_{2}^{-8}x_{3}^{-7}-4x_{2}^{-8}x_{3}^{-8}-x_{2}^{-9}x_{3}^{-7}-x_{2}^{-9}x_{3}^{-8}\end{array}~~~~}{(1-x_{1} )^7 (1-x_{2} ) (1-x_{1} x_{2}^2x_{3}^2) }\\&
	+\displaystyle \frac{~~~~\begin{array}{l}-x_{1} x_{2}^{-2}x_{3}^2-x_{1} x_{2}^{-2}-x_{2}^{-3}x_{3}^2+x_{1} x_{2}^{-2}x_{3}^{-1}-x_{1} x_{2}^{-3}x_{3}^{-1}+x_{1} x_{2}^{-2}x_{3}^{-3} \\ \hline  +x_{1} x_{2}^{-3}x_{3}^{-2}-x_{2}^{-4}-x_{2}^{-3}x_{3}^{-2}+x_{2}^{-4}x_{3}^{-1}+x_{2}^{-3}x_{3}^{-3}-x_{2}^{-4}x_{3}^{-2} \\ \hline  -x_{2}^{-5}x_{3}^{-2}+x_{2}^{-3}x_{3}^{-5}+x_{2}^{-4}x_{3}^{-4}+x_{2}^{-5}x_{3}^{-4}\end{array}~~~~}{(1-x_{1} )^6 (1-x_{1} x_{2}^2x_{3}^2) (1-x_{3} )^2 }\\&
	+\displaystyle \frac{x_{1} x_{2}^{-3}}{(1-x_{1} )^6 (1-x_{1} x_{2} ) (1-x_{3} )^2 }\\&
	+\displaystyle \frac{5x_{2}^{-7}x_{3}^{-8}+4x_{2}^{-8}x_{3}^{-8}+x_{2}^{-9}x_{3}^{-8}}{(1-x_{1} )^7 (1-x_{2} ) (1-x_{3} ) }\\&
	+\displaystyle \frac{~~~~\begin{array}{l}-x_{2}^{-3}x_{3}^{-5}+x_{1}^{-1}x_{2}^{-2}x_{3}^{-5}-x_{2}^{-4}x_{3}^{-6}-x_{2}^{-5}x_{3}^{-5}+x_{1}^{-1}x_{2}^{-3}x_{3}^{-6} \\ \hline  -x_{2}^{-6}x_{3}^{-6}+x_{2}^{-6}x_{3}^{-7}+x_{2}^{-6}x_{3}^{-8}+x_{2}^{-7}x_{3}^{-7}+2x_{1}^{-1}x_{2}^{-6}x_{3}^{-7}+x_{2}^{-7}x_{3}^{-8} \\ \hline  +2x_{1}^{-1}x_{2}^{-6}x_{3}^{-8}\end{array}~~~~}{(1-x_{1} )^6 (1-x_{2} )^2 (1-x_{1} x_{2}^2x_{3}^2) }\\&
	+\displaystyle \frac{-3x_{2}^{-6}x_{3}^{-8}-x_{2}^{-7}x_{3}^{-8}}{(1-x_{1} )^6 (1-x_{2} )^2 (1-x_{3} ) }\\&
	+\displaystyle \frac{~~~~\begin{array}{l}-x_{2}^{-3}x_{3}^{-5}+x_{1}^{-1}x_{2}^{-2}x_{3}^{-5}+x_{2}^{-4}x_{3}^{-5}-x_{2}^{-4}x_{3}^{-6}+x_{1}^{-1}x_{2}^{-3}x_{3}^{-6} \\ \hline  +x_{2}^{-5}x_{3}^{-6}-x_{1}^{-1}x_{2}^{-5}x_{3}^{-7}-x_{1}^{-1}x_{2}^{-5}x_{3}^{-8}\end{array}~~~~}{(1-x_{1} )^5 (1-x_{2} )^3 (1-x_{1} x_{2}^2x_{3}^2) }\\&
	+\displaystyle \frac{~~~~\begin{array}{l}-x_{1} x_{2}^{-2}x_{3}^4-2x_{1} x_{2}^{-2}x_{3}^3-x_{2}^{-3}x_{3}^4+2x_{1} x_{2}^{-2}x_{3} -2x_{2}^{-3}x_{3}^3+x_{1} x_{2}^{-2} \\ \hline  -x_{2}^{-3}x_{3}^2+x_{2}^{-3}+2x_{2}^{-3}x_{3}^{-1}+x_{2}^{-3}x_{3}^{-2}\end{array}~~~~}{(1-x_{1} )^5 (1-x_{1} x_{2}^2x_{3}^2) (1-x_{3}^2)^3 }\\&
	+\displaystyle \frac{x_{2}^{-5}x_{3}^{-8}}{(1-x_{1} )^5 (1-x_{2} )^3 (1-x_{3} ) }\\=&
	-x_{1}^3x_{2}^{-1}x_{3}^{-1}\cdot \frac{1}{120}\left(x_{1} \partial_{1} -x_{2} \partial_{2} \right)^5\left(x_{2} \partial_{2} -x_{3} \partial_{3} \right)\cdot\frac{1}{(1-x_{1} ) (1-x_{1} x_{2} ) (1-x_{1} x_{2} x_{3} ) }\\&
	\displaystyle +\left(x_{1} x_{2}^{-3}x_{3}^2+2x_{1} x_{2}^{-3}x_{3} +x_{1} x_{2}^{-3}\right)\cdot \frac{1}{48}\left(x_{1} \partial_{1} -x_{2} \partial_{2} \right)^4\left(\frac{x_{3} \partial_{3} }{2}\right)^2\cdot\frac{1}{(1-x_{1} ) (1-x_{1} x_{2} ) (1-x_{3}^2) }\\&
	\displaystyle +\left(\begin{array}{l}-x_{1}^9x_{2}^{-3}x_{3}^{-4}-6x_{1}^8x_{2}^{-3}x_{3}^{-4} \\ -x_{1}^9x_{2}^{-4}x_{3}^{-5}-15x_{1}^7x_{2}^{-3}x_{3}^{-4} \\ -6x_{1}^8x_{2}^{-4}x_{3}^{-5}-20x_{1}^6x_{2}^{-3}x_{3}^{-4} \\ -15x_{1}^7x_{2}^{-4}x_{3}^{-5}-15x_{1}^5x_{2}^{-3}x_{3}^{-4} \\ -20x_{1}^6x_{2}^{-4}x_{3}^{-5}-6x_{1}^4x_{2}^{-3}x_{3}^{-4} \\ -15x_{1}^5x_{2}^{-4}x_{3}^{-5}-x_{1}^3x_{2}^{-3}x_{3}^{-4} \\ -6x_{1}^4x_{2}^{-4}x_{3}^{-5}-x_{1}^3x_{2}^{-4}x_{3}^{-5}\end{array}\right)\begin{array}{l}\cdot \frac{1}{720}\left(\frac{x_{1} \partial_{1} }{2}-\frac{x_{2} \partial_{2} }{2}\right)^6\\
		~~\cdot\frac{1}{(1-x_{1}^2) (1-x_{1} x_{2} ) (1-x_{1} x_{2} x_{3}^2) }\end{array}\\&
	\displaystyle +\left(\begin{array}{l}x_{1}^9x_{2}^{-2}x_{3}^{-4}+6x_{1}^8x_{2}^{-2}x_{3}^{-4} \\ +x_{1}^9x_{2}^{-3}x_{3}^{-5}+15x_{1}^7x_{2}^{-2}x_{3}^{-4} \\ +6x_{1}^8x_{2}^{-3}x_{3}^{-5}+20x_{1}^6x_{2}^{-2}x_{3}^{-4} \\ +15x_{1}^7x_{2}^{-3}x_{3}^{-5}+15x_{1}^5x_{2}^{-2}x_{3}^{-4} \\ +20x_{1}^6x_{2}^{-3}x_{3}^{-5}+6x_{1}^4x_{2}^{-2}x_{3}^{-4} \\ +x_{1}^8x_{2}^{-5}x_{3}^{-6}+15x_{1}^5x_{2}^{-3}x_{3}^{-5} \\ +x_{1}^3x_{2}^{-2}x_{3}^{-4}+6x_{1}^7x_{2}^{-5}x_{3}^{-6} \\ +6x_{1}^4x_{2}^{-3}x_{3}^{-5}+x_{1}^8x_{2}^{-6}x_{3}^{-7} \\ +15x_{1}^6x_{2}^{-5}x_{3}^{-6}+x_{1}^3x_{2}^{-3}x_{3}^{-5} \\ +6x_{1}^7x_{2}^{-6}x_{3}^{-7}+20x_{1}^5x_{2}^{-5}x_{3}^{-6} \\ +15x_{1}^6x_{2}^{-6}x_{3}^{-7}+15x_{1}^4x_{2}^{-5}x_{3}^{-6} \\ +20x_{1}^5x_{2}^{-6}x_{3}^{-7}+6x_{1}^3x_{2}^{-5}x_{3}^{-6} \\ +15x_{1}^4x_{2}^{-6}x_{3}^{-7}+x_{1}^2x_{2}^{-5}x_{3}^{-6} \\ +6x_{1}^3x_{2}^{-6}x_{3}^{-7}+x_{1}^2x_{2}^{-6}x_{3}^{-7}\end{array}\right)\begin{array}{l}\cdot \frac{1}{720}\left(\begin{array}{l}\frac{x_{1} \partial_{1} }{2}-\frac{x_{2} \partial_{2} }{2} \\ +\frac{x_{3} \partial_{3} }{4}\end{array}\right)^6\\
		~~\cdot\frac{1}{(1-x_{1}^2) (1-x_{1} x_{2}^2x_{3}^2) (1-x_{1} x_{2} ) }\end{array}\\&
	\displaystyle -x_{1} x_{2}^{-3}x_{3}^{-1}\cdot \frac{1}{48}\left(x_{1} \partial_{1} -x_{2} \partial_{2} \right)^4\left(-x_{2} \partial_{2} +x_{3} \partial_{3} \right)^2\cdot\frac{1}{(1-x_{1} ) (1-x_{1} x_{2} x_{3} ) (1-x_{3} ) }\\&
	\displaystyle +\left(\begin{array}{l}-x_{1}^8x_{2}^{-6}x_{3}^{-7}-6x_{1}^7x_{2}^{-6}x_{3}^{-7} \\ -15x_{1}^6x_{2}^{-6}x_{3}^{-7}-20x_{1}^5x_{2}^{-6}x_{3}^{-7} \\ -15x_{1}^4x_{2}^{-6}x_{3}^{-7}-6x_{1}^3x_{2}^{-6}x_{3}^{-7} \\ -x_{1}^2x_{2}^{-6}x_{3}^{-7}\end{array}\right)\begin{array}{l}\cdot \frac{1}{720}\left(\begin{array}{l}\frac{x_{1} \partial_{1} }{2}-\frac{x_{2} \partial_{2} }{2} \\ +\frac{x_{3} \partial_{3} }{2}\end{array}\right)^6\\
		~~\cdot\frac{1}{(1-x_{1}^2) (1-x_{1} x_{2} ) (1-x_{2} x_{3} ) }\end{array}\\&
	\displaystyle +\left(\begin{array}{l}x_{1} x_{2}^{-2}-x_{1} x_{2}^{-2}x_{3}^{-1}+x_{2}^{-3}x_{3}^{-1} \\ -x_{2}^{-3}x_{3}^{-3}\end{array}\right)\begin{array}{l}\cdot \frac{1}{48}\left(x_{1} \partial_{1} -\frac{x_{2} \partial_{2} }{2}\right)^4\left(-x_{2} \partial_{2} +x_{3} \partial_{3} \right)^2\\
		~~\cdot\frac{1}{(1-x_{1} ) (1-x_{1} x_{2}^2x_{3}^2) (1-x_{3} ) }\end{array}\\&
	\displaystyle -x_{2}^{-5}x_{3}^{-5}\cdot \frac{1}{720}\left(x_{1} \partial_{1} \right)^6\cdot\frac{1}{(1-x_{1} ) (1-x_{2} x_{3}^2) (1-x_{2} x_{3} ) }\\&
	\displaystyle +x_{2}^{-6}x_{3}^{-6}\cdot \frac{1}{720}\left(x_{1} \partial_{1} \right)^6\cdot\frac{1}{(1-x_{1} ) (1-x_{2} x_{3} ) (1-x_{3} ) }\\&
	\displaystyle +\left(x_{1}^3x_{2}^{-1}+x_{1}^2x_{2}^{-2}x_{3}^{-1}\right)\cdot \frac{1}{120}\left(x_{1} \partial_{1} -x_{2} \partial_{2} \right)^5\left(x_{2} \partial_{2} -\frac{x_{3} \partial_{3} }{2}\right)\cdot\frac{1}{(1-x_{1} ) (1-x_{1} x_{2} ) (1-x_{1} x_{2} x_{3}^2) }\\&
	\displaystyle -x_{2}^{-6}x_{3}^{-5}\cdot \frac{1}{120}\left(x_{1} \partial_{1} \right)^5\left(-x_{2} \partial_{2} +x_{3} \partial_{3} \right)\cdot\frac{1}{(1-x_{1} ) (1-x_{2} x_{3} ) (1-x_{3} ) }\\&
	\displaystyle +\left(\begin{array}{l}-x_{1}^7x_{2}^{-5}x_{3}^{-5}-6x_{1}^6x_{2}^{-5}x_{3}^{-5} \\ -15x_{1}^5x_{2}^{-5}x_{3}^{-5}-20x_{1}^4x_{2}^{-5}x_{3}^{-5} \\ -15x_{1}^3x_{2}^{-5}x_{3}^{-5}-6x_{1}^2x_{2}^{-5}x_{3}^{-5} \\ -x_{1} x_{2}^{-5}x_{3}^{-5}\end{array}\right)\begin{array}{l}\cdot \frac{1}{720}\left(\begin{array}{l}\frac{x_{1} \partial_{1} }{2}+\frac{x_{2} \partial_{2} }{2} \\ -\frac{x_{3} \partial_{3} }{2}\end{array}\right)^6\\
		~~\cdot\frac{1}{(1-x_{1}^2) (1-x_{2} x_{3} ) (1-x_{1} x_{2} x_{3}^2) }\end{array}\\&
	\displaystyle +x_{1} x_{2}^{-3}x_{3}^{-1}\cdot \frac{1}{48}\left(x_{1} \partial_{1} -x_{2} \partial_{2} \right)^4\left(-2x_{2} \partial_{2} +x_{3} \partial_{3} \right)^2\cdot\frac{1}{(1-x_{1} ) (1-x_{1} x_{2} x_{3}^2) (1-x_{3} ) }\\&
	\displaystyle +x_{2}^{-4}x_{3}^{-6}\cdot \frac{1}{48}\left(x_{1} \partial_{1} \right)^4\left(x_{2} \partial_{2} -x_{3} \partial_{3} \right)^2\cdot\frac{1}{(1-x_{1} ) (1-x_{2} ) (1-x_{2} x_{3} ) }\\&
	\displaystyle +x_{2}^{-3}x_{3}^{-4}\cdot \frac{1}{720}\left(x_{1} \partial_{1} +x_{2} \partial_{2} -x_{3} \partial_{3} \right)^6\cdot\frac{1}{(1-x_{1} ) (1-x_{1} x_{2} x_{3}^2) (1-x_{2} x_{3} ) }\\&
	\displaystyle +x_{2}^{-4}x_{3}^{-6}\cdot \frac{1}{120}\left(x_{1} \partial_{1} \right)^5\left(x_{2} \partial_{2} -x_{3} \partial_{3} \right)\cdot\frac{1}{(1-x_{1} ) (1-x_{2} ) (1-x_{2} x_{3} ) }\\&
	\displaystyle +\left(\begin{array}{l}-x_{1}^7x_{2}^{-4}x_{3}^{-3}-x_{1}^7x_{2}^{-5}x_{3}^{-3} \\ -6x_{1}^6x_{2}^{-4}x_{3}^{-3}-6x_{1}^6x_{2}^{-5}x_{3}^{-3} \\ -15x_{1}^5x_{2}^{-4}x_{3}^{-3}-x_{1}^6x_{2}^{-5}x_{3}^{-4} \\ -15x_{1}^5x_{2}^{-5}x_{3}^{-3}-20x_{1}^4x_{2}^{-4}x_{3}^{-3} \\ -x_{1}^6x_{2}^{-6}x_{3}^{-4}-6x_{1}^5x_{2}^{-5}x_{3}^{-4} \\ -20x_{1}^4x_{2}^{-5}x_{3}^{-3}-15x_{1}^3x_{2}^{-4}x_{3}^{-3} \\ -6x_{1}^5x_{2}^{-6}x_{3}^{-4}-15x_{1}^4x_{2}^{-5}x_{3}^{-4} \\ -15x_{1}^3x_{2}^{-5}x_{3}^{-3}-6x_{1}^2x_{2}^{-4}x_{3}^{-3} \\ -15x_{1}^4x_{2}^{-6}x_{3}^{-4}-20x_{1}^3x_{2}^{-5}x_{3}^{-4} \\ -6x_{1}^2x_{2}^{-5}x_{3}^{-3}-x_{1} x_{2}^{-4}x_{3}^{-3} \\ -20x_{1}^3x_{2}^{-6}x_{3}^{-4}-15x_{1}^2x_{2}^{-5}x_{3}^{-4} \\ -x_{1} x_{2}^{-5}x_{3}^{-3}-15x_{1}^2x_{2}^{-6}x_{3}^{-4} \\ -6x_{1} x_{2}^{-5}x_{3}^{-4}-6x_{1} x_{2}^{-6}x_{3}^{-4} \\ -x_{2}^{-5}x_{3}^{-4}-x_{2}^{-6}x_{3}^{-4}\end{array}\right)\begin{array}{l}\cdot \frac{1}{720}\left(\begin{array}{l}\frac{x_{1} \partial_{1} }{2}-\frac{x_{2} \partial_{2} }{2} \\ +\frac{x_{3} \partial_{3} }{4}\end{array}\right)^6\\
		~~\cdot\frac{1}{(1-x_{1}^2) (1-x_{1} x_{2}^2x_{3}^2) (1-x_{2} x_{3}^2) }\end{array}\\&
	\displaystyle +\left(\begin{array}{l}x_{1}^7x_{2}^{-5}x_{3}^{-4}+6x_{1}^6x_{2}^{-5}x_{3}^{-4} \\ +15x_{1}^5x_{2}^{-5}x_{3}^{-4}+20x_{1}^4x_{2}^{-5}x_{3}^{-4} \\ +15x_{1}^3x_{2}^{-5}x_{3}^{-4}+6x_{1}^2x_{2}^{-5}x_{3}^{-4} \\ +x_{1} x_{2}^{-5}x_{3}^{-4}\end{array}\right)\begin{array}{l}\cdot \frac{1}{720}\left(\frac{x_{1} \partial_{1} }{2}-x_{2} \partial_{2} +\frac{x_{3} \partial_{3} }{2}\right)^6\\
		~~\cdot\frac{1}{(1-x_{1}^2) (1-x_{1} x_{2} x_{3} ) (1-x_{2} x_{3}^2) }\end{array}\\&
	\displaystyle +\left(x_{1} x_{2}^{-4}x_{3}^{-3}-x_{2}^{-4}x_{3}^{-3}\right)\cdot \frac{1}{720}\left(x_{1} \partial_{1} -x_{2} \partial_{2} \right)^6\cdot\frac{1}{(1-x_{1} ) (1-x_{1} x_{2} x_{3}^2) (1-x_{1} x_{2} x_{3} ) }\\&
	\displaystyle +\left(\begin{array}{l}x_{1}^7x_{2}^{-6}x_{3}^{-5}+6x_{1}^6x_{2}^{-6}x_{3}^{-5} \\ +15x_{1}^5x_{2}^{-6}x_{3}^{-5}+x_{1}^6x_{2}^{-7}x_{3}^{-6} \\ +20x_{1}^4x_{2}^{-6}x_{3}^{-5}+6x_{1}^5x_{2}^{-7}x_{3}^{-6} \\ +15x_{1}^3x_{2}^{-6}x_{3}^{-5}+15x_{1}^4x_{2}^{-7}x_{3}^{-6} \\ +6x_{1}^2x_{2}^{-6}x_{3}^{-5}+20x_{1}^3x_{2}^{-7}x_{3}^{-6} \\ +x_{1} x_{2}^{-6}x_{3}^{-5}+15x_{1}^2x_{2}^{-7}x_{3}^{-6} \\ +6x_{1} x_{2}^{-7}x_{3}^{-6}+x_{2}^{-7}x_{3}^{-6}\end{array}\right)\begin{array}{l}\cdot \frac{1}{720}\left(\frac{x_{1} \partial_{1} }{2}\right)^6\\
		~~\cdot\frac{1}{(1-x_{1}^2) (1-x_{2} ) (1-x_{2} x_{3}^2) }\end{array}\\&
	\displaystyle +\left(-x_{1} x_{2}^{-4}x_{3}^{-3}+x_{2}^{-4}x_{3}^{-3}\right)\cdot \frac{1}{720}\left(x_{1} \partial_{1} -2x_{2} \partial_{2} +x_{3} \partial_{3} \right)^6\cdot\frac{1}{(1-x_{1} ) (1-x_{2} x_{3}^2) (1-x_{1} x_{2} x_{3} ) }\\&
	\displaystyle +x_{1} x_{2}^{-4}x_{3}^{-2}\cdot \frac{1}{120}\left(x_{1} \partial_{1} -x_{2} \partial_{2} \right)^5\left(-x_{2} \partial_{2} +x_{3} \partial_{3} \right)\cdot\frac{1}{(1-x_{1} ) (1-x_{1} x_{2} x_{3} ) (1-x_{3} ) }\\&
	\displaystyle +\left(\begin{array}{l}x_{1}^9x_{2}^{-3}x_{3}^{-1}+6x_{1}^8x_{2}^{-3}x_{3}^{-1} \\ +15x_{1}^7x_{2}^{-3}x_{3}^{-1}+x_{1}^8x_{2}^{-4}x_{3}^{-2} \\ +20x_{1}^6x_{2}^{-3}x_{3}^{-1}+x_{1}^8x_{2}^{-3}x_{3}^{-4} \\ +6x_{1}^7x_{2}^{-4}x_{3}^{-2}+15x_{1}^5x_{2}^{-3}x_{3}^{-1} \\ +6x_{1}^7x_{2}^{-3}x_{3}^{-4}+x_{1}^7x_{2}^{-4}x_{3}^{-3} \\ +15x_{1}^6x_{2}^{-4}x_{3}^{-2}+6x_{1}^4x_{2}^{-3}x_{3}^{-1} \\ +x_{1}^8x_{2}^{-4}x_{3}^{-5}+x_{1}^7x_{2}^{-4}x_{3}^{-4} \\ +15x_{1}^6x_{2}^{-3}x_{3}^{-4}+6x_{1}^6x_{2}^{-4}x_{3}^{-3} \\ +20x_{1}^5x_{2}^{-4}x_{3}^{-2}+x_{1}^3x_{2}^{-3}x_{3}^{-1} \\ +6x_{1}^7x_{2}^{-4}x_{3}^{-5}+6x_{1}^6x_{2}^{-4}x_{3}^{-4} \\ +20x_{1}^5x_{2}^{-3}x_{3}^{-4}+15x_{1}^5x_{2}^{-4}x_{3}^{-3} \\ +15x_{1}^4x_{2}^{-4}x_{3}^{-2}+x_{1}^7x_{2}^{-5}x_{3}^{-5} \\ +15x_{1}^6x_{2}^{-4}x_{3}^{-5}+x_{1}^6x_{2}^{-5}x_{3}^{-4} \\ +15x_{1}^5x_{2}^{-4}x_{3}^{-4}+15x_{1}^4x_{2}^{-3}x_{3}^{-4} \\ +20x_{1}^4x_{2}^{-4}x_{3}^{-3}+6x_{1}^3x_{2}^{-4}x_{3}^{-2} \\ +6x_{1}^6x_{2}^{-5}x_{3}^{-5}+20x_{1}^5x_{2}^{-4}x_{3}^{-5} \\ +6x_{1}^5x_{2}^{-5}x_{3}^{-4}+20x_{1}^4x_{2}^{-4}x_{3}^{-4} \\ +6x_{1}^3x_{2}^{-3}x_{3}^{-4}+15x_{1}^3x_{2}^{-4}x_{3}^{-3} \\ +x_{1}^2x_{2}^{-4}x_{3}^{-2}+15x_{1}^5x_{2}^{-5}x_{3}^{-5} \\ +15x_{1}^4x_{2}^{-4}x_{3}^{-5}+15x_{1}^4x_{2}^{-5}x_{3}^{-4} \\ +15x_{1}^3x_{2}^{-4}x_{3}^{-4}+x_{1}^2x_{2}^{-3}x_{3}^{-4} \\ +6x_{1}^2x_{2}^{-4}x_{3}^{-3}+20x_{1}^4x_{2}^{-5}x_{3}^{-5} \\ +6x_{1}^3x_{2}^{-4}x_{3}^{-5}+20x_{1}^3x_{2}^{-5}x_{3}^{-4} \\ +6x_{1}^2x_{2}^{-4}x_{3}^{-4}+x_{1} x_{2}^{-4}x_{3}^{-3} \\ +15x_{1}^3x_{2}^{-5}x_{3}^{-5}+x_{1}^2x_{2}^{-4}x_{3}^{-5} \\ +15x_{1}^2x_{2}^{-5}x_{3}^{-4}+x_{1} x_{2}^{-4}x_{3}^{-4} \\ +6x_{1}^2x_{2}^{-5}x_{3}^{-5}+6x_{1} x_{2}^{-5}x_{3}^{-4} \\ +x_{1} x_{2}^{-5}x_{3}^{-5}+x_{2}^{-5}x_{3}^{-4}\end{array}\right)\begin{array}{l}\cdot \frac{1}{720}\left(\frac{x_{1} \partial_{1} }{2}-\frac{x_{3} \partial_{3} }{4}\right)^6\\
		~~\cdot\frac{1}{(1-x_{1}^2) (1-x_{1} x_{2} x_{3}^2) (1-x_{1} x_{2}^2x_{3}^2) }\end{array}\\&
	\displaystyle +\left(\begin{array}{l}-x_{1}^9x_{2}^{-4}x_{3}^{-2}-6x_{1}^8x_{2}^{-4}x_{3}^{-2} \\ -15x_{1}^7x_{2}^{-4}x_{3}^{-2}-20x_{1}^6x_{2}^{-4}x_{3}^{-2} \\ -15x_{1}^5x_{2}^{-4}x_{3}^{-2}-x_{1}^7x_{2}^{-5}x_{3}^{-4} \\ -6x_{1}^4x_{2}^{-4}x_{3}^{-2}-6x_{1}^6x_{2}^{-5}x_{3}^{-4} \\ -x_{1}^3x_{2}^{-4}x_{3}^{-2}-15x_{1}^5x_{2}^{-5}x_{3}^{-4} \\ -20x_{1}^4x_{2}^{-5}x_{3}^{-4}-15x_{1}^3x_{2}^{-5}x_{3}^{-4} \\ -6x_{1}^2x_{2}^{-5}x_{3}^{-4}-x_{1} x_{2}^{-5}x_{3}^{-4}\end{array}\right)\begin{array}{l}\cdot \frac{1}{720}\left(\frac{x_{1} \partial_{1} }{2}-\frac{x_{2} \partial_{2} }{2}\right)^6\\
		~~\cdot\frac{1}{(1-x_{1}^2) (1-x_{1} x_{2} x_{3}^2) (1-x_{1} x_{2} x_{3} ) }\end{array}\\&
	\displaystyle +\left(\begin{array}{l}x_{1}^{10}x_{2}^{-2}x_{3}^{-1}+6x_{1}^9x_{2}^{-2}x_{3}^{-1} \\ +15x_{1}^8x_{2}^{-2}x_{3}^{-1}+x_{1}^9x_{2}^{-3}x_{3}^{-2} \\ +20x_{1}^7x_{2}^{-2}x_{3}^{-1}+6x_{1}^8x_{2}^{-3}x_{3}^{-2} \\ +15x_{1}^6x_{2}^{-2}x_{3}^{-1}+15x_{1}^7x_{2}^{-3}x_{3}^{-2} \\ +6x_{1}^5x_{2}^{-2}x_{3}^{-1}+20x_{1}^6x_{2}^{-3}x_{3}^{-2} \\ +x_{1}^4x_{2}^{-2}x_{3}^{-1}+15x_{1}^5x_{2}^{-3}x_{3}^{-2} \\ +6x_{1}^4x_{2}^{-3}x_{3}^{-2}+x_{1}^3x_{2}^{-3}x_{3}^{-2} \\ -x_{1}^7x_{2}^{-6}x_{3}^{-5}-6x_{1}^6x_{2}^{-6}x_{3}^{-5} \\ -15x_{1}^5x_{2}^{-6}x_{3}^{-5}-x_{1}^6x_{2}^{-7}x_{3}^{-6} \\ -20x_{1}^4x_{2}^{-6}x_{3}^{-5}-6x_{1}^5x_{2}^{-7}x_{3}^{-6} \\ -15x_{1}^3x_{2}^{-6}x_{3}^{-5}-15x_{1}^4x_{2}^{-7}x_{3}^{-6} \\ -6x_{1}^2x_{2}^{-6}x_{3}^{-5}-20x_{1}^3x_{2}^{-7}x_{3}^{-6} \\ -x_{1} x_{2}^{-6}x_{3}^{-5}-15x_{1}^2x_{2}^{-7}x_{3}^{-6} \\ -6x_{1} x_{2}^{-7}x_{3}^{-6}-x_{2}^{-7}x_{3}^{-6}\end{array}\right)\begin{array}{l}\cdot \frac{1}{720}\left(\frac{x_{1} \partial_{1} }{2}-\frac{x_{3} \partial_{3} }{4}\right)^6\\
		~~\cdot\frac{1}{(1-x_{1}^2) (1-x_{2} ) (1-x_{1} x_{2}^2x_{3}^2) }\end{array}\\&
	\displaystyle +\left(\begin{array}{l}-x_{1}^7x_{2}^{-4}x_{3}^{-3}-5x_{1}^6x_{2}^{-4}x_{3}^{-3} \\ -10x_{1}^5x_{2}^{-4}x_{3}^{-3}-x_{1}^6x_{2}^{-5}x_{3}^{-4} \\ -10x_{1}^4x_{2}^{-4}x_{3}^{-3}-5x_{1}^5x_{2}^{-5}x_{3}^{-4} \\ -5x_{1}^3x_{2}^{-4}x_{3}^{-3}-x_{1}^5x_{2}^{-5}x_{3}^{-5} \\ -10x_{1}^4x_{2}^{-5}x_{3}^{-4}-x_{1}^2x_{2}^{-4}x_{3}^{-3} \\ -5x_{1}^4x_{2}^{-5}x_{3}^{-5}-10x_{1}^3x_{2}^{-5}x_{3}^{-4} \\ -10x_{1}^3x_{2}^{-5}x_{3}^{-5}-5x_{1}^2x_{2}^{-5}x_{3}^{-4} \\ -x_{1}^4x_{2}^{-6}x_{3}^{-6}-10x_{1}^2x_{2}^{-5}x_{3}^{-5} \\ -x_{1} x_{2}^{-5}x_{3}^{-4}-5x_{1}^3x_{2}^{-6}x_{3}^{-6} \\ -5x_{1} x_{2}^{-5}x_{3}^{-5}-10x_{1}^2x_{2}^{-6}x_{3}^{-6} \\ -x_{2}^{-5}x_{3}^{-5}-10x_{1} x_{2}^{-6}x_{3}^{-6}-5x_{2}^{-6}x_{3}^{-6} \\ -x_{1}^{-1}x_{2}^{-6}x_{3}^{-6}\end{array}\right)\begin{array}{l}\cdot \frac{1}{120}\left(\frac{x_{1} \partial_{1} }{2}-\frac{x_{3} \partial_{3} }{4}\right)^5\left(x_{2} \partial_{2} -\frac{x_{3} \partial_{3} }{2}\right)\\
		~~\cdot\frac{1}{(1-x_{1}^2) (1-x_{2} ) (1-x_{1} x_{2} x_{3}^2) }\end{array}\\&
	\displaystyle +\left(\begin{array}{l}-x_{1}^{-1}x_{2}^{-2}x_{3}^{-5}+x_{1}^{-1}x_{2}^{-3}x_{3}^{-5} \\ -x_{1}^{-1}x_{2}^{-3}x_{3}^{-6}+x_{1}^{-1}x_{2}^{-4}x_{3}^{-6}\end{array}\right)\begin{array}{l}\cdot \frac{1}{36}\left(x_{1} \partial_{1} -\frac{x_{3} \partial_{3} }{2}\right)^3\left(x_{2} \partial_{2} -\frac{x_{3} \partial_{3} }{2}\right)^3\\
		~~\cdot\frac{1}{(1-x_{1} ) (1-x_{2} ) (1-x_{1} x_{2} x_{3}^2) }\end{array}\\&
	\displaystyle +\left(\begin{array}{l}-x_{1}^{10}x_{2}^{-3}x_{3}^{-2}-6x_{1}^9x_{2}^{-3}x_{3}^{-2} \\ -15x_{1}^8x_{2}^{-3}x_{3}^{-2}-20x_{1}^7x_{2}^{-3}x_{3}^{-2} \\ -15x_{1}^6x_{2}^{-3}x_{3}^{-2}-6x_{1}^5x_{2}^{-3}x_{3}^{-2} \\ -x_{1}^4x_{2}^{-3}x_{3}^{-2}\end{array}\right)\begin{array}{l}\cdot \frac{1}{720}\left(\frac{x_{1} \partial_{1} }{2}-\frac{x_{3} \partial_{3} }{2}\right)^6\\
		~~\cdot\frac{1}{(1-x_{1}^2) (1-x_{2} ) (1-x_{1} x_{2} x_{3} ) }\end{array}\\&
	\displaystyle +\left(\begin{array}{l}x_{1}^{-1}x_{2}^{-2}x_{3}^{-5}-x_{1}^{-1}x_{2}^{-3}x_{3}^{-5} \\ +x_{1}^{-1}x_{2}^{-3}x_{3}^{-6}-x_{1}^{-1}x_{2}^{-4}x_{3}^{-6}\end{array}\right)\begin{array}{l}\cdot \frac{1}{36}\left(x_{1} \partial_{1} -\frac{x_{3} \partial_{3} }{2}\right)^3\left(x_{2} \partial_{2} -x_{3} \partial_{3} \right)^3\\
		~~\cdot\frac{1}{(1-x_{1} ) (1-x_{2} ) (1-x_{1} x_{2}^2x_{3}^2) }\end{array}\\&
	\displaystyle +x_{2}^{-4}x_{3}^{-6}\cdot \frac{1}{720}\left(x_{1} \partial_{1} \right)^6\cdot\frac{1}{(1-x_{1} ) (1-x_{2} ) (1-x_{2} x_{3} ) }\\&
	\displaystyle +\left(\begin{array}{l}x_{1}^7x_{2}^{-4}x_{3}^{-3}+5x_{1}^6x_{2}^{-4}x_{3}^{-3} \\ +10x_{1}^5x_{2}^{-4}x_{3}^{-3}+x_{1}^6x_{2}^{-5}x_{3}^{-4} \\ +10x_{1}^4x_{2}^{-4}x_{3}^{-3}+5x_{1}^5x_{2}^{-5}x_{3}^{-4} \\ +5x_{1}^3x_{2}^{-4}x_{3}^{-3}+x_{1}^5x_{2}^{-5}x_{3}^{-5} \\ +10x_{1}^4x_{2}^{-5}x_{3}^{-4}+x_{1}^2x_{2}^{-4}x_{3}^{-3} \\ +5x_{1}^4x_{2}^{-5}x_{3}^{-5}+10x_{1}^3x_{2}^{-5}x_{3}^{-4} \\ +10x_{1}^3x_{2}^{-5}x_{3}^{-5}+5x_{1}^2x_{2}^{-5}x_{3}^{-4} \\ +x_{1}^4x_{2}^{-6}x_{3}^{-6}+10x_{1}^2x_{2}^{-5}x_{3}^{-5} \\ +x_{1} x_{2}^{-5}x_{3}^{-4}+5x_{1}^3x_{2}^{-6}x_{3}^{-6} \\ +5x_{1} x_{2}^{-5}x_{3}^{-5}+10x_{1}^2x_{2}^{-6}x_{3}^{-6} \\ +x_{2}^{-5}x_{3}^{-5}+10x_{1} x_{2}^{-6}x_{3}^{-6}+5x_{2}^{-6}x_{3}^{-6} \\ +x_{1}^{-1}x_{2}^{-6}x_{3}^{-6}\end{array}\right)\begin{array}{l}\cdot \frac{1}{120}\left(\frac{x_{1} \partial_{1} }{2}-\frac{x_{3} \partial_{3} }{4}\right)^5\left(x_{2} \partial_{2} -x_{3} \partial_{3} \right)\\
		~~\cdot\frac{1}{(1-x_{1}^2) (1-x_{2} ) (1-x_{1} x_{2}^2x_{3}^2) }\end{array}\\&
	\displaystyle +\left(\begin{array}{l}x_{1}^3x_{2}^{-1}+x_{1}^2x_{2}^{-2}x_{3}^{-1}+x_{1} x_{2}^{-3} \\ +2x_{1} x_{2}^{-3}x_{3}^{-1}+2x_{1} x_{2}^{-3}x_{3}^{-2} \\ +x_{1}^2x_{2}^{-3}x_{3}^{-4}+x_{1}^2x_{2}^{-4}x_{3}^{-5} \\ +x_{2}^{-4}x_{3}^{-3}\end{array}\right)\begin{array}{l}\cdot \frac{1}{720}\left(x_{1} \partial_{1} -x_{2} \partial_{2} \right)^6\\
		~~\cdot\frac{1}{(1-x_{1} ) (1-x_{1} x_{2} ) (1-x_{1} x_{2} x_{3}^2) }\end{array}\\&
	\displaystyle +\left(\begin{array}{l}-x_{2}^{-4}x_{3}^{-2}-2x_{2}^{-4}x_{3}^{-3}-2x_{2}^{-4}x_{3}^{-4} \\ -x_{1} x_{2}^{-4}x_{3}^{-6}-x_{1} x_{2}^{-5}x_{3}^{-7} \\ -x_{1}^{-1}x_{2}^{-5}x_{3}^{-5}\end{array}\right)\begin{array}{l}\cdot \frac{1}{120}\left(x_{1} \partial_{1} -x_{2} \partial_{2} +\frac{x_{3} \partial_{3} }{2}\right)^5\left(x_{2} \partial_{2} -x_{3} \partial_{3} \right)\\
		~~\cdot\frac{1}{(1-x_{1} ) (1-x_{1} x_{2}^2x_{3}^2) (1-x_{1} x_{2} ) }\end{array}\\&
	\displaystyle +\left(\begin{array}{l}x_{2}^{-4}x_{3}^{-2}+2x_{2}^{-4}x_{3}^{-3}+2x_{2}^{-4}x_{3}^{-4} \\ +x_{1} x_{2}^{-4}x_{3}^{-6}+x_{1} x_{2}^{-5}x_{3}^{-7} \\ +x_{1}^{-1}x_{2}^{-5}x_{3}^{-5}\end{array}\right)\begin{array}{l}\cdot \frac{1}{120}\left(x_{1} \partial_{1} -x_{2} \partial_{2} +\frac{x_{3} \partial_{3} }{2}\right)^5\left(x_{2} \partial_{2} -\frac{x_{3} \partial_{3} }{2}\right)\\
		~~\cdot\frac{1}{(1-x_{1} ) (1-x_{1} x_{2} ) (1-x_{2} x_{3}^2) }\end{array}\\&
	\displaystyle +\left(\begin{array}{l}2x_{2}^{-5}x_{3}^{-1}+x_{2}^{-5}x_{3}^{-2}-x_{2}^{-5}x_{3}^{-3} \\ +2x_{2}^{-6}x_{3}^{-4}+x_{2}^{-6}x_{3}^{-5}\end{array}\right)\begin{array}{l}\cdot \frac{1}{720}\left(x_{1} \partial_{1} \right)^6\\
		~~\cdot\frac{1}{(1-x_{1} ) (1-x_{2} x_{3}^2) (1-x_{3} ) }\end{array}\\&
	\displaystyle +\left(\begin{array}{l}-x_{1} x_{2}^{-2}x_{3}^4-2x_{1} x_{2}^{-2}x_{3}^3-x_{2}^{-3}x_{3}^4 \\ +2x_{1} x_{2}^{-2}x_{3} -2x_{2}^{-3}x_{3}^3+x_{1} x_{2}^{-2} \\ -x_{2}^{-3}x_{3}^2+x_{2}^{-3}+2x_{2}^{-3}x_{3}^{-1} \\ +x_{2}^{-3}x_{3}^{-2}\end{array}\right)\begin{array}{l}\cdot \frac{1}{120}\left(x_{1} \partial_{1} -\frac{x_{2} \partial_{2} }{2}\right)^5\left(-\frac{x_{2} \partial_{2} }{2}+\frac{x_{3} \partial_{3} }{2}\right)\\
		~~\cdot\frac{1}{(1-x_{1} ) (1-x_{1} x_{2}^2x_{3}^2) (1-x_{3}^2) }\end{array}\\&
	\displaystyle +\left(-x_{1} x_{2}^{-3}x_{3}^2-2x_{1} x_{2}^{-3}x_{3} -x_{1} x_{2}^{-3}\right)\cdot \frac{1}{120}\left(x_{1} \partial_{1} -x_{2} \partial_{2} \right)^5\left(-x_{2} \partial_{2} +\frac{x_{3} \partial_{3} }{2}\right)\cdot\frac{1}{(1-x_{1} ) (1-x_{1} x_{2} x_{3}^2) (1-x_{3}^2) }\\&
	\displaystyle +\left(x_{2}^{-5}x_{3}^{-2}-x_{2}^{-5}x_{3}^{-3}+x_{2}^{-6}x_{3}^{-5}\right)\cdot \frac{1}{120}\left(x_{1} \partial_{1} \right)^5\left(-2x_{2} \partial_{2} +x_{3} \partial_{3} \right)\cdot\frac{1}{(1-x_{1} ) (1-x_{2} x_{3}^2) (1-x_{3} ) }\\&
	\displaystyle +\left(x_{1} x_{2}^{-3}x_{3}^2+2x_{1} x_{2}^{-3}x_{3} +x_{1} x_{2}^{-3}\right)\cdot \frac{1}{120}\left(x_{1} \partial_{1} -x_{2} \partial_{2} \right)^5\left(\frac{x_{3} \partial_{3} }{2}\right)\cdot\frac{1}{(1-x_{1} ) (1-x_{1} x_{2} ) (1-x_{3}^2) }\\&
	\displaystyle +\left(
}
		\caption{$23$ combinatorial chambers of \(C_3\):(1, 0, 0), (0, 1, 0), (0, 0, 1), (1, 1, 0), (0, 1, 1), (1, 1, 1), (0, 2, 1), (1, 2, 1), (2, 2, 1)}
	\end{center}
\end{table}
\begin{longtable}{|c|c|c|}\caption{\footnotesize V.p.f. of \(C_3\):(1, 0, 0), (0, 1, 0), (0, 0, 1), (1, 1, 0), (0, 1, 1), (1, 1, 1), (0, 2, 1), (1, 2, 1), (2, 2, 1)}\\\hline N & Polynomial/Lattice & Shift(s)\\ \hline
	\endfirsthead\multicolumn{3}{c}{{\bfseries \tablename\ \thetable{} -- continued from previous page}} \\
	\hline  N &  Polynomial/Lattice & Shift(s) \\ \hline
	\endhead\hline \multicolumn{3}{|c|}{{Continued on next page}} \\ \hline\endfoot\endlastfoot\multirow{5}{*}{1}&\(\Lambda=\langle((2, 0, 0),(0, 2, 0),(0, 0, 1))\rangle\)&\\
\hline\end{longtable}

\begin{longtable}{|ccccc|}\caption{\footnotesize V.p.f. of \(C_3\):(1, 0, 0), (0, 1, 0), (0, 0, 1), (1, 1, 0), (0, 1, 1), (1, 1, 1), (0, 2, 1), (1, 2, 1), (2, 2, 1)}\\\hline N & Defining inequalities & Vertices& Int. Pt.& Neighbors \\ \hline
	\endfirsthead\multicolumn{5}{c}{{\bfseries \tablename\ \thetable{} -- continued from previous page}} \\
	\hline N & Defining inequalities & Vertices&Int. Pt. & Neighbors\\ \hline
	\endhead\hline \multicolumn{5}{|c|}{{Continued on next page}} \\ \hline\endfoot\endlastfoot1&
	\(\begin{array}{rcl}-x_{1}+x_{2}-x_{3}&\geq& 0\\
		x_{1}+x_{2}-2x_{3}&\geq& 0\\
		-x_{2}+2x_{3}&\geq& 0\\
	\end{array}\)&\(\begin{array}{l}(1, 3, 2): 14\checkmark\\(1, 2, 1): 7\checkmark\\(0, 2, 1): 3\checkmark\end{array}\)&(2, 7, 4): 100\checkmark& [45, 27, 36]\\\hline
	2&
	\(\begin{array}{rcl}-x_{1}+x_{2}-x_{3}&\geq& 0\\
		x_{1}-x_{2}+2x_{3}&\geq& 0\\
		x_{1}-2x_{3}&\geq& 0\\
	\end{array}\)&\(\begin{array}{l}(1, 1, 0): 2\checkmark\\(2, 3, 1): 13\checkmark\\(2, 4, 1): 14\checkmark\end{array}\)&(5, 8, 2): 97\checkmark& [57, 49, 54]\\\hline
	3&
	\(\begin{array}{rcl}x_{1}-x_{3}&\geq& 0\\
		-x_{1}+x_{2}-2x_{3}&\geq& 0\\
		-x_{1}+2x_{3}&\geq& 0\\
	\end{array}\)&\(\begin{array}{l}(1, 3, 1): 8\checkmark\\(0, 1, 0): 1\checkmark\\(2, 4, 1): 14\checkmark\end{array}\)&(3, 8, 2): 56\checkmark& [33, 54, 49]\\\hline
	4&
	\(\begin{array}{rcl}x_{1}-x_{3}&\geq& 0\\
		-x_{2}+2x_{3}&\geq& 0\\
		-x_{1}+2x_{2}-2x_{3}&\geq& 0\\
	\end{array}\)&\(\begin{array}{l}(1, 2, 1): 7\checkmark\\(2, 2, 1): 10\checkmark\\(2, 3, 2): 23\checkmark\end{array}\)&(5, 7, 4): 290& [60, 58, 63]\\\hline
	5&
	\(\begin{array}{rcl}x_{2}-x_{3}&\geq& 0\\
		-2x_{1}+x_{2}&\geq& 0\\
		x_{1}-2x_{2}+2x_{3}&\geq& 0\\
	\end{array}\)&\(\begin{array}{l}(1, 2, 2): 9\checkmark\\(0, 1, 1): 2\checkmark\\(2, 4, 3): 43\checkmark\end{array}\)&(3, 7, 6): 220& [37, 42, 48]\\\hline
	6&
	\(\begin{array}{rcl}x_{2}-x_{3}&\geq& 0\\
		-x_{1}-x_{2}+2x_{3}&\geq& 0\\
		2x_{1}-x_{2}&\geq& 0\\
	\end{array}\)&\(\begin{array}{l}(1, 1, 1): 4\checkmark\\(1, 2, 2): 9\checkmark\\(2, 4, 3): 43\checkmark\end{array}\)&(4, 7, 6): 307& [38, 59, 47]\\\hline
	7&
	\(\begin{array}{rcl}x_{1}-x_{2}&\geq& 0\\
		x_{2}-x_{3}&\geq& 0\\
		-x_{2}+2x_{3}&\geq& 0\\
	\end{array}\)&\(\begin{array}{l}(1, 1, 1): 4\checkmark\\(1, 0, 0): 1\checkmark\\(2, 2, 1): 10\checkmark\end{array}\)&(4, 3, 2): 28\checkmark& [63, 3, 17]\\
\hline
\end{longtable}

\allowdisplaybreaks\begin{align*}&~~~
	\frac{1}{(1-x_{1} ) (1-x_{2} ) (1-x_{3} ) (1-x_{1} x_{2} ) (1-x_{2} x_{3} ) (1-x_{1} x_{2} x_{3} ) (1-x_{2}^2x_{3} ) (1-x_{1} x_{2}^2x_{3} ) (1-x_{1}^2x_{2}^2x_{3} ) }\\=&~~~
	\displaystyle \frac{x_{1} x_{2}^{-4}}{(1-x_{1} )^4(1-x_{1}^2) (1-x_{1} x_{2} ) (1-x_{3} )^3 }\\&
	+\displaystyle \frac{x_{2}^{-6}x_{3}^{-5}}{(1-x_{1} )^4(1-x_{1}^2) (1-x_{2} )^2(1-x_{2}^2) (1-x_{3} ) }\\&
	+\displaystyle \frac{-x_{2}^{-6}x_{3}^{-5}}{(1-x_{1} )^4(1-x_{1}^2) (1-x_{2} )^2(1-x_{2}^2) (1-x_{2}^2x_{3} ) }\\&
	+\displaystyle \frac{-x_{1}^4x_{2}^{-1}x_{3} -x_{1}^3x_{2}^{-2}-x_{1}^2x_{2}^{-3}-x_{1} x_{2}^{-4}}{(1-x_{1} )^4(1-x_{1}^2) (1-x_{1}^2x_{2}^2x_{3} ) (1-x_{3} )^3 }\\&
	+\displaystyle \frac{-x_{2}^{-8}x_{3}^{-5}}{(1-x_{1} )^4(1-x_{1}^2)^2 (1-x_{2} )^2 (1-x_{3} ) }\\&
	+\displaystyle \frac{x_{2}^{-8}x_{3}^{-5}}{(1-x_{1} )^4(1-x_{1}^2)^2 (1-x_{2} )^2 (1-x_{2} x_{3} ) }\\&
	+\displaystyle \frac{x_{2}^{-9}x_{3}^{-5}+x_{2}^{-10}x_{3}^{-5}}{(1-x_{1} )^4(1-x_{1}^2)^3 (1-x_{2} ) (1-x_{3} ) }\\&
	+\displaystyle \frac{x_{1}^{-1}x_{2}^{-5}+x_{1}^{-1}x_{2}^{-6}x_{3}^{-1}+2x_{1}^{-1}x_{2}^{-8}x_{3}^{-4}+3x_{1}^{-1}x_{2}^{-9}x_{3}^{-4}+x_{1}^{-1}x_{2}^{-10}x_{3}^{-4}}{(1-x_{1} )^7 (1-x_{1}^2x_{2}^2x_{3} ) (1-x_{3} ) }\\&
	+\displaystyle \frac{~~~~\begin{array}{l}-x_{1}^2x_{2}^{-3}-x_{1}^2x_{2}^{-3}x_{3}^{-1}-x_{1}^2x_{2}^{-3}x_{3}^{-2}-x_{1} x_{2}^{-4}-x_{1}^2x_{2}^{-4}x_{3}^{-2}-x_{1} x_{2}^{-4}x_{3}^{-2} \\ \hline  -x_{1} x_{2}^{-4}x_{3}^{-3}-x_{1} x_{2}^{-5}x_{3}^{-3}-x_{2}^{-5}x_{3}^{-2}-x_{2}^{-6}x_{3}^{-3}-x_{1}^{-1}x_{2}^{-6}x_{3}^{-3} \\ \hline  -x_{1}^{-1}x_{2}^{-7}x_{3}^{-4}\end{array}~~~~}{(1-x_{1} )^6(1-x_{1}^2) (1-x_{1} x_{2}^2x_{3} ) (1-x_{3} ) }\\&
	+\displaystyle \frac{-x_{2}^{-10}x_{3}^{-5}}{(1-x_{1} )^4(1-x_{1}^2)^3 (1-x_{2} ) (1-x_{1}^2x_{2}^2x_{3} ) }\\&
	+\displaystyle \frac{~~~~\begin{array}{l}-x_{1} x_{2}^{-5}x_{3}^{-1}-x_{1} x_{2}^{-5}x_{3}^{-3}-x_{2}^{-6}x_{3}^{-1}-x_{2}^{-6}x_{3}^{-4}-x_{2}^{-8}x_{3}^{-4}-2x_{2}^{-9}x_{3}^{-4} \\ \hline  -x_{2}^{-10}x_{3}^{-4}\end{array}~~~~}{(1-x_{1} )^4(1-x_{1}^2)^3 (1-x_{1}^2x_{2}^2x_{3} ) (1-x_{3} ) }\\&
	+\displaystyle \frac{~~~~\begin{array}{l}-x_{1}^2x_{2}^{-3}+x_{1}^2x_{2}^{-3}x_{3}^{-1}-x_{1} x_{2}^{-4}+x_{1} x_{2}^{-4}x_{3}^{-2}-x_{2}^{-5}x_{3}^{-1}+x_{2}^{-5}x_{3}^{-2} \\ \hline  -x_{1}^{-1}x_{2}^{-6}x_{3}^{-1}+x_{1}^{-1}x_{2}^{-6}x_{3}^{-3}\end{array}~~~~}{(1-x_{1} )^5(1-x_{1}^2) (1-x_{1} x_{2}^2x_{3} ) (1-x_{3} )^2 }\\&
	+\displaystyle \frac{-x_{2}^{-8}x_{3}^{-5}}{(1-x_{1} )^5(1-x_{1}^2) (1-x_{2} )^2 (1-x_{3} ) }\\&
	+\displaystyle \frac{x_{1}^{-1}x_{2}^{-4}+x_{1}^{-1}x_{2}^{-6}x_{3}^{-3}+x_{1}^{-1}x_{2}^{-8}x_{3}^{-4}}{(1-x_{1} )^7 (1-x_{2} x_{3} ) (1-x_{1} x_{2}^2x_{3} ) }\\&
	+\displaystyle \frac{x_{1} x_{2}^{-4}}{(1-x_{1} )^5(1-x_{1}^2) (1-x_{1} x_{2} ) (1-x_{3} )^2 }\\&
	+\displaystyle \frac{2x_{2}^{-9}x_{3}^{-5}+x_{2}^{-10}x_{3}^{-5}}{(1-x_{1} )^5(1-x_{1}^2)^2 (1-x_{2} ) (1-x_{3} ) }\\&
	+\displaystyle \frac{-x_{2}^{-9}x_{3}^{-5}-x_{2}^{-10}x_{3}^{-5}}{(1-x_{1} )^5(1-x_{1}^2)^2 (1-x_{2} ) (1-x_{1}^2x_{2}^2x_{3} ) }\\&
	+\displaystyle \frac{~~~~\begin{array}{l}x_{1}^2x_{2}^{-3}+x_{1}^2x_{2}^{-3}x_{3}^{-2}+x_{1} x_{2}^{-4}-x_{2}^{-5}+x_{1} x_{2}^{-4}x_{3}^{-3}-x_{1} x_{2}^{-5}x_{3}^{-2} \\ \hline  +x_{2}^{-5}x_{3}^{-1}-x_{2}^{-6}x_{3}^{-1}-x_{1} x_{2}^{-6}x_{3}^{-3}+x_{2}^{-5}x_{3}^{-3}+x_{1}^{-1}x_{2}^{-6}x_{3}^{-1} \\ \hline  -x_{2}^{-6}x_{3}^{-3}-x_{2}^{-7}x_{3}^{-4}+x_{1}^{-1}x_{2}^{-6}x_{3}^{-4}-2x_{2}^{-8}x_{3}^{-4}-3x_{2}^{-9}x_{3}^{-4}-x_{2}^{-10}x_{3}^{-4}\end{array}~~~~}{(1-x_{1} )^5(1-x_{1}^2)^2 (1-x_{1}^2x_{2}^2x_{3} ) (1-x_{3} ) }\\&
	+\displaystyle \frac{x_{1}^{-1}x_{2}^{-6}x_{3}^{-4}}{(1-x_{1} )^6 (1-x_{2} )^2 (1-x_{1} x_{2}^2x_{3} ) }\\&
	+\displaystyle \frac{2x_{2}^{-9}x_{3}^{-5}+x_{2}^{-10}x_{3}^{-5}}{(1-x_{1} )^6(1-x_{1}^2) (1-x_{2} ) (1-x_{3} ) }\\&
	+\displaystyle \frac{-x_{1}^{-1}x_{2}^{-6}x_{3}^{-4}-x_{1}^{-1}x_{2}^{-9}x_{3}^{-5}-x_{1}^{-1}x_{2}^{-10}x_{3}^{-5}}{(1-x_{1} )^7 (1-x_{2} ) (1-x_{1} x_{2}^2x_{3} ) }\\&
	+\displaystyle \frac{x_{1}^{-1}x_{2}^{-6}x_{3}^{-4}-x_{1}^{-1}x_{2}^{-7}x_{3}^{-4}+x_{1}^{-1}x_{2}^{-9}x_{3}^{-5}+x_{1}^{-1}x_{2}^{-10}x_{3}^{-5}}{(1-x_{1} )^7 (1-x_{2} ) (1-x_{1}^2x_{2}^2x_{3} ) }\\&
	+\displaystyle \frac{-x_{2}^{-6}x_{3}^{-4}-x_{2}^{-9}x_{3}^{-5}-x_{2}^{-10}x_{3}^{-5}}{(1-x_{1} )^6(1-x_{1}^2) (1-x_{2} ) (1-x_{1}^2x_{2}^2x_{3} ) }\\&
	+\displaystyle \frac{~~~~\begin{array}{l}x_{1}^2x_{2}^{-3}+x_{1}^2x_{2}^{-3}x_{3}^{-1}+x_{1}^2x_{2}^{-3}x_{3}^{-2}+x_{1} x_{2}^{-4}+x_{1}^2x_{2}^{-4}x_{3}^{-2}+x_{1} x_{2}^{-4}x_{3}^{-2} \\ \hline  -x_{2}^{-5}+x_{1} x_{2}^{-4}x_{3}^{-3}+x_{1} x_{2}^{-5}x_{3}^{-3}+x_{2}^{-5}x_{3}^{-2}-x_{2}^{-6}x_{3}^{-1}+x_{2}^{-6}x_{3}^{-3} \\ \hline  +x_{1}^{-1}x_{2}^{-6}x_{3}^{-3}-2x_{2}^{-8}x_{3}^{-4}+x_{1}^{-1}x_{2}^{-7}x_{3}^{-4}-3x_{2}^{-9}x_{3}^{-4}-x_{2}^{-10}x_{3}^{-4}\end{array}~~~~}{(1-x_{1} )^6(1-x_{1}^2) (1-x_{1}^2x_{2}^2x_{3} ) (1-x_{3} ) }\\&
	+\displaystyle \frac{-x_{1}^{-1}x_{2}^{-4}+x_{1}^{-1}x_{2}^{-5}x_{3}^{-2}-x_{1}^{-1}x_{2}^{-6}x_{3}^{-3}-x_{1}^{-1}x_{2}^{-8}x_{3}^{-4}}{(1-x_{1} )^7 (1-x_{2} x_{3} ) (1-x_{1}^2x_{2}^2x_{3} ) }\\&
	+\displaystyle \frac{-x_{1}^{-1}x_{2}^{-5}-x_{1}^{-1}x_{2}^{-6}x_{3}^{-1}-2x_{1}^{-1}x_{2}^{-8}x_{3}^{-4}-3x_{1}^{-1}x_{2}^{-9}x_{3}^{-4}-x_{1}^{-1}x_{2}^{-10}x_{3}^{-4}}{(1-x_{1} )^7 (1-x_{1} x_{2}^2x_{3} ) (1-x_{3} ) }\\&
	+\displaystyle \frac{-x_{1}^2x_{2}^{-3}x_{3}^{-2}-x_{1} x_{2}^{-4}x_{3}^{-1}-x_{2}^{-5}x_{3}^{-3}}{(1-x_{1} )^6(1-x_{1}^2) (1-x_{1} x_{2} ) (1-x_{1} x_{2}^2x_{3} ) }\\&
	+\displaystyle \frac{x_{1} x_{2}^{-5}x_{3}^{-1}-x_{1} x_{2}^{-5}x_{3}^{-2}+x_{2}^{-6}x_{3}^{-1}-x_{2}^{-6}x_{3}^{-3}}{(1-x_{1} )^4(1-x_{1}^2)^2 (1-x_{2}^2x_{3} ) (1-x_{3} )^2 }\\&
	+\displaystyle \frac{~~~~\begin{array}{l}x_{1}^4x_{2}^{-1}x_{3} +x_{1}^4x_{2}^{-1}+x_{1}^4x_{2}^{-2}+x_{1}^3x_{2}^{-2}+x_{1}^3x_{2}^{-2}x_{3}^{-1}+x_{1}^3x_{2}^{-3}x_{3}^{-1} \\ \hline  -x_{1}^2x_{2}^{-3}x_{3}^{-1}-x_{1} x_{2}^{-4}x_{3}^{-2}-x_{2}^{-5}+2x_{2}^{-5}x_{3}^{-1}-x_{2}^{-5}x_{3}^{-2}-x_{2}^{-6}x_{3}^{-1} \\ \hline  +x_{2}^{-6}x_{3}^{-2}+x_{1}^{-1}x_{2}^{-6}x_{3}^{-1}-x_{1}^{-1}x_{2}^{-6}x_{3}^{-3}\end{array}~~~~}{(1-x_{1} )^5(1-x_{1}^2) (1-x_{1}^2x_{2}^2x_{3} ) (1-x_{3} )^2 }\\&
	+\displaystyle \frac{~~~~\begin{array}{l}x_{1}^4x_{2}^{-1}+x_{1}^3x_{2}^{-2}x_{3}^{-1}-x_{1}^2x_{2}^{-3}-x_{1} x_{2}^{-4}-x_{1} x_{2}^{-5}x_{3}^{-1}+x_{1} x_{2}^{-5}x_{3}^{-2} \\ \hline  -x_{2}^{-6}x_{3}^{-1}+x_{2}^{-6}x_{3}^{-3}\end{array}~~~~}{(1-x_{1} )^4(1-x_{1}^2)^2 (1-x_{1}^2x_{2}^2x_{3} ) (1-x_{3} )^2 }\\&
	+\displaystyle \frac{x_{1}^2x_{2}^{-3}x_{3}^{-2}+x_{1} x_{2}^{-4}x_{3}^{-1}+x_{1} x_{2}^{-5}x_{3}^{-3}+x_{2}^{-5}x_{3}^{-3}}{(1-x_{1} )^6(1-x_{1}^2) (1-x_{1}^2x_{2}^2x_{3} ) (1-x_{1} x_{2} ) }\\&
	+\displaystyle \frac{x_{1} x_{2}^{-4}}{(1-x_{1} )^4(1-x_{1}^2)^2 (1-x_{1} x_{2} ) (1-x_{3} )^2 }\\&
	+\displaystyle \frac{-x_{2}^{-7}x_{3}^{-5}}{(1-x_{1} )^4(1-x_{1}^2) (1-x_{2} )^3 (1-x_{2} x_{3} ) }\\&
	+\displaystyle \frac{x_{2}^{-7}x_{3}^{-5}}{(1-x_{1} )^4(1-x_{1}^2) (1-x_{2} )^3 (1-x_{2}^2x_{3} ) }\\&
	+\displaystyle \frac{-x_{1}^{-1}x_{2}^{-6}x_{3}^{-4}}{(1-x_{1} )^6 (1-x_{2} )^2 (1-x_{1} x_{2} x_{3} ) }\\&
	+\displaystyle \frac{x_{1}^{-1}x_{2}^{-7}x_{3}^{-4}}{(1-x_{1} )^7 (1-x_{2} ) (1-x_{1} x_{2} x_{3} ) }\\&
	+\displaystyle \frac{-x_{2}^{-6}x_{3}^{-3}+x_{1}^{-1}x_{2}^{-7}x_{3}^{-4}}{(1-x_{1} )^7 (1-x_{1}^2x_{2}^2x_{3} ) (1-x_{1} x_{2} x_{3} ) }\\&
	+\displaystyle \frac{~~~~\begin{array}{l}-x_{1}^2x_{2}^{-3}-x_{1}^2x_{2}^{-3}x_{3}^{-2}-x_{1} x_{2}^{-4}-x_{1} x_{2}^{-4}x_{3}^{-3}-x_{2}^{-5}x_{3}^{-1}-x_{2}^{-5}x_{3}^{-3} \\ \hline  -x_{1}^{-1}x_{2}^{-6}x_{3}^{-1}-x_{1}^{-1}x_{2}^{-6}x_{3}^{-4}\end{array}~~~~}{(1-x_{1} )^5(1-x_{1}^2)^2 (1-x_{1} x_{2}^2x_{3} ) (1-x_{3} ) }\\&
	+\displaystyle \frac{-x_{1} x_{2}^{-4}x_{3}^{-1}-x_{1}^{-1}x_{2}^{-6}x_{3}^{-2}}{(1-x_{1} )^5(1-x_{1}^2)^2 (1-x_{1} x_{2} ) (1-x_{1} x_{2}^2x_{3} ) }\\&
	+\displaystyle \frac{x_{1} x_{2}^{-4}x_{3}^{-1}-x_{1} x_{2}^{-5}x_{3}^{-3}+x_{1}^{-1}x_{2}^{-6}x_{3}^{-2}}{(1-x_{1} )^5(1-x_{1}^2)^2 (1-x_{1}^2x_{2}^2x_{3} ) (1-x_{1} x_{2} ) }\\&
	+\displaystyle \frac{x_{2}^{-4}+x_{2}^{-8}x_{3}^{-4}}{(1-x_{1} )^5(1-x_{1}^2)^2 (1-x_{1}^2x_{2}^2x_{3} ) (1-x_{2} x_{3} ) }\\&
	+\displaystyle \frac{x_{2}^{-6}x_{3}^{-4}}{(1-x_{1} )^6(1-x_{1}^2) (1-x_{2} ) (1-x_{2}^2x_{3} ) }\\&
	+\displaystyle \frac{x_{2}^{-7}x_{3}^{-4}}{(1-x_{1} )^6(1-x_{1}^2) (1-x_{1} x_{2} x_{3} ) (1-x_{2}^2x_{3} ) }\\&
	+\displaystyle \frac{-x_{1}^{-1}x_{2}^{-7}x_{3}^{-4}}{(1-x_{1} )^7 (1-x_{1} x_{2} x_{3} ) (1-x_{1} x_{2}^2x_{3} ) }\\&
	+\displaystyle \frac{-x_{1} x_{2}^{-4}x_{3}^{-3}-x_{1} x_{2}^{-5}x_{3}^{-3}-x_{2}^{-7}x_{3}^{-4}-x_{1}^{-1}x_{2}^{-7}x_{3}^{-4}}{(1-x_{1} )^6(1-x_{1}^2) (1-x_{1} x_{2} x_{3} ) (1-x_{1}^2x_{2}^2x_{3} ) }\\&
	+\displaystyle \frac{x_{1} x_{2}^{-5}x_{3}^{-1}+x_{1} x_{2}^{-5}x_{3}^{-3}+x_{2}^{-6}x_{3}^{-1}+x_{2}^{-6}x_{3}^{-4}}{(1-x_{1} )^4(1-x_{1}^2)^3 (1-x_{2}^2x_{3} ) (1-x_{3} ) }\\&
	+\displaystyle \frac{x_{2}^{-6}x_{3}^{-2}}{(1-x_{1} )^4(1-x_{1}^2)^3 (1-x_{1} x_{2} ) (1-x_{2}^2x_{3} ) }\\&
	+\displaystyle \frac{-x_{2}^{-6}x_{3}^{-2}}{(1-x_{1} )^4(1-x_{1}^2)^3 (1-x_{1}^2x_{2}^2x_{3} ) (1-x_{1} x_{2} ) }\\&
	+\displaystyle \frac{x_{1}^3x_{2}^{-2}x_{3} }{(1-x_{1} )^4(1-x_{1}^2) (1-x_{1} x_{2} x_{3} ) (1-x_{3} )^3 }\\&
	+\displaystyle \frac{-x_{2}^{-7}x_{3}^{-5}}{(1-x_{1} )^5(1-x_{1}^2) (1-x_{2} )(1-x_{2}^2) (1-x_{3} ) }\\&
	+\displaystyle \frac{x_{2}^{-7}x_{3}^{-5}}{(1-x_{1} )^5(1-x_{1}^2) (1-x_{2} )(1-x_{2}^2) (1-x_{2}^2x_{3} ) }\\&
	+\displaystyle \frac{x_{1} x_{2}^{-4}x_{3}^{-3}+x_{1} x_{2}^{-5}x_{3}^{-3}+x_{1}^{-1}x_{2}^{-7}x_{3}^{-4}}{(1-x_{1} )^6(1-x_{1}^2) (1-x_{1} x_{2} x_{3} ) (1-x_{1} x_{2}^2x_{3} ) }\\&
	+\displaystyle \frac{-x_{1}^3x_{2}^{-2}x_{3} -x_{1}^3x_{2}^{-2}-x_{1}^3x_{2}^{-3}}{(1-x_{1} )^5(1-x_{1}^2) (1-x_{1} x_{2} x_{3} ) (1-x_{3} )^2 }\\&
	+\displaystyle \frac{-x_{1}^3x_{2}^{-2}}{(1-x_{1} )^4(1-x_{1}^2)^2 (1-x_{1} x_{2} x_{3} ) (1-x_{3} )^2 }\\&
	+\displaystyle \frac{x_{1} x_{2}^{-4}x_{3}^{-3}+x_{1}^{-1}x_{2}^{-6}x_{3}^{-4}}{(1-x_{1} )^5(1-x_{1}^2)^2 (1-x_{1} x_{2} x_{3} ) (1-x_{1} x_{2}^2x_{3} ) }\\&
	+\displaystyle \frac{-x_{1} x_{2}^{-4}x_{3}^{-3}+x_{2}^{-7}x_{3}^{-4}-x_{1}^{-1}x_{2}^{-6}x_{3}^{-4}}{(1-x_{1} )^5(1-x_{1}^2)^2 (1-x_{1}^2x_{2}^2x_{3} ) (1-x_{1} x_{2} x_{3} ) }\\&
	+\displaystyle \frac{-x_{2}^{-6}x_{3}^{-4}}{(1-x_{1} )^4(1-x_{1}^2)^3 (1-x_{2}^2x_{3} ) (1-x_{1} x_{2} x_{3} ) }\\&
	+\displaystyle \frac{x_{2}^{-6}x_{3}^{-4}}{(1-x_{1} )^4(1-x_{1}^2)^3 (1-x_{1}^2x_{2}^2x_{3} ) (1-x_{1} x_{2} x_{3} ) }\\&
	+\displaystyle \frac{x_{1} x_{2}^{-5}x_{3}^{-2}+x_{1} x_{2}^{-6}x_{3}^{-3}+x_{2}^{-6}x_{3}^{-3}+x_{2}^{-7}x_{3}^{-4}}{(1-x_{1} )^5(1-x_{1}^2)^2 (1-x_{2}^2x_{3} ) (1-x_{3} ) }\\&
	+\displaystyle \frac{-x_{2}^{-7}x_{3}^{-4}}{(1-x_{1} )^5(1-x_{1}^2)^2 (1-x_{2}^2x_{3} ) (1-x_{1} x_{2} x_{3} ) }\\&
	+\displaystyle \frac{x_{2}^{-4}+x_{2}^{-6}x_{3}^{-3}+x_{2}^{-8}x_{3}^{-4}}{(1-x_{1} )^6(1-x_{1}^2) (1-x_{1}^2x_{2}^2x_{3} ) (1-x_{2} x_{3} ) }\\&
	+\displaystyle \frac{-x_{2}^{-9}x_{3}^{-5}}{(1-x_{1} )^6(1-x_{1}^2) (1-x_{2} ) (1-x_{2} x_{3} ) }\\&
	+\displaystyle \frac{-x_{2}^{-9}x_{3}^{-5}}{(1-x_{1} )^4(1-x_{1}^2)^3 (1-x_{2} ) (1-x_{2} x_{3} ) }\\&
	+\displaystyle \frac{x_{2}^{-8}x_{3}^{-4}}{(1-x_{1} )^4(1-x_{1}^2)^3 (1-x_{1}^2x_{2}^2x_{3} ) (1-x_{2} x_{3} ) }\\&
	+\displaystyle \frac{x_{2}^{-8}x_{3}^{-5}}{(1-x_{1} )^5(1-x_{1}^2) (1-x_{2} )^2 (1-x_{2} x_{3} ) }\\&
	+\displaystyle \frac{-x_{2}^{-9}x_{3}^{-5}}{(1-x_{1} )^5(1-x_{1}^2)^2 (1-x_{2} ) (1-x_{2} x_{3} ) }\\&
	+\displaystyle \frac{x_{2}^{-6}x_{3}^{-1}}{(1-x_{1} )^5(1-x_{1}^2) (1-x_{2} x_{3} ) (1-x_{3} )^2 }\\&
	+\displaystyle \frac{x_{2}^{-5}x_{3}^{-3}}{(1-x_{1} )^6 (1-x_{1} x_{2} )^2 (1-x_{2} x_{3} ) }\\&
	+\displaystyle \frac{-x_{2}^{-5}x_{3}^{-3}}{(1-x_{1} )^6 (1-x_{1} x_{2} )^2 (1-x_{1} x_{2}^2x_{3} ) }\\&
	+\displaystyle \frac{-x_{1}^{-1}x_{2}^{-6}x_{3}^{-3}}{(1-x_{1} )^7 (1-x_{1} x_{2} ) (1-x_{2} x_{3} ) }\\&
	+\displaystyle \frac{-x_{2}^{-5}x_{3}^{-3}+x_{1}^{-1}x_{2}^{-6}x_{3}^{-3}}{(1-x_{1} )^7 (1-x_{1}^2x_{2}^2x_{3} ) (1-x_{1} x_{2} ) }\\&
	+\displaystyle \frac{x_{2}^{-6}x_{3}^{-1}}{(1-x_{1} )^5(1-x_{1}^2)^2 (1-x_{2} x_{3} ) (1-x_{3} ) }\\&
	+\displaystyle \frac{-x_{1}^{-1}x_{2}^{-5}+x_{1}^{-1}x_{2}^{-5}x_{3}^{-1}-x_{1}^{-1}x_{2}^{-6}x_{3}^{-1}+x_{1}^{-1}x_{2}^{-6}x_{3}^{-2}}{(1-x_{1} )^6 (1-x_{1} x_{2}^2x_{3} ) (1-x_{3} )^2 }\\&
	+\displaystyle \frac{x_{2}^{-6}x_{3}^{-1}}{(1-x_{1} )^6(1-x_{1}^2) (1-x_{2} x_{3} ) (1-x_{3} ) }\\&
	+\displaystyle \frac{x_{1}^{-1}x_{2}^{-5}-x_{1}^{-1}x_{2}^{-5}x_{3}^{-1}+x_{1}^{-1}x_{2}^{-6}x_{3}^{-1}-x_{1}^{-1}x_{2}^{-6}x_{3}^{-2}}{(1-x_{1} )^6 (1-x_{1}^2x_{2}^2x_{3} ) (1-x_{3} )^2 }\\&
	+\displaystyle \frac{-x_{1} x_{2}^{-5}x_{3}^{-3}}{(1-x_{1} )^6(1-x_{1}^2) (1-x_{1} x_{2} ) (1-x_{2}^2x_{3} ) }\\&
	+\displaystyle \frac{-x_{2}^{-6}x_{3}^{-3}}{(1-x_{1} )^6(1-x_{1}^2) (1-x_{2}^2x_{3} ) (1-x_{2} x_{3} ) }\\&
	+\displaystyle \frac{x_{2}^{-5}x_{3}^{-3}}{(1-x_{1} )^7 (1-x_{1} x_{2} ) (1-x_{1} x_{2}^2x_{3} ) }\\&
	+\displaystyle \frac{-x_{2}^{-5}x_{3}^{-1}-x_{2}^{-6}x_{3}^{-2}}{(1-x_{1} )^5(1-x_{1}^2) (1-x_{2}^2x_{3} ) (1-x_{3} )^2 }\\&
	+\displaystyle \frac{x_{1}^2x_{2}^{-3}x_{3}^{-2}}{(1-x_{1} )^5(1-x_{1}^2) (1-x_{1} x_{2} )^2 (1-x_{1} x_{2} x_{3} ) }\\&
	+\displaystyle \frac{-x_{1}^2x_{2}^{-3}x_{3}^{-2}}{(1-x_{1} )^5(1-x_{1}^2) (1-x_{1}^2x_{2}^2x_{3} ) (1-x_{1} x_{2} )^2 }\\&
	+\displaystyle \frac{x_{1} x_{2}^{-5}x_{3}^{-3}}{(1-x_{1} )^5(1-x_{1}^2)^2 (1-x_{1} x_{2} ) (1-x_{2}^2x_{3} ) }\\=&~~~
	\displaystyle \frac{x_{1}^5x_{2}^{-4}+4x_{1}^4x_{2}^{-4}+6x_{1}^3x_{2}^{-4}+4x_{1}^2x_{2}^{-4}+x_{1} x_{2}^{-4}}{(1-x_{1}^2)^5 (1-x_{1} x_{2} ) (1-x_{3} )^3 }\\&
	+\displaystyle \frac{~~~~\begin{array}{l}-x_{1}^7x_{2}^{-5}x_{3}^{-3}-5x_{1}^6x_{2}^{-5}x_{3}^{-3}-10x_{1}^5x_{2}^{-5}x_{3}^{-3}-10x_{1}^4x_{2}^{-5}x_{3}^{-3}+x_{1}^4x_{2}^{-6}x_{3}^{-2} \\ \hline  -5x_{1}^3x_{2}^{-5}x_{3}^{-3}+4x_{1}^3x_{2}^{-6}x_{3}^{-2}-x_{1}^2x_{2}^{-5}x_{3}^{-3}+6x_{1}^2x_{2}^{-6}x_{3}^{-2}+4x_{1} x_{2}^{-6}x_{3}^{-2} \\ \hline  +x_{2}^{-6}x_{3}^{-2}\end{array}~~~~}{(1-x_{1}^2)^7 (1-x_{1} x_{2} ) (1-x_{2}^2x_{3} ) }\\&
	+\displaystyle \frac{~~~~\begin{array}{l}-x_{1}^7x_{2}^{-3}x_{3}^{-2}-5x_{1}^6x_{2}^{-3}x_{3}^{-2}-10x_{1}^5x_{2}^{-3}x_{3}^{-2}-10x_{1}^4x_{2}^{-3}x_{3}^{-2}-5x_{1}^3x_{2}^{-3}x_{3}^{-2} \\ \hline  -x_{1}^2x_{2}^{-3}x_{3}^{-2}\end{array}~~~~}{(1-x_{1}^2)^6 (1-x_{1}^2x_{2}^2x_{3} ) (1-x_{1} x_{2} )^2 }\\&
	+\displaystyle \frac{~~~~\begin{array}{l}x_{1}^7x_{2}^{-3}x_{3}^{-2}+5x_{1}^6x_{2}^{-3}x_{3}^{-2}+10x_{1}^5x_{2}^{-3}x_{3}^{-2}+10x_{1}^4x_{2}^{-3}x_{3}^{-2}+5x_{1}^3x_{2}^{-3}x_{3}^{-2} \\ \hline  +x_{1}^2x_{2}^{-3}x_{3}^{-2}\end{array}~~~~}{(1-x_{1}^2)^6 (1-x_{1} x_{2} )^2 (1-x_{1} x_{2} x_{3} ) }\\&
	+\displaystyle \frac{x_{2}^{-5}x_{3}^{-3}}{(1-x_{1} )^7 (1-x_{1} x_{2} ) (1-x_{1} x_{2}^2x_{3} ) }\\&
	+\displaystyle \frac{~~~~\begin{array}{l}-x_{1}^5x_{2}^{-5}x_{3}^{-2}-x_{1}^4x_{2}^{-5}x_{3}^{-1}-x_{1}^5x_{2}^{-6}x_{3}^{-2}-4x_{1}^4x_{2}^{-5}x_{3}^{-2}+x_{1}^4x_{2}^{-6}x_{3}^{-1} \\ \hline  -4x_{1}^3x_{2}^{-5}x_{3}^{-1}-5x_{1}^4x_{2}^{-6}x_{3}^{-2}-6x_{1}^3x_{2}^{-5}x_{3}^{-2}+4x_{1}^3x_{2}^{-6}x_{3}^{-1}-6x_{1}^2x_{2}^{-5}x_{3}^{-1} \\ \hline  -x_{1}^4x_{2}^{-6}x_{3}^{-3}-10x_{1}^3x_{2}^{-6}x_{3}^{-2}-4x_{1}^2x_{2}^{-5}x_{3}^{-2}+6x_{1}^2x_{2}^{-6}x_{3}^{-1}-4x_{1} x_{2}^{-5}x_{3}^{-1} \\ \hline  -4x_{1}^3x_{2}^{-6}x_{3}^{-3}-10x_{1}^2x_{2}^{-6}x_{3}^{-2}-x_{1} x_{2}^{-5}x_{3}^{-2}+4x_{1} x_{2}^{-6}x_{3}^{-1}-x_{2}^{-5}x_{3}^{-1} \\ \hline  -6x_{1}^2x_{2}^{-6}x_{3}^{-3}-5x_{1} x_{2}^{-6}x_{3}^{-2}+x_{2}^{-6}x_{3}^{-1}-4x_{1} x_{2}^{-6}x_{3}^{-3}-x_{2}^{-6}x_{3}^{-2} \\ \hline  -x_{2}^{-6}x_{3}^{-3}\end{array}~~~~}{(1-x_{1}^2)^6 (1-x_{2}^2x_{3} ) (1-x_{3} )^2 }\\&
	+\displaystyle \frac{~~~~\begin{array}{l}-x_{1}^6x_{2}^{-6}x_{3}^{-3}-6x_{1}^5x_{2}^{-6}x_{3}^{-3}-15x_{1}^4x_{2}^{-6}x_{3}^{-3}-20x_{1}^3x_{2}^{-6}x_{3}^{-3}-15x_{1}^2x_{2}^{-6}x_{3}^{-3} \\ \hline  -6x_{1} x_{2}^{-6}x_{3}^{-3}-x_{2}^{-6}x_{3}^{-3}\end{array}~~~~}{(1-x_{1}^2)^7 (1-x_{2}^2x_{3} ) (1-x_{2} x_{3} ) }\\&
	+\displaystyle \frac{x_{1}^{-1}x_{2}^{-5}-x_{1}^{-1}x_{2}^{-5}x_{3}^{-1}+x_{1}^{-1}x_{2}^{-6}x_{3}^{-1}-x_{1}^{-1}x_{2}^{-6}x_{3}^{-2}}{(1-x_{1} )^6 (1-x_{1}^2x_{2}^2x_{3} ) (1-x_{3} )^2 }\\&
	+\displaystyle \frac{-x_{1}^{-1}x_{2}^{-5}+x_{1}^{-1}x_{2}^{-5}x_{3}^{-1}-x_{1}^{-1}x_{2}^{-6}x_{3}^{-1}+x_{1}^{-1}x_{2}^{-6}x_{3}^{-2}}{(1-x_{1} )^6 (1-x_{1} x_{2}^2x_{3} ) (1-x_{3} )^2 }\\&
	+\displaystyle \frac{~~~~\begin{array}{l}x_{1}^6x_{2}^{-6}x_{3}^{-1}+7x_{1}^5x_{2}^{-6}x_{3}^{-1}+20x_{1}^4x_{2}^{-6}x_{3}^{-1}+30x_{1}^3x_{2}^{-6}x_{3}^{-1}+25x_{1}^2x_{2}^{-6}x_{3}^{-1} \\ \hline  +11x_{1} x_{2}^{-6}x_{3}^{-1}+2x_{2}^{-6}x_{3}^{-1}\end{array}~~~~}{(1-x_{1}^2)^7 (1-x_{2} x_{3} ) (1-x_{3} ) }\\&
	+\displaystyle \frac{-x_{2}^{-5}x_{3}^{-3}+x_{1}^{-1}x_{2}^{-6}x_{3}^{-3}}{(1-x_{1} )^7 (1-x_{1}^2x_{2}^2x_{3} ) (1-x_{1} x_{2} ) }\\&
	+\displaystyle \frac{-x_{1}^{-1}x_{2}^{-6}x_{3}^{-3}}{(1-x_{1} )^7 (1-x_{1} x_{2} ) (1-x_{2} x_{3} ) }\\&
	+\displaystyle \frac{-x_{2}^{-5}x_{3}^{-3}}{(1-x_{1} )^6 (1-x_{1} x_{2} )^2 (1-x_{1} x_{2}^2x_{3} ) }\\&
	+\displaystyle \frac{x_{2}^{-5}x_{3}^{-3}}{(1-x_{1} )^6 (1-x_{1} x_{2} )^2 (1-x_{2} x_{3} ) }\\&
	+\displaystyle \frac{~~~~\begin{array}{l}x_{1}^5x_{2}^{-6}x_{3}^{-1}+5x_{1}^4x_{2}^{-6}x_{3}^{-1}+10x_{1}^3x_{2}^{-6}x_{3}^{-1}+10x_{1}^2x_{2}^{-6}x_{3}^{-1}+5x_{1} x_{2}^{-6}x_{3}^{-1} \\ \hline  +x_{2}^{-6}x_{3}^{-1}\end{array}~~~~}{(1-x_{1}^2)^6 (1-x_{2} x_{3} ) (1-x_{3} )^2 }\\&
	+\displaystyle \frac{~~~~\begin{array}{l}-x_{1}^6x_{2}^{-9}x_{3}^{-5}-7x_{1}^5x_{2}^{-9}x_{3}^{-5}-21x_{1}^4x_{2}^{-9}x_{3}^{-5}-34x_{1}^3x_{2}^{-9}x_{3}^{-5}-31x_{1}^2x_{2}^{-9}x_{3}^{-5} \\ \hline  -15x_{1} x_{2}^{-9}x_{3}^{-5}-3x_{2}^{-9}x_{3}^{-5}\end{array}~~~~}{(1-x_{1}^2)^7 (1-x_{2} ) (1-x_{2} x_{3} ) }\\&
	+\displaystyle \frac{~~~~\begin{array}{l}x_{1}^5x_{2}^{-8}x_{3}^{-5}+6x_{1}^4x_{2}^{-8}x_{3}^{-5}+14x_{1}^3x_{2}^{-8}x_{3}^{-5}+16x_{1}^2x_{2}^{-8}x_{3}^{-5}+9x_{1} x_{2}^{-8}x_{3}^{-5} \\ \hline  +2x_{2}^{-8}x_{3}^{-5}\end{array}~~~~}{(1-x_{1}^2)^6 (1-x_{2} )^2 (1-x_{2} x_{3} ) }\\&
	+\displaystyle \frac{~~~~\begin{array}{l}x_{1}^6x_{2}^{-4}+7x_{1}^5x_{2}^{-4}+20x_{1}^4x_{2}^{-4}+30x_{1}^3x_{2}^{-4}+25x_{1}^2x_{2}^{-4}+x_{1}^6x_{2}^{-6}x_{3}^{-3} \\ \hline  +11x_{1} x_{2}^{-4}+6x_{1}^5x_{2}^{-6}x_{3}^{-3}+2x_{2}^{-4}+15x_{1}^4x_{2}^{-6}x_{3}^{-3}+x_{1}^6x_{2}^{-8}x_{3}^{-4}+20x_{1}^3x_{2}^{-6}x_{3}^{-3} \\ \hline  +7x_{1}^5x_{2}^{-8}x_{3}^{-4}+15x_{1}^2x_{2}^{-6}x_{3}^{-3}+21x_{1}^4x_{2}^{-8}x_{3}^{-4}+6x_{1} x_{2}^{-6}x_{3}^{-3}+34x_{1}^3x_{2}^{-8}x_{3}^{-4} \\ \hline  +x_{2}^{-6}x_{3}^{-3}+31x_{1}^2x_{2}^{-8}x_{3}^{-4}+15x_{1} x_{2}^{-8}x_{3}^{-4}+3x_{2}^{-8}x_{3}^{-4}\end{array}~~~~}{(1-x_{1}^2)^7 (1-x_{1}^2x_{2}^2x_{3} ) (1-x_{2} x_{3} ) }\\&
	+\displaystyle \frac{~~~~\begin{array}{l}x_{1}^6x_{2}^{-7}x_{3}^{-4}+5x_{1}^5x_{2}^{-7}x_{3}^{-4}-x_{1}^4x_{2}^{-6}x_{3}^{-4}+10x_{1}^4x_{2}^{-7}x_{3}^{-4}-4x_{1}^3x_{2}^{-6}x_{3}^{-4} \\ \hline  +10x_{1}^3x_{2}^{-7}x_{3}^{-4}-6x_{1}^2x_{2}^{-6}x_{3}^{-4}+5x_{1}^2x_{2}^{-7}x_{3}^{-4}-4x_{1} x_{2}^{-6}x_{3}^{-4}+x_{1} x_{2}^{-7}x_{3}^{-4} \\ \hline  -x_{2}^{-6}x_{3}^{-4}\end{array}~~~~}{(1-x_{1}^2)^7 (1-x_{2}^2x_{3} ) (1-x_{1} x_{2} x_{3} ) }\\&
	+\displaystyle \frac{~~~~\begin{array}{l}x_{1}^6x_{2}^{-5}x_{3}^{-2}+x_{1}^5x_{2}^{-5}x_{3}^{-1}+5x_{1}^5x_{2}^{-5}x_{3}^{-2}+4x_{1}^4x_{2}^{-5}x_{3}^{-1}+x_{1}^6x_{2}^{-6}x_{3}^{-3} \\ \hline  +x_{1}^5x_{2}^{-5}x_{3}^{-3}+10x_{1}^4x_{2}^{-5}x_{3}^{-2}+x_{1}^4x_{2}^{-6}x_{3}^{-1}+6x_{1}^3x_{2}^{-5}x_{3}^{-1}+6x_{1}^5x_{2}^{-6}x_{3}^{-3} \\ \hline  +4x_{1}^4x_{2}^{-5}x_{3}^{-3}+10x_{1}^3x_{2}^{-5}x_{3}^{-2}+4x_{1}^3x_{2}^{-6}x_{3}^{-1}+4x_{1}^2x_{2}^{-5}x_{3}^{-1}+15x_{1}^4x_{2}^{-6}x_{3}^{-3} \\ \hline  +6x_{1}^3x_{2}^{-5}x_{3}^{-3}+5x_{1}^2x_{2}^{-5}x_{3}^{-2}+6x_{1}^2x_{2}^{-6}x_{3}^{-1}+x_{1} x_{2}^{-5}x_{3}^{-1}+x_{1}^5x_{2}^{-7}x_{3}^{-4} \\ \hline  +x_{1}^4x_{2}^{-6}x_{3}^{-4}+20x_{1}^3x_{2}^{-6}x_{3}^{-3}+4x_{1}^2x_{2}^{-5}x_{3}^{-3}+x_{1} x_{2}^{-5}x_{3}^{-2}+4x_{1} x_{2}^{-6}x_{3}^{-1} \\ \hline  +5x_{1}^4x_{2}^{-7}x_{3}^{-4}+4x_{1}^3x_{2}^{-6}x_{3}^{-4}+15x_{1}^2x_{2}^{-6}x_{3}^{-3}+x_{1} x_{2}^{-5}x_{3}^{-3}+x_{2}^{-6}x_{3}^{-1} \\ \hline  +10x_{1}^3x_{2}^{-7}x_{3}^{-4}+6x_{1}^2x_{2}^{-6}x_{3}^{-4}+6x_{1} x_{2}^{-6}x_{3}^{-3}+10x_{1}^2x_{2}^{-7}x_{3}^{-4}+4x_{1} x_{2}^{-6}x_{3}^{-4} \\ \hline  +x_{2}^{-6}x_{3}^{-3}+5x_{1} x_{2}^{-7}x_{3}^{-4}+x_{2}^{-6}x_{3}^{-4}+x_{2}^{-7}x_{3}^{-4}\end{array}~~~~}{(1-x_{1}^2)^7 (1-x_{2}^2x_{3} ) (1-x_{3} ) }\\&
	+\displaystyle \frac{~~~~\begin{array}{l}-x_{1}^7x_{2}^{-4}x_{3}^{-3}-x_{1}^7x_{2}^{-5}x_{3}^{-3}-7x_{1}^6x_{2}^{-4}x_{3}^{-3}-6x_{1}^6x_{2}^{-5}x_{3}^{-3}-20x_{1}^5x_{2}^{-4}x_{3}^{-3} \\ \hline  -15x_{1}^5x_{2}^{-5}x_{3}^{-3}-30x_{1}^4x_{2}^{-4}x_{3}^{-3}-20x_{1}^4x_{2}^{-5}x_{3}^{-3}-25x_{1}^3x_{2}^{-4}x_{3}^{-3}-x_{1}^6x_{2}^{-7}x_{3}^{-4} \\ \hline  -15x_{1}^3x_{2}^{-5}x_{3}^{-3}-11x_{1}^2x_{2}^{-4}x_{3}^{-3}-6x_{1}^5x_{2}^{-7}x_{3}^{-4}-6x_{1}^2x_{2}^{-5}x_{3}^{-3}-2x_{1} x_{2}^{-4}x_{3}^{-3} \\ \hline  -16x_{1}^4x_{2}^{-7}x_{3}^{-4}-x_{1}^3x_{2}^{-6}x_{3}^{-4}-x_{1} x_{2}^{-5}x_{3}^{-3}-25x_{1}^3x_{2}^{-7}x_{3}^{-4}-4x_{1}^2x_{2}^{-6}x_{3}^{-4} \\ \hline  -25x_{1}^2x_{2}^{-7}x_{3}^{-4}-6x_{1} x_{2}^{-6}x_{3}^{-4}-16x_{1} x_{2}^{-7}x_{3}^{-4}-4x_{2}^{-6}x_{3}^{-4}-6x_{2}^{-7}x_{3}^{-4} \\ \hline  -x_{1}^{-1}x_{2}^{-6}x_{3}^{-4}-x_{1}^{-1}x_{2}^{-7}x_{3}^{-4}\end{array}~~~~}{(1-x_{1}^2)^7 (1-x_{1}^2x_{2}^2x_{3} ) (1-x_{1} x_{2} x_{3} ) }\\&
	+\displaystyle \frac{~~~~\begin{array}{l}x_{1}^7x_{2}^{-4}x_{3}^{-3}+x_{1}^7x_{2}^{-5}x_{3}^{-3}+7x_{1}^6x_{2}^{-4}x_{3}^{-3}+6x_{1}^6x_{2}^{-5}x_{3}^{-3}+20x_{1}^5x_{2}^{-4}x_{3}^{-3} \\ \hline  +15x_{1}^5x_{2}^{-5}x_{3}^{-3}+30x_{1}^4x_{2}^{-4}x_{3}^{-3}+20x_{1}^4x_{2}^{-5}x_{3}^{-3}+25x_{1}^3x_{2}^{-4}x_{3}^{-3}+15x_{1}^3x_{2}^{-5}x_{3}^{-3} \\ \hline  +11x_{1}^2x_{2}^{-4}x_{3}^{-3}+x_{1}^5x_{2}^{-7}x_{3}^{-4}+x_{1}^4x_{2}^{-6}x_{3}^{-4}+6x_{1}^2x_{2}^{-5}x_{3}^{-3}+2x_{1} x_{2}^{-4}x_{3}^{-3} \\ \hline  +6x_{1}^4x_{2}^{-7}x_{3}^{-4}+5x_{1}^3x_{2}^{-6}x_{3}^{-4}+x_{1} x_{2}^{-5}x_{3}^{-3}+15x_{1}^3x_{2}^{-7}x_{3}^{-4}+10x_{1}^2x_{2}^{-6}x_{3}^{-4} \\ \hline  +20x_{1}^2x_{2}^{-7}x_{3}^{-4}+10x_{1} x_{2}^{-6}x_{3}^{-4}+15x_{1} x_{2}^{-7}x_{3}^{-4}+5x_{2}^{-6}x_{3}^{-4}+6x_{2}^{-7}x_{3}^{-4} \\ \hline  +x_{1}^{-1}x_{2}^{-6}x_{3}^{-4}+x_{1}^{-1}x_{2}^{-7}x_{3}^{-4}\end{array}~~~~}{(1-x_{1}^2)^7 (1-x_{1} x_{2} x_{3} ) (1-x_{1} x_{2}^2x_{3} ) }\\&
	+\displaystyle \frac{~~~~\begin{array}{l}-x_{1}^8x_{2}^{-2}x_{3} -x_{1}^8x_{2}^{-2}-5x_{1}^7x_{2}^{-2}x_{3} -x_{1}^8x_{2}^{-3}-6x_{1}^7x_{2}^{-2}-10x_{1}^6x_{2}^{-2}x_{3}  \\ \hline  -5x_{1}^7x_{2}^{-3}-14x_{1}^6x_{2}^{-2}-10x_{1}^5x_{2}^{-2}x_{3} -10x_{1}^6x_{2}^{-3}-16x_{1}^5x_{2}^{-2}-5x_{1}^4x_{2}^{-2}x_{3}  \\ \hline  -10x_{1}^5x_{2}^{-3}-9x_{1}^4x_{2}^{-2}-x_{1}^3x_{2}^{-2}x_{3} -5x_{1}^4x_{2}^{-3}-2x_{1}^3x_{2}^{-2}-x_{1}^3x_{2}^{-3}\end{array}~~~~}{(1-x_{1}^2)^6 (1-x_{1} x_{2} x_{3} ) (1-x_{3} )^2 }\\&
	+\displaystyle \frac{~~~~\begin{array}{l}x_{1}^5x_{2}^{-6}x_{3}^{-5}+x_{1}^5x_{2}^{-7}x_{3}^{-5}+5x_{1}^4x_{2}^{-6}x_{3}^{-5}+5x_{1}^4x_{2}^{-7}x_{3}^{-5}+10x_{1}^3x_{2}^{-6}x_{3}^{-5} \\ \hline  +10x_{1}^3x_{2}^{-7}x_{3}^{-5}+10x_{1}^2x_{2}^{-6}x_{3}^{-5}+10x_{1}^2x_{2}^{-7}x_{3}^{-5}+5x_{1} x_{2}^{-6}x_{3}^{-5}+5x_{1} x_{2}^{-7}x_{3}^{-5} \\ \hline  +x_{2}^{-6}x_{3}^{-5}+x_{2}^{-7}x_{3}^{-5}\end{array}~~~~}{(1-x_{1}^2)^6 (1-x_{2}^2)^2 (1-x_{2}^2x_{3} ) }\\&
	+\displaystyle \frac{~~~~\begin{array}{l}-x_{1}^5x_{2}^{-6}x_{3}^{-5}-x_{1}^5x_{2}^{-7}x_{3}^{-5}-5x_{1}^4x_{2}^{-6}x_{3}^{-5}-5x_{1}^4x_{2}^{-7}x_{3}^{-5}-10x_{1}^3x_{2}^{-6}x_{3}^{-5} \\ \hline  -10x_{1}^3x_{2}^{-7}x_{3}^{-5}-10x_{1}^2x_{2}^{-6}x_{3}^{-5}-10x_{1}^2x_{2}^{-7}x_{3}^{-5}-5x_{1} x_{2}^{-6}x_{3}^{-5}-5x_{1} x_{2}^{-7}x_{3}^{-5} \\ \hline  -x_{2}^{-6}x_{3}^{-5}-x_{2}^{-7}x_{3}^{-5}\end{array}~~~~}{(1-x_{1}^2)^6 (1-x_{2}^2)^2 (1-x_{3} ) }\\&
	+\displaystyle \frac{x_{1}^7x_{2}^{-2}x_{3} +4x_{1}^6x_{2}^{-2}x_{3} +6x_{1}^5x_{2}^{-2}x_{3} +4x_{1}^4x_{2}^{-2}x_{3} +x_{1}^3x_{2}^{-2}x_{3} }{(1-x_{1}^2)^5 (1-x_{1} x_{2} x_{3} ) (1-x_{3} )^3 }\\&
	+\displaystyle \frac{~~~~\begin{array}{l}x_{1}^8x_{2}^{-3}x_{3}^{-2}+6x_{1}^7x_{2}^{-3}x_{3}^{-2}+x_{1}^7x_{2}^{-4}x_{3}^{-1}+15x_{1}^6x_{2}^{-3}x_{3}^{-2}+7x_{1}^6x_{2}^{-4}x_{3}^{-1} \\ \hline  +20x_{1}^5x_{2}^{-3}x_{3}^{-2}+20x_{1}^5x_{2}^{-4}x_{3}^{-1}+x_{1}^7x_{2}^{-5}x_{3}^{-3}+15x_{1}^4x_{2}^{-3}x_{3}^{-2}+30x_{1}^4x_{2}^{-4}x_{3}^{-1} \\ \hline  +6x_{1}^6x_{2}^{-5}x_{3}^{-3}+6x_{1}^3x_{2}^{-3}x_{3}^{-2}+25x_{1}^3x_{2}^{-4}x_{3}^{-1}+16x_{1}^5x_{2}^{-5}x_{3}^{-3}+x_{1}^2x_{2}^{-3}x_{3}^{-2} \\ \hline  +11x_{1}^2x_{2}^{-4}x_{3}^{-1}+25x_{1}^4x_{2}^{-5}x_{3}^{-3}+2x_{1} x_{2}^{-4}x_{3}^{-1}+25x_{1}^3x_{2}^{-5}x_{3}^{-3}+x_{1}^3x_{2}^{-6}x_{3}^{-2} \\ \hline  +16x_{1}^2x_{2}^{-5}x_{3}^{-3}+4x_{1}^2x_{2}^{-6}x_{3}^{-2}+6x_{1} x_{2}^{-5}x_{3}^{-3}+6x_{1} x_{2}^{-6}x_{3}^{-2}+x_{2}^{-5}x_{3}^{-3} \\ \hline  +4x_{2}^{-6}x_{3}^{-2}+x_{1}^{-1}x_{2}^{-6}x_{3}^{-2}\end{array}~~~~}{(1-x_{1}^2)^7 (1-x_{1}^2x_{2}^2x_{3} ) (1-x_{1} x_{2} ) }\\&
	+\displaystyle \frac{-x_{1}^{-1}x_{2}^{-7}x_{3}^{-4}}{(1-x_{1} )^7 (1-x_{1} x_{2} x_{3} ) (1-x_{1} x_{2}^2x_{3} ) }\\&
	+\displaystyle \frac{~~~~\begin{array}{l}x_{1}^6x_{2}^{-6}x_{3}^{-4}+6x_{1}^5x_{2}^{-6}x_{3}^{-4}+15x_{1}^4x_{2}^{-6}x_{3}^{-4}+20x_{1}^3x_{2}^{-6}x_{3}^{-4}+15x_{1}^2x_{2}^{-6}x_{3}^{-4} \\ \hline  +6x_{1} x_{2}^{-6}x_{3}^{-4}+x_{2}^{-6}x_{3}^{-4}\end{array}~~~~}{(1-x_{1}^2)^7 (1-x_{2} ) (1-x_{2}^2x_{3} ) }\\&
	+\displaystyle \frac{~~~~\begin{array}{l}-x_{1}^8x_{2}^{-3}x_{3}^{-2}-6x_{1}^7x_{2}^{-3}x_{3}^{-2}-x_{1}^7x_{2}^{-4}x_{3}^{-1}-15x_{1}^6x_{2}^{-3}x_{3}^{-2}-7x_{1}^6x_{2}^{-4}x_{3}^{-1} \\ \hline  -20x_{1}^5x_{2}^{-3}x_{3}^{-2}-20x_{1}^5x_{2}^{-4}x_{3}^{-1}-15x_{1}^4x_{2}^{-3}x_{3}^{-2}-30x_{1}^4x_{2}^{-4}x_{3}^{-1}-x_{1}^6x_{2}^{-5}x_{3}^{-3} \\ \hline  -6x_{1}^3x_{2}^{-3}x_{3}^{-2}-25x_{1}^3x_{2}^{-4}x_{3}^{-1}-6x_{1}^5x_{2}^{-5}x_{3}^{-3}-x_{1}^2x_{2}^{-3}x_{3}^{-2}-11x_{1}^2x_{2}^{-4}x_{3}^{-1} \\ \hline  -15x_{1}^4x_{2}^{-5}x_{3}^{-3}-x_{1}^4x_{2}^{-6}x_{3}^{-2}-2x_{1} x_{2}^{-4}x_{3}^{-1}-20x_{1}^3x_{2}^{-5}x_{3}^{-3}-5x_{1}^3x_{2}^{-6}x_{3}^{-2} \\ \hline  -15x_{1}^2x_{2}^{-5}x_{3}^{-3}-10x_{1}^2x_{2}^{-6}x_{3}^{-2}-6x_{1} x_{2}^{-5}x_{3}^{-3}-10x_{1} x_{2}^{-6}x_{3}^{-2}-x_{2}^{-5}x_{3}^{-3} \\ \hline  -5x_{2}^{-6}x_{3}^{-2}-x_{1}^{-1}x_{2}^{-6}x_{3}^{-2}\end{array}~~~~}{(1-x_{1}^2)^7 (1-x_{1} x_{2} ) (1-x_{1} x_{2}^2x_{3} ) }\\&
	+\displaystyle \frac{~~~~\begin{array}{l}-x_{1}^8x_{2}^{-3}-x_{1}^8x_{2}^{-3}x_{3}^{-1}-7x_{1}^7x_{2}^{-3}-x_{1}^8x_{2}^{-3}x_{3}^{-2}-6x_{1}^7x_{2}^{-3}x_{3}^{-1} \\ \hline  -x_{1}^7x_{2}^{-4}-20x_{1}^6x_{2}^{-3}-x_{1}^8x_{2}^{-4}x_{3}^{-2}-7x_{1}^7x_{2}^{-3}x_{3}^{-2}-15x_{1}^6x_{2}^{-3}x_{3}^{-1} \\ \hline  -7x_{1}^6x_{2}^{-4}-30x_{1}^5x_{2}^{-3}-7x_{1}^7x_{2}^{-4}x_{3}^{-2}-20x_{1}^6x_{2}^{-3}x_{3}^{-2}-20x_{1}^5x_{2}^{-3}x_{3}^{-1} \\ \hline  -20x_{1}^5x_{2}^{-4}-25x_{1}^4x_{2}^{-3}-x_{1}^7x_{2}^{-4}x_{3}^{-3}-21x_{1}^6x_{2}^{-4}x_{3}^{-2}-30x_{1}^5x_{2}^{-3}x_{3}^{-2} \\ \hline  -15x_{1}^4x_{2}^{-3}x_{3}^{-1}-30x_{1}^4x_{2}^{-4}-11x_{1}^3x_{2}^{-3}-x_{1}^7x_{2}^{-5}x_{3}^{-3}-7x_{1}^6x_{2}^{-4}x_{3}^{-3} \\ \hline  -x_{1}^6x_{2}^{-5}x_{3}^{-2}-35x_{1}^5x_{2}^{-4}x_{3}^{-2}-x_{1}^5x_{2}^{-5}x_{3}^{-1}-25x_{1}^4x_{2}^{-3}x_{3}^{-2}-6x_{1}^3x_{2}^{-3}x_{3}^{-1} \\ \hline  -25x_{1}^3x_{2}^{-4}-2x_{1}^2x_{2}^{-3}-6x_{1}^6x_{2}^{-5}x_{3}^{-3}-20x_{1}^5x_{2}^{-4}x_{3}^{-3}-6x_{1}^5x_{2}^{-5}x_{3}^{-2} \\ \hline  -35x_{1}^4x_{2}^{-4}x_{3}^{-2}-5x_{1}^4x_{2}^{-5}x_{3}^{-1}-11x_{1}^3x_{2}^{-3}x_{3}^{-2}-x_{1}^2x_{2}^{-3}x_{3}^{-1}-11x_{1}^2x_{2}^{-4} \\ \hline  -x_{1}^6x_{2}^{-6}x_{3}^{-3}-16x_{1}^5x_{2}^{-5}x_{3}^{-3}-30x_{1}^4x_{2}^{-4}x_{3}^{-3}-15x_{1}^4x_{2}^{-5}x_{3}^{-2}-x_{1}^4x_{2}^{-6}x_{3}^{-1} \\ \hline  -21x_{1}^3x_{2}^{-4}x_{3}^{-2}-10x_{1}^3x_{2}^{-5}x_{3}^{-1}-2x_{1}^2x_{2}^{-3}x_{3}^{-2}-2x_{1} x_{2}^{-4}-7x_{1}^5x_{2}^{-6}x_{3}^{-3} \\ \hline  -25x_{1}^4x_{2}^{-5}x_{3}^{-3}-25x_{1}^3x_{2}^{-4}x_{3}^{-3}-20x_{1}^3x_{2}^{-5}x_{3}^{-2}-5x_{1}^3x_{2}^{-6}x_{3}^{-1}-7x_{1}^2x_{2}^{-4}x_{3}^{-2} \\ \hline  -10x_{1}^2x_{2}^{-5}x_{3}^{-1}-21x_{1}^4x_{2}^{-6}x_{3}^{-3}-25x_{1}^3x_{2}^{-5}x_{3}^{-3}-11x_{1}^2x_{2}^{-4}x_{3}^{-3}-15x_{1}^2x_{2}^{-5}x_{3}^{-2} \\ \hline  -10x_{1}^2x_{2}^{-6}x_{3}^{-1}-x_{1} x_{2}^{-4}x_{3}^{-2}-5x_{1} x_{2}^{-5}x_{3}^{-1}-x_{1}^5x_{2}^{-7}x_{3}^{-4}-x_{1}^4x_{2}^{-6}x_{3}^{-4} \\ \hline  -35x_{1}^3x_{2}^{-6}x_{3}^{-3}-16x_{1}^2x_{2}^{-5}x_{3}^{-3}-2x_{1} x_{2}^{-4}x_{3}^{-3}-6x_{1} x_{2}^{-5}x_{3}^{-2}-10x_{1} x_{2}^{-6}x_{3}^{-1} \\ \hline  -x_{2}^{-5}x_{3}^{-1}-6x_{1}^4x_{2}^{-7}x_{3}^{-4}-5x_{1}^3x_{2}^{-6}x_{3}^{-4}-35x_{1}^2x_{2}^{-6}x_{3}^{-3}-6x_{1} x_{2}^{-5}x_{3}^{-3} \\ \hline  -x_{2}^{-5}x_{3}^{-2}-5x_{2}^{-6}x_{3}^{-1}-15x_{1}^3x_{2}^{-7}x_{3}^{-4}-10x_{1}^2x_{2}^{-6}x_{3}^{-4}-21x_{1} x_{2}^{-6}x_{3}^{-3} \\ \hline  -x_{2}^{-5}x_{3}^{-3}-x_{1}^{-1}x_{2}^{-6}x_{3}^{-1}-20x_{1}^2x_{2}^{-7}x_{3}^{-4}-10x_{1} x_{2}^{-6}x_{3}^{-4}-7x_{2}^{-6}x_{3}^{-3} \\ \hline  -15x_{1} x_{2}^{-7}x_{3}^{-4}-5x_{2}^{-6}x_{3}^{-4}-x_{1}^{-1}x_{2}^{-6}x_{3}^{-3}-6x_{2}^{-7}x_{3}^{-4}-x_{1}^{-1}x_{2}^{-6}x_{3}^{-4} \\ \hline  -x_{1}^{-1}x_{2}^{-7}x_{3}^{-4}\end{array}~~~~}{(1-x_{1}^2)^7 (1-x_{1} x_{2}^2x_{3} ) (1-x_{3} ) }\\&
	+\displaystyle \frac{-x_{2}^{-6}x_{3}^{-3}+x_{1}^{-1}x_{2}^{-7}x_{3}^{-4}}{(1-x_{1} )^7 (1-x_{1}^2x_{2}^2x_{3} ) (1-x_{1} x_{2} x_{3} ) }\\&
	+\displaystyle \frac{x_{1}^{-1}x_{2}^{-7}x_{3}^{-4}}{(1-x_{1} )^7 (1-x_{2} ) (1-x_{1} x_{2} x_{3} ) }\\&
	+\displaystyle \frac{-x_{1}^{-1}x_{2}^{-6}x_{3}^{-4}}{(1-x_{1} )^6 (1-x_{2} )^2 (1-x_{1} x_{2} x_{3} ) }\\&
	+\displaystyle \frac{x_{1}^4x_{2}^{-7}x_{3}^{-5}+4x_{1}^3x_{2}^{-7}x_{3}^{-5}+6x_{1}^2x_{2}^{-7}x_{3}^{-5}+4x_{1} x_{2}^{-7}x_{3}^{-5}+x_{2}^{-7}x_{3}^{-5}}{(1-x_{1}^2)^5 (1-x_{2} )^3 (1-x_{2}^2x_{3} ) }\\&
	+\displaystyle \frac{-x_{1}^4x_{2}^{-7}x_{3}^{-5}-4x_{1}^3x_{2}^{-7}x_{3}^{-5}-6x_{1}^2x_{2}^{-7}x_{3}^{-5}-4x_{1} x_{2}^{-7}x_{3}^{-5}-x_{2}^{-7}x_{3}^{-5}}{(1-x_{1}^2)^5 (1-x_{2} )^3 (1-x_{2} x_{3} ) }\\&
	+\displaystyle \frac{x_{1}^6x_{2}^{-4}+6x_{1}^5x_{2}^{-4}+14x_{1}^4x_{2}^{-4}+16x_{1}^3x_{2}^{-4}+9x_{1}^2x_{2}^{-4}+2x_{1} x_{2}^{-4}}{(1-x_{1}^2)^6 (1-x_{1} x_{2} ) (1-x_{3} )^2 }\\&
	+\displaystyle \frac{~~~~\begin{array}{l}x_{1}^9x_{2}^{-1}x_{3} +x_{1}^9x_{2}^{-1}+5x_{1}^8x_{2}^{-1}x_{3} +x_{1}^9x_{2}^{-2}+6x_{1}^8x_{2}^{-1}+10x_{1}^7x_{2}^{-1}x_{3}  \\ \hline  +6x_{1}^8x_{2}^{-2}+14x_{1}^7x_{2}^{-1}+10x_{1}^6x_{2}^{-1}x_{3} +x_{1}^8x_{2}^{-2}x_{3}^{-1}+15x_{1}^7x_{2}^{-2}+16x_{1}^6x_{2}^{-1} \\ \hline  +5x_{1}^5x_{2}^{-1}x_{3} +x_{1}^8x_{2}^{-3}x_{3}^{-1}+6x_{1}^7x_{2}^{-2}x_{3}^{-1}+20x_{1}^6x_{2}^{-2}+9x_{1}^5x_{2}^{-1} \\ \hline  +x_{1}^4x_{2}^{-1}x_{3} +4x_{1}^7x_{2}^{-3}x_{3}^{-1}+14x_{1}^6x_{2}^{-2}x_{3}^{-1}-x_{1}^6x_{2}^{-3}+15x_{1}^5x_{2}^{-2} \\ \hline  +2x_{1}^4x_{2}^{-1}+5x_{1}^6x_{2}^{-3}x_{3}^{-1}+16x_{1}^5x_{2}^{-2}x_{3}^{-1}-4x_{1}^5x_{2}^{-3}+6x_{1}^4x_{2}^{-2}-x_{1}^5x_{2}^{-4} \\ \hline  +9x_{1}^4x_{2}^{-2}x_{3}^{-1}-6x_{1}^4x_{2}^{-3}+x_{1}^3x_{2}^{-2}-x_{1}^6x_{2}^{-4}x_{3}^{-2}-x_{1}^5x_{2}^{-5}-5x_{1}^4x_{2}^{-3}x_{3}^{-1} \\ \hline  -4x_{1}^4x_{2}^{-4}+2x_{1}^3x_{2}^{-2}x_{3}^{-1}-4x_{1}^3x_{2}^{-3}-5x_{1}^5x_{2}^{-4}x_{3}^{-2}+x_{1}^5x_{2}^{-5}x_{3}^{-1} \\ \hline  -5x_{1}^4x_{2}^{-5}-4x_{1}^3x_{2}^{-3}x_{3}^{-1}-6x_{1}^3x_{2}^{-4}-x_{1}^2x_{2}^{-3}-x_{1}^5x_{2}^{-6}x_{3}^{-1}-10x_{1}^4x_{2}^{-4}x_{3}^{-2} \\ \hline  +6x_{1}^4x_{2}^{-5}x_{3}^{-1}-10x_{1}^3x_{2}^{-5}-x_{1}^2x_{2}^{-3}x_{3}^{-1}-4x_{1}^2x_{2}^{-4}+x_{1}^5x_{2}^{-6}x_{3}^{-2} \\ \hline  -x_{1}^4x_{2}^{-5}x_{3}^{-2}-5x_{1}^4x_{2}^{-6}x_{3}^{-1}-10x_{1}^3x_{2}^{-4}x_{3}^{-2}+14x_{1}^3x_{2}^{-5}x_{3}^{-1}-10x_{1}^2x_{2}^{-5} \\ \hline  -x_{1} x_{2}^{-4}+5x_{1}^4x_{2}^{-6}x_{3}^{-2}-4x_{1}^3x_{2}^{-5}x_{3}^{-2}-9x_{1}^3x_{2}^{-6}x_{3}^{-1}-5x_{1}^2x_{2}^{-4}x_{3}^{-2} \\ \hline  +16x_{1}^2x_{2}^{-5}x_{3}^{-1}-5x_{1} x_{2}^{-5}+10x_{1}^3x_{2}^{-6}x_{3}^{-2}-6x_{1}^2x_{2}^{-5}x_{3}^{-2}-6x_{1}^2x_{2}^{-6}x_{3}^{-1} \\ \hline  -x_{1} x_{2}^{-4}x_{3}^{-2}+9x_{1} x_{2}^{-5}x_{3}^{-1}-x_{2}^{-5}-x_{1}^3x_{2}^{-6}x_{3}^{-3}+10x_{1}^2x_{2}^{-6}x_{3}^{-2} \\ \hline  -4x_{1} x_{2}^{-5}x_{3}^{-2}+x_{1} x_{2}^{-6}x_{3}^{-1}+2x_{2}^{-5}x_{3}^{-1}-4x_{1}^2x_{2}^{-6}x_{3}^{-3}+5x_{1} x_{2}^{-6}x_{3}^{-2} \\ \hline  -x_{2}^{-5}x_{3}^{-2}+3x_{2}^{-6}x_{3}^{-1}-6x_{1} x_{2}^{-6}x_{3}^{-3}+x_{2}^{-6}x_{3}^{-2}+x_{1}^{-1}x_{2}^{-6}x_{3}^{-1} \\ \hline  -4x_{2}^{-6}x_{3}^{-3}-x_{1}^{-1}x_{2}^{-6}x_{3}^{-3}\end{array}~~~~}{(1-x_{1}^2)^6 (1-x_{1}^2x_{2}^2x_{3} ) (1-x_{3} )^2 }\\&
	+\displaystyle \frac{-x_{1}^{-1}x_{2}^{-5}-x_{1}^{-1}x_{2}^{-6}x_{3}^{-1}-2x_{1}^{-1}x_{2}^{-8}x_{3}^{-4}-3x_{1}^{-1}x_{2}^{-9}x_{3}^{-4}-x_{1}^{-1}x_{2}^{-10}x_{3}^{-4}}{(1-x_{1} )^7 (1-x_{1} x_{2}^2x_{3} ) (1-x_{3} ) }\\&
	+\displaystyle \frac{-x_{1}^{-1}x_{2}^{-4}+x_{1}^{-1}x_{2}^{-5}x_{3}^{-2}-x_{1}^{-1}x_{2}^{-6}x_{3}^{-3}-x_{1}^{-1}x_{2}^{-8}x_{3}^{-4}}{(1-x_{1} )^7 (1-x_{2} x_{3} ) (1-x_{1}^2x_{2}^2x_{3} ) }\\&
	+\displaystyle \frac{~~~~\begin{array}{l}x_{1}^8x_{2}^{-3}+x_{1}^8x_{2}^{-3}x_{3}^{-1}+7x_{1}^7x_{2}^{-3}+x_{1}^8x_{2}^{-3}x_{3}^{-2}+6x_{1}^7x_{2}^{-3}x_{3}^{-1} \\ \hline  +x_{1}^7x_{2}^{-4}+20x_{1}^6x_{2}^{-3}+x_{1}^8x_{2}^{-4}x_{3}^{-2}+7x_{1}^7x_{2}^{-3}x_{3}^{-2}+15x_{1}^6x_{2}^{-3}x_{3}^{-1} \\ \hline  +7x_{1}^6x_{2}^{-4}+30x_{1}^5x_{2}^{-3}+7x_{1}^7x_{2}^{-4}x_{3}^{-2}+20x_{1}^6x_{2}^{-3}x_{3}^{-2}-x_{1}^6x_{2}^{-5}+20x_{1}^5x_{2}^{-3}x_{3}^{-1} \\ \hline  +20x_{1}^5x_{2}^{-4}+25x_{1}^4x_{2}^{-3}+x_{1}^7x_{2}^{-4}x_{3}^{-3}+21x_{1}^6x_{2}^{-4}x_{3}^{-2}+30x_{1}^5x_{2}^{-3}x_{3}^{-2} \\ \hline  -7x_{1}^5x_{2}^{-5}+15x_{1}^4x_{2}^{-3}x_{3}^{-1}+30x_{1}^4x_{2}^{-4}+11x_{1}^3x_{2}^{-3}+x_{1}^7x_{2}^{-5}x_{3}^{-3}+7x_{1}^6x_{2}^{-4}x_{3}^{-3} \\ \hline  -x_{1}^6x_{2}^{-6}x_{3}^{-1}+35x_{1}^5x_{2}^{-4}x_{3}^{-2}+25x_{1}^4x_{2}^{-3}x_{3}^{-2}-20x_{1}^4x_{2}^{-5}+6x_{1}^3x_{2}^{-3}x_{3}^{-1} \\ \hline  +25x_{1}^3x_{2}^{-4}+2x_{1}^2x_{2}^{-3}+6x_{1}^6x_{2}^{-5}x_{3}^{-3}+20x_{1}^5x_{2}^{-4}x_{3}^{-3}+x_{1}^5x_{2}^{-5}x_{3}^{-2} \\ \hline  -7x_{1}^5x_{2}^{-6}x_{3}^{-1}+35x_{1}^4x_{2}^{-4}x_{3}^{-2}+x_{1}^4x_{2}^{-5}x_{3}^{-1}+11x_{1}^3x_{2}^{-3}x_{3}^{-2}-30x_{1}^3x_{2}^{-5} \\ \hline  +x_{1}^2x_{2}^{-3}x_{3}^{-1}+11x_{1}^2x_{2}^{-4}+15x_{1}^5x_{2}^{-5}x_{3}^{-3}+30x_{1}^4x_{2}^{-4}x_{3}^{-3}+5x_{1}^4x_{2}^{-5}x_{3}^{-2} \\ \hline  -20x_{1}^4x_{2}^{-6}x_{3}^{-1}+21x_{1}^3x_{2}^{-4}x_{3}^{-2}+4x_{1}^3x_{2}^{-5}x_{3}^{-1}+2x_{1}^2x_{2}^{-3}x_{3}^{-2}-25x_{1}^2x_{2}^{-5} \\ \hline  +2x_{1} x_{2}^{-4}+x_{1}^5x_{2}^{-6}x_{3}^{-3}+21x_{1}^4x_{2}^{-5}x_{3}^{-3}+25x_{1}^3x_{2}^{-4}x_{3}^{-3}+10x_{1}^3x_{2}^{-5}x_{3}^{-2} \\ \hline  -29x_{1}^3x_{2}^{-6}x_{3}^{-1}+7x_{1}^2x_{2}^{-4}x_{3}^{-2}+6x_{1}^2x_{2}^{-5}x_{3}^{-1}-11x_{1} x_{2}^{-5}+6x_{1}^4x_{2}^{-6}x_{3}^{-3} \\ \hline  +19x_{1}^3x_{2}^{-5}x_{3}^{-3}+11x_{1}^2x_{2}^{-4}x_{3}^{-3}+10x_{1}^2x_{2}^{-5}x_{3}^{-2}-21x_{1}^2x_{2}^{-6}x_{3}^{-1}+x_{1} x_{2}^{-4}x_{3}^{-2} \\ \hline  +4x_{1} x_{2}^{-5}x_{3}^{-1}-2x_{2}^{-5}-2x_{1}^6x_{2}^{-8}x_{3}^{-4}+15x_{1}^3x_{2}^{-6}x_{3}^{-3}+12x_{1}^2x_{2}^{-5}x_{3}^{-3} \\ \hline  +2x_{1} x_{2}^{-4}x_{3}^{-3}+5x_{1} x_{2}^{-5}x_{3}^{-2}-5x_{1} x_{2}^{-6}x_{3}^{-1}+x_{2}^{-5}x_{3}^{-1}-3x_{1}^6x_{2}^{-9}x_{3}^{-4} \\ \hline  -14x_{1}^5x_{2}^{-8}x_{3}^{-4}+x_{1}^4x_{2}^{-7}x_{3}^{-4}+x_{1}^3x_{2}^{-6}x_{3}^{-4}+20x_{1}^2x_{2}^{-6}x_{3}^{-3}+5x_{1} x_{2}^{-5}x_{3}^{-3} \\ \hline  +x_{2}^{-5}x_{3}^{-2}+2x_{2}^{-6}x_{3}^{-1}-x_{1}^6x_{2}^{-10}x_{3}^{-4}-21x_{1}^5x_{2}^{-9}x_{3}^{-4}-41x_{1}^4x_{2}^{-8}x_{3}^{-4} \\ \hline  +5x_{1}^3x_{2}^{-7}x_{3}^{-4}+4x_{1}^2x_{2}^{-6}x_{3}^{-4}+15x_{1} x_{2}^{-6}x_{3}^{-3}+x_{2}^{-5}x_{3}^{-3}+x_{1}^{-1}x_{2}^{-6}x_{3}^{-1} \\ \hline  -7x_{1}^5x_{2}^{-10}x_{3}^{-4}-62x_{1}^4x_{2}^{-9}x_{3}^{-4}-64x_{1}^3x_{2}^{-8}x_{3}^{-4}+10x_{1}^2x_{2}^{-7}x_{3}^{-4}+6x_{1} x_{2}^{-6}x_{3}^{-4} \\ \hline  +6x_{2}^{-6}x_{3}^{-3}-21x_{1}^4x_{2}^{-10}x_{3}^{-4}-98x_{1}^3x_{2}^{-9}x_{3}^{-4}-56x_{1}^2x_{2}^{-8}x_{3}^{-4}+10x_{1} x_{2}^{-7}x_{3}^{-4} \\ \hline  +4x_{2}^{-6}x_{3}^{-4}+x_{1}^{-1}x_{2}^{-6}x_{3}^{-3}-34x_{1}^3x_{2}^{-10}x_{3}^{-4}-87x_{1}^2x_{2}^{-9}x_{3}^{-4}-26x_{1} x_{2}^{-8}x_{3}^{-4} \\ \hline  +5x_{2}^{-7}x_{3}^{-4}+x_{1}^{-1}x_{2}^{-6}x_{3}^{-4}-31x_{1}^2x_{2}^{-10}x_{3}^{-4}-41x_{1} x_{2}^{-9}x_{3}^{-4}-5x_{2}^{-8}x_{3}^{-4} \\ \hline  +x_{1}^{-1}x_{2}^{-7}x_{3}^{-4}-15x_{1} x_{2}^{-10}x_{3}^{-4}-8x_{2}^{-9}x_{3}^{-4}-3x_{2}^{-10}x_{3}^{-4}\end{array}~~~~}{(1-x_{1}^2)^7 (1-x_{1}^2x_{2}^2x_{3} ) (1-x_{3} ) }\\&
	+\displaystyle \frac{~~~~\begin{array}{l}-x_{1}^6x_{2}^{-6}x_{3}^{-4}-6x_{1}^5x_{2}^{-6}x_{3}^{-4}-15x_{1}^4x_{2}^{-6}x_{3}^{-4}-20x_{1}^3x_{2}^{-6}x_{3}^{-4}-x_{1}^6x_{2}^{-9}x_{3}^{-5} \\ \hline  -15x_{1}^2x_{2}^{-6}x_{3}^{-4}-x_{1}^6x_{2}^{-10}x_{3}^{-5}-7x_{1}^5x_{2}^{-9}x_{3}^{-5}-6x_{1} x_{2}^{-6}x_{3}^{-4}-7x_{1}^5x_{2}^{-10}x_{3}^{-5} \\ \hline  -20x_{1}^4x_{2}^{-9}x_{3}^{-5}-x_{2}^{-6}x_{3}^{-4}-21x_{1}^4x_{2}^{-10}x_{3}^{-5}-30x_{1}^3x_{2}^{-9}x_{3}^{-5}-34x_{1}^3x_{2}^{-10}x_{3}^{-5} \\ \hline  -25x_{1}^2x_{2}^{-9}x_{3}^{-5}-31x_{1}^2x_{2}^{-10}x_{3}^{-5}-11x_{1} x_{2}^{-9}x_{3}^{-5}-15x_{1} x_{2}^{-10}x_{3}^{-5}-2x_{2}^{-9}x_{3}^{-5} \\ \hline  -3x_{2}^{-10}x_{3}^{-5}\end{array}~~~~}{(1-x_{1}^2)^7 (1-x_{2} ) (1-x_{1}^2x_{2}^2x_{3} ) }\\&
	+\displaystyle \frac{x_{1}^{-1}x_{2}^{-6}x_{3}^{-4}-x_{1}^{-1}x_{2}^{-7}x_{3}^{-4}+x_{1}^{-1}x_{2}^{-9}x_{3}^{-5}+x_{1}^{-1}x_{2}^{-10}x_{3}^{-5}}{(1-x_{1} )^7 (1-x_{2} ) (1-x_{1}^2x_{2}^2x_{3} ) }\\&
	+\displaystyle \frac{-x_{1}^{-1}x_{2}^{-6}x_{3}^{-4}-x_{1}^{-1}x_{2}^{-9}x_{3}^{-5}-x_{1}^{-1}x_{2}^{-10}x_{3}^{-5}}{(1-x_{1} )^7 (1-x_{2} ) (1-x_{1} x_{2}^2x_{3} ) }\\&
	+\displaystyle \frac{x_{1}^{-1}x_{2}^{-6}x_{3}^{-4}}{(1-x_{1} )^6 (1-x_{2} )^2 (1-x_{1} x_{2}^2x_{3} ) }\\&
	+\displaystyle \frac{~~~~\begin{array}{l}2x_{1}^6x_{2}^{-9}x_{3}^{-5}+x_{1}^6x_{2}^{-10}x_{3}^{-5}+14x_{1}^5x_{2}^{-9}x_{3}^{-5}+7x_{1}^5x_{2}^{-10}x_{3}^{-5}+41x_{1}^4x_{2}^{-9}x_{3}^{-5} \\ \hline  +21x_{1}^4x_{2}^{-10}x_{3}^{-5}+64x_{1}^3x_{2}^{-9}x_{3}^{-5}+34x_{1}^3x_{2}^{-10}x_{3}^{-5}+56x_{1}^2x_{2}^{-9}x_{3}^{-5} \\ \hline  +31x_{1}^2x_{2}^{-10}x_{3}^{-5}+26x_{1} x_{2}^{-9}x_{3}^{-5}+15x_{1} x_{2}^{-10}x_{3}^{-5}+5x_{2}^{-9}x_{3}^{-5}+3x_{2}^{-10}x_{3}^{-5}\end{array}~~~~}{(1-x_{1}^2)^7 (1-x_{2} ) (1-x_{3} ) }\\&
	+\displaystyle \frac{x_{1}^{-1}x_{2}^{-4}+x_{1}^{-1}x_{2}^{-6}x_{3}^{-3}+x_{1}^{-1}x_{2}^{-8}x_{3}^{-4}}{(1-x_{1} )^7 (1-x_{2} x_{3} ) (1-x_{1} x_{2}^2x_{3} ) }\\&
	+\displaystyle \frac{~~~~\begin{array}{l}-x_{1}^5x_{2}^{-8}x_{3}^{-5}-6x_{1}^4x_{2}^{-8}x_{3}^{-5}-14x_{1}^3x_{2}^{-8}x_{3}^{-5}-16x_{1}^2x_{2}^{-8}x_{3}^{-5}-9x_{1} x_{2}^{-8}x_{3}^{-5} \\ \hline  -2x_{2}^{-8}x_{3}^{-5}\end{array}~~~~}{(1-x_{1}^2)^6 (1-x_{2} )^2 (1-x_{3} ) }\\&
	+\displaystyle \frac{~~~~\begin{array}{l}-x_{1}^7x_{2}^{-3}+x_{1}^7x_{2}^{-3}x_{3}^{-1}-5x_{1}^6x_{2}^{-3}+5x_{1}^6x_{2}^{-3}x_{3}^{-1}-x_{1}^6x_{2}^{-4}-10x_{1}^5x_{2}^{-3} \\ \hline  +10x_{1}^5x_{2}^{-3}x_{3}^{-1}-5x_{1}^5x_{2}^{-4}-10x_{1}^4x_{2}^{-3}+x_{1}^6x_{2}^{-4}x_{3}^{-2}+10x_{1}^4x_{2}^{-3}x_{3}^{-1} \\ \hline  -10x_{1}^4x_{2}^{-4}-5x_{1}^3x_{2}^{-3}+5x_{1}^5x_{2}^{-4}x_{3}^{-2}-x_{1}^5x_{2}^{-5}x_{3}^{-1}+5x_{1}^3x_{2}^{-3}x_{3}^{-1} \\ \hline  -10x_{1}^3x_{2}^{-4}-x_{1}^2x_{2}^{-3}+x_{1}^5x_{2}^{-5}x_{3}^{-2}+10x_{1}^4x_{2}^{-4}x_{3}^{-2}-5x_{1}^4x_{2}^{-5}x_{3}^{-1} \\ \hline  +x_{1}^2x_{2}^{-3}x_{3}^{-1}-5x_{1}^2x_{2}^{-4}+5x_{1}^4x_{2}^{-5}x_{3}^{-2}-x_{1}^4x_{2}^{-6}x_{3}^{-1}+10x_{1}^3x_{2}^{-4}x_{3}^{-2} \\ \hline  -10x_{1}^3x_{2}^{-5}x_{3}^{-1}-x_{1} x_{2}^{-4}+10x_{1}^3x_{2}^{-5}x_{3}^{-2}-5x_{1}^3x_{2}^{-6}x_{3}^{-1}+5x_{1}^2x_{2}^{-4}x_{3}^{-2} \\ \hline  -10x_{1}^2x_{2}^{-5}x_{3}^{-1}+x_{1}^4x_{2}^{-6}x_{3}^{-3}+10x_{1}^2x_{2}^{-5}x_{3}^{-2}-10x_{1}^2x_{2}^{-6}x_{3}^{-1}+x_{1} x_{2}^{-4}x_{3}^{-2} \\ \hline  -5x_{1} x_{2}^{-5}x_{3}^{-1}+5x_{1}^3x_{2}^{-6}x_{3}^{-3}+5x_{1} x_{2}^{-5}x_{3}^{-2}-10x_{1} x_{2}^{-6}x_{3}^{-1}-x_{2}^{-5}x_{3}^{-1} \\ \hline  +10x_{1}^2x_{2}^{-6}x_{3}^{-3}+x_{2}^{-5}x_{3}^{-2}-5x_{2}^{-6}x_{3}^{-1}+10x_{1} x_{2}^{-6}x_{3}^{-3}-x_{1}^{-1}x_{2}^{-6}x_{3}^{-1} \\ \hline  +5x_{2}^{-6}x_{3}^{-3}+x_{1}^{-1}x_{2}^{-6}x_{3}^{-3}\end{array}~~~~}{(1-x_{1}^2)^6 (1-x_{1} x_{2}^2x_{3} ) (1-x_{3} )^2 }\\&
	+\displaystyle \frac{x_{1}^{-1}x_{2}^{-5}+x_{1}^{-1}x_{2}^{-6}x_{3}^{-1}+2x_{1}^{-1}x_{2}^{-8}x_{3}^{-4}+3x_{1}^{-1}x_{2}^{-9}x_{3}^{-4}+x_{1}^{-1}x_{2}^{-10}x_{3}^{-4}}{(1-x_{1} )^7 (1-x_{1}^2x_{2}^2x_{3} ) (1-x_{3} ) }\\&
	+\displaystyle \frac{~~~~\begin{array}{l}-x_{1}^8x_{2}^{-1}x_{3} -4x_{1}^7x_{2}^{-1}x_{3} -6x_{1}^6x_{2}^{-1}x_{3} -x_{1}^7x_{2}^{-2}-4x_{1}^5x_{2}^{-1}x_{3} -4x_{1}^6x_{2}^{-2} \\ \hline  -x_{1}^4x_{2}^{-1}x_{3} -x_{1}^6x_{2}^{-3}-6x_{1}^5x_{2}^{-2}-4x_{1}^5x_{2}^{-3}-4x_{1}^4x_{2}^{-2}-x_{1}^5x_{2}^{-4}-6x_{1}^4x_{2}^{-3} \\ \hline  -x_{1}^3x_{2}^{-2}-4x_{1}^4x_{2}^{-4}-4x_{1}^3x_{2}^{-3}-6x_{1}^3x_{2}^{-4}-x_{1}^2x_{2}^{-3}-4x_{1}^2x_{2}^{-4}-x_{1} x_{2}^{-4}\end{array}~~~~}{(1-x_{1}^2)^5 (1-x_{1}^2x_{2}^2x_{3} ) (1-x_{3} )^3 }\\&
	+\displaystyle \frac{~~~~\begin{array}{l}-x_{1}^4x_{2}^{-4}x_{3}^{-5}-2x_{1}^4x_{2}^{-5}x_{3}^{-5}-4x_{1}^3x_{2}^{-4}x_{3}^{-5}-x_{1}^4x_{2}^{-6}x_{3}^{-5}-8x_{1}^3x_{2}^{-5}x_{3}^{-5} \\ \hline  -6x_{1}^2x_{2}^{-4}x_{3}^{-5}-4x_{1}^3x_{2}^{-6}x_{3}^{-5}-12x_{1}^2x_{2}^{-5}x_{3}^{-5}-4x_{1} x_{2}^{-4}x_{3}^{-5}-6x_{1}^2x_{2}^{-6}x_{3}^{-5} \\ \hline  -8x_{1} x_{2}^{-5}x_{3}^{-5}-x_{2}^{-4}x_{3}^{-5}-4x_{1} x_{2}^{-6}x_{3}^{-5}-2x_{2}^{-5}x_{3}^{-5}-x_{2}^{-6}x_{3}^{-5}\end{array}~~~~}{(1-x_{1}^2)^5 (1-x_{2}^2)^3 (1-x_{2}^2x_{3} ) }\\&
	+\displaystyle \frac{~~~~\begin{array}{l}x_{1}^4x_{2}^{-4}x_{3}^{-5}+2x_{1}^4x_{2}^{-5}x_{3}^{-5}+4x_{1}^3x_{2}^{-4}x_{3}^{-5}+x_{1}^4x_{2}^{-6}x_{3}^{-5}+8x_{1}^3x_{2}^{-5}x_{3}^{-5} \\ \hline  +6x_{1}^2x_{2}^{-4}x_{3}^{-5}+4x_{1}^3x_{2}^{-6}x_{3}^{-5}+12x_{1}^2x_{2}^{-5}x_{3}^{-5}+4x_{1} x_{2}^{-4}x_{3}^{-5}+6x_{1}^2x_{2}^{-6}x_{3}^{-5} \\ \hline  +8x_{1} x_{2}^{-5}x_{3}^{-5}+x_{2}^{-4}x_{3}^{-5}+4x_{1} x_{2}^{-6}x_{3}^{-5}+2x_{2}^{-5}x_{3}^{-5}+x_{2}^{-6}x_{3}^{-5}\end{array}~~~~}{(1-x_{1}^2)^5 (1-x_{2}^2)^3 (1-x_{3} ) }\\=&
	\left(\begin{array}{l}x_{1}^5x_{2}^{-4}+4x_{1}^4x_{2}^{-4}+6x_{1}^3x_{2}^{-4} \\ +4x_{1}^2x_{2}^{-4}+x_{1} x_{2}^{-4}\end{array}\right)\begin{array}{l}\cdot \frac{1}{48}\left(\frac{x_{1} \partial_{1} }{2}-\frac{x_{2} \partial_{2} }{2}\right)^4\left(x_{3} \partial_{3} \right)^2\\
		~~\cdot\frac{1}{(1-x_{1}^2) (1-x_{1} x_{2} ) (1-x_{3} ) }\end{array}\\&
	\displaystyle +\left(\begin{array}{l}-x_{1}^7x_{2}^{-5}x_{3}^{-3}-5x_{1}^6x_{2}^{-5}x_{3}^{-3} \\ -10x_{1}^5x_{2}^{-5}x_{3}^{-3}-10x_{1}^4x_{2}^{-5}x_{3}^{-3} \\ +x_{1}^4x_{2}^{-6}x_{3}^{-2}-5x_{1}^3x_{2}^{-5}x_{3}^{-3} \\ +4x_{1}^3x_{2}^{-6}x_{3}^{-2}-x_{1}^2x_{2}^{-5}x_{3}^{-3} \\ +6x_{1}^2x_{2}^{-6}x_{3}^{-2}+4x_{1} x_{2}^{-6}x_{3}^{-2} \\ +x_{2}^{-6}x_{3}^{-2}\end{array}\right)\begin{array}{l}\cdot \frac{1}{720}\left(\begin{array}{l}\frac{x_{1} \partial_{1} }{2}-\frac{x_{2} \partial_{2} }{2} \\ +x_{3} \partial_{3} \end{array}\right)^6\\
		~~\cdot\frac{1}{(1-x_{1}^2) (1-x_{1} x_{2} ) (1-x_{2}^2x_{3} ) }\end{array}\\&
	\displaystyle +\left(\begin{array}{l}-x_{1}^7x_{2}^{-3}x_{3}^{-2}-5x_{1}^6x_{2}^{-3}x_{3}^{-2} \\ -10x_{1}^5x_{2}^{-3}x_{3}^{-2}-10x_{1}^4x_{2}^{-3}x_{3}^{-2} \\ -5x_{1}^3x_{2}^{-3}x_{3}^{-2}-x_{1}^2x_{2}^{-3}x_{3}^{-2}\end{array}\right)\begin{array}{l}\cdot \frac{1}{120}\left(\frac{x_{1} \partial_{1} }{2}-\frac{x_{2} \partial_{2} }{2}\right)^5\left(x_{2} \partial_{2} -2x_{3} \partial_{3} \right)\\
		~~\cdot\frac{1}{(1-x_{1}^2) (1-x_{1}^2x_{2}^2x_{3} ) (1-x_{1} x_{2} ) }\end{array}\\&
	\displaystyle +\left(\begin{array}{l}x_{1}^7x_{2}^{-3}x_{3}^{-2}+5x_{1}^6x_{2}^{-3}x_{3}^{-2} \\ +10x_{1}^5x_{2}^{-3}x_{3}^{-2}+10x_{1}^4x_{2}^{-3}x_{3}^{-2} \\ +5x_{1}^3x_{2}^{-3}x_{3}^{-2}+x_{1}^2x_{2}^{-3}x_{3}^{-2}\end{array}\right)\begin{array}{l}\cdot \frac{1}{120}\left(\frac{x_{1} \partial_{1} }{2}-\frac{x_{2} \partial_{2} }{2}\right)^5\left(x_{2} \partial_{2} -x_{3} \partial_{3} \right)\\
		~~\cdot\frac{1}{(1-x_{1}^2) (1-x_{1} x_{2} ) (1-x_{1} x_{2} x_{3} ) }\end{array}\\&
	\displaystyle +x_{2}^{-5}x_{3}^{-3}\cdot \frac{1}{720}\left(x_{1} \partial_{1} -x_{2} \partial_{2} +x_{3} \partial_{3} \right)^6\cdot\frac{1}{(1-x_{1} ) (1-x_{1} x_{2} ) (1-x_{1} x_{2}^2x_{3} ) }\\&
	\displaystyle +\left(\begin{array}{l}-x_{1}^5x_{2}^{-5}x_{3}^{-2}-x_{1}^4x_{2}^{-5}x_{3}^{-1} \\ -x_{1}^5x_{2}^{-6}x_{3}^{-2}-4x_{1}^4x_{2}^{-5}x_{3}^{-2} \\ +x_{1}^4x_{2}^{-6}x_{3}^{-1}-4x_{1}^3x_{2}^{-5}x_{3}^{-1} \\ -5x_{1}^4x_{2}^{-6}x_{3}^{-2}-6x_{1}^3x_{2}^{-5}x_{3}^{-2} \\ +4x_{1}^3x_{2}^{-6}x_{3}^{-1}-6x_{1}^2x_{2}^{-5}x_{3}^{-1} \\ -x_{1}^4x_{2}^{-6}x_{3}^{-3}-10x_{1}^3x_{2}^{-6}x_{3}^{-2} \\ -4x_{1}^2x_{2}^{-5}x_{3}^{-2}+6x_{1}^2x_{2}^{-6}x_{3}^{-1} \\ -4x_{1} x_{2}^{-5}x_{3}^{-1}-4x_{1}^3x_{2}^{-6}x_{3}^{-3} \\ -10x_{1}^2x_{2}^{-6}x_{3}^{-2}-x_{1} x_{2}^{-5}x_{3}^{-2} \\ +4x_{1} x_{2}^{-6}x_{3}^{-1}-x_{2}^{-5}x_{3}^{-1}-6x_{1}^2x_{2}^{-6}x_{3}^{-3} \\ -5x_{1} x_{2}^{-6}x_{3}^{-2}+x_{2}^{-6}x_{3}^{-1}-4x_{1} x_{2}^{-6}x_{3}^{-3} \\ -x_{2}^{-6}x_{3}^{-2}-x_{2}^{-6}x_{3}^{-3}\end{array}\right)\begin{array}{l}\cdot \frac{1}{120}\left(\frac{x_{1} \partial_{1} }{2}\right)^5\left(-\frac{x_{2} \partial_{2} }{2}+x_{3} \partial_{3} \right)\\
		~~\cdot\frac{1}{(1-x_{1}^2) (1-x_{2}^2x_{3} ) (1-x_{3} ) }\end{array}\\&
	\displaystyle +\left(\begin{array}{l}-x_{1}^6x_{2}^{-6}x_{3}^{-3}-6x_{1}^5x_{2}^{-6}x_{3}^{-3} \\ -15x_{1}^4x_{2}^{-6}x_{3}^{-3}-20x_{1}^3x_{2}^{-6}x_{3}^{-3} \\ -15x_{1}^2x_{2}^{-6}x_{3}^{-3}-6x_{1} x_{2}^{-6}x_{3}^{-3} \\ -x_{2}^{-6}x_{3}^{-3}\end{array}\right)\begin{array}{l}\cdot \frac{1}{720}\left(\frac{x_{1} \partial_{1} }{2}\right)^6\\
		~~\cdot\frac{1}{(1-x_{1}^2) (1-x_{2}^2x_{3} ) (1-x_{2} x_{3} ) }\end{array}\\&
	\displaystyle +\left(\begin{array}{l}x_{1}^{-1}x_{2}^{-5}-x_{1}^{-1}x_{2}^{-5}x_{3}^{-1} \\ +x_{1}^{-1}x_{2}^{-6}x_{3}^{-1}-x_{1}^{-1}x_{2}^{-6}x_{3}^{-2}\end{array}\right)\begin{array}{l}\cdot \frac{1}{120}\left(x_{1} \partial_{1} -x_{2} \partial_{2} \right)^5\left(-\frac{x_{2} \partial_{2} }{2}+x_{3} \partial_{3} \right)\\
		~~\cdot\frac{1}{(1-x_{1} ) (1-x_{1}^2x_{2}^2x_{3} ) (1-x_{3} ) }\end{array}\\&
	\displaystyle +\left(\begin{array}{l}-x_{1}^{-1}x_{2}^{-5}+x_{1}^{-1}x_{2}^{-5}x_{3}^{-1} \\ -x_{1}^{-1}x_{2}^{-6}x_{3}^{-1}+x_{1}^{-1}x_{2}^{-6}x_{3}^{-2}\end{array}\right)\begin{array}{l}\cdot \frac{1}{120}\left(x_{1} \partial_{1} -\frac{x_{2} \partial_{2} }{2}\right)^5\left(-\frac{x_{2} \partial_{2} }{2}+x_{3} \partial_{3} \right)\\
		~~\cdot\frac{1}{(1-x_{1} ) (1-x_{1} x_{2}^2x_{3} ) (1-x_{3} ) }\end{array}\\&
	\displaystyle +\left(\begin{array}{l}x_{1}^6x_{2}^{-6}x_{3}^{-1}+7x_{1}^5x_{2}^{-6}x_{3}^{-1} \\ +20x_{1}^4x_{2}^{-6}x_{3}^{-1}+30x_{1}^3x_{2}^{-6}x_{3}^{-1} \\ +25x_{1}^2x_{2}^{-6}x_{3}^{-1}+11x_{1} x_{2}^{-6}x_{3}^{-1} \\ +2x_{2}^{-6}x_{3}^{-1}\end{array}\right)\begin{array}{l}\cdot \frac{1}{720}\left(\frac{x_{1} \partial_{1} }{2}\right)^6\\
		~~\cdot\frac{1}{(1-x_{1}^2) (1-x_{2} x_{3} ) (1-x_{3} ) }\end{array}\\&
	\displaystyle +\left(-x_{2}^{-5}x_{3}^{-3}+x_{1}^{-1}x_{2}^{-6}x_{3}^{-3}\right)\cdot \frac{1}{720}\left(x_{1} \partial_{1} -x_{2} \partial_{2} \right)^6\cdot\frac{1}{(1-x_{1} ) (1-x_{1}^2x_{2}^2x_{3} ) (1-x_{1} x_{2} ) }\\&
	\displaystyle -x_{1}^{-1}x_{2}^{-6}x_{3}^{-3}\cdot \frac{1}{720}\left(x_{1} \partial_{1} -x_{2} \partial_{2} +x_{3} \partial_{3} \right)^6\cdot\frac{1}{(1-x_{1} ) (1-x_{1} x_{2} ) (1-x_{2} x_{3} ) }\\&
	\displaystyle -x_{2}^{-5}x_{3}^{-3}\cdot \frac{1}{120}\left(x_{1} \partial_{1} -x_{2} \partial_{2} +x_{3} \partial_{3} \right)^5\left(x_{2} \partial_{2} -2x_{3} \partial_{3} \right)\cdot\frac{1}{(1-x_{1} ) (1-x_{1} x_{2} ) (1-x_{1} x_{2}^2x_{3} ) }\\&
	\displaystyle +x_{2}^{-5}x_{3}^{-3}\cdot \frac{1}{120}\left(x_{1} \partial_{1} -x_{2} \partial_{2} +x_{3} \partial_{3} \right)^5\left(x_{2} \partial_{2} -x_{3} \partial_{3} \right)\cdot\frac{1}{(1-x_{1} ) (1-x_{1} x_{2} ) (1-x_{2} x_{3} ) }\\&
	\displaystyle +\left(\begin{array}{l}x_{1}^5x_{2}^{-6}x_{3}^{-1}+5x_{1}^4x_{2}^{-6}x_{3}^{-1} \\ +10x_{1}^3x_{2}^{-6}x_{3}^{-1}+10x_{1}^2x_{2}^{-6}x_{3}^{-1} \\ +5x_{1} x_{2}^{-6}x_{3}^{-1}+x_{2}^{-6}x_{3}^{-1}\end{array}\right)\begin{array}{l}\cdot \frac{1}{120}\left(\frac{x_{1} \partial_{1} }{2}\right)^5\left(-x_{2} \partial_{2} +x_{3} \partial_{3} \right)\\
		~~\cdot\frac{1}{(1-x_{1}^2) (1-x_{2} x_{3} ) (1-x_{3} ) }\end{array}\\&
	\displaystyle +\left(\begin{array}{l}-x_{1}^6x_{2}^{-9}x_{3}^{-5}-7x_{1}^5x_{2}^{-9}x_{3}^{-5} \\ -21x_{1}^4x_{2}^{-9}x_{3}^{-5}-34x_{1}^3x_{2}^{-9}x_{3}^{-5} \\ -31x_{1}^2x_{2}^{-9}x_{3}^{-5}-15x_{1} x_{2}^{-9}x_{3}^{-5} \\ -3x_{2}^{-9}x_{3}^{-5}\end{array}\right)\begin{array}{l}\cdot \frac{1}{720}\left(\frac{x_{1} \partial_{1} }{2}\right)^6\\
		~~\cdot\frac{1}{(1-x_{1}^2) (1-x_{2} ) (1-x_{2} x_{3} ) }\end{array}\\&
	\displaystyle +\left(\begin{array}{l}x_{1}^5x_{2}^{-8}x_{3}^{-5}+6x_{1}^4x_{2}^{-8}x_{3}^{-5} \\ +14x_{1}^3x_{2}^{-8}x_{3}^{-5}+16x_{1}^2x_{2}^{-8}x_{3}^{-5} \\ +9x_{1} x_{2}^{-8}x_{3}^{-5}+2x_{2}^{-8}x_{3}^{-5}\end{array}\right)\begin{array}{l}\cdot \frac{1}{120}\left(\frac{x_{1} \partial_{1} }{2}\right)^5\left(x_{2} \partial_{2} -x_{3} \partial_{3} \right)\\
		~~\cdot\frac{1}{(1-x_{1}^2) (1-x_{2} ) (1-x_{2} x_{3} ) }\end{array}\\&
	\displaystyle +\left(\begin{array}{l}x_{1}^6x_{2}^{-4}+7x_{1}^5x_{2}^{-4}+20x_{1}^4x_{2}^{-4} \\ +30x_{1}^3x_{2}^{-4}+25x_{1}^2x_{2}^{-4}+x_{1}^6x_{2}^{-6}x_{3}^{-3} \\ +11x_{1} x_{2}^{-4}+6x_{1}^5x_{2}^{-6}x_{3}^{-3}+2x_{2}^{-4} \\ +15x_{1}^4x_{2}^{-6}x_{3}^{-3}+x_{1}^6x_{2}^{-8}x_{3}^{-4} \\ +20x_{1}^3x_{2}^{-6}x_{3}^{-3}+7x_{1}^5x_{2}^{-8}x_{3}^{-4} \\ +15x_{1}^2x_{2}^{-6}x_{3}^{-3}+21x_{1}^4x_{2}^{-8}x_{3}^{-4} \\ +6x_{1} x_{2}^{-6}x_{3}^{-3}+34x_{1}^3x_{2}^{-8}x_{3}^{-4} \\ +x_{2}^{-6}x_{3}^{-3}+31x_{1}^2x_{2}^{-8}x_{3}^{-4} \\ +15x_{1} x_{2}^{-8}x_{3}^{-4}+3x_{2}^{-8}x_{3}^{-4}\end{array}\right)\begin{array}{l}\cdot \frac{1}{720}\left(\frac{x_{1} \partial_{1} }{2}-x_{2} \partial_{2} +x_{3} \partial_{3} \right)^6\\
		~~\cdot\frac{1}{(1-x_{1}^2) (1-x_{1}^2x_{2}^2x_{3} ) (1-x_{2} x_{3} ) }\end{array}\\&
	\displaystyle +\left(\begin{array}{l}x_{1}^6x_{2}^{-7}x_{3}^{-4}+5x_{1}^5x_{2}^{-7}x_{3}^{-4} \\ -x_{1}^4x_{2}^{-6}x_{3}^{-4}+10x_{1}^4x_{2}^{-7}x_{3}^{-4} \\ -4x_{1}^3x_{2}^{-6}x_{3}^{-4}+10x_{1}^3x_{2}^{-7}x_{3}^{-4} \\ -6x_{1}^2x_{2}^{-6}x_{3}^{-4}+5x_{1}^2x_{2}^{-7}x_{3}^{-4} \\ -4x_{1} x_{2}^{-6}x_{3}^{-4}+x_{1} x_{2}^{-7}x_{3}^{-4} \\ -x_{2}^{-6}x_{3}^{-4}\end{array}\right)\begin{array}{l}\cdot \frac{1}{720}\left(\begin{array}{l}\frac{x_{1} \partial_{1} }{2}+\frac{x_{2} \partial_{2} }{2} \\ -x_{3} \partial_{3} \end{array}\right)^6\\
		~~\cdot\frac{1}{(1-x_{1}^2) (1-x_{2}^2x_{3} ) (1-x_{1} x_{2} x_{3} ) }\end{array}\\&
	\displaystyle +\left(\begin{array}{l}x_{1}^6x_{2}^{-5}x_{3}^{-2}+x_{1}^5x_{2}^{-5}x_{3}^{-1} \\ +5x_{1}^5x_{2}^{-5}x_{3}^{-2}+4x_{1}^4x_{2}^{-5}x_{3}^{-1} \\ +x_{1}^6x_{2}^{-6}x_{3}^{-3}+x_{1}^5x_{2}^{-5}x_{3}^{-3} \\ +10x_{1}^4x_{2}^{-5}x_{3}^{-2}+x_{1}^4x_{2}^{-6}x_{3}^{-1} \\ +6x_{1}^3x_{2}^{-5}x_{3}^{-1}+6x_{1}^5x_{2}^{-6}x_{3}^{-3} \\ +4x_{1}^4x_{2}^{-5}x_{3}^{-3}+10x_{1}^3x_{2}^{-5}x_{3}^{-2} \\ +4x_{1}^3x_{2}^{-6}x_{3}^{-1}+4x_{1}^2x_{2}^{-5}x_{3}^{-1} \\ +15x_{1}^4x_{2}^{-6}x_{3}^{-3}+6x_{1}^3x_{2}^{-5}x_{3}^{-3} \\ +5x_{1}^2x_{2}^{-5}x_{3}^{-2}+6x_{1}^2x_{2}^{-6}x_{3}^{-1} \\ +x_{1} x_{2}^{-5}x_{3}^{-1}+x_{1}^5x_{2}^{-7}x_{3}^{-4} \\ +x_{1}^4x_{2}^{-6}x_{3}^{-4}+20x_{1}^3x_{2}^{-6}x_{3}^{-3} \\ +4x_{1}^2x_{2}^{-5}x_{3}^{-3}+x_{1} x_{2}^{-5}x_{3}^{-2} \\ +4x_{1} x_{2}^{-6}x_{3}^{-1}+5x_{1}^4x_{2}^{-7}x_{3}^{-4} \\ +4x_{1}^3x_{2}^{-6}x_{3}^{-4}+15x_{1}^2x_{2}^{-6}x_{3}^{-3} \\ +x_{1} x_{2}^{-5}x_{3}^{-3}+x_{2}^{-6}x_{3}^{-1}+10x_{1}^3x_{2}^{-7}x_{3}^{-4} \\ +6x_{1}^2x_{2}^{-6}x_{3}^{-4}+6x_{1} x_{2}^{-6}x_{3}^{-3} \\ +10x_{1}^2x_{2}^{-7}x_{3}^{-4}+4x_{1} x_{2}^{-6}x_{3}^{-4} \\ +x_{2}^{-6}x_{3}^{-3}+5x_{1} x_{2}^{-7}x_{3}^{-4}+x_{2}^{-6}x_{3}^{-4} \\ +x_{2}^{-7}x_{3}^{-4}\end{array}\right)\begin{array}{l}\cdot \frac{1}{720}\left(\frac{x_{1} \partial_{1} }{2}\right)^6\\
		~~\cdot\frac{1}{(1-x_{1}^2) (1-x_{2}^2x_{3} ) (1-x_{3} ) }\end{array}\\&
	\displaystyle +\left(\begin{array}{l}-x_{1}^7x_{2}^{-4}x_{3}^{-3}-x_{1}^7x_{2}^{-5}x_{3}^{-3} \\ -7x_{1}^6x_{2}^{-4}x_{3}^{-3}-6x_{1}^6x_{2}^{-5}x_{3}^{-3} \\ -20x_{1}^5x_{2}^{-4}x_{3}^{-3}-15x_{1}^5x_{2}^{-5}x_{3}^{-3} \\ -30x_{1}^4x_{2}^{-4}x_{3}^{-3}-20x_{1}^4x_{2}^{-5}x_{3}^{-3} \\ -25x_{1}^3x_{2}^{-4}x_{3}^{-3}-x_{1}^6x_{2}^{-7}x_{3}^{-4} \\ -15x_{1}^3x_{2}^{-5}x_{3}^{-3}-11x_{1}^2x_{2}^{-4}x_{3}^{-3} \\ -6x_{1}^5x_{2}^{-7}x_{3}^{-4}-6x_{1}^2x_{2}^{-5}x_{3}^{-3} \\ -2x_{1} x_{2}^{-4}x_{3}^{-3}-16x_{1}^4x_{2}^{-7}x_{3}^{-4} \\ -x_{1}^3x_{2}^{-6}x_{3}^{-4}-x_{1} x_{2}^{-5}x_{3}^{-3} \\ -25x_{1}^3x_{2}^{-7}x_{3}^{-4}-4x_{1}^2x_{2}^{-6}x_{3}^{-4} \\ -25x_{1}^2x_{2}^{-7}x_{3}^{-4}-6x_{1} x_{2}^{-6}x_{3}^{-4} \\ -16x_{1} x_{2}^{-7}x_{3}^{-4}-4x_{2}^{-6}x_{3}^{-4} \\ -6x_{2}^{-7}x_{3}^{-4}-x_{1}^{-1}x_{2}^{-6}x_{3}^{-4} \\ -x_{1}^{-1}x_{2}^{-7}x_{3}^{-4}\end{array}\right)\begin{array}{l}\cdot \frac{1}{720}\left(\frac{x_{1} \partial_{1} }{2}-\frac{x_{2} \partial_{2} }{2}\right)^6\\
		~~\cdot\frac{1}{(1-x_{1}^2) (1-x_{1}^2x_{2}^2x_{3} ) (1-x_{1} x_{2} x_{3} ) }\end{array}\\&
	\displaystyle +\left(\begin{array}{l}x_{1}^7x_{2}^{-4}x_{3}^{-3}+x_{1}^7x_{2}^{-5}x_{3}^{-3} \\ +7x_{1}^6x_{2}^{-4}x_{3}^{-3}+6x_{1}^6x_{2}^{-5}x_{3}^{-3} \\ +20x_{1}^5x_{2}^{-4}x_{3}^{-3}+15x_{1}^5x_{2}^{-5}x_{3}^{-3} \\ +30x_{1}^4x_{2}^{-4}x_{3}^{-3}+20x_{1}^4x_{2}^{-5}x_{3}^{-3} \\ +25x_{1}^3x_{2}^{-4}x_{3}^{-3}+15x_{1}^3x_{2}^{-5}x_{3}^{-3} \\ +11x_{1}^2x_{2}^{-4}x_{3}^{-3}+x_{1}^5x_{2}^{-7}x_{3}^{-4} \\ +x_{1}^4x_{2}^{-6}x_{3}^{-4}+6x_{1}^2x_{2}^{-5}x_{3}^{-3} \\ +2x_{1} x_{2}^{-4}x_{3}^{-3}+6x_{1}^4x_{2}^{-7}x_{3}^{-4} \\ +5x_{1}^3x_{2}^{-6}x_{3}^{-4}+x_{1} x_{2}^{-5}x_{3}^{-3} \\ +15x_{1}^3x_{2}^{-7}x_{3}^{-4}+10x_{1}^2x_{2}^{-6}x_{3}^{-4} \\ +20x_{1}^2x_{2}^{-7}x_{3}^{-4}+10x_{1} x_{2}^{-6}x_{3}^{-4} \\ +15x_{1} x_{2}^{-7}x_{3}^{-4}+5x_{2}^{-6}x_{3}^{-4} \\ +6x_{2}^{-7}x_{3}^{-4}+x_{1}^{-1}x_{2}^{-6}x_{3}^{-4} \\ +x_{1}^{-1}x_{2}^{-7}x_{3}^{-4}\end{array}\right)\begin{array}{l}\cdot \frac{1}{720}\left(\frac{x_{1} \partial_{1} }{2}-\frac{x_{3} \partial_{3} }{2}\right)^6\\
		~~\cdot\frac{1}{(1-x_{1}^2) (1-x_{1} x_{2} x_{3} ) (1-x_{1} x_{2}^2x_{3} ) }\end{array}\\&
	\displaystyle +\left(\begin{array}{l}-x_{1}^8x_{2}^{-2}x_{3} -x_{1}^8x_{2}^{-2}-5x_{1}^7x_{2}^{-2}x_{3}  \\ -x_{1}^8x_{2}^{-3}-6x_{1}^7x_{2}^{-2}-10x_{1}^6x_{2}^{-2}x_{3}  \\ -5x_{1}^7x_{2}^{-3}-14x_{1}^6x_{2}^{-2}-10x_{1}^5x_{2}^{-2}x_{3}  \\ -10x_{1}^6x_{2}^{-3}-16x_{1}^5x_{2}^{-2}-5x_{1}^4x_{2}^{-2}x_{3}  \\ -10x_{1}^5x_{2}^{-3}-9x_{1}^4x_{2}^{-2}-x_{1}^3x_{2}^{-2}x_{3}  \\ -5x_{1}^4x_{2}^{-3}-2x_{1}^3x_{2}^{-2}-x_{1}^3x_{2}^{-3}\end{array}\right)\begin{array}{l}\cdot \frac{1}{120}\left(\frac{x_{1} \partial_{1} }{2}-\frac{x_{2} \partial_{2} }{2}\right)^5\left(-x_{2} \partial_{2} +x_{3} \partial_{3} \right)\\
		~~\cdot\frac{1}{(1-x_{1}^2) (1-x_{1} x_{2} x_{3} ) (1-x_{3} ) }\end{array}\\&
	\displaystyle +\left(\begin{array}{l}x_{1}^5x_{2}^{-6}x_{3}^{-5}+x_{1}^5x_{2}^{-7}x_{3}^{-5} \\ +5x_{1}^4x_{2}^{-6}x_{3}^{-5}+5x_{1}^4x_{2}^{-7}x_{3}^{-5} \\ +10x_{1}^3x_{2}^{-6}x_{3}^{-5}+10x_{1}^3x_{2}^{-7}x_{3}^{-5} \\ +10x_{1}^2x_{2}^{-6}x_{3}^{-5}+10x_{1}^2x_{2}^{-7}x_{3}^{-5} \\ +5x_{1} x_{2}^{-6}x_{3}^{-5}+5x_{1} x_{2}^{-7}x_{3}^{-5} \\ +x_{2}^{-6}x_{3}^{-5}+x_{2}^{-7}x_{3}^{-5}\end{array}\right)\begin{array}{l}\cdot \frac{1}{120}\left(\frac{x_{1} \partial_{1} }{2}\right)^5\left(\frac{x_{2} \partial_{2} }{2}-x_{3} \partial_{3} \right)\\
		~~\cdot\frac{1}{(1-x_{1}^2) (1-x_{2}^2) (1-x_{2}^2x_{3} ) }\end{array}\\&
	\displaystyle +\left(\begin{array}{l}-x_{1}^5x_{2}^{-6}x_{3}^{-5}-x_{1}^5x_{2}^{-7}x_{3}^{-5} \\ -5x_{1}^4x_{2}^{-6}x_{3}^{-5}-5x_{1}^4x_{2}^{-7}x_{3}^{-5} \\ -10x_{1}^3x_{2}^{-6}x_{3}^{-5}-10x_{1}^3x_{2}^{-7}x_{3}^{-5} \\ -10x_{1}^2x_{2}^{-6}x_{3}^{-5}-10x_{1}^2x_{2}^{-7}x_{3}^{-5} \\ -5x_{1} x_{2}^{-6}x_{3}^{-5}-5x_{1} x_{2}^{-7}x_{3}^{-5} \\ -x_{2}^{-6}x_{3}^{-5}-x_{2}^{-7}x_{3}^{-5}\end{array}\right)\begin{array}{l}\cdot \frac{1}{120}\left(\frac{x_{1} \partial_{1} }{2}\right)^5\left(\frac{x_{2} \partial_{2} }{2}\right)\\
		~~\cdot\frac{1}{(1-x_{1}^2) (1-x_{2}^2) (1-x_{3} ) }\end{array}\\&
	\displaystyle +\left(\begin{array}{l}x_{1}^7x_{2}^{-2}x_{3} +4x_{1}^6x_{2}^{-2}x_{3} +6x_{1}^5x_{2}^{-2}x_{3}  \\ +4x_{1}^4x_{2}^{-2}x_{3} +x_{1}^3x_{2}^{-2}x_{3} \end{array}\right)\begin{array}{l}\cdot \frac{1}{48}\left(\frac{x_{1} \partial_{1} }{2}-\frac{x_{2} \partial_{2} }{2}\right)^4\left(-x_{2} \partial_{2} +x_{3} \partial_{3} \right)^2\\
		~~\cdot\frac{1}{(1-x_{1}^2) (1-x_{1} x_{2} x_{3} ) (1-x_{3} ) }\end{array}\\&
	\displaystyle +\left(\begin{array}{l}x_{1}^8x_{2}^{-3}x_{3}^{-2}+6x_{1}^7x_{2}^{-3}x_{3}^{-2} \\ +x_{1}^7x_{2}^{-4}x_{3}^{-1}+15x_{1}^6x_{2}^{-3}x_{3}^{-2} \\ +7x_{1}^6x_{2}^{-4}x_{3}^{-1}+20x_{1}^5x_{2}^{-3}x_{3}^{-2} \\ +20x_{1}^5x_{2}^{-4}x_{3}^{-1}+x_{1}^7x_{2}^{-5}x_{3}^{-3} \\ +15x_{1}^4x_{2}^{-3}x_{3}^{-2}+30x_{1}^4x_{2}^{-4}x_{3}^{-1} \\ +6x_{1}^6x_{2}^{-5}x_{3}^{-3}+6x_{1}^3x_{2}^{-3}x_{3}^{-2} \\ +25x_{1}^3x_{2}^{-4}x_{3}^{-1}+16x_{1}^5x_{2}^{-5}x_{3}^{-3} \\ +x_{1}^2x_{2}^{-3}x_{3}^{-2}+11x_{1}^2x_{2}^{-4}x_{3}^{-1} \\ +25x_{1}^4x_{2}^{-5}x_{3}^{-3}+2x_{1} x_{2}^{-4}x_{3}^{-1} \\ +25x_{1}^3x_{2}^{-5}x_{3}^{-3}+x_{1}^3x_{2}^{-6}x_{3}^{-2} \\ +16x_{1}^2x_{2}^{-5}x_{3}^{-3}+4x_{1}^2x_{2}^{-6}x_{3}^{-2} \\ +6x_{1} x_{2}^{-5}x_{3}^{-3}+6x_{1} x_{2}^{-6}x_{3}^{-2} \\ +x_{2}^{-5}x_{3}^{-3}+4x_{2}^{-6}x_{3}^{-2}+x_{1}^{-1}x_{2}^{-6}x_{3}^{-2}\end{array}\right)\begin{array}{l}\cdot \frac{1}{720}\left(\frac{x_{1} \partial_{1} }{2}-\frac{x_{2} \partial_{2} }{2}\right)^6\\
		~~\cdot\frac{1}{(1-x_{1}^2) (1-x_{1}^2x_{2}^2x_{3} ) (1-x_{1} x_{2} ) }\end{array}\\&
	\displaystyle -x_{1}^{-1}x_{2}^{-7}x_{3}^{-4}\cdot \frac{1}{720}\left(x_{1} \partial_{1} -x_{3} \partial_{3} \right)^6\cdot\frac{1}{(1-x_{1} ) (1-x_{1} x_{2} x_{3} ) (1-x_{1} x_{2}^2x_{3} ) }\\&
	\displaystyle +\left(\begin{array}{l}x_{1}^6x_{2}^{-6}x_{3}^{-4}+6x_{1}^5x_{2}^{-6}x_{3}^{-4} \\ +15x_{1}^4x_{2}^{-6}x_{3}^{-4}+20x_{1}^3x_{2}^{-6}x_{3}^{-4} \\ +15x_{1}^2x_{2}^{-6}x_{3}^{-4}+6x_{1} x_{2}^{-6}x_{3}^{-4} \\ +x_{2}^{-6}x_{3}^{-4}\end{array}\right)\begin{array}{l}\cdot \frac{1}{720}\left(\frac{x_{1} \partial_{1} }{2}\right)^6\\
		~~\cdot\frac{1}{(1-x_{1}^2) (1-x_{2} ) (1-x_{2}^2x_{3} ) }\end{array}\\&
	\displaystyle +\left(\begin{array}{l}-x_{1}^8x_{2}^{-3}x_{3}^{-2}-6x_{1}^7x_{2}^{-3}x_{3}^{-2} \\ -x_{1}^7x_{2}^{-4}x_{3}^{-1}-15x_{1}^6x_{2}^{-3}x_{3}^{-2} \\ -7x_{1}^6x_{2}^{-4}x_{3}^{-1}-20x_{1}^5x_{2}^{-3}x_{3}^{-2} \\ -20x_{1}^5x_{2}^{-4}x_{3}^{-1}-15x_{1}^4x_{2}^{-3}x_{3}^{-2} \\ -30x_{1}^4x_{2}^{-4}x_{3}^{-1}-x_{1}^6x_{2}^{-5}x_{3}^{-3} \\ -6x_{1}^3x_{2}^{-3}x_{3}^{-2}-25x_{1}^3x_{2}^{-4}x_{3}^{-1} \\ -6x_{1}^5x_{2}^{-5}x_{3}^{-3}-x_{1}^2x_{2}^{-3}x_{3}^{-2} \\ -11x_{1}^2x_{2}^{-4}x_{3}^{-1}-15x_{1}^4x_{2}^{-5}x_{3}^{-3} \\ -x_{1}^4x_{2}^{-6}x_{3}^{-2}-2x_{1} x_{2}^{-4}x_{3}^{-1} \\ -20x_{1}^3x_{2}^{-5}x_{3}^{-3}-5x_{1}^3x_{2}^{-6}x_{3}^{-2} \\ -15x_{1}^2x_{2}^{-5}x_{3}^{-3}-10x_{1}^2x_{2}^{-6}x_{3}^{-2} \\ -6x_{1} x_{2}^{-5}x_{3}^{-3}-10x_{1} x_{2}^{-6}x_{3}^{-2} \\ -x_{2}^{-5}x_{3}^{-3}-5x_{2}^{-6}x_{3}^{-2}-x_{1}^{-1}x_{2}^{-6}x_{3}^{-2}\end{array}\right)\begin{array}{l}\cdot \frac{1}{720}\left(\begin{array}{l}\frac{x_{1} \partial_{1} }{2}-\frac{x_{2} \partial_{2} }{2} \\ +\frac{x_{3} \partial_{3} }{2}\end{array}\right)^6\\
		~~\cdot\frac{1}{(1-x_{1}^2) (1-x_{1} x_{2} ) (1-x_{1} x_{2}^2x_{3} ) }\end{array}\\&
	\displaystyle +\left(\begin{array}{l}-x_{1}^8x_{2}^{-3}-x_{1}^8x_{2}^{-3}x_{3}^{-1}-7x_{1}^7x_{2}^{-3} \\ -x_{1}^8x_{2}^{-3}x_{3}^{-2}-6x_{1}^7x_{2}^{-3}x_{3}^{-1} \\ -x_{1}^7x_{2}^{-4}-20x_{1}^6x_{2}^{-3}-x_{1}^8x_{2}^{-4}x_{3}^{-2} \\ -7x_{1}^7x_{2}^{-3}x_{3}^{-2}-15x_{1}^6x_{2}^{-3}x_{3}^{-1} \\ -7x_{1}^6x_{2}^{-4}-30x_{1}^5x_{2}^{-3}-7x_{1}^7x_{2}^{-4}x_{3}^{-2} \\ -20x_{1}^6x_{2}^{-3}x_{3}^{-2}-20x_{1}^5x_{2}^{-3}x_{3}^{-1} \\ -20x_{1}^5x_{2}^{-4}-25x_{1}^4x_{2}^{-3}-x_{1}^7x_{2}^{-4}x_{3}^{-3} \\ -21x_{1}^6x_{2}^{-4}x_{3}^{-2}-30x_{1}^5x_{2}^{-3}x_{3}^{-2} \\ -15x_{1}^4x_{2}^{-3}x_{3}^{-1}-30x_{1}^4x_{2}^{-4}-11x_{1}^3x_{2}^{-3} \\ -x_{1}^7x_{2}^{-5}x_{3}^{-3}-7x_{1}^6x_{2}^{-4}x_{3}^{-3} \\ -x_{1}^6x_{2}^{-5}x_{3}^{-2}-35x_{1}^5x_{2}^{-4}x_{3}^{-2} \\ -x_{1}^5x_{2}^{-5}x_{3}^{-1}-25x_{1}^4x_{2}^{-3}x_{3}^{-2} \\ -6x_{1}^3x_{2}^{-3}x_{3}^{-1}-25x_{1}^3x_{2}^{-4}-2x_{1}^2x_{2}^{-3} \\ -6x_{1}^6x_{2}^{-5}x_{3}^{-3}-20x_{1}^5x_{2}^{-4}x_{3}^{-3} \\ -6x_{1}^5x_{2}^{-5}x_{3}^{-2}-35x_{1}^4x_{2}^{-4}x_{3}^{-2} \\ -5x_{1}^4x_{2}^{-5}x_{3}^{-1}-11x_{1}^3x_{2}^{-3}x_{3}^{-2} \\ -x_{1}^2x_{2}^{-3}x_{3}^{-1}-11x_{1}^2x_{2}^{-4}-x_{1}^6x_{2}^{-6}x_{3}^{-3} \\ -16x_{1}^5x_{2}^{-5}x_{3}^{-3}-30x_{1}^4x_{2}^{-4}x_{3}^{-3} \\ -15x_{1}^4x_{2}^{-5}x_{3}^{-2}-x_{1}^4x_{2}^{-6}x_{3}^{-1} \\ -21x_{1}^3x_{2}^{-4}x_{3}^{-2}-10x_{1}^3x_{2}^{-5}x_{3}^{-1} \\ -2x_{1}^2x_{2}^{-3}x_{3}^{-2}-2x_{1} x_{2}^{-4}-7x_{1}^5x_{2}^{-6}x_{3}^{-3} \\ -25x_{1}^4x_{2}^{-5}x_{3}^{-3}-25x_{1}^3x_{2}^{-4}x_{3}^{-3} \\ -20x_{1}^3x_{2}^{-5}x_{3}^{-2}-5x_{1}^3x_{2}^{-6}x_{3}^{-1} \\ -7x_{1}^2x_{2}^{-4}x_{3}^{-2}-10x_{1}^2x_{2}^{-5}x_{3}^{-1} \\ -21x_{1}^4x_{2}^{-6}x_{3}^{-3}-25x_{1}^3x_{2}^{-5}x_{3}^{-3} \\ -11x_{1}^2x_{2}^{-4}x_{3}^{-3}-15x_{1}^2x_{2}^{-5}x_{3}^{-2} \\ -10x_{1}^2x_{2}^{-6}x_{3}^{-1}-x_{1} x_{2}^{-4}x_{3}^{-2} \\ -5x_{1} x_{2}^{-5}x_{3}^{-1}-x_{1}^5x_{2}^{-7}x_{3}^{-4} \\ -x_{1}^4x_{2}^{-6}x_{3}^{-4}-35x_{1}^3x_{2}^{-6}x_{3}^{-3} \\ -16x_{1}^2x_{2}^{-5}x_{3}^{-3}-2x_{1} x_{2}^{-4}x_{3}^{-3} \\ -6x_{1} x_{2}^{-5}x_{3}^{-2}-10x_{1} x_{2}^{-6}x_{3}^{-1} \\ -x_{2}^{-5}x_{3}^{-1}-6x_{1}^4x_{2}^{-7}x_{3}^{-4}-5x_{1}^3x_{2}^{-6}x_{3}^{-4} \\ -35x_{1}^2x_{2}^{-6}x_{3}^{-3}-6x_{1} x_{2}^{-5}x_{3}^{-3} \\ -x_{2}^{-5}x_{3}^{-2}-5x_{2}^{-6}x_{3}^{-1}-15x_{1}^3x_{2}^{-7}x_{3}^{-4} \\ -10x_{1}^2x_{2}^{-6}x_{3}^{-4}-21x_{1} x_{2}^{-6}x_{3}^{-3} \\ -x_{2}^{-5}x_{3}^{-3}-x_{1}^{-1}x_{2}^{-6}x_{3}^{-1} \\ -20x_{1}^2x_{2}^{-7}x_{3}^{-4}-10x_{1} x_{2}^{-6}x_{3}^{-4} \\ -7x_{2}^{-6}x_{3}^{-3}-15x_{1} x_{2}^{-7}x_{3}^{-4} \\ -5x_{2}^{-6}x_{3}^{-4}-x_{1}^{-1}x_{2}^{-6}x_{3}^{-3} \\ -6x_{2}^{-7}x_{3}^{-4}-x_{1}^{-1}x_{2}^{-6}x_{3}^{-4} \\ -x_{1}^{-1}x_{2}^{-7}x_{3}^{-4}\end{array}\right)\begin{array}{l}\cdot \frac{1}{720}\left(\frac{x_{1} \partial_{1} }{2}-\frac{x_{2} \partial_{2} }{4}\right)^6\\
		~~\cdot\frac{1}{(1-x_{1}^2) (1-x_{1} x_{2}^2x_{3} ) (1-x_{3} ) }\end{array}\\&
	\displaystyle +\left(-x_{2}^{-6}x_{3}^{-3}+x_{1}^{-1}x_{2}^{-7}x_{3}^{-4}\right)\cdot \frac{1}{720}\left(x_{1} \partial_{1} -x_{2} \partial_{2} \right)^6\cdot\frac{1}{(1-x_{1} ) (1-x_{1}^2x_{2}^2x_{3} ) (1-x_{1} x_{2} x_{3} ) }\\&
	\displaystyle +x_{1}^{-1}x_{2}^{-7}x_{3}^{-4}\cdot \frac{1}{720}\left(x_{1} \partial_{1} -x_{3} \partial_{3} \right)^6\cdot\frac{1}{(1-x_{1} ) (1-x_{2} ) (1-x_{1} x_{2} x_{3} ) }\\&
	\displaystyle -x_{1}^{-1}x_{2}^{-6}x_{3}^{-4}\cdot \frac{1}{120}\left(x_{1} \partial_{1} -x_{3} \partial_{3} \right)^5\left(x_{2} \partial_{2} -x_{3} \partial_{3} \right)\cdot\frac{1}{(1-x_{1} ) (1-x_{2} ) (1-x_{1} x_{2} x_{3} ) }\\&
	\displaystyle +\left(
}
		\caption{$48$ combinatorial chambers of \(A_4\):(1, 0, 0, 0), (0, 1, 0, 0), (0, 0, 1, 0), (0, 0, 0, 1), (1, 1, 0, 0), (0, 1, 1, 0), (0, 0, 1, 1), (1, 1, 1, 0), (0, 1, 1, 1), (1, 1, 1, 1)}
	\end{center}
\end{table}
\begin{longtable}{|c|c|c|}\caption{\footnotesize V.p.f. of \(A_4\):(1, 0, 0, 0), (0, 1, 0, 0), (0, 0, 1, 0), (0, 0, 0, 1), (1, 1, 0, 0), (0, 1, 1, 0), (0, 0, 1, 1), (1, 1, 1, 0), (0, 1, 1, 1), (1, 1, 1, 1)}\\\hline N & Polynomial/Lattice & Shift(s)\\ \hline
	\endfirsthead\multicolumn{3}{c}{{\bfseries \tablename\ \thetable{} -- continued from previous page}} \\
	\hline  N &  Polynomial/Lattice & Shift(s) \\ \hline
	\endhead\hline \multicolumn{3}{|c|}{{Continued on next page}} \\ \hline\endfoot\endlastfoot1&\(\begin{array}{l}-\frac{x_{1}^5x_{2} }{60}+\frac{x_{1}^5x_{4} }{20}+\frac{x_{1}^4x_{2}^2}{24}-\frac{x_{1}^4x_{2} x_{4} }{8}-\frac{x_{1}^3x_{2}^3}{36}+\frac{x_{1}^3x_{2}^2x_{4} }{12} \\ +\frac{x_{1} x_{2}^5}{120}-\frac{x_{1} x_{2}^4x_{4} }{24}+\frac{x_{1} x_{2}^3x_{4}^2}{12}-\frac{x_{1} x_{2}^2x_{4}^3}{12}+\frac{x_{1} x_{2} x_{4}^4}{24} \\ -\frac{x_{1} x_{4}^5}{120}-\frac{x_{2}^6}{360}+\frac{x_{2}^5x_{4} }{60}-\frac{x_{2}^4x_{4}^2}{24}+\frac{x_{2}^3x_{4}^3}{18}-\frac{x_{2}^2x_{4}^4}{24}+\frac{x_{2} x_{4}^5}{60} \\ -\frac{x_{4}^6}{360}+\frac{x_{1}^5}{20}-\frac{x_{1}^4x_{2} }{6}+\frac{x_{1}^4x_{4} }{8}+\frac{x_{1}^3x_{2}^2}{4}-\frac{x_{1}^3x_{2} x_{4} }{2}-\frac{x_{1}^2x_{2}^3}{6} \\ +\frac{x_{1}^2x_{2}^2x_{4} }{2}+\frac{x_{1} x_{2}^4}{24}-\frac{x_{1} x_{2}^3x_{4} }{6}+\frac{x_{1} x_{2}^2x_{4}^2}{4}-\frac{x_{1} x_{2} x_{4}^3}{6}+\frac{x_{1} x_{4}^4}{24} \\ +\frac{x_{1}^4}{8}-\frac{7}{18}x_{1}^3x_{2} -\frac{x_{1}^3x_{4} }{3}+\frac{11}{24}x_{1}^2x_{2}^2+\frac{x_{1}^2x_{2} x_{4} }{8}-\frac{19}{72}x_{1} x_{2}^3 \\ +\frac{19}{24}x_{1} x_{2}^2x_{4} +\frac{x_{1} x_{2} x_{4}^2}{8}-\frac{x_{1} x_{4}^3}{24}+\frac{x_{2}^4}{18}-\frac{2}{9}x_{2}^3x_{4} +\frac{x_{2}^2x_{4}^2}{3} \\ -\frac{2}{9}x_{2} x_{4}^3+\frac{x_{4}^4}{18}-\frac{x_{1}^3}{3}+\frac{x_{1}^2x_{2} }{3}-\frac{5}{8}x_{1}^2x_{4} +\frac{5}{24}x_{1} x_{2}^2+\frac{25}{12}x_{1} x_{2} x_{4}  \\ -\frac{x_{1} x_{4}^2}{24}-\frac{x_{2}^3}{12}+\frac{x_{2}^2x_{4} }{4}+\frac{x_{2} x_{4}^2}{4}-\frac{x_{4}^3}{12}-\frac{5}{8}x_{1}^2+\frac{76}{45}x_{1} x_{2}  \\ +\frac{5}{6}x_{1} x_{4} -\frac{19}{360}x_{2}^2+\frac{289}{180}x_{2} x_{4} -\frac{19}{360}x_{4}^2+\frac{47}{60}x_{1} +\frac{13}{12}x_{2} +\frac{13}{12}x_{4}  \\ +1\end{array}\)&-\\
\hline\end{longtable}

\begin{longtable}{|ccccc|}\caption{\footnotesize V.p.f. of \(A_4\):(1, 0, 0, 0), (0, 1, 0, 0), (0, 0, 1, 0), (0, 0, 0, 1), (1, 1, 0, 0), (0, 1, 1, 0), (0, 0, 1, 1), (1, 1, 1, 0), (0, 1, 1, 1), (1, 1, 1, 1)}\\\hline N & Defining inequalities & Vertices& Int. Pt.& Neighbors \\ \hline
	\endfirsthead\multicolumn{5}{c}{{\bfseries \tablename\ \thetable{} -- continued from previous page}} \\
	\hline N & Defining inequalities & Vertices&Int. Pt. & Neighbors\\ \hline
	\endhead\hline \multicolumn{5}{|c|}{{Continued on next page}} \\ \hline\endfoot\endlastfoot1&
	\(\begin{array}{rcl}-x_{2}+x_{3}-x_{4}&\geq& 0\\
		x_{2}-x_{4}&\geq& 0\\
		x_{1}-x_{2}+x_{4}&\geq& 0\\
		-x_{1}+x_{4}&\geq& 0\\
	\end{array}\)&\(\begin{array}{l}(0, 1, 2, 1): 5\checkmark\\(1, 1, 2, 1): 10\checkmark\\(0, 0, 1, 0): 1\checkmark\\(1, 2, 3, 1): 20\checkmark\end{array}\)&(2, 4, 8, 3): 216& [223, 169, 211, 218]\\\hline
	2&
	\(\begin{array}{rcl}-x_{1}+x_{2}-x_{3}&\geq& 0\\
		x_{3}-x_{4}&\geq& 0\\
		x_{1}-x_{3}+x_{4}&\geq& 0\\
		-x_{1}+x_{4}&\geq& 0\\
	\end{array}\)&\(\begin{array}{l}(0, 1, 1, 1): 4\checkmark\\(1, 2, 1, 1): 10\checkmark\\(0, 1, 0, 0): 1\checkmark\\(1, 3, 2, 1): 20\checkmark\end{array}\)&(2, 7, 4, 3): 189& [197, 55, 185, 192]\\\hline
	3&
	\(\begin{array}{rcl}-x_{1}+x_{2}-x_{3}&\geq& 0\\
		-x_{1}+x_{3}&\geq& 0\\
		x_{1}-x_{3}+x_{4}&\geq& 0\\
		x_{1}-x_{4}&\geq& 0\\
	\end{array}\)&\(\begin{array}{l}(1, 2, 1, 0): 5\checkmark\\(1, 2, 1, 1): 10\checkmark\\(0, 1, 0, 0): 1\checkmark\\(1, 3, 2, 1): 20\checkmark\end{array}\)&(3, 8, 4, 2): 216& [198, 62, 186, 191]\\\hline
	4&
	\(\begin{array}{rcl}x_{1}-x_{2}+x_{3}-x_{4}&\geq& 0\\
		-x_{1}+x_{3}&\geq& 0\\
		x_{1}-x_{3}+x_{4}&\geq& 0\\
		-x_{1}+x_{2}-x_{4}&\geq& 0\\
	\end{array}\)&\(\begin{array}{l}(1, 2, 1, 0): 5\checkmark\\(1, 1, 1, 0): 4\checkmark\\(1, 2, 2, 1): 18\checkmark\\(2, 3, 2, 1): 33\end{array}\)&(5, 8, 6, 2): 675& [198, 131, 97, 126]\\\hline
	5&
	\(\begin{array}{rcl}x_{1}-x_{2}+x_{3}-x_{4}&\geq& 0\\
		x_{1}-x_{3}+x_{4}&\geq& 0\\
		x_{2}-x_{4}&\geq& 0\\
		-x_{1}+x_{4}&\geq& 0\\
	\end{array}\)&\(\begin{array}{l}(1, 2, 3, 2): 33\\(1, 2, 2, 1): 18\checkmark\\(1, 1, 1, 1): 8\checkmark\\(1, 1, 2, 1): 10\checkmark\end{array}\)&(4, 6, 8, 5): 1393& [164, 223, 145, 234]\\\hline
	6&
	\(\begin{array}{rcl}x_{3}-x_{4}&\geq& 0\\
		-x_{1}+x_{2}-x_{3}+x_{4}&\geq& 0\\
		x_{1}-x_{3}&\geq& 0\\
		x_{1}-x_{2}+x_{4}&\geq& 0\\
	\end{array}\)&\(\begin{array}{l}(1, 1, 1, 1): 8\checkmark\\(1, 1, 0, 0): 2\checkmark\\(1, 2, 1, 1): 10\checkmark\\(2, 3, 2, 1): 33\end{array}\)&(5, 7, 4, 3): 496& [84, 132, 204, 111]\\
\hline
\end{longtable}

\allowdisplaybreaks\begin{align*}&~~~
	\frac{1}{(1-x_{1} ) (1-x_{2} ) (1-x_{3} ) (1-x_{4} ) (1-x_{1} x_{2} ) (1-x_{2} x_{3} ) (1-x_{3} x_{4} ) (1-x_{1} x_{2} x_{3} ) (1-x_{2} x_{3} x_{4} ) (1-x_{1} x_{2} x_{3} x_{4} ) }\\=&~~~
	\displaystyle \frac{-x_{2}^{-3}x_{4}^{-2}}{(1-x_{1} )^4 (1-x_{1} x_{2} ) (1-x_{3} )^4 (1-x_{4} ) }\\&
	+\displaystyle \frac{x_{2}^{-4}x_{3}^{-3}x_{4}^{-2}}{(1-x_{1} )^6 (1-x_{3} x_{4} ) (1-x_{2} x_{3} x_{4} ) (1-x_{4} )^2 }\\&
	+\displaystyle \frac{x_{2}^{-4}x_{3}^{-3}x_{4}^{-3}}{(1-x_{1} )^7 (1-x_{3} x_{4} ) (1-x_{2} x_{3} ) (1-x_{4} ) }\\&
	+\displaystyle \frac{-x_{2}^{-4}x_{3}^{-3}x_{4}^{-2}}{(1-x_{1} )^6 (1-x_{3} x_{4} ) (1-x_{2} x_{3} ) (1-x_{4} )^2 }\\&
	+\displaystyle \frac{x_{1} x_{2}^{-2}x_{3}^{-2}x_{4}^{-1}}{(1-x_{1} )^5 (1-x_{1} x_{2} x_{3} x_{4} ) (1-x_{3} x_{4} ) (1-x_{4} )^3 }\\&
	+\displaystyle \frac{-x_{1} x_{2}^{-2}x_{3}^{-2}x_{4}^{-1}}{(1-x_{1} )^5 (1-x_{3} x_{4} ) (1-x_{1} x_{2} x_{3} ) (1-x_{4} )^3 }\\&
	+\displaystyle \frac{-x_{1}^2x_{2}^{-1}x_{3}^{-2}x_{4}^{-1}}{(1-x_{1} )^5 (1-x_{1} x_{2} ) (1-x_{1} x_{2} x_{3} ) (1-x_{4} )^3 }\\&
	+\displaystyle \frac{x_{1}^2x_{2}^{-1}x_{3}^{-2}x_{4}^{-1}}{(1-x_{1} )^5 (1-x_{1} x_{2} x_{3} x_{4} ) (1-x_{1} x_{2} ) (1-x_{4} )^3 }\\&
	+\displaystyle \frac{-x_{1} x_{2}^{-2}x_{3}^{-2}x_{4}^{-2}}{(1-x_{1} )^6 (1-x_{1} x_{2} x_{3} x_{4} ) (1-x_{1} x_{2} )^2 (1-x_{1} x_{2} x_{3} ) }\\&
	+\displaystyle \frac{-x_{1}^2x_{2}^{-1}x_{3}^{-2}x_{4}^{-1}}{(1-x_{1} )^5 (1-x_{1} x_{2} x_{3} x_{4} ) (1-x_{1} x_{2} )^2 (1-x_{4} )^2 }\\&
	+\displaystyle \frac{x_{1}^2x_{2}^{-1}x_{3}^{-2}x_{4}^{-1}}{(1-x_{1} )^5 (1-x_{1} x_{2} )^2 (1-x_{3} x_{4} ) (1-x_{4} )^2 }\\&
	+\displaystyle \frac{x_{1} x_{2}^{-2}x_{3}^{-2}x_{4}^{-2}}{(1-x_{1} )^6 (1-x_{1} x_{2} )^2 (1-x_{3} x_{4} ) (1-x_{1} x_{2} x_{3} ) }\\&
	+\displaystyle \frac{-x_{1}^2x_{2}^{-1}x_{3}^{-2}x_{4}^{-2}+x_{1} x_{2}^{-2}x_{3}^{-2}x_{4}^{-2}-x_{1} x_{2}^{-3}x_{3}^{-3}x_{4}^{-2}}{(1-x_{1} )^7 (1-x_{1} x_{2} x_{3} x_{4} ) (1-x_{1} x_{2} ) (1-x_{1} x_{2} x_{3} ) }\\&
	+\displaystyle \frac{x_{1}^2x_{2}^{-1}x_{3}^{-2}x_{4}^{-2}-x_{1} x_{2}^{-2}x_{3}^{-2}x_{4}^{-2}+x_{1} x_{2}^{-3}x_{3}^{-3}x_{4}^{-2}}{(1-x_{1} )^7 (1-x_{1} x_{2} ) (1-x_{2} x_{3} x_{4} ) (1-x_{1} x_{2} x_{3} ) }\\&
	+\displaystyle \frac{x_{1}^2x_{2}^{-1}x_{3}^{-2}x_{4}^{-1}+x_{1} x_{2}^{-3}x_{3}^{-3}x_{4}^{-1}}{(1-x_{1} )^6 (1-x_{1} x_{2} x_{3} x_{4} ) (1-x_{1} x_{2} ) (1-x_{4} )^2 }\\&
	+\displaystyle \frac{-x_{1}^2x_{2}^{-1}x_{3}^{-2}x_{4}^{-1}-x_{1} x_{2}^{-3}x_{3}^{-3}x_{4}^{-1}}{(1-x_{1} )^6 (1-x_{1} x_{2} ) (1-x_{1} x_{2} x_{3} ) (1-x_{4} )^2 }\\&
	+\displaystyle \frac{x_{2}^{-3}x_{3}^{-2}x_{4}^{-2}}{(1-x_{1} )^7 (1-x_{2} x_{3} ) (1-x_{2} x_{3} x_{4} ) (1-x_{3} x_{4} ) }\\&
	+\displaystyle \frac{x_{1} x_{2}^{-4}x_{3}^{-4}x_{4}^{-2}}{(1-x_{1} )^6 (1-x_{1} x_{2} ) (1-x_{2} x_{3} x_{4} ) (1-x_{4} )^2 }\\&
	+\displaystyle \frac{-x_{1} x_{2}^{-2}x_{3}^{-2}x_{4}^{-2}-x_{1} x_{2}^{-3}x_{3}^{-3}x_{4}^{-3}}{(1-x_{1} )^7 (1-x_{1} x_{2} x_{3} x_{4} ) (1-x_{2} x_{3} ) (1-x_{1} x_{2} ) }\\&
	+\displaystyle \frac{x_{1} x_{2}^{-4}x_{3}^{-4}x_{4}^{-3}}{(1-x_{1} )^7 (1-x_{1} x_{2} ) (1-x_{2} x_{3} ) (1-x_{4} ) }\\&
	+\displaystyle \frac{x_{1} x_{2}^{-2}x_{3}^{-2}x_{4}^{-2}}{(1-x_{1} )^6 (1-x_{1} x_{2} x_{3} x_{4} ) (1-x_{2} x_{3} ) (1-x_{1} x_{2} )^2 }\\&
	+\displaystyle \frac{x_{1} x_{2}^{-2}x_{3}^{-2}x_{4}^{-2}}{(1-x_{1} )^7 (1-x_{1} x_{2} ) (1-x_{2} x_{3} x_{4} ) (1-x_{2} x_{3} ) }\\&
	+\displaystyle \frac{-x_{1} x_{2}^{-4}x_{3}^{-4}x_{4}^{-2}}{(1-x_{1} )^6 (1-x_{1} x_{2} ) (1-x_{2} x_{3} ) (1-x_{4} )^2 }\\&
	+\displaystyle \frac{-x_{1} x_{2}^{-2}x_{3}^{-2}x_{4}^{-2}}{(1-x_{1} )^6 (1-x_{1} x_{2} )^2 (1-x_{2} x_{3} ) (1-x_{3} x_{4} ) }\\&
	+\displaystyle \frac{x_{2}^{-4}x_{3}^{-2}x_{4}^{-2}}{(1-x_{1} )^5 (1-x_{2} x_{3} x_{4} ) (1-x_{3} )^2 (1-x_{4} )^2 }\\&
	+\displaystyle \frac{x_{2}^{-2}x_{4}^{-2}}{(1-x_{1} )^5 (1-x_{2} x_{3} x_{4} ) (1-x_{3} )^3 (1-x_{2} x_{3} ) }\\&
	+\displaystyle \frac{x_{2}^{-3}x_{4}^{-2}}{(1-x_{1} )^7 (1-x_{1} x_{2} x_{3} x_{4} ) (1-x_{3} ) (1-x_{2} x_{3} ) }\\&
	+\displaystyle \frac{x_{2}^{-4}x_{3}^{-2}x_{4}^{-3}}{(1-x_{1} )^6 (1-x_{2} x_{3} ) (1-x_{3} )^2 (1-x_{4} ) }\\&
	+\displaystyle \frac{-x_{2}^{-3}x_{4}^{-2}}{(1-x_{1} )^7 (1-x_{2} x_{3} ) (1-x_{3} ) (1-x_{3} x_{4} ) }\\&
	+\displaystyle \frac{x_{2}^{-3}x_{4}^{-2}-x_{2}^{-3}x_{3}^{-1}x_{4}^{-3}}{(1-x_{1} )^6 (1-x_{1} x_{2} x_{3} x_{4} ) (1-x_{3} )^2 (1-x_{2} x_{3} ) }\\&
	+\displaystyle \frac{-x_{2}^{-3}x_{4}^{-2}}{(1-x_{1} )^6 (1-x_{2} x_{3} ) (1-x_{3} )^2 (1-x_{3} x_{4} ) }\\&
	+\displaystyle \frac{-x_{2}^{-2}x_{4}^{-2}+x_{2}^{-3}x_{4}^{-2}}{(1-x_{1} )^5 (1-x_{1} x_{2} x_{3} x_{4} ) (1-x_{3} )^3 (1-x_{2} x_{3} ) }\\&
	+\displaystyle \frac{-x_{2}^{-3}x_{4}^{-2}}{(1-x_{1} )^5 (1-x_{2} x_{3} ) (1-x_{3} )^3 (1-x_{3} x_{4} ) }\\&
	+\displaystyle \frac{-x_{2}^{-4}x_{3}^{-2}x_{4}^{-2}}{(1-x_{1} )^5 (1-x_{2} x_{3} ) (1-x_{3} )^2 (1-x_{4} )^2 }\\&
	+\displaystyle \frac{-x_{2}^{-2}x_{3} x_{4}^{-1}-x_{2}^{-3}x_{3}^{-2}x_{4}^{-2}-x_{2}^{-3}x_{3}^{-2}x_{4}^{-3}-x_{2}^{-6}x_{3}^{-3}x_{4}^{-1}}{(1-x_{1} )^7 (1-x_{1} x_{2} x_{3} x_{4} ) (1-x_{2} x_{3} ) (1-x_{3} x_{4} ) }\\&
	+\displaystyle \frac{-x_{2}^{-6}x_{3}^{-4}x_{4}^{-2}}{(1-x_{1} )^7 (1-x_{2} ) (1-x_{2} x_{3} ) (1-x_{1} x_{2} x_{3} x_{4} ) }\\&
	+\displaystyle \frac{x_{2}^{-6}x_{3}^{-4}x_{4}^{-2}}{(1-x_{1} )^7 (1-x_{2} ) (1-x_{2} x_{3} ) (1-x_{3} x_{4} ) }\\&
	+\displaystyle \frac{x_{2}^{-5}x_{3}^{-4}x_{4}^{-2}}{(1-x_{1} )^6 (1-x_{2} )^2 (1-x_{2} x_{3} ) (1-x_{1} x_{2} x_{3} x_{4} ) }\\&
	+\displaystyle \frac{-x_{2}^{-5}x_{3}^{-4}x_{4}^{-2}}{(1-x_{1} )^6 (1-x_{2} )^2 (1-x_{2} x_{3} ) (1-x_{3} x_{4} ) }\\&
	+\displaystyle \frac{x_{2}^{-4}x_{3}^{-4}x_{4}^{-2}}{(1-x_{1} )^5 (1-x_{2} )^3 (1-x_{2} x_{3} ) (1-x_{3} x_{4} ) }\\&
	+\displaystyle \frac{x_{2}^{-3}x_{3}^{-4}x_{4}^{-2}}{(1-x_{1} )^4 (1-x_{2} )^4 (1-x_{2} x_{3} x_{4} ) (1-x_{2} x_{3} ) }\\&
	+\displaystyle \frac{-x_{2}^{-3}x_{3}^{-4}x_{4}^{-2}}{(1-x_{1} )^4 (1-x_{2} )^4 (1-x_{2} x_{3} ) (1-x_{3} x_{4} ) }\\&
	+\displaystyle \frac{x_{1} x_{2}^{-2}x_{3}^{-2}x_{4}^{-1}+x_{2}^{-3}x_{3}^{-2}x_{4}^{-1}+x_{2}^{-5}x_{3}^{-2}x_{4} }{(1-x_{1} )^7 (1-x_{3} x_{4} ) (1-x_{2} x_{3} x_{4} ) (1-x_{4} ) }\\&
	+\displaystyle \frac{x_{1} x_{2}^{-2}x_{3}^{-2}x_{4}^{-1}+x_{2}^{-3}x_{3}^{-2}x_{4}^{-1}+x_{2}^{-5}x_{3}^{-2}x_{4} }{(1-x_{1} )^6 (1-x_{3} x_{4} ) (1-x_{1} x_{2} x_{3} x_{4} ) (1-x_{4} )^2 }\\&
	+\displaystyle \frac{-x_{1} x_{2}^{-2}x_{3}^{-2}x_{4}^{-1}-x_{2}^{-3}x_{3}^{-2}x_{4}^{-1}-x_{2}^{-5}x_{3}^{-2}x_{4} }{(1-x_{1} )^6 (1-x_{3} x_{4} ) (1-x_{1} x_{2} x_{3} ) (1-x_{4} )^2 }\\&
	+\displaystyle \frac{x_{2}^{-3}x_{3}^{-3}x_{4}^{-1}-x_{2}^{-3}x_{3}^{-4}x_{4}^{-2}}{(1-x_{1} )^4 (1-x_{2} )^4 (1-x_{2} x_{3} x_{4} ) (1-x_{1} x_{2} x_{3} ) }\\&
	+\displaystyle \frac{-x_{2}^{-3}x_{3}^{-3}x_{4}^{-1}+x_{2}^{-3}x_{3}^{-4}x_{4}^{-2}}{(1-x_{1} )^4 (1-x_{2} )^4 (1-x_{1} x_{2} x_{3} ) (1-x_{3} x_{4} ) }\\&
	+\displaystyle \frac{-x_{1} x_{2}^{-2}x_{3}^{-2}x_{4}^{-1}-x_{2}^{-3}x_{3}^{-2}x_{4}^{-1}-x_{2}^{-5}x_{3}^{-2}x_{4} -3x_{2}^{-6}x_{3}^{-3}-x_{2}^{-4}x_{3}^{-3}x_{4}^{-3}}{(1-x_{1} )^7 (1-x_{1} x_{2} x_{3} x_{4} ) (1-x_{3} x_{4} ) (1-x_{4} ) }\\&
	+\displaystyle \frac{-3x_{2}^{-5}x_{3}^{-4}x_{4}^{-1}}{(1-x_{1} )^6 (1-x_{2} )^2 (1-x_{3} x_{4} ) (1-x_{4} ) }\\&
	+\displaystyle \frac{2x_{2}^{-4}x_{3}^{-4}x_{4}^{-1}}{(1-x_{1} )^5 (1-x_{2} )^3 (1-x_{3} x_{4} ) (1-x_{4} ) }\\&
	+\displaystyle \frac{3x_{2}^{-6}x_{3}^{-4}x_{4}^{-1}}{(1-x_{1} )^7 (1-x_{2} ) (1-x_{3} x_{4} ) (1-x_{4} ) }\\&
	+\displaystyle \frac{-x_{2}^{-3}x_{3}^{-4}x_{4}^{-1}}{(1-x_{1} )^4 (1-x_{2} )^4 (1-x_{3} x_{4} ) (1-x_{4} ) }\\&
	+\displaystyle \frac{-x_{2}^{-2}x_{3}^{-3}+x_{2}^{-3}x_{3}^{-3}x_{4} }{(1-x_{1} )^4 (1-x_{2} )^3 (1-x_{1} x_{2} x_{3} x_{4} ) (1-x_{4} )^2 }\\&
	+\displaystyle \frac{-x_{2}^{-2}x_{3}^{-3}+x_{2}^{-3}x_{3}^{-3}x_{4} -x_{2}^{-4}x_{3}^{-4}x_{4}^{-1}}{(1-x_{1} )^5 (1-x_{2} )^3 (1-x_{2} x_{3} x_{4} ) (1-x_{4} ) }\\&
	+\displaystyle \frac{x_{2}^{-2}x_{3}^{-3}-x_{2}^{-3}x_{3}^{-3}x_{4} }{(1-x_{1} )^4 (1-x_{2} )^3 (1-x_{1} x_{2} x_{3} ) (1-x_{4} )^2 }\\&
	+\displaystyle \frac{x_{2}^{-2}x_{3}^{-3}x_{4}^{-2}-x_{2}^{-4}x_{3}^{-4}x_{4}^{-2}}{(1-x_{1} )^5 (1-x_{2} )^3 (1-x_{2} x_{3} ) (1-x_{1} x_{2} x_{3} x_{4} ) }\\&
	+\displaystyle \frac{x_{2}^{-2}x_{3}^{-3}-x_{2}^{-3}x_{3}^{-3}x_{4} +x_{2}^{-3}x_{3}^{-4}x_{4}^{-2}-x_{2}^{-4}x_{3}^{-4}x_{4}^{-1}}{(1-x_{1} )^5 (1-x_{2} )^3 (1-x_{1} x_{2} x_{3} x_{4} ) (1-x_{4} ) }\\&
	+\displaystyle \frac{-x_{2}^{-3}x_{3}^{-4}x_{4}^{-2}}{(1-x_{1} )^5 (1-x_{2} )^3 (1-x_{2} x_{3} ) (1-x_{4} ) }\\&
	+\displaystyle \frac{-x_{2}^{-3}x_{3}^{-4}x_{4}^{-1}}{(1-x_{1} )^4 (1-x_{2} )^3 (1-x_{2} x_{3} x_{4} ) (1-x_{4} )^2 }\\&
	+\displaystyle \frac{x_{2}^{-3}x_{3}^{-4}x_{4}^{-1}}{(1-x_{1} )^4 (1-x_{2} )^3 (1-x_{2} x_{3} ) (1-x_{4} )^2 }\\&
	+\displaystyle \frac{3x_{2}^{-6}x_{3}^{-3}x_{4}^{-1}-x_{2}^{-6}x_{3}^{-4}x_{4}^{-2}-3x_{2}^{-6}x_{3}^{-5}x_{4}^{-3}}{(1-x_{1} )^7 (1-x_{2} ) (1-x_{3} x_{4} ) (1-x_{1} x_{2} x_{3} ) }\\&
	+\displaystyle \frac{-3x_{2}^{-5}x_{3}^{-3}x_{4}^{-1}+x_{2}^{-5}x_{3}^{-4}x_{4}^{-2}+3x_{2}^{-5}x_{3}^{-5}x_{4}^{-3}}{(1-x_{1} )^6 (1-x_{2} )^2 (1-x_{3} x_{4} ) (1-x_{1} x_{2} x_{3} ) }\\&
	+\displaystyle \frac{-x_{2}^{-2}x_{3}^{-3}x_{4}^{-1}+x_{2}^{-3}x_{3}^{-3}-x_{2}^{-4}x_{3}^{-3}x_{4}^{-1}+x_{2}^{-4}x_{3}^{-5}x_{4}^{-3}}{(1-x_{1} )^5 (1-x_{2} )^3 (1-x_{2} x_{3} x_{4} ) (1-x_{1} x_{2} x_{3} ) }\\&
	+\displaystyle \frac{2x_{2}^{-4}x_{3}^{-3}x_{4}^{-1}-x_{2}^{-4}x_{3}^{-4}x_{4}^{-2}-x_{2}^{-4}x_{3}^{-5}x_{4}^{-3}}{(1-x_{1} )^5 (1-x_{2} )^3 (1-x_{3} x_{4} ) (1-x_{1} x_{2} x_{3} ) }\\&
	+\displaystyle \frac{3x_{2}^{-6}x_{3}^{-5}x_{4}^{-3}}{(1-x_{1} )^7 (1-x_{2} ) (1-x_{3} ) (1-x_{3} x_{4} ) }\\&
	+\displaystyle \frac{x_{2}^{-4}x_{3}^{-5}x_{4}^{-3}}{(1-x_{1} )^5 (1-x_{2} )^3 (1-x_{3} ) (1-x_{3} x_{4} ) }\\&
	+\displaystyle \frac{~~~~\begin{array}{l}x_{1} x_{2}^{-2}x_{3}^{-2}x_{4}^{-2}+x_{2}^{-5}x_{3}^{-2}-x_{2}^{-3}x_{3}^{-3}x_{4}^{-2}-x_{2}^{-3}x_{3}^{-4}x_{4}^{-2}-x_{2}^{-4}x_{3}^{-4}x_{4}^{-2} \\ \hline  -2x_{2}^{-5}x_{3}^{-4}x_{4}^{-2}\end{array}~~~~}{(1-x_{1} )^7 (1-x_{3} x_{4} ) (1-x_{2} x_{3} x_{4} ) (1-x_{1} x_{2} x_{3} ) }\\&
	+\displaystyle \frac{~~~~\begin{array}{l}x_{2}^{-3}x_{3} x_{4}^{-1}-x_{2}^{-3}x_{3}^{-2}x_{4}^{-1}-x_{2}^{-3}x_{3}^{-3}x_{4}^{-1}-x_{2}^{-4}x_{3}^{-3}x_{4}^{-1}-2x_{2}^{-5}x_{3}^{-3}x_{4}^{-1} \\ \hline  -3x_{2}^{-6}x_{3}^{-4}x_{4}^{-2}\end{array}~~~~}{(1-x_{1} )^7 (1-x_{1} x_{2} x_{3} x_{4} ) (1-x_{3} ) (1-x_{3} x_{4} ) }\\&
	+\displaystyle \frac{x_{2}^{-3}x_{3}^{-2}x_{4}^{-1}+x_{2}^{-3}x_{3}^{-3}x_{4}^{-1}+x_{2}^{-4}x_{3}^{-3}x_{4}^{-1}+2x_{2}^{-5}x_{3}^{-3}x_{4}^{-1}}{(1-x_{1} )^7 (1-x_{3} x_{4} ) (1-x_{3} ) (1-x_{2} x_{3} x_{4} ) }\\&
	+\displaystyle \frac{x_{1} x_{2}^{-2}x_{4}^{-1}+x_{1} x_{2}^{-2}x_{3}^{-1}x_{4}^{-1}+x_{1} x_{2}^{-3}x_{3}^{-1}x_{4}^{-1}+2x_{1} x_{2}^{-4}x_{3}^{-1}x_{4}^{-1}}{(1-x_{1} )^6 (1-x_{3} x_{4} ) (1-x_{3} )^2 (1-x_{1} x_{2} x_{3} ) }\\&
	+\displaystyle \frac{-x_{2}^{-3}x_{3}^{-1}x_{4}^{-1}-x_{2}^{-3}x_{3}^{-2}x_{4}^{-1}-x_{2}^{-4}x_{3}^{-2}x_{4}^{-1}}{(1-x_{1} )^6 (1-x_{3} x_{4} ) (1-x_{3} )^2 (1-x_{2} x_{3} x_{4} ) }\\&
	+\displaystyle \frac{-x_{1} x_{2}^{-2}x_{3} x_{4}^{-1}-x_{1} x_{2}^{-2}x_{4}^{-1}-x_{1} x_{2}^{-3}x_{4}^{-1}}{(1-x_{1} )^5 (1-x_{3} x_{4} ) (1-x_{3} )^3 (1-x_{1} x_{2} x_{3} ) }\\&
	+\displaystyle \frac{-x_{1} x_{2}^{-2}x_{4}^{-1}}{(1-x_{1} )^4 (1-x_{1} x_{2} x_{3} x_{4} ) (1-x_{3} )^4 (1-x_{1} x_{2} x_{3} ) }\\&
	+\displaystyle \frac{-x_{1} x_{2}^{-2}x_{3} }{(1-x_{1} )^4 (1-x_{1} x_{2} x_{3} x_{4} ) (1-x_{3} )^4 (1-x_{3} x_{4} ) }\\&
	+\displaystyle \frac{x_{1} x_{2}^{-2}x_{3} x_{4}^{-1}}{(1-x_{1} )^4 (1-x_{3} x_{4} ) (1-x_{3} )^4 (1-x_{1} x_{2} x_{3} ) }\\&
	+\displaystyle \frac{x_{1} x_{2}^{-2}x_{3} +x_{1} x_{2}^{-2}+x_{1} x_{2}^{-3}+x_{2}^{-3}x_{3} x_{4}^{-1}-x_{2}^{-3}x_{3}^{-1}x_{4}^{-1}}{(1-x_{1} )^5 (1-x_{1} x_{2} x_{3} x_{4} ) (1-x_{3} )^3 (1-x_{3} x_{4} ) }\\&
	+\displaystyle \frac{x_{2}^{-3}x_{3}^{-1}x_{4}^{-1}}{(1-x_{1} )^5 (1-x_{2} x_{3} x_{4} ) (1-x_{3} )^3 (1-x_{3} x_{4} ) }\\&
	+\displaystyle \frac{-2x_{2}^{-5}x_{3}^{-4}x_{4}^{-3}}{(1-x_{1} )^6 (1-x_{2} ) (1-x_{3} )^2 (1-x_{3} x_{4} ) }\\&
	+\displaystyle \frac{~~~~\begin{array}{l}-x_{1} x_{2}^{-2}x_{3}^{-2}x_{4}^{-2}-x_{2}^{-5}x_{3}^{-2}+x_{2}^{-3}x_{3}^{-3}x_{4}^{-2}-3x_{2}^{-6}x_{3}^{-2}+x_{2}^{-3}x_{3}^{-4}x_{4}^{-2} \\ \hline  +x_{2}^{-4}x_{3}^{-4}x_{4}^{-2}+x_{2}^{-6}x_{3}^{-3}x_{4}^{-1}+2x_{2}^{-5}x_{3}^{-4}x_{4}^{-2}+3x_{2}^{-6}x_{3}^{-4}x_{4}^{-2}\end{array}~~~~}{(1-x_{1} )^7 (1-x_{1} x_{2} x_{3} x_{4} ) (1-x_{3} x_{4} ) (1-x_{1} x_{2} x_{3} ) }\\&
	+\displaystyle \frac{-x_{1} x_{2}^{-2}x_{3}^{-1}x_{4}^{-2}}{(1-x_{1} )^4 (1-x_{1} x_{2} x_{3} x_{4} ) (1-x_{3} )^3 (1-x_{1} x_{2} )^2 }\\&
	+\displaystyle \frac{x_{1} x_{2}^{-2}x_{3}^{-1}x_{4}^{-2}}{(1-x_{1} )^4 (1-x_{1} x_{2} )^2 (1-x_{3} )^3 (1-x_{3} x_{4} ) }\\&
	+\displaystyle \frac{-x_{2}^{-3}x_{3}^{-4}x_{4}^{-3}}{(1-x_{1} )^4 (1-x_{2} )^3 (1-x_{3} )^2 (1-x_{3} x_{4} ) }\\&
	+\displaystyle \frac{-x_{2}^{-5}x_{3}^{-3}-3x_{2}^{-6}x_{3}^{-3}x_{4}^{-1}+2x_{2}^{-5}x_{3}^{-4}x_{4}^{-2}+x_{2}^{-6}x_{3}^{-4}x_{4}^{-2}+3x_{2}^{-6}x_{3}^{-5}x_{4}^{-3}}{(1-x_{1} )^7 (1-x_{2} ) (1-x_{1} x_{2} x_{3} ) (1-x_{1} x_{2} x_{3} x_{4} ) }\\&
	+\displaystyle \frac{x_{2}^{-5}x_{3}^{-3}-2x_{2}^{-5}x_{3}^{-4}x_{4}^{-2}}{(1-x_{1} )^7 (1-x_{2} ) (1-x_{1} x_{2} x_{3} ) (1-x_{2} x_{3} x_{4} ) }\\&
	+\displaystyle \frac{x_{2}^{-5}x_{3}^{-3}x_{4} -2x_{2}^{-5}x_{3}^{-4}x_{4}^{-1}}{(1-x_{1} )^7 (1-x_{2} ) (1-x_{2} x_{3} x_{4} ) (1-x_{4} ) }\\&
	+\displaystyle \frac{~~~~\begin{array}{l}-x_{2}^{-3}x_{3}^{-3}x_{4}^{-2}-x_{2}^{-3}x_{3}^{-4}x_{4}^{-2}-x_{2}^{-4}x_{3}^{-4}x_{4}^{-2}+2x_{2}^{-5}x_{3}^{-3}x_{4}^{-2} \\ \hline  -2x_{2}^{-5}x_{3}^{-4}x_{4}^{-2}\end{array}~~~~}{(1-x_{1} )^7 (1-x_{1} x_{2} x_{3} ) (1-x_{3} ) (1-x_{1} x_{2} x_{3} x_{4} ) }\\&
	+\displaystyle \frac{~~~~\begin{array}{l}x_{2}^{-3}x_{3}^{-3}x_{4}^{-2}+x_{2}^{-3}x_{3}^{-4}x_{4}^{-2}+x_{2}^{-4}x_{3}^{-4}x_{4}^{-2}-2x_{2}^{-5}x_{3}^{-3}x_{4}^{-2} \\ \hline  +2x_{2}^{-5}x_{3}^{-4}x_{4}^{-2}\end{array}~~~~}{(1-x_{1} )^7 (1-x_{1} x_{2} x_{3} ) (1-x_{3} ) (1-x_{2} x_{3} x_{4} ) }\\&
	+\displaystyle \frac{x_{2}^{-5}x_{3}^{-3}x_{4} -2x_{2}^{-5}x_{3}^{-4}x_{4}^{-1}}{(1-x_{1} )^6 (1-x_{2} ) (1-x_{1} x_{2} x_{3} x_{4} ) (1-x_{4} )^2 }\\&
	+\displaystyle \frac{-2x_{2}^{-5}x_{3}^{-3}x_{4}^{-1}}{(1-x_{1} )^6 (1-x_{1} x_{2} x_{3} x_{4} ) (1-x_{3} ) (1-x_{4} )^2 }\\&
	+\displaystyle \frac{2x_{2}^{-5}x_{3}^{-3}x_{4}^{-1}}{(1-x_{1} )^6 (1-x_{1} x_{2} x_{3} ) (1-x_{3} ) (1-x_{4} )^2 }\\&
	+\displaystyle \frac{-x_{2}^{-5}x_{3}^{-3}x_{4} +2x_{2}^{-5}x_{3}^{-4}x_{4}^{-1}}{(1-x_{1} )^6 (1-x_{2} ) (1-x_{1} x_{2} x_{3} ) (1-x_{4} )^2 }\\&
	+\displaystyle \frac{~~~~\begin{array}{l}x_{1} x_{2}^{-2}x_{4}^{-1}+x_{1} x_{2}^{-2}x_{3}^{-1}x_{4}^{-1}+x_{1} x_{2}^{-3}x_{3}^{-1}x_{4}^{-1}+x_{2}^{-3}x_{3}^{-1}x_{4}^{-2} \\ \hline  -x_{2}^{-3}x_{3}^{-2}x_{4}^{-2}\end{array}~~~~}{(1-x_{1} )^5 (1-x_{1} x_{2} x_{3} x_{4} ) (1-x_{3} )^3 (1-x_{1} x_{2} x_{3} ) }\\&
	+\displaystyle \frac{-x_{2}^{-3}x_{3}^{-1}x_{4}^{-2}+x_{2}^{-3}x_{3}^{-2}x_{4}^{-2}}{(1-x_{1} )^5 (1-x_{1} x_{2} x_{3} ) (1-x_{3} )^3 (1-x_{2} x_{3} x_{4} ) }\\&
	+\displaystyle \frac{x_{2}^{-2}x_{3}^{-3}x_{4}^{-1}-x_{2}^{-3}x_{3}^{-3}-x_{2}^{-4}x_{3}^{-3}x_{4}^{-1}+x_{2}^{-4}x_{3}^{-4}x_{4}^{-2}}{(1-x_{1} )^5 (1-x_{2} )^3 (1-x_{1} x_{2} x_{3} x_{4} ) (1-x_{1} x_{2} x_{3} ) }\\&
	+\displaystyle \frac{-x_{2}^{-3}x_{3}^{-1}x_{4}^{-1}}{(1-x_{1} )^4 (1-x_{1} x_{2} x_{3} x_{4} ) (1-x_{3} )^3 (1-x_{4} )^2 }\\&
	+\displaystyle \frac{x_{2}^{-3}x_{3}^{-1}x_{4}^{-1}}{(1-x_{1} )^4 (1-x_{1} x_{2} x_{3} ) (1-x_{3} )^3 (1-x_{4} )^2 }\\&
	+\displaystyle \frac{~~~~\begin{array}{l}-x_{1} x_{2}^{-2}x_{3}^{-1}x_{4}^{-1}-x_{1} x_{2}^{-2}x_{3}^{-2}x_{4}^{-1}-x_{1} x_{2}^{-3}x_{3}^{-2}x_{4}^{-1}-2x_{1} x_{2}^{-4}x_{3}^{-2}x_{4}^{-1} \\ \hline  -x_{2}^{-3}x_{3}^{-1}x_{4}^{-2}+x_{2}^{-3}x_{3}^{-2}x_{4}^{-2}+x_{2}^{-3}x_{3}^{-3}x_{4}^{-2}-x_{2}^{-4}x_{3}^{-2}x_{4}^{-2} \\ \hline  +x_{2}^{-4}x_{3}^{-3}x_{4}^{-2}\end{array}~~~~}{(1-x_{1} )^6 (1-x_{1} x_{2} x_{3} ) (1-x_{3} )^2 (1-x_{1} x_{2} x_{3} x_{4} ) }\\&
	+\displaystyle \frac{~~~~\begin{array}{l}x_{2}^{-3}x_{3}^{-1}x_{4}^{-2}-x_{2}^{-3}x_{3}^{-2}x_{4}^{-2}-x_{2}^{-3}x_{3}^{-3}x_{4}^{-2}+x_{2}^{-4}x_{3}^{-2}x_{4}^{-2} \\ \hline  -x_{2}^{-4}x_{3}^{-3}x_{4}^{-2}\end{array}~~~~}{(1-x_{1} )^6 (1-x_{1} x_{2} x_{3} ) (1-x_{3} )^2 (1-x_{2} x_{3} x_{4} ) }\\&
	+\displaystyle \frac{x_{2}^{-4}x_{3}^{-3}+2x_{2}^{-5}x_{3}^{-3}x_{4}^{-1}-x_{2}^{-4}x_{3}^{-4}x_{4}^{-2}-x_{2}^{-5}x_{3}^{-4}x_{4}^{-2}-x_{2}^{-5}x_{3}^{-5}x_{4}^{-3}}{(1-x_{1} )^6 (1-x_{2} )^2 (1-x_{1} x_{2} x_{3} ) (1-x_{1} x_{2} x_{3} x_{4} ) }\\&
	+\displaystyle \frac{-x_{2}^{-4}x_{3}^{-3}+x_{2}^{-5}x_{3}^{-3}x_{4}^{-1}+x_{2}^{-4}x_{3}^{-4}x_{4}^{-2}-2x_{2}^{-5}x_{3}^{-5}x_{4}^{-3}}{(1-x_{1} )^6 (1-x_{2} )^2 (1-x_{1} x_{2} x_{3} ) (1-x_{2} x_{3} x_{4} ) }\\&
	+\displaystyle \frac{-x_{2}^{-4}x_{3}^{-3}x_{4} +x_{2}^{-4}x_{3}^{-4}x_{4}^{-1}+x_{2}^{-5}x_{3}^{-4}x_{4}^{-1}-2x_{2}^{-5}x_{3}^{-5}x_{4}^{-2}}{(1-x_{1} )^6 (1-x_{2} )^2 (1-x_{2} x_{3} x_{4} ) (1-x_{4} ) }\\&
	+\displaystyle \frac{x_{2}^{-3}x_{3}^{-1}x_{4}^{-1}+x_{2}^{-4}x_{3}^{-2}x_{4}^{-1}}{(1-x_{1} )^5 (1-x_{1} x_{2} x_{3} x_{4} ) (1-x_{3} )^2 (1-x_{4} )^2 }\\&
	+\displaystyle \frac{-x_{2}^{-4}x_{3}^{-3}x_{4} +x_{2}^{-4}x_{3}^{-4}x_{4}^{-1}}{(1-x_{1} )^5 (1-x_{2} )^2 (1-x_{1} x_{2} x_{3} x_{4} ) (1-x_{4} )^2 }\\&
	+\displaystyle \frac{x_{2}^{-4}x_{3}^{-3}x_{4} -x_{2}^{-4}x_{3}^{-4}x_{4}^{-1}}{(1-x_{1} )^5 (1-x_{2} )^2 (1-x_{1} x_{2} x_{3} ) (1-x_{4} )^2 }\\&
	+\displaystyle \frac{-x_{2}^{-3}x_{3}^{-1}x_{4}^{-1}-x_{2}^{-4}x_{3}^{-2}x_{4}^{-1}}{(1-x_{1} )^5 (1-x_{1} x_{2} x_{3} ) (1-x_{3} )^2 (1-x_{4} )^2 }\\&
	+\displaystyle \frac{-3x_{2}^{-5}x_{3}^{-5}x_{4}^{-3}}{(1-x_{1} )^6 (1-x_{2} )^2 (1-x_{3} ) (1-x_{3} x_{4} ) }\\&
	+\displaystyle \frac{x_{1} x_{2}^{-2}x_{4}^{-2}+x_{1} x_{2}^{-2}x_{3}^{-1}x_{4}^{-2}-x_{2}^{-3}x_{4}^{-2}}{(1-x_{1} )^5 (1-x_{1} x_{2} x_{3} x_{4} ) (1-x_{3} )^3 (1-x_{1} x_{2} ) }\\&
	+\displaystyle \frac{x_{1} x_{2}^{-2}x_{4}^{-2}-x_{2}^{-3}x_{4}^{-2}}{(1-x_{1} )^6 (1-x_{1} x_{2} x_{3} x_{4} ) (1-x_{3} )^2 (1-x_{1} x_{2} ) }\\&
	+\displaystyle \frac{-x_{2}^{-3}x_{4}^{-2}}{(1-x_{1} )^7 (1-x_{1} x_{2} x_{3} x_{4} ) (1-x_{3} ) (1-x_{1} x_{2} ) }\\&
	+\displaystyle \frac{-x_{1}^2x_{2}^{-1}x_{3}^{-2}x_{4}^{-1}-x_{1} x_{2}^{-3}x_{3}^{-3}x_{4}^{-1}-x_{2}^{-3}x_{3}^{-1}x_{4}^{-2}-x_{1} x_{2}^{-4}x_{3}^{-4}x_{4}^{-3}}{(1-x_{1} )^7 (1-x_{1} x_{2} x_{3} x_{4} ) (1-x_{1} x_{2} ) (1-x_{4} ) }\\&
	+\displaystyle \frac{x_{1}^2x_{2}^{-1}x_{3}^{-2}x_{4}^{-1}+x_{1} x_{2}^{-3}x_{3}^{-3}x_{4}^{-1}+x_{2}^{-3}x_{3}^{-1}x_{4}^{-2}}{(1-x_{1} )^7 (1-x_{1} x_{2} ) (1-x_{2} x_{3} x_{4} ) (1-x_{4} ) }\\&
	+\displaystyle \frac{x_{2}^{-3}x_{4}^{-2}}{(1-x_{1} )^7 (1-x_{1} x_{2} ) (1-x_{3} ) (1-x_{2} x_{3} x_{4} ) }\\&
	+\displaystyle \frac{x_{2}^{-3}x_{4}^{-1}-2x_{2}^{-5}x_{3}^{-3}x_{4}^{-1}}{(1-x_{1} )^7 (1-x_{2} x_{3} x_{4} ) (1-x_{3} ) (1-x_{4} ) }\\&
	+\displaystyle \frac{-x_{2}^{-3}x_{4}^{-2}}{(1-x_{1} )^6 (1-x_{1} x_{2} ) (1-x_{3} )^2 (1-x_{4} ) }\\&
	+\displaystyle \frac{x_{2}^{-3}x_{4}^{-2}}{(1-x_{1} )^6 (1-x_{1} x_{2} ) (1-x_{3} )^2 (1-x_{2} x_{3} x_{4} ) }\\&
	+\displaystyle \frac{x_{2}^{-3}x_{4}^{-1}+x_{2}^{-3}x_{3}^{-1}x_{4}^{-1}+x_{2}^{-4}x_{3}^{-2}x_{4}^{-1}}{(1-x_{1} )^6 (1-x_{2} x_{3} x_{4} ) (1-x_{3} )^2 (1-x_{4} ) }\\&
	+\displaystyle \frac{-x_{2}^{-3}x_{4}^{-2}}{(1-x_{1} )^5 (1-x_{1} x_{2} ) (1-x_{3} )^3 (1-x_{4} ) }\\&
	+\displaystyle \frac{-x_{1} x_{2}^{-2}x_{3}^{-1}x_{4}^{-2}+x_{2}^{-3}x_{4}^{-2}}{(1-x_{1} )^5 (1-x_{1} x_{2} ) (1-x_{3} )^3 (1-x_{2} x_{3} x_{4} ) }\\&
	+\displaystyle \frac{-x_{2}^{-3}x_{4}^{-1}+x_{2}^{-3}x_{4}^{-2}+x_{2}^{-3}x_{3}^{-1}x_{4}^{-1}}{(1-x_{1} )^5 (1-x_{1} x_{2} x_{3} x_{4} ) (1-x_{3} )^3 (1-x_{4} ) }\\&
	+\displaystyle \frac{x_{2}^{-3}x_{4}^{-1}-x_{2}^{-3}x_{3}^{-1}x_{4}^{-1}}{(1-x_{1} )^5 (1-x_{2} x_{3} x_{4} ) (1-x_{3} )^3 (1-x_{4} ) }\\&
	+\displaystyle \frac{x_{1} x_{2}^{-2}x_{4}^{-2}}{(1-x_{1} )^4 (1-x_{1} x_{2} x_{3} x_{4} ) (1-x_{3} )^4 (1-x_{1} x_{2} ) }\\&
	+\displaystyle \frac{x_{2}^{-4}x_{3}^{-3}x_{4} -x_{2}^{-4}x_{3}^{-4}x_{4}^{-1}+2x_{2}^{-5}x_{3}^{-4}x_{4}^{-1}-2x_{2}^{-5}x_{3}^{-5}x_{4}^{-2}}{(1-x_{1} )^6 (1-x_{2} )^2 (1-x_{1} x_{2} x_{3} x_{4} ) (1-x_{4} ) }\\&
	+\displaystyle \frac{-x_{2}^{-5}x_{3}^{-5}x_{4}^{-3}}{(1-x_{1} )^6 (1-x_{2} )^2 (1-x_{3} ) (1-x_{1} x_{2} x_{3} x_{4} ) }\\&
	+\displaystyle \frac{4x_{2}^{-5}x_{3}^{-5}x_{4}^{-3}}{(1-x_{1} )^6 (1-x_{2} )^2 (1-x_{3} ) (1-x_{4} ) }\\&
	+\displaystyle \frac{x_{2}^{-3}x_{4}^{-2}}{(1-x_{1} )^4 (1-x_{1} x_{2} x_{3} x_{4} ) (1-x_{3} )^4 (1-x_{4} ) }\\&
	+\displaystyle \frac{-x_{2}^{-3}x_{3}^{-5}x_{4}^{-3}}{(1-x_{1} )^4 (1-x_{2} )^4 (1-x_{3} ) (1-x_{4} ) }\\&
	+\displaystyle \frac{-x_{2}^{-3}x_{4}^{-1}+2x_{2}^{-5}x_{3}^{-3}x_{4}^{-1}+6x_{2}^{-6}x_{3}^{-4}x_{4}^{-2}}{(1-x_{1} )^7 (1-x_{1} x_{2} x_{3} x_{4} ) (1-x_{3} ) (1-x_{4} ) }\\&
	+\displaystyle \frac{-x_{2}^{-5}x_{3}^{-3}x_{4} +2x_{2}^{-5}x_{3}^{-4}x_{4}^{-1}-3x_{2}^{-6}x_{3}^{-4}x_{4}^{-1}+6x_{2}^{-6}x_{3}^{-5}x_{4}^{-2}}{(1-x_{1} )^7 (1-x_{2} ) (1-x_{1} x_{2} x_{3} x_{4} ) (1-x_{4} ) }\\&
	+\displaystyle \frac{3x_{2}^{-6}x_{3}^{-5}x_{4}^{-3}}{(1-x_{1} )^7 (1-x_{2} ) (1-x_{3} ) (1-x_{1} x_{2} x_{3} x_{4} ) }\\&
	+\displaystyle \frac{-6x_{2}^{-6}x_{3}^{-5}x_{4}^{-3}}{(1-x_{1} )^7 (1-x_{2} ) (1-x_{3} ) (1-x_{4} ) }\\&
	+\displaystyle \frac{~~~~\begin{array}{l}-x_{2}^{-3}x_{4}^{-1}+x_{2}^{-3}x_{4}^{-2}-x_{2}^{-3}x_{3}^{-1}x_{4}^{-1}-x_{2}^{-4}x_{3}^{-2}x_{4}^{-1}-x_{2}^{-4}x_{3}^{-2}x_{4}^{-3} \\ \hline  -2x_{2}^{-5}x_{3}^{-3}x_{4}^{-2}\end{array}~~~~}{(1-x_{1} )^6 (1-x_{1} x_{2} x_{3} x_{4} ) (1-x_{3} )^2 (1-x_{4} ) }\\&
	+\displaystyle \frac{-x_{2}^{-4}x_{3}^{-5}x_{4}^{-3}}{(1-x_{1} )^5 (1-x_{2} )^3 (1-x_{3} ) (1-x_{2} x_{3} x_{4} ) }\\&
	+\displaystyle \frac{2x_{2}^{-5}x_{3}^{-4}x_{4}^{-3}}{(1-x_{1} )^6 (1-x_{2} ) (1-x_{3} )^2 (1-x_{4} ) }\\&
	+\displaystyle \frac{x_{2}^{-3}x_{3}^{-5}x_{4}^{-3}}{(1-x_{1} )^4 (1-x_{2} )^4 (1-x_{3} ) (1-x_{2} x_{3} x_{4} ) }\\&
	+\displaystyle \frac{2x_{2}^{-4}x_{3}^{-4}x_{4}^{-3}}{(1-x_{1} )^5 (1-x_{2} )^2 (1-x_{3} )^2 (1-x_{3} x_{4} ) }\\&
	+\displaystyle \frac{-2x_{2}^{-4}x_{3}^{-4}x_{4}^{-3}}{(1-x_{1} )^5 (1-x_{2} )^2 (1-x_{3} )^2 (1-x_{4} ) }\\&
	+\displaystyle \frac{x_{2}^{-3}x_{3}^{-4}x_{4}^{-1}+x_{2}^{-3}x_{3}^{-5}x_{4}^{-2}}{(1-x_{1} )^4 (1-x_{2} )^4 (1-x_{2} x_{3} x_{4} ) (1-x_{4} ) }\\&
	+\displaystyle \frac{~~~~\begin{array}{l}-x_{1} x_{2}^{-2}-x_{1} x_{2}^{-2}x_{3}^{-1}-x_{1} x_{2}^{-3}x_{3}^{-1}+x_{2}^{-3}x_{3} x_{4}^{-1}-2x_{1} x_{2}^{-4}x_{3}^{-1} \\ \hline  +x_{2}^{-3}x_{3}^{-1}x_{4}^{-1}+x_{2}^{-3}x_{3}^{-2}x_{4}^{-1}+x_{2}^{-4}x_{3}^{-2}x_{4}^{-1}+2x_{2}^{-5}x_{3}^{-3}x_{4}^{-2}\end{array}~~~~}{(1-x_{1} )^6 (1-x_{1} x_{2} x_{3} x_{4} ) (1-x_{3} )^2 (1-x_{3} x_{4} ) }\\&
	+\displaystyle \frac{x_{2}^{-3}x_{3}^{-4}x_{4}^{-3}}{(1-x_{1} )^4 (1-x_{2} )^3 (1-x_{3} )^2 (1-x_{4} ) }\\=&
	-x_{2}^{-3}x_{4}^{-2}\cdot \frac{1}{36}\left(x_{1} \partial_{1} -x_{2} \partial_{2} \right)^3\left(x_{3} \partial_{3} \right)^3\cdot\frac{1}{(1-x_{1} ) (1-x_{1} x_{2} ) (1-x_{3} ) (1-x_{4} ) }\\&
	\displaystyle +x_{2}^{-4}x_{3}^{-3}x_{4}^{-2}\cdot \frac{1}{120}\left(x_{1} \partial_{1} \right)^5\left(-x_{3} \partial_{3} +x_{4} \partial_{4} \right)\cdot\frac{1}{(1-x_{1} ) (1-x_{3} x_{4} ) (1-x_{2} x_{3} x_{4} ) (1-x_{4} ) }\\&
	\displaystyle +x_{2}^{-4}x_{3}^{-3}x_{4}^{-3}\cdot \frac{1}{720}\left(x_{1} \partial_{1} \right)^6\cdot\frac{1}{(1-x_{1} ) (1-x_{3} x_{4} ) (1-x_{2} x_{3} ) (1-x_{4} ) }\\&
	\displaystyle -x_{2}^{-4}x_{3}^{-3}x_{4}^{-2}\cdot \frac{1}{120}\left(x_{1} \partial_{1} \right)^5\left(x_{2} \partial_{2} -x_{3} \partial_{3} +x_{4} \partial_{4} \right)\cdot\frac{1}{(1-x_{1} ) (1-x_{3} x_{4} ) (1-x_{2} x_{3} ) (1-x_{4} ) }\\&
	\displaystyle +x_{1} x_{2}^{-2}x_{3}^{-2}x_{4}^{-1}\cdot \frac{1}{48}\left(x_{1} \partial_{1} -x_{2} \partial_{2} \right)^4\left(-x_{3} \partial_{3} +x_{4} \partial_{4} \right)^2\cdot\frac{1}{(1-x_{1} ) (1-x_{1} x_{2} x_{3} x_{4} ) (1-x_{3} x_{4} ) (1-x_{4} ) }\\&
	\displaystyle -x_{1} x_{2}^{-2}x_{3}^{-2}x_{4}^{-1}\cdot \frac{1}{48}\left(x_{1} \partial_{1} -x_{2} \partial_{2} \right)^4\left(x_{2} \partial_{2} -x_{3} \partial_{3} +x_{4} \partial_{4} \right)^2\cdot\frac{1}{(1-x_{1} ) (1-x_{3} x_{4} ) (1-x_{1} x_{2} x_{3} ) (1-x_{4} ) }\\&
	\displaystyle -x_{1}^2x_{2}^{-1}x_{3}^{-2}x_{4}^{-1}\cdot \frac{1}{48}\left(x_{1} \partial_{1} -x_{2} \partial_{2} \right)^4\left(x_{4} \partial_{4} \right)^2\cdot\frac{1}{(1-x_{1} ) (1-x_{1} x_{2} ) (1-x_{1} x_{2} x_{3} ) (1-x_{4} ) }\\&
	\displaystyle +x_{1}^2x_{2}^{-1}x_{3}^{-2}x_{4}^{-1}\cdot \frac{1}{48}\left(x_{1} \partial_{1} -x_{2} \partial_{2} \right)^4\left(-x_{3} \partial_{3} +x_{4} \partial_{4} \right)^2\cdot\frac{1}{(1-x_{1} ) (1-x_{1} x_{2} x_{3} x_{4} ) (1-x_{1} x_{2} ) (1-x_{4} ) }\\&
	\displaystyle -x_{1} x_{2}^{-2}x_{3}^{-2}x_{4}^{-2}\cdot \frac{1}{120}\left(x_{1} \partial_{1} -x_{2} \partial_{2} \right)^5\left(x_{2} \partial_{2} -x_{3} \partial_{3} \right)\cdot\frac{1}{(1-x_{1} ) (1-x_{1} x_{2} x_{3} x_{4} ) (1-x_{1} x_{2} ) (1-x_{1} x_{2} x_{3} ) }\\&
	\displaystyle -x_{1}^2x_{2}^{-1}x_{3}^{-2}x_{4}^{-1}\cdot \frac{1}{24}\left(x_{1} \partial_{1} -x_{2} \partial_{2} \right)^4\left(-x_{3} \partial_{3} +x_{4} \partial_{4} \right)\left(x_{2} \partial_{2} -x_{3} \partial_{3} \right)\cdot\frac{1}{(1-x_{1} ) (1-x_{1} x_{2} x_{3} x_{4} ) (1-x_{1} x_{2} ) (1-x_{4} ) }\\&
	\displaystyle +x_{1}^2x_{2}^{-1}x_{3}^{-2}x_{4}^{-1}\cdot \frac{1}{24}\left(x_{1} \partial_{1} -x_{2} \partial_{2} \right)^4\left(-x_{3} \partial_{3} +x_{4} \partial_{4} \right)\left(x_{2} \partial_{2} \right)\cdot\frac{1}{(1-x_{1} ) (1-x_{1} x_{2} ) (1-x_{3} x_{4} ) (1-x_{4} ) }\\&
	\displaystyle +x_{1} x_{2}^{-2}x_{3}^{-2}x_{4}^{-2}\cdot \frac{1}{120}\left(x_{1} \partial_{1} -x_{2} \partial_{2} \right)^5\left(x_{2} \partial_{2} -x_{3} \partial_{3} +x_{4} \partial_{4} \right)\cdot\frac{1}{(1-x_{1} ) (1-x_{1} x_{2} ) (1-x_{3} x_{4} ) (1-x_{1} x_{2} x_{3} ) }\\&
	\displaystyle +\left(\begin{array}{l}-x_{1}^2x_{2}^{-1}x_{3}^{-2}x_{4}^{-2}+x_{1} x_{2}^{-2}x_{3}^{-2}x_{4}^{-2} \\ -x_{1} x_{2}^{-3}x_{3}^{-3}x_{4}^{-2}\end{array}\right)\begin{array}{l}\cdot \frac{1}{720}\left(x_{1} \partial_{1} -x_{2} \partial_{2} \right)^6\\
		~~\cdot\frac{1}{(1-x_{1} ) (1-x_{1} x_{2} x_{3} x_{4} ) (1-x_{1} x_{2} ) (1-x_{1} x_{2} x_{3} ) }\end{array}\\&
	\displaystyle +\left(\begin{array}{l}x_{1}^2x_{2}^{-1}x_{3}^{-2}x_{4}^{-2}-x_{1} x_{2}^{-2}x_{3}^{-2}x_{4}^{-2} \\ +x_{1} x_{2}^{-3}x_{3}^{-3}x_{4}^{-2}\end{array}\right)\begin{array}{l}\cdot \frac{1}{720}\left(x_{1} \partial_{1} -x_{2} \partial_{2} +x_{4} \partial_{4} \right)^6\\
		~~\cdot\frac{1}{(1-x_{1} ) (1-x_{1} x_{2} ) (1-x_{2} x_{3} x_{4} ) (1-x_{1} x_{2} x_{3} ) }\end{array}\\&
	\displaystyle +\left(x_{1}^2x_{2}^{-1}x_{3}^{-2}x_{4}^{-1}+x_{1} x_{2}^{-3}x_{3}^{-3}x_{4}^{-1}\right)\cdot \frac{1}{120}\left(x_{1} \partial_{1} -x_{2} \partial_{2} \right)^5\left(-x_{3} \partial_{3} +x_{4} \partial_{4} \right)\cdot\frac{1}{(1-x_{1} ) (1-x_{1} x_{2} x_{3} x_{4} ) (1-x_{1} x_{2} ) (1-x_{4} ) }\\&
	\displaystyle +\left(-x_{1}^2x_{2}^{-1}x_{3}^{-2}x_{4}^{-1}-x_{1} x_{2}^{-3}x_{3}^{-3}x_{4}^{-1}\right)\cdot \frac{1}{120}\left(x_{1} \partial_{1} -x_{2} \partial_{2} \right)^5\left(x_{4} \partial_{4} \right)\cdot\frac{1}{(1-x_{1} ) (1-x_{1} x_{2} ) (1-x_{1} x_{2} x_{3} ) (1-x_{4} ) }\\&
	\displaystyle +x_{2}^{-3}x_{3}^{-2}x_{4}^{-2}\cdot \frac{1}{720}\left(x_{1} \partial_{1} \right)^6\cdot\frac{1}{(1-x_{1} ) (1-x_{2} x_{3} ) (1-x_{2} x_{3} x_{4} ) (1-x_{3} x_{4} ) }\\&
	\displaystyle +x_{1} x_{2}^{-4}x_{3}^{-4}x_{4}^{-2}\cdot \frac{1}{120}\left(x_{1} \partial_{1} -x_{2} \partial_{2} +x_{3} \partial_{3} \right)^5\left(-x_{3} \partial_{3} +x_{4} \partial_{4} \right)\cdot\frac{1}{(1-x_{1} ) (1-x_{1} x_{2} ) (1-x_{2} x_{3} x_{4} ) (1-x_{4} ) }\\&
	\displaystyle +\left(-x_{1} x_{2}^{-2}x_{3}^{-2}x_{4}^{-2}-x_{1} x_{2}^{-3}x_{3}^{-3}x_{4}^{-3}\right)\cdot \frac{1}{720}\left(\begin{array}{l}x_{1} \partial_{1} -x_{2} \partial_{2} +x_{3} \partial_{3}  \\ -x_{4} \partial_{4} \end{array}\right)^6\cdot\frac{1}{(1-x_{1} ) (1-x_{1} x_{2} x_{3} x_{4} ) (1-x_{2} x_{3} ) (1-x_{1} x_{2} ) }\\&
	\displaystyle +x_{1} x_{2}^{-4}x_{3}^{-4}x_{4}^{-3}\cdot \frac{1}{720}\left(x_{1} \partial_{1} -x_{2} \partial_{2} +x_{3} \partial_{3} \right)^6\cdot\frac{1}{(1-x_{1} ) (1-x_{1} x_{2} ) (1-x_{2} x_{3} ) (1-x_{4} ) }\\&
	\displaystyle +x_{1} x_{2}^{-2}x_{3}^{-2}x_{4}^{-2}\cdot \frac{1}{120}\left(\begin{array}{l}x_{1} \partial_{1} -x_{2} \partial_{2} +x_{3} \partial_{3}  \\ -x_{4} \partial_{4} \end{array}\right)^5\left(x_{2} \partial_{2} -x_{3} \partial_{3} \right)\cdot\frac{1}{(1-x_{1} ) (1-x_{1} x_{2} x_{3} x_{4} ) (1-x_{2} x_{3} ) (1-x_{1} x_{2} ) }\\&
	\displaystyle +x_{1} x_{2}^{-2}x_{3}^{-2}x_{4}^{-2}\cdot \frac{1}{720}\left(x_{1} \partial_{1} -x_{2} \partial_{2} +x_{3} \partial_{3} \right)^6\cdot\frac{1}{(1-x_{1} ) (1-x_{1} x_{2} ) (1-x_{2} x_{3} x_{4} ) (1-x_{2} x_{3} ) }\\&
	\displaystyle -x_{1} x_{2}^{-4}x_{3}^{-4}x_{4}^{-2}\cdot \frac{1}{120}\left(x_{1} \partial_{1} -x_{2} \partial_{2} +x_{3} \partial_{3} \right)^5\left(x_{4} \partial_{4} \right)\cdot\frac{1}{(1-x_{1} ) (1-x_{1} x_{2} ) (1-x_{2} x_{3} ) (1-x_{4} ) }\\&
	\displaystyle -x_{1} x_{2}^{-2}x_{3}^{-2}x_{4}^{-2}\cdot \frac{1}{120}\left(\begin{array}{l}x_{1} \partial_{1} -x_{2} \partial_{2} +x_{3} \partial_{3}  \\ -x_{4} \partial_{4} \end{array}\right)^5\left(x_{2} \partial_{2} -x_{3} \partial_{3} +x_{4} \partial_{4} \right)\cdot\frac{1}{(1-x_{1} ) (1-x_{1} x_{2} ) (1-x_{2} x_{3} ) (1-x_{3} x_{4} ) }\\&
	\displaystyle +x_{2}^{-4}x_{3}^{-2}x_{4}^{-2}\cdot \frac{1}{24}\left(x_{1} \partial_{1} \right)^4\left(-x_{2} \partial_{2} +x_{3} \partial_{3} \right)\left(-x_{2} \partial_{2} +x_{4} \partial_{4} \right)\cdot\frac{1}{(1-x_{1} ) (1-x_{2} x_{3} x_{4} ) (1-x_{3} ) (1-x_{4} ) }\\&
	\displaystyle +x_{2}^{-2}x_{4}^{-2}\cdot \frac{1}{48}\left(x_{1} \partial_{1} \right)^4\left(-x_{2} \partial_{2} +x_{3} \partial_{3} \right)^2\cdot\frac{1}{(1-x_{1} ) (1-x_{2} x_{3} x_{4} ) (1-x_{3} ) (1-x_{2} x_{3} ) }\\&
	\displaystyle +x_{2}^{-3}x_{4}^{-2}\cdot \frac{1}{720}\left(x_{1} \partial_{1} -x_{4} \partial_{4} \right)^6\cdot\frac{1}{(1-x_{1} ) (1-x_{1} x_{2} x_{3} x_{4} ) (1-x_{3} ) (1-x_{2} x_{3} ) }\\&
	\displaystyle +x_{2}^{-4}x_{3}^{-2}x_{4}^{-3}\cdot \frac{1}{120}\left(x_{1} \partial_{1} \right)^5\left(-x_{2} \partial_{2} +x_{3} \partial_{3} \right)\cdot\frac{1}{(1-x_{1} ) (1-x_{2} x_{3} ) (1-x_{3} ) (1-x_{4} ) }\\&
	\displaystyle -x_{2}^{-3}x_{4}^{-2}\cdot \frac{1}{720}\left(x_{1} \partial_{1} \right)^6\cdot\frac{1}{(1-x_{1} ) (1-x_{2} x_{3} ) (1-x_{3} ) (1-x_{3} x_{4} ) }\\&
	\displaystyle +\left(x_{2}^{-3}x_{4}^{-2}-x_{2}^{-3}x_{3}^{-1}x_{4}^{-3}\right)\cdot \frac{1}{120}\left(x_{1} \partial_{1} -x_{4} \partial_{4} \right)^5\left(-x_{2} \partial_{2} +x_{3} \partial_{3} \right)\cdot\frac{1}{(1-x_{1} ) (1-x_{1} x_{2} x_{3} x_{4} ) (1-x_{3} ) (1-x_{2} x_{3} ) }\\&
	\displaystyle -x_{2}^{-3}x_{4}^{-2}\cdot \frac{1}{120}\left(x_{1} \partial_{1} \right)^5\left(-x_{2} \partial_{2} +x_{3} \partial_{3} -x_{4} \partial_{4} \right)\cdot\frac{1}{(1-x_{1} ) (1-x_{2} x_{3} ) (1-x_{3} ) (1-x_{3} x_{4} ) }\\&
	\displaystyle +\left(-x_{2}^{-2}x_{4}^{-2}+x_{2}^{-3}x_{4}^{-2}\right)\cdot \frac{1}{48}\left(x_{1} \partial_{1} -x_{4} \partial_{4} \right)^4\left(-x_{2} \partial_{2} +x_{3} \partial_{3} \right)^2\cdot\frac{1}{(1-x_{1} ) (1-x_{1} x_{2} x_{3} x_{4} ) (1-x_{3} ) (1-x_{2} x_{3} ) }\\&
	\displaystyle -x_{2}^{-3}x_{4}^{-2}\cdot \frac{1}{48}\left(x_{1} \partial_{1} \right)^4\left(-x_{2} \partial_{2} +x_{3} \partial_{3} -x_{4} \partial_{4} \right)^2\cdot\frac{1}{(1-x_{1} ) (1-x_{2} x_{3} ) (1-x_{3} ) (1-x_{3} x_{4} ) }\\&
	\displaystyle -x_{2}^{-4}x_{3}^{-2}x_{4}^{-2}\cdot \frac{1}{24}\left(x_{1} \partial_{1} \right)^4\left(-x_{2} \partial_{2} +x_{3} \partial_{3} \right)\left(x_{4} \partial_{4} \right)\cdot\frac{1}{(1-x_{1} ) (1-x_{2} x_{3} ) (1-x_{3} ) (1-x_{4} ) }\\&
	\displaystyle +\left(\begin{array}{l}-x_{2}^{-2}x_{3} x_{4}^{-1}-x_{2}^{-3}x_{3}^{-2}x_{4}^{-2} \\ -x_{2}^{-3}x_{3}^{-2}x_{4}^{-3}-x_{2}^{-6}x_{3}^{-3}x_{4}^{-1}\end{array}\right)\begin{array}{l}\cdot \frac{1}{720}\left(\begin{array}{l}x_{1} \partial_{1} -x_{2} \partial_{2} +x_{3} \partial_{3}  \\ -x_{4} \partial_{4} \end{array}\right)^6\\
		~~\cdot\frac{1}{(1-x_{1} ) (1-x_{1} x_{2} x_{3} x_{4} ) (1-x_{2} x_{3} ) (1-x_{3} x_{4} ) }\end{array}\\&
	\displaystyle -x_{2}^{-6}x_{3}^{-4}x_{4}^{-2}\cdot \frac{1}{720}\left(x_{1} \partial_{1} -x_{4} \partial_{4} \right)^6\cdot\frac{1}{(1-x_{1} ) (1-x_{2} ) (1-x_{2} x_{3} ) (1-x_{1} x_{2} x_{3} x_{4} ) }\\&
	\displaystyle +x_{2}^{-6}x_{3}^{-4}x_{4}^{-2}\cdot \frac{1}{720}\left(x_{1} \partial_{1} \right)^6\cdot\frac{1}{(1-x_{1} ) (1-x_{2} ) (1-x_{2} x_{3} ) (1-x_{3} x_{4} ) }\\&
	\displaystyle +x_{2}^{-5}x_{3}^{-4}x_{4}^{-2}\cdot \frac{1}{120}\left(x_{1} \partial_{1} -x_{4} \partial_{4} \right)^5\left(x_{2} \partial_{2} -x_{3} \partial_{3} \right)\cdot\frac{1}{(1-x_{1} ) (1-x_{2} ) (1-x_{2} x_{3} ) (1-x_{1} x_{2} x_{3} x_{4} ) }\\&
	\displaystyle -x_{2}^{-5}x_{3}^{-4}x_{4}^{-2}\cdot \frac{1}{120}\left(x_{1} \partial_{1} \right)^5\left(x_{2} \partial_{2} -x_{3} \partial_{3} +x_{4} \partial_{4} \right)\cdot\frac{1}{(1-x_{1} ) (1-x_{2} ) (1-x_{2} x_{3} ) (1-x_{3} x_{4} ) }\\&
	\displaystyle +x_{2}^{-4}x_{3}^{-4}x_{4}^{-2}\cdot \frac{1}{48}\left(x_{1} \partial_{1} \right)^4\left(x_{2} \partial_{2} -x_{3} \partial_{3} +x_{4} \partial_{4} \right)^2\cdot\frac{1}{(1-x_{1} ) (1-x_{2} ) (1-x_{2} x_{3} ) (1-x_{3} x_{4} ) }\\&
	\displaystyle +x_{2}^{-3}x_{3}^{-4}x_{4}^{-2}\cdot \frac{1}{36}\left(x_{1} \partial_{1} \right)^3\left(x_{2} \partial_{2} -x_{3} \partial_{3} \right)^3\cdot\frac{1}{(1-x_{1} ) (1-x_{2} ) (1-x_{2} x_{3} x_{4} ) (1-x_{2} x_{3} ) }\\&
	\displaystyle -x_{2}^{-3}x_{3}^{-4}x_{4}^{-2}\cdot \frac{1}{36}\left(x_{1} \partial_{1} \right)^3\left(x_{2} \partial_{2} -x_{3} \partial_{3} +x_{4} \partial_{4} \right)^3\cdot\frac{1}{(1-x_{1} ) (1-x_{2} ) (1-x_{2} x_{3} ) (1-x_{3} x_{4} ) }\\&
	\displaystyle +\left(\begin{array}{l}x_{1} x_{2}^{-2}x_{3}^{-2}x_{4}^{-1}+x_{2}^{-3}x_{3}^{-2}x_{4}^{-1} \\ +x_{2}^{-5}x_{3}^{-2}x_{4} \end{array}\right)\begin{array}{l}\cdot \frac{1}{720}\left(x_{1} \partial_{1} \right)^6\\
		~~\cdot\frac{1}{(1-x_{1} ) (1-x_{3} x_{4} ) (1-x_{2} x_{3} x_{4} ) (1-x_{4} ) }\end{array}\\&
	\displaystyle +\left(\begin{array}{l}x_{1} x_{2}^{-2}x_{3}^{-2}x_{4}^{-1}+x_{2}^{-3}x_{3}^{-2}x_{4}^{-1} \\ +x_{2}^{-5}x_{3}^{-2}x_{4} \end{array}\right)\begin{array}{l}\cdot \frac{1}{120}\left(x_{1} \partial_{1} -x_{2} \partial_{2} \right)^5\left(-x_{3} \partial_{3} +x_{4} \partial_{4} \right)\\
		~~\cdot\frac{1}{(1-x_{1} ) (1-x_{3} x_{4} ) (1-x_{1} x_{2} x_{3} x_{4} ) (1-x_{4} ) }\end{array}\\&
	\displaystyle +\left(\begin{array}{l}-x_{1} x_{2}^{-2}x_{3}^{-2}x_{4}^{-1}-x_{2}^{-3}x_{3}^{-2}x_{4}^{-1} \\ -x_{2}^{-5}x_{3}^{-2}x_{4} \end{array}\right)\begin{array}{l}\cdot \frac{1}{120}\left(x_{1} \partial_{1} -x_{2} \partial_{2} \right)^5\left(x_{2} \partial_{2} -x_{3} \partial_{3} +x_{4} \partial_{4} \right)\\
		~~\cdot\frac{1}{(1-x_{1} ) (1-x_{3} x_{4} ) (1-x_{1} x_{2} x_{3} ) (1-x_{4} ) }\end{array}\\&
	\displaystyle +\left(x_{2}^{-3}x_{3}^{-3}x_{4}^{-1}-x_{2}^{-3}x_{3}^{-4}x_{4}^{-2}\right)\cdot \frac{1}{36}\left(x_{1} \partial_{1} -x_{3} \partial_{3} +x_{4} \partial_{4} \right)^3\left(x_{2} \partial_{2} -x_{3} \partial_{3} \right)^3\cdot\frac{1}{(1-x_{1} ) (1-x_{2} ) (1-x_{2} x_{3} x_{4} ) (1-x_{1} x_{2} x_{3} ) }\\&
	\displaystyle +\left(-x_{2}^{-3}x_{3}^{-3}x_{4}^{-1}+x_{2}^{-3}x_{3}^{-4}x_{4}^{-2}\right)\cdot \frac{1}{36}\left(x_{1} \partial_{1} -x_{3} \partial_{3} +x_{4} \partial_{4} \right)^3\left(x_{2} \partial_{2} -x_{3} \partial_{3} +x_{4} \partial_{4} \right)^3\cdot\frac{1}{(1-x_{1} ) (1-x_{2} ) (1-x_{1} x_{2} x_{3} ) (1-x_{3} x_{4} ) }\\&
	\displaystyle +\left(\begin{array}{l}-x_{1} x_{2}^{-2}x_{3}^{-2}x_{4}^{-1}-x_{2}^{-3}x_{3}^{-2}x_{4}^{-1} \\ -x_{2}^{-5}x_{3}^{-2}x_{4} -3x_{2}^{-6}x_{3}^{-3}-x_{2}^{-4}x_{3}^{-3}x_{4}^{-3}\end{array}\right)\begin{array}{l}\cdot \frac{1}{720}\left(x_{1} \partial_{1} -x_{2} \partial_{2} \right)^6\\
		~~\cdot\frac{1}{(1-x_{1} ) (1-x_{1} x_{2} x_{3} x_{4} ) (1-x_{3} x_{4} ) (1-x_{4} ) }\end{array}\\&
	\displaystyle -3x_{2}^{-5}x_{3}^{-4}x_{4}^{-1}\cdot \frac{1}{120}\left(x_{1} \partial_{1} \right)^5\left(x_{2} \partial_{2} \right)\cdot\frac{1}{(1-x_{1} ) (1-x_{2} ) (1-x_{3} x_{4} ) (1-x_{4} ) }\\&
	\displaystyle +2x_{2}^{-4}x_{3}^{-4}x_{4}^{-1}\cdot \frac{1}{48}\left(x_{1} \partial_{1} \right)^4\left(x_{2} \partial_{2} \right)^2\cdot\frac{1}{(1-x_{1} ) (1-x_{2} ) (1-x_{3} x_{4} ) (1-x_{4} ) }\\&
	\displaystyle +3x_{2}^{-6}x_{3}^{-4}x_{4}^{-1}\cdot \frac{1}{720}\left(x_{1} \partial_{1} \right)^6\cdot\frac{1}{(1-x_{1} ) (1-x_{2} ) (1-x_{3} x_{4} ) (1-x_{4} ) }\\&
	\displaystyle -x_{2}^{-3}x_{3}^{-4}x_{4}^{-1}\cdot \frac{1}{36}\left(x_{1} \partial_{1} \right)^3\left(x_{2} \partial_{2} \right)^3\cdot\frac{1}{(1-x_{1} ) (1-x_{2} ) (1-x_{3} x_{4} ) (1-x_{4} ) }\\&
	\displaystyle +\left(-x_{2}^{-2}x_{3}^{-3}+x_{2}^{-3}x_{3}^{-3}x_{4} \right)\cdot \frac{1}{12}\left(x_{1} \partial_{1} -x_{3} \partial_{3} \right)^3\left(x_{2} \partial_{2} -x_{3} \partial_{3} \right)^2\left(-x_{3} \partial_{3} +x_{4} \partial_{4} \right)\cdot\frac{1}{(1-x_{1} ) (1-x_{2} ) (1-x_{1} x_{2} x_{3} x_{4} ) (1-x_{4} ) }\\&
	\displaystyle +\left(-x_{2}^{-2}x_{3}^{-3}+x_{2}^{-3}x_{3}^{-3}x_{4} -x_{2}^{-4}x_{3}^{-4}x_{4}^{-1}\right)\cdot \frac{1}{48}\left(x_{1} \partial_{1} \right)^4\left(x_{2} \partial_{2} -x_{3} \partial_{3} \right)^2\cdot\frac{1}{(1-x_{1} ) (1-x_{2} ) (1-x_{2} x_{3} x_{4} ) (1-x_{4} ) }\\&
	\displaystyle +\left(x_{2}^{-2}x_{3}^{-3}-x_{2}^{-3}x_{3}^{-3}x_{4} \right)\cdot \frac{1}{12}\left(x_{1} \partial_{1} -x_{3} \partial_{3} \right)^3\left(x_{2} \partial_{2} -x_{3} \partial_{3} \right)^2\left(x_{4} \partial_{4} \right)\cdot\frac{1}{(1-x_{1} ) (1-x_{2} ) (1-x_{1} x_{2} x_{3} ) (1-x_{4} ) }\\&
	\displaystyle +\left(x_{2}^{-2}x_{3}^{-3}x_{4}^{-2}-x_{2}^{-4}x_{3}^{-4}x_{4}^{-2}\right)\cdot \frac{1}{48}\left(x_{1} \partial_{1} -x_{4} \partial_{4} \right)^4\left(x_{2} \partial_{2} -x_{3} \partial_{3} \right)^2\cdot\frac{1}{(1-x_{1} ) (1-x_{2} ) (1-x_{2} x_{3} ) (1-x_{1} x_{2} x_{3} x_{4} ) }\\&
	\displaystyle +\left(\begin{array}{l}x_{2}^{-2}x_{3}^{-3}-x_{2}^{-3}x_{3}^{-3}x_{4} +x_{2}^{-3}x_{3}^{-4}x_{4}^{-2} \\ -x_{2}^{-4}x_{3}^{-4}x_{4}^{-1}\end{array}\right)\begin{array}{l}\cdot \frac{1}{48}\left(x_{1} \partial_{1} -x_{3} \partial_{3} \right)^4\left(x_{2} \partial_{2} -x_{3} \partial_{3} \right)^2\\
		~~\cdot\frac{1}{(1-x_{1} ) (1-x_{2} ) (1-x_{1} x_{2} x_{3} x_{4} ) (1-x_{4} ) }\end{array}\\&
	\displaystyle -x_{2}^{-3}x_{3}^{-4}x_{4}^{-2}\cdot \frac{1}{48}\left(x_{1} \partial_{1} \right)^4\left(x_{2} \partial_{2} -x_{3} \partial_{3} \right)^2\cdot\frac{1}{(1-x_{1} ) (1-x_{2} ) (1-x_{2} x_{3} ) (1-x_{4} ) }\\&
	\displaystyle -x_{2}^{-3}x_{3}^{-4}x_{4}^{-1}\cdot \frac{1}{12}\left(x_{1} \partial_{1} \right)^3\left(x_{2} \partial_{2} -x_{3} \partial_{3} \right)^2\left(-x_{3} \partial_{3} +x_{4} \partial_{4} \right)\cdot\frac{1}{(1-x_{1} ) (1-x_{2} ) (1-x_{2} x_{3} x_{4} ) (1-x_{4} ) }\\&
	\displaystyle +x_{2}^{-3}x_{3}^{-4}x_{4}^{-1}\cdot \frac{1}{12}\left(x_{1} \partial_{1} \right)^3\left(x_{2} \partial_{2} -x_{3} \partial_{3} \right)^2\left(x_{4} \partial_{4} \right)\cdot\frac{1}{(1-x_{1} ) (1-x_{2} ) (1-x_{2} x_{3} ) (1-x_{4} ) }\\&
	\displaystyle +\left(\begin{array}{l}3x_{2}^{-6}x_{3}^{-3}x_{4}^{-1}-x_{2}^{-6}x_{3}^{-4}x_{4}^{-2} \\ -3x_{2}^{-6}x_{3}^{-5}x_{4}^{-3}\end{array}\right)\begin{array}{l}\cdot \frac{1}{720}\left(x_{1} \partial_{1} -x_{3} \partial_{3} +x_{4} \partial_{4} \right)^6\\
		~~\cdot\frac{1}{(1-x_{1} ) (1-x_{2} ) (1-x_{3} x_{4} ) (1-x_{1} x_{2} x_{3} ) }\end{array}\\&
	\displaystyle +\left(\begin{array}{l}-3x_{2}^{-5}x_{3}^{-3}x_{4}^{-1}+x_{2}^{-5}x_{3}^{-4}x_{4}^{-2} \\ +3x_{2}^{-5}x_{3}^{-5}x_{4}^{-3}\end{array}\right)\begin{array}{l}\cdot \frac{1}{120}\left(x_{1} \partial_{1} -x_{3} \partial_{3} +x_{4} \partial_{4} \right)^5\left(x_{2} \partial_{2} -x_{3} \partial_{3} +x_{4} \partial_{4} \right)\\
		~~\cdot\frac{1}{(1-x_{1} ) (1-x_{2} ) (1-x_{3} x_{4} ) (1-x_{1} x_{2} x_{3} ) }\end{array}\\&
	\displaystyle +\left(\begin{array}{l}-x_{2}^{-2}x_{3}^{-3}x_{4}^{-1}+x_{2}^{-3}x_{3}^{-3} \\ -x_{2}^{-4}x_{3}^{-3}x_{4}^{-1}+x_{2}^{-4}x_{3}^{-5}x_{4}^{-3}\end{array}\right)\begin{array}{l}\cdot \frac{1}{48}\left(x_{1} \partial_{1} -x_{3} \partial_{3} +x_{4} \partial_{4} \right)^4\left(x_{2} \partial_{2} -x_{3} \partial_{3} \right)^2\\
		~~\cdot\frac{1}{(1-x_{1} ) (1-x_{2} ) (1-x_{2} x_{3} x_{4} ) (1-x_{1} x_{2} x_{3} ) }\end{array}\\&
	\displaystyle +\left(\begin{array}{l}2x_{2}^{-4}x_{3}^{-3}x_{4}^{-1}-x_{2}^{-4}x_{3}^{-4}x_{4}^{-2} \\ -x_{2}^{-4}x_{3}^{-5}x_{4}^{-3}\end{array}\right)\begin{array}{l}\cdot \frac{1}{48}\left(x_{1} \partial_{1} -x_{3} \partial_{3} +x_{4} \partial_{4} \right)^4\left(x_{2} \partial_{2} -x_{3} \partial_{3} +x_{4} \partial_{4} \right)^2\\
		~~\cdot\frac{1}{(1-x_{1} ) (1-x_{2} ) (1-x_{3} x_{4} ) (1-x_{1} x_{2} x_{3} ) }\end{array}\\&
	\displaystyle +3x_{2}^{-6}x_{3}^{-5}x_{4}^{-3}\cdot \frac{1}{720}\left(x_{1} \partial_{1} \right)^6\cdot\frac{1}{(1-x_{1} ) (1-x_{2} ) (1-x_{3} ) (1-x_{3} x_{4} ) }\\&
	\displaystyle +x_{2}^{-4}x_{3}^{-5}x_{4}^{-3}\cdot \frac{1}{48}\left(x_{1} \partial_{1} \right)^4\left(x_{2} \partial_{2} \right)^2\cdot\frac{1}{(1-x_{1} ) (1-x_{2} ) (1-x_{3} ) (1-x_{3} x_{4} ) }\\&
	\displaystyle +\left(\begin{array}{l}x_{1} x_{2}^{-2}x_{3}^{-2}x_{4}^{-2}+x_{2}^{-5}x_{3}^{-2} \\ -x_{2}^{-3}x_{3}^{-3}x_{4}^{-2}-x_{2}^{-3}x_{3}^{-4}x_{4}^{-2} \\ -x_{2}^{-4}x_{3}^{-4}x_{4}^{-2}-2x_{2}^{-5}x_{3}^{-4}x_{4}^{-2}\end{array}\right)\begin{array}{l}\cdot \frac{1}{720}\left(x_{1} \partial_{1} -x_{3} \partial_{3} +x_{4} \partial_{4} \right)^6\\
		~~\cdot\frac{1}{(1-x_{1} ) (1-x_{3} x_{4} ) (1-x_{2} x_{3} x_{4} ) (1-x_{1} x_{2} x_{3} ) }\end{array}\\&
	\displaystyle +\left(\begin{array}{l}x_{2}^{-3}x_{3} x_{4}^{-1}-x_{2}^{-3}x_{3}^{-2}x_{4}^{-1} \\ -x_{2}^{-3}x_{3}^{-3}x_{4}^{-1}-x_{2}^{-4}x_{3}^{-3}x_{4}^{-1} \\ -2x_{2}^{-5}x_{3}^{-3}x_{4}^{-1}-3x_{2}^{-6}x_{3}^{-4}x_{4}^{-2}\end{array}\right)\begin{array}{l}\cdot \frac{1}{720}\left(x_{1} \partial_{1} -x_{2} \partial_{2} \right)^6\\
		~~\cdot\frac{1}{(1-x_{1} ) (1-x_{1} x_{2} x_{3} x_{4} ) (1-x_{3} ) (1-x_{3} x_{4} ) }\end{array}\\&
	\displaystyle +\left(\begin{array}{l}x_{2}^{-3}x_{3}^{-2}x_{4}^{-1}+x_{2}^{-3}x_{3}^{-3}x_{4}^{-1} \\ +x_{2}^{-4}x_{3}^{-3}x_{4}^{-1}+2x_{2}^{-5}x_{3}^{-3}x_{4}^{-1}\end{array}\right)\begin{array}{l}\cdot \frac{1}{720}\left(x_{1} \partial_{1} \right)^6\\
		~~\cdot\frac{1}{(1-x_{1} ) (1-x_{3} x_{4} ) (1-x_{3} ) (1-x_{2} x_{3} x_{4} ) }\end{array}\\&
	\displaystyle +\left(\begin{array}{l}x_{1} x_{2}^{-2}x_{4}^{-1}+x_{1} x_{2}^{-2}x_{3}^{-1}x_{4}^{-1} \\ +x_{1} x_{2}^{-3}x_{3}^{-1}x_{4}^{-1}+2x_{1} x_{2}^{-4}x_{3}^{-1}x_{4}^{-1}\end{array}\right)\begin{array}{l}\cdot \frac{1}{120}\left(x_{1} \partial_{1} -x_{2} \partial_{2} \right)^5\left(-x_{2} \partial_{2} +x_{3} \partial_{3} -x_{4} \partial_{4} \right)\\
		~~\cdot\frac{1}{(1-x_{1} ) (1-x_{3} x_{4} ) (1-x_{3} ) (1-x_{1} x_{2} x_{3} ) }\end{array}\\&
	\displaystyle +\left(\begin{array}{l}-x_{2}^{-3}x_{3}^{-1}x_{4}^{-1}-x_{2}^{-3}x_{3}^{-2}x_{4}^{-1} \\ -x_{2}^{-4}x_{3}^{-2}x_{4}^{-1}\end{array}\right)\begin{array}{l}\cdot \frac{1}{120}\left(x_{1} \partial_{1} \right)^5\left(x_{3} \partial_{3} -x_{4} \partial_{4} \right)\\
		~~\cdot\frac{1}{(1-x_{1} ) (1-x_{3} x_{4} ) (1-x_{3} ) (1-x_{2} x_{3} x_{4} ) }\end{array}\\&
	\displaystyle +\left(\begin{array}{l}-x_{1} x_{2}^{-2}x_{3} x_{4}^{-1}-x_{1} x_{2}^{-2}x_{4}^{-1} \\ -x_{1} x_{2}^{-3}x_{4}^{-1}\end{array}\right)\begin{array}{l}\cdot \frac{1}{48}\left(x_{1} \partial_{1} -x_{2} \partial_{2} \right)^4\left(-x_{2} \partial_{2} +x_{3} \partial_{3} -x_{4} \partial_{4} \right)^2\\
		~~\cdot\frac{1}{(1-x_{1} ) (1-x_{3} x_{4} ) (1-x_{3} ) (1-x_{1} x_{2} x_{3} ) }\end{array}\\&
	\displaystyle -x_{1} x_{2}^{-2}x_{4}^{-1}\cdot \frac{1}{36}\left(x_{1} \partial_{1} -x_{2} \partial_{2} \right)^3\left(-x_{2} \partial_{2} +x_{3} \partial_{3} \right)^3\cdot\frac{1}{(1-x_{1} ) (1-x_{1} x_{2} x_{3} x_{4} ) (1-x_{3} ) (1-x_{1} x_{2} x_{3} ) }\\&
	\displaystyle -x_{1} x_{2}^{-2}x_{3} \cdot \frac{1}{36}\left(x_{1} \partial_{1} -x_{2} \partial_{2} \right)^3\left(x_{3} \partial_{3} -x_{4} \partial_{4} \right)^3\cdot\frac{1}{(1-x_{1} ) (1-x_{1} x_{2} x_{3} x_{4} ) (1-x_{3} ) (1-x_{3} x_{4} ) }\\&
	\displaystyle +x_{1} x_{2}^{-2}x_{3} x_{4}^{-1}\cdot \frac{1}{36}\left(x_{1} \partial_{1} -x_{2} \partial_{2} \right)^3\left(-x_{2} \partial_{2} +x_{3} \partial_{3} -x_{4} \partial_{4} \right)^3\cdot\frac{1}{(1-x_{1} ) (1-x_{3} x_{4} ) (1-x_{3} ) (1-x_{1} x_{2} x_{3} ) }\\&
	\displaystyle +\left(\begin{array}{l}x_{1} x_{2}^{-2}x_{3} +x_{1} x_{2}^{-2}+x_{1} x_{2}^{-3} \\ +x_{2}^{-3}x_{3} x_{4}^{-1}-x_{2}^{-3}x_{3}^{-1}x_{4}^{-1}\end{array}\right)\begin{array}{l}\cdot \frac{1}{48}\left(x_{1} \partial_{1} -x_{2} \partial_{2} \right)^4\left(x_{3} \partial_{3} -x_{4} \partial_{4} \right)^2\\
		~~\cdot\frac{1}{(1-x_{1} ) (1-x_{1} x_{2} x_{3} x_{4} ) (1-x_{3} ) (1-x_{3} x_{4} ) }\end{array}\\&
	\displaystyle +x_{2}^{-3}x_{3}^{-1}x_{4}^{-1}\cdot \frac{1}{48}\left(x_{1} \partial_{1} \right)^4\left(x_{3} \partial_{3} -x_{4} \partial_{4} \right)^2\cdot\frac{1}{(1-x_{1} ) (1-x_{2} x_{3} x_{4} ) (1-x_{3} ) (1-x_{3} x_{4} ) }\\&
	\displaystyle -2x_{2}^{-5}x_{3}^{-4}x_{4}^{-3}\cdot \frac{1}{120}\left(x_{1} \partial_{1} \right)^5\left(x_{3} \partial_{3} -x_{4} \partial_{4} \right)\cdot\frac{1}{(1-x_{1} ) (1-x_{2} ) (1-x_{3} ) (1-x_{3} x_{4} ) }\\&
	\displaystyle +\left(\begin{array}{l}-x_{1} x_{2}^{-2}x_{3}^{-2}x_{4}^{-2}-x_{2}^{-5}x_{3}^{-2} \\ +x_{2}^{-3}x_{3}^{-3}x_{4}^{-2}-3x_{2}^{-6}x_{3}^{-2} \\ +x_{2}^{-3}x_{3}^{-4}x_{4}^{-2}+x_{2}^{-4}x_{3}^{-4}x_{4}^{-2} \\ +x_{2}^{-6}x_{3}^{-3}x_{4}^{-1}+2x_{2}^{-5}x_{3}^{-4}x_{4}^{-2} \\ +3x_{2}^{-6}x_{3}^{-4}x_{4}^{-2}\end{array}\right)\begin{array}{l}\cdot \frac{1}{720}\left(x_{1} \partial_{1} -x_{2} \partial_{2} \right)^6\\
		~~\cdot\frac{1}{(1-x_{1} ) (1-x_{1} x_{2} x_{3} x_{4} ) (1-x_{3} x_{4} ) (1-x_{1} x_{2} x_{3} ) }\end{array}\\&
	\displaystyle -x_{1} x_{2}^{-2}x_{3}^{-1}x_{4}^{-2}\cdot \frac{1}{12}\left(x_{1} \partial_{1} -x_{2} \partial_{2} \right)^3\left(x_{3} \partial_{3} -x_{4} \partial_{4} \right)^2\left(x_{2} \partial_{2} -x_{4} \partial_{4} \right)\cdot\frac{1}{(1-x_{1} ) (1-x_{1} x_{2} x_{3} x_{4} ) (1-x_{3} ) (1-x_{1} x_{2} ) }\\&
	\displaystyle +x_{1} x_{2}^{-2}x_{3}^{-1}x_{4}^{-2}\cdot \frac{1}{12}\left(x_{1} \partial_{1} -x_{2} \partial_{2} \right)^3\left(x_{3} \partial_{3} -x_{4} \partial_{4} \right)^2\left(x_{2} \partial_{2} \right)\cdot\frac{1}{(1-x_{1} ) (1-x_{1} x_{2} ) (1-x_{3} ) (1-x_{3} x_{4} ) }\\&
	\displaystyle -x_{2}^{-3}x_{3}^{-4}x_{4}^{-3}\cdot \frac{1}{12}\left(x_{1} \partial_{1} \right)^3\left(x_{2} \partial_{2} \right)^2\left(x_{3} \partial_{3} -x_{4} \partial_{4} \right)\cdot\frac{1}{(1-x_{1} ) (1-x_{2} ) (1-x_{3} ) (1-x_{3} x_{4} ) }\\&
	\displaystyle +\left(\begin{array}{l}-x_{2}^{-5}x_{3}^{-3}-3x_{2}^{-6}x_{3}^{-3}x_{4}^{-1} \\ +2x_{2}^{-5}x_{3}^{-4}x_{4}^{-2}+x_{2}^{-6}x_{3}^{-4}x_{4}^{-2} \\ +3x_{2}^{-6}x_{3}^{-5}x_{4}^{-3}\end{array}\right)\begin{array}{l}\cdot \frac{1}{720}\left(x_{1} \partial_{1} -x_{3} \partial_{3} \right)^6\\
		~~\cdot\frac{1}{(1-x_{1} ) (1-x_{2} ) (1-x_{1} x_{2} x_{3} ) (1-x_{1} x_{2} x_{3} x_{4} ) }\end{array}\\&
	\displaystyle +\left(x_{2}^{-5}x_{3}^{-3}-2x_{2}^{-5}x_{3}^{-4}x_{4}^{-2}\right)\cdot \frac{1}{720}\left(x_{1} \partial_{1} -x_{3} \partial_{3} +x_{4} \partial_{4} \right)^6\cdot\frac{1}{(1-x_{1} ) (1-x_{2} ) (1-x_{1} x_{2} x_{3} ) (1-x_{2} x_{3} x_{4} ) }\\&
	\displaystyle +\left(x_{2}^{-5}x_{3}^{-3}x_{4} -2x_{2}^{-5}x_{3}^{-4}x_{4}^{-1}\right)\cdot \frac{1}{720}\left(x_{1} \partial_{1} \right)^6\cdot\frac{1}{(1-x_{1} ) (1-x_{2} ) (1-x_{2} x_{3} x_{4} ) (1-x_{4} ) }\\&
	\displaystyle +\left(\begin{array}{l}-x_{2}^{-3}x_{3}^{-3}x_{4}^{-2}-x_{2}^{-3}x_{3}^{-4}x_{4}^{-2} \\ -x_{2}^{-4}x_{3}^{-4}x_{4}^{-2}+2x_{2}^{-5}x_{3}^{-3}x_{4}^{-2} \\ -2x_{2}^{-5}x_{3}^{-4}x_{4}^{-2}\end{array}\right)\begin{array}{l}\cdot \frac{1}{720}\left(x_{1} \partial_{1} -x_{2} \partial_{2} \right)^6\\
		~~\cdot\frac{1}{(1-x_{1} ) (1-x_{1} x_{2} x_{3} ) (1-x_{3} ) (1-x_{1} x_{2} x_{3} x_{4} ) }\end{array}\\&
	\displaystyle +\left(\begin{array}{l}x_{2}^{-3}x_{3}^{-3}x_{4}^{-2}+x_{2}^{-3}x_{3}^{-4}x_{4}^{-2} \\ +x_{2}^{-4}x_{3}^{-4}x_{4}^{-2}-2x_{2}^{-5}x_{3}^{-3}x_{4}^{-2} \\ +2x_{2}^{-5}x_{3}^{-4}x_{4}^{-2}\end{array}\right)\begin{array}{l}\cdot \frac{1}{720}\left(x_{1} \partial_{1} -x_{2} \partial_{2} +x_{4} \partial_{4} \right)^6\\
		~~\cdot\frac{1}{(1-x_{1} ) (1-x_{1} x_{2} x_{3} ) (1-x_{3} ) (1-x_{2} x_{3} x_{4} ) }\end{array}\\&
	\displaystyle +\left(x_{2}^{-5}x_{3}^{-3}x_{4} -2x_{2}^{-5}x_{3}^{-4}x_{4}^{-1}\right)\cdot \frac{1}{120}\left(x_{1} \partial_{1} -x_{3} \partial_{3} \right)^5\left(-x_{3} \partial_{3} +x_{4} \partial_{4} \right)\cdot\frac{1}{(1-x_{1} ) (1-x_{2} ) (1-x_{1} x_{2} x_{3} x_{4} ) (1-x_{4} ) }\\&
	\displaystyle -2x_{2}^{-5}x_{3}^{-3}x_{4}^{-1}\cdot \frac{1}{120}\left(x_{1} \partial_{1} -x_{2} \partial_{2} \right)^5\left(-x_{2} \partial_{2} +x_{4} \partial_{4} \right)\cdot\frac{1}{(1-x_{1} ) (1-x_{1} x_{2} x_{3} x_{4} ) (1-x_{3} ) (1-x_{4} ) }\\&
	\displaystyle +2x_{2}^{-5}x_{3}^{-3}x_{4}^{-1}\cdot \frac{1}{120}\left(x_{1} \partial_{1} -x_{2} \partial_{2} \right)^5\left(x_{4} \partial_{4} \right)\cdot\frac{1}{(1-x_{1} ) (1-x_{1} x_{2} x_{3} ) (1-x_{3} ) (1-x_{4} ) }\\&
	\displaystyle +\left(-x_{2}^{-5}x_{3}^{-3}x_{4} +2x_{2}^{-5}x_{3}^{-4}x_{4}^{-1}\right)\cdot \frac{1}{120}\left(x_{1} \partial_{1} -x_{3} \partial_{3} \right)^5\left(x_{4} \partial_{4} \right)\cdot\frac{1}{(1-x_{1} ) (1-x_{2} ) (1-x_{1} x_{2} x_{3} ) (1-x_{4} ) }\\&
	\displaystyle +\left(\begin{array}{l}x_{1} x_{2}^{-2}x_{4}^{-1}+x_{1} x_{2}^{-2}x_{3}^{-1}x_{4}^{-1} \\ +x_{1} x_{2}^{-3}x_{3}^{-1}x_{4}^{-1}+x_{2}^{-3}x_{3}^{-1}x_{4}^{-2} \\ -x_{2}^{-3}x_{3}^{-2}x_{4}^{-2}\end{array}\right)\begin{array}{l}\cdot \frac{1}{48}\left(x_{1} \partial_{1} -x_{2} \partial_{2} \right)^4\left(-x_{2} \partial_{2} +x_{3} \partial_{3} \right)^2\\
		~~\cdot\frac{1}{(1-x_{1} ) (1-x_{1} x_{2} x_{3} x_{4} ) (1-x_{3} ) (1-x_{1} x_{2} x_{3} ) }\end{array}\\&
	\displaystyle +\left(-x_{2}^{-3}x_{3}^{-1}x_{4}^{-2}+x_{2}^{-3}x_{3}^{-2}x_{4}^{-2}\right)\cdot \frac{1}{48}\left(x_{1} \partial_{1} -x_{2} \partial_{2} +x_{4} \partial_{4} \right)^4\left(-x_{2} \partial_{2} +x_{3} \partial_{3} \right)^2\cdot\frac{1}{(1-x_{1} ) (1-x_{1} x_{2} x_{3} ) (1-x_{3} ) (1-x_{2} x_{3} x_{4} ) }\\&
	\displaystyle +\left(\begin{array}{l}x_{2}^{-2}x_{3}^{-3}x_{4}^{-1}-x_{2}^{-3}x_{3}^{-3} \\ -x_{2}^{-4}x_{3}^{-3}x_{4}^{-1}+x_{2}^{-4}x_{3}^{-4}x_{4}^{-2}\end{array}\right)\begin{array}{l}\cdot \frac{1}{48}\left(x_{1} \partial_{1} -x_{3} \partial_{3} \right)^4\left(x_{2} \partial_{2} -x_{3} \partial_{3} \right)^2\\
		~~\cdot\frac{1}{(1-x_{1} ) (1-x_{2} ) (1-x_{1} x_{2} x_{3} x_{4} ) (1-x_{1} x_{2} x_{3} ) }\end{array}\\&
	\displaystyle -x_{2}^{-3}x_{3}^{-1}x_{4}^{-1}\cdot \frac{1}{12}\left(x_{1} \partial_{1} -x_{2} \partial_{2} \right)^3\left(-x_{2} \partial_{2} +x_{3} \partial_{3} \right)^2\left(-x_{2} \partial_{2} +x_{4} \partial_{4} \right)\cdot\frac{1}{(1-x_{1} ) (1-x_{1} x_{2} x_{3} x_{4} ) (1-x_{3} ) (1-x_{4} ) }\\&
	\displaystyle +x_{2}^{-3}x_{3}^{-1}x_{4}^{-1}\cdot \frac{1}{12}\left(x_{1} \partial_{1} -x_{2} \partial_{2} \right)^3\left(-x_{2} \partial_{2} +x_{3} \partial_{3} \right)^2\left(x_{4} \partial_{4} \right)\cdot\frac{1}{(1-x_{1} ) (1-x_{1} x_{2} x_{3} ) (1-x_{3} ) (1-x_{4} ) }\\&
	\displaystyle +\left(\begin{array}{l}-x_{1} x_{2}^{-2}x_{3}^{-1}x_{4}^{-1}-x_{1} x_{2}^{-2}x_{3}^{-2}x_{4}^{-1} \\ -x_{1} x_{2}^{-3}x_{3}^{-2}x_{4}^{-1}-2x_{1} x_{2}^{-4}x_{3}^{-2}x_{4}^{-1} \\ -x_{2}^{-3}x_{3}^{-1}x_{4}^{-2}+x_{2}^{-3}x_{3}^{-2}x_{4}^{-2} \\ +x_{2}^{-3}x_{3}^{-3}x_{4}^{-2}-x_{2}^{-4}x_{3}^{-2}x_{4}^{-2} \\ +x_{2}^{-4}x_{3}^{-3}x_{4}^{-2}\end{array}\right)\begin{array}{l}\cdot \frac{1}{120}\left(x_{1} \partial_{1} -x_{2} \partial_{2} \right)^5\left(-x_{2} \partial_{2} +x_{3} \partial_{3} \right)\\
		~~\cdot\frac{1}{(1-x_{1} ) (1-x_{1} x_{2} x_{3} ) (1-x_{3} ) (1-x_{1} x_{2} x_{3} x_{4} ) }\end{array}\\&
	\displaystyle +\left(\begin{array}{l}x_{2}^{-3}x_{3}^{-1}x_{4}^{-2}-x_{2}^{-3}x_{3}^{-2}x_{4}^{-2} \\ -x_{2}^{-3}x_{3}^{-3}x_{4}^{-2}+x_{2}^{-4}x_{3}^{-2}x_{4}^{-2} \\ -x_{2}^{-4}x_{3}^{-3}x_{4}^{-2}\end{array}\right)\begin{array}{l}\cdot \frac{1}{120}\left(x_{1} \partial_{1} -x_{2} \partial_{2} +x_{4} \partial_{4} \right)^5\left(-x_{2} \partial_{2} +x_{3} \partial_{3} \right)\\
		~~\cdot\frac{1}{(1-x_{1} ) (1-x_{1} x_{2} x_{3} ) (1-x_{3} ) (1-x_{2} x_{3} x_{4} ) }\end{array}\\&
	\displaystyle +\left(\begin{array}{l}x_{2}^{-4}x_{3}^{-3}+2x_{2}^{-5}x_{3}^{-3}x_{4}^{-1} \\ -x_{2}^{-4}x_{3}^{-4}x_{4}^{-2}-x_{2}^{-5}x_{3}^{-4}x_{4}^{-2} \\ -x_{2}^{-5}x_{3}^{-5}x_{4}^{-3}\end{array}\right)\begin{array}{l}\cdot \frac{1}{120}\left(x_{1} \partial_{1} -x_{3} \partial_{3} \right)^5\left(x_{2} \partial_{2} -x_{3} \partial_{3} \right)\\
		~~\cdot\frac{1}{(1-x_{1} ) (1-x_{2} ) (1-x_{1} x_{2} x_{3} ) (1-x_{1} x_{2} x_{3} x_{4} ) }\end{array}\\&
	\displaystyle +\left(\begin{array}{l}-x_{2}^{-4}x_{3}^{-3}+x_{2}^{-5}x_{3}^{-3}x_{4}^{-1} \\ +x_{2}^{-4}x_{3}^{-4}x_{4}^{-2}-2x_{2}^{-5}x_{3}^{-5}x_{4}^{-3}\end{array}\right)\begin{array}{l}\cdot \frac{1}{120}\left(x_{1} \partial_{1} -x_{3} \partial_{3} +x_{4} \partial_{4} \right)^5\left(x_{2} \partial_{2} -x_{3} \partial_{3} \right)\\
		~~\cdot\frac{1}{(1-x_{1} ) (1-x_{2} ) (1-x_{1} x_{2} x_{3} ) (1-x_{2} x_{3} x_{4} ) }\end{array}\\&
	\displaystyle +\left(\begin{array}{l}-x_{2}^{-4}x_{3}^{-3}x_{4} +x_{2}^{-4}x_{3}^{-4}x_{4}^{-1} \\ +x_{2}^{-5}x_{3}^{-4}x_{4}^{-1}-2x_{2}^{-5}x_{3}^{-5}x_{4}^{-2}\end{array}\right)\begin{array}{l}\cdot \frac{1}{120}\left(x_{1} \partial_{1} \right)^5\left(x_{2} \partial_{2} -x_{3} \partial_{3} \right)\\
		~~\cdot\frac{1}{(1-x_{1} ) (1-x_{2} ) (1-x_{2} x_{3} x_{4} ) (1-x_{4} ) }\end{array}\\&
	\displaystyle +\left(x_{2}^{-3}x_{3}^{-1}x_{4}^{-1}+x_{2}^{-4}x_{3}^{-2}x_{4}^{-1}\right)\cdot \frac{1}{24}\left(x_{1} \partial_{1} -x_{2} \partial_{2} \right)^4\left(-x_{2} \partial_{2} +x_{3} \partial_{3} \right)\left(-x_{2} \partial_{2} +x_{4} \partial_{4} \right)\cdot\frac{1}{(1-x_{1} ) (1-x_{1} x_{2} x_{3} x_{4} ) (1-x_{3} ) (1-x_{4} ) }\\&
	\displaystyle +\left(-x_{2}^{-4}x_{3}^{-3}x_{4} +x_{2}^{-4}x_{3}^{-4}x_{4}^{-1}\right)\cdot \frac{1}{24}\left(x_{1} \partial_{1} -x_{3} \partial_{3} \right)^4\left(x_{2} \partial_{2} -x_{3} \partial_{3} \right)\left(-x_{3} \partial_{3} +x_{4} \partial_{4} \right)\cdot\frac{1}{(1-x_{1} ) (1-x_{2} ) (1-x_{1} x_{2} x_{3} x_{4} ) (1-x_{4} ) }\\&
	\displaystyle +\left(x_{2}^{-4}x_{3}^{-3}x_{4} -x_{2}^{-4}x_{3}^{-4}x_{4}^{-1}\right)\cdot \frac{1}{24}\left(x_{1} \partial_{1} -x_{3} \partial_{3} \right)^4\left(x_{2} \partial_{2} -x_{3} \partial_{3} \right)\left(x_{4} \partial_{4} \right)\cdot\frac{1}{(1-x_{1} ) (1-x_{2} ) (1-x_{1} x_{2} x_{3} ) (1-x_{4} ) }\\&
	\displaystyle +\left(-x_{2}^{-3}x_{3}^{-1}x_{4}^{-1}-x_{2}^{-4}x_{3}^{-2}x_{4}^{-1}\right)\cdot \frac{1}{24}\left(x_{1} \partial_{1} -x_{2} \partial_{2} \right)^4\left(-x_{2} \partial_{2} +x_{3} \partial_{3} \right)\left(x_{4} \partial_{4} \right)\cdot\frac{1}{(1-x_{1} ) (1-x_{1} x_{2} x_{3} ) (1-x_{3} ) (1-x_{4} ) }\\&
	\displaystyle -3x_{2}^{-5}x_{3}^{-5}x_{4}^{-3}\cdot \frac{1}{120}\left(x_{1} \partial_{1} \right)^5\left(x_{2} \partial_{2} \right)\cdot\frac{1}{(1-x_{1} ) (1-x_{2} ) (1-x_{3} ) (1-x_{3} x_{4} ) }\\&
	\displaystyle +\left(\begin{array}{l}x_{1} x_{2}^{-2}x_{4}^{-2}+x_{1} x_{2}^{-2}x_{3}^{-1}x_{4}^{-2} \\ -x_{2}^{-3}x_{4}^{-2}\end{array}\right)\begin{array}{l}\cdot \frac{1}{48}\left(x_{1} \partial_{1} -x_{2} \partial_{2} \right)^4\left(x_{3} \partial_{3} -x_{4} \partial_{4} \right)^2\\
		~~\cdot\frac{1}{(1-x_{1} ) (1-x_{1} x_{2} x_{3} x_{4} ) (1-x_{3} ) (1-x_{1} x_{2} ) }\end{array}\\&
	\displaystyle +\left(x_{1} x_{2}^{-2}x_{4}^{-2}-x_{2}^{-3}x_{4}^{-2}\right)\cdot \frac{1}{120}\left(x_{1} \partial_{1} -x_{2} \partial_{2} \right)^5\left(x_{3} \partial_{3} -x_{4} \partial_{4} \right)\cdot\frac{1}{(1-x_{1} ) (1-x_{1} x_{2} x_{3} x_{4} ) (1-x_{3} ) (1-x_{1} x_{2} ) }\\&
	\displaystyle -x_{2}^{-3}x_{4}^{-2}\cdot \frac{1}{720}\left(x_{1} \partial_{1} -x_{2} \partial_{2} \right)^6\cdot\frac{1}{(1-x_{1} ) (1-x_{1} x_{2} x_{3} x_{4} ) (1-x_{3} ) (1-x_{1} x_{2} ) }\\&
	\displaystyle +\left(\begin{array}{l}-x_{1}^2x_{2}^{-1}x_{3}^{-2}x_{4}^{-1}-x_{1} x_{2}^{-3}x_{3}^{-3}x_{4}^{-1} \\ -x_{2}^{-3}x_{3}^{-1}x_{4}^{-2}-x_{1} x_{2}^{-4}x_{3}^{-4}x_{4}^{-3}\end{array}\right)\begin{array}{l}\cdot \frac{1}{720}\left(x_{1} \partial_{1} -x_{2} \partial_{2} \right)^6\\
		~~\cdot\frac{1}{(1-x_{1} ) (1-x_{1} x_{2} x_{3} x_{4} ) (1-x_{1} x_{2} ) (1-x_{4} ) }\end{array}\\&
	\displaystyle +\left(\begin{array}{l}x_{1}^2x_{2}^{-1}x_{3}^{-2}x_{4}^{-1}+x_{1} x_{2}^{-3}x_{3}^{-3}x_{4}^{-1} \\ +x_{2}^{-3}x_{3}^{-1}x_{4}^{-2}\end{array}\right)\begin{array}{l}\cdot \frac{1}{720}\left(x_{1} \partial_{1} -x_{2} \partial_{2} +x_{3} \partial_{3} \right)^6\\
		~~\cdot\frac{1}{(1-x_{1} ) (1-x_{1} x_{2} ) (1-x_{2} x_{3} x_{4} ) (1-x_{4} ) }\end{array}\\&
	\displaystyle +x_{2}^{-3}x_{4}^{-2}\cdot \frac{1}{720}\left(x_{1} \partial_{1} -x_{2} \partial_{2} +x_{4} \partial_{4} \right)^6\cdot\frac{1}{(1-x_{1} ) (1-x_{1} x_{2} ) (1-x_{3} ) (1-x_{2} x_{3} x_{4} ) }\\&
	\displaystyle +\left(x_{2}^{-3}x_{4}^{-1}-2x_{2}^{-5}x_{3}^{-3}x_{4}^{-1}\right)\cdot \frac{1}{720}\left(x_{1} \partial_{1} \right)^6\cdot\frac{1}{(1-x_{1} ) (1-x_{2} x_{3} x_{4} ) (1-x_{3} ) (1-x_{4} ) }\\&
	\displaystyle -x_{2}^{-3}x_{4}^{-2}\cdot \frac{1}{120}\left(x_{1} \partial_{1} -x_{2} \partial_{2} \right)^5\left(x_{3} \partial_{3} \right)\cdot\frac{1}{(1-x_{1} ) (1-x_{1} x_{2} ) (1-x_{3} ) (1-x_{4} ) }\\&
	\displaystyle +x_{2}^{-3}x_{4}^{-2}\cdot \frac{1}{120}\left(x_{1} \partial_{1} -x_{2} \partial_{2} +x_{4} \partial_{4} \right)^5\left(x_{3} \partial_{3} -x_{4} \partial_{4} \right)\cdot\frac{1}{(1-x_{1} ) (1-x_{1} x_{2} ) (1-x_{3} ) (1-x_{2} x_{3} x_{4} ) }\\&
	\displaystyle +\left(\begin{array}{l}x_{2}^{-3}x_{4}^{-1}+x_{2}^{-3}x_{3}^{-1}x_{4}^{-1} \\ +x_{2}^{-4}x_{3}^{-2}x_{4}^{-1}\end{array}\right)\begin{array}{l}\cdot \frac{1}{120}\left(x_{1} \partial_{1} \right)^5\left(-x_{2} \partial_{2} +x_{3} \partial_{3} \right)\\
		~~\cdot\frac{1}{(1-x_{1} ) (1-x_{2} x_{3} x_{4} ) (1-x_{3} ) (1-x_{4} ) }\end{array}\\&
	\displaystyle -x_{2}^{-3}x_{4}^{-2}\cdot \frac{1}{48}\left(x_{1} \partial_{1} -x_{2} \partial_{2} \right)^4\left(x_{3} \partial_{3} \right)^2\cdot\frac{1}{(1-x_{1} ) (1-x_{1} x_{2} ) (1-x_{3} ) (1-x_{4} ) }\\&
	\displaystyle +\left(-x_{1} x_{2}^{-2}x_{3}^{-1}x_{4}^{-2}+x_{2}^{-3}x_{4}^{-2}\right)\cdot \frac{1}{48}\left(x_{1} \partial_{1} -x_{2} \partial_{2} +x_{4} \partial_{4} \right)^4\left(x_{3} \partial_{3} -x_{4} \partial_{4} \right)^2\cdot\frac{1}{(1-x_{1} ) (1-x_{1} x_{2} ) (1-x_{3} ) (1-x_{2} x_{3} x_{4} ) }\\&
	\displaystyle +\left(-x_{2}^{-3}x_{4}^{-1}+x_{2}^{-3}x_{4}^{-2}+x_{2}^{-3}x_{3}^{-1}x_{4}^{-1}\right)\cdot \frac{1}{48}\left(x_{1} \partial_{1} -x_{2} \partial_{2} \right)^4\left(-x_{2} \partial_{2} +x_{3} \partial_{3} \right)^2\cdot\frac{1}{(1-x_{1} ) (1-x_{1} x_{2} x_{3} x_{4} ) (1-x_{3} ) (1-x_{4} ) }\\&
	\displaystyle +\left(x_{2}^{-3}x_{4}^{-1}-x_{2}^{-3}x_{3}^{-1}x_{4}^{-1}\right)\cdot \frac{1}{48}\left(x_{1} \partial_{1} \right)^4\left(-x_{2} \partial_{2} +x_{3} \partial_{3} \right)^2\cdot\frac{1}{(1-x_{1} ) (1-x_{2} x_{3} x_{4} ) (1-x_{3} ) (1-x_{4} ) }\\&
	\displaystyle +x_{1} x_{2}^{-2}x_{4}^{-2}\cdot \frac{1}{36}\left(x_{1} \partial_{1} -x_{2} \partial_{2} \right)^3\left(x_{3} \partial_{3} -x_{4} \partial_{4} \right)^3\cdot\frac{1}{(1-x_{1} ) (1-x_{1} x_{2} x_{3} x_{4} ) (1-x_{3} ) (1-x_{1} x_{2} ) }\\&
	\displaystyle +\left(\begin{array}{l}x_{2}^{-4}x_{3}^{-3}x_{4} -x_{2}^{-4}x_{3}^{-4}x_{4}^{-1} \\ +2x_{2}^{-5}x_{3}^{-4}x_{4}^{-1}-2x_{2}^{-5}x_{3}^{-5}x_{4}^{-2}\end{array}\right)\begin{array}{l}\cdot \frac{1}{120}\left(x_{1} \partial_{1} -x_{3} \partial_{3} \right)^5\left(x_{2} \partial_{2} -x_{3} \partial_{3} \right)\\
		~~\cdot\frac{1}{(1-x_{1} ) (1-x_{2} ) (1-x_{1} x_{2} x_{3} x_{4} ) (1-x_{4} ) }\end{array}\\&
	\displaystyle -x_{2}^{-5}x_{3}^{-5}x_{4}^{-3}\cdot \frac{1}{120}\left(x_{1} \partial_{1} -x_{4} \partial_{4} \right)^5\left(x_{2} \partial_{2} -x_{4} \partial_{4} \right)\cdot\frac{1}{(1-x_{1} ) (1-x_{2} ) (1-x_{3} ) (1-x_{1} x_{2} x_{3} x_{4} ) }\\&
	\displaystyle +4x_{2}^{-5}x_{3}^{-5}x_{4}^{-3}\cdot \frac{1}{120}\left(x_{1} \partial_{1} \right)^5\left(x_{2} \partial_{2} \right)\cdot\frac{1}{(1-x_{1} ) (1-x_{2} ) (1-x_{3} ) (1-x_{4} ) }\\&
	\displaystyle +x_{2}^{-3}x_{4}^{-2}\cdot \frac{1}{36}\left(x_{1} \partial_{1} -x_{2} \partial_{2} \right)^3\left(-x_{2} \partial_{2} +x_{3} \partial_{3} \right)^3\cdot\frac{1}{(1-x_{1} ) (1-x_{1} x_{2} x_{3} x_{4} ) (1-x_{3} ) (1-x_{4} ) }\\&
	\displaystyle -x_{2}^{-3}x_{3}^{-5}x_{4}^{-3}\cdot \frac{1}{36}\left(x_{1} \partial_{1} \right)^3\left(x_{2} \partial_{2} \right)^3\cdot\frac{1}{(1-x_{1} ) (1-x_{2} ) (1-x_{3} ) (1-x_{4} ) }\\&
	\displaystyle +\left(\begin{array}{l}-x_{2}^{-3}x_{4}^{-1}+2x_{2}^{-5}x_{3}^{-3}x_{4}^{-1} \\ +6x_{2}^{-6}x_{3}^{-4}x_{4}^{-2}\end{array}\right)\begin{array}{l}\cdot \frac{1}{720}\left(x_{1} \partial_{1} -x_{2} \partial_{2} \right)^6\\
		~~\cdot\frac{1}{(1-x_{1} ) (1-x_{1} x_{2} x_{3} x_{4} ) (1-x_{3} ) (1-x_{4} ) }\end{array}\\&
	\displaystyle +\left(\begin{array}{l}-x_{2}^{-5}x_{3}^{-3}x_{4} +2x_{2}^{-5}x_{3}^{-4}x_{4}^{-1} \\ -3x_{2}^{-6}x_{3}^{-4}x_{4}^{-1}+6x_{2}^{-6}x_{3}^{-5}x_{4}^{-2}\end{array}\right)\begin{array}{l}\cdot \frac{1}{720}\left(x_{1} \partial_{1} -x_{3} \partial_{3} \right)^6\\
		~~\cdot\frac{1}{(1-x_{1} ) (1-x_{2} ) (1-x_{1} x_{2} x_{3} x_{4} ) (1-x_{4} ) }\end{array}\\&
	\displaystyle +3x_{2}^{-6}x_{3}^{-5}x_{4}^{-3}\cdot \frac{1}{720}\left(x_{1} \partial_{1} -x_{4} \partial_{4} \right)^6\cdot\frac{1}{(1-x_{1} ) (1-x_{2} ) (1-x_{3} ) (1-x_{1} x_{2} x_{3} x_{4} ) }\\&
	\displaystyle -6x_{2}^{-6}x_{3}^{-5}x_{4}^{-3}\cdot \frac{1}{720}\left(x_{1} \partial_{1} \right)^6\cdot\frac{1}{(1-x_{1} ) (1-x_{2} ) (1-x_{3} ) (1-x_{4} ) }\\&
	\displaystyle +\left(\begin{array}{l}-x_{2}^{-3}x_{4}^{-1}+x_{2}^{-3}x_{4}^{-2}-x_{2}^{-3}x_{3}^{-1}x_{4}^{-1} \\ -x_{2}^{-4}x_{3}^{-2}x_{4}^{-1}-x_{2}^{-4}x_{3}^{-2}x_{4}^{-3} \\ -2x_{2}^{-5}x_{3}^{-3}x_{4}^{-2}\end{array}\right)\begin{array}{l}\cdot \frac{1}{120}\left(x_{1} \partial_{1} -x_{2} \partial_{2} \right)^5\left(-x_{2} \partial_{2} +x_{3} \partial_{3} \right)\\
		~~\cdot\frac{1}{(1-x_{1} ) (1-x_{1} x_{2} x_{3} x_{4} ) (1-x_{3} ) (1-x_{4} ) }\end{array}\\&
	\displaystyle -x_{2}^{-4}x_{3}^{-5}x_{4}^{-3}\cdot \frac{1}{48}\left(x_{1} \partial_{1} \right)^4\left(x_{2} \partial_{2} -x_{4} \partial_{4} \right)^2\cdot\frac{1}{(1-x_{1} ) (1-x_{2} ) (1-x_{3} ) (1-x_{2} x_{3} x_{4} ) }\\&
	\displaystyle +2x_{2}^{-5}x_{3}^{-4}x_{4}^{-3}\cdot \frac{1}{120}\left(x_{1} \partial_{1} \right)^5\left(x_{3} \partial_{3} \right)\cdot\frac{1}{(1-x_{1} ) (1-x_{2} ) (1-x_{3} ) (1-x_{4} ) }\\&
	\displaystyle +x_{2}^{-3}x_{3}^{-5}x_{4}^{-3}\cdot \frac{1}{36}\left(x_{1} \partial_{1} \right)^3\left(x_{2} \partial_{2} -x_{4} \partial_{4} \right)^3\cdot\frac{1}{(1-x_{1} ) (1-x_{2} ) (1-x_{3} ) (1-x_{2} x_{3} x_{4} ) }\\&
	\displaystyle +2x_{2}^{-4}x_{3}^{-4}x_{4}^{-3}\cdot \frac{1}{24}\left(x_{1} \partial_{1} \right)^4\left(x_{2} \partial_{2} \right)\left(x_{3} \partial_{3} -x_{4} \partial_{4} \right)\cdot\frac{1}{(1-x_{1} ) (1-x_{2} ) (1-x_{3} ) (1-x_{3} x_{4} ) }\\&
	\displaystyle -2x_{2}^{-4}x_{3}^{-4}x_{4}^{-3}\cdot \frac{1}{24}\left(x_{1} \partial_{1} \right)^4\left(x_{2} \partial_{2} \right)\left(x_{3} \partial_{3} \right)\cdot\frac{1}{(1-x_{1} ) (1-x_{2} ) (1-x_{3} ) (1-x_{4} ) }\\&
	\displaystyle +\left(x_{2}^{-3}x_{3}^{-4}x_{4}^{-1}+x_{2}^{-3}x_{3}^{-5}x_{4}^{-2}\right)\cdot \frac{1}{36}\left(x_{1} \partial_{1} \right)^3\left(x_{2} \partial_{2} -x_{3} \partial_{3} \right)^3\cdot\frac{1}{(1-x_{1} ) (1-x_{2} ) (1-x_{2} x_{3} x_{4} ) (1-x_{4} ) }\\&
	\displaystyle +\left(\begin{array}{l}-x_{1} x_{2}^{-2}-x_{1} x_{2}^{-2}x_{3}^{-1}-x_{1} x_{2}^{-3}x_{3}^{-1} \\ +x_{2}^{-3}x_{3} x_{4}^{-1}-2x_{1} x_{2}^{-4}x_{3}^{-1} \\ +x_{2}^{-3}x_{3}^{-1}x_{4}^{-1}+x_{2}^{-3}x_{3}^{-2}x_{4}^{-1} \\ +x_{2}^{-4}x_{3}^{-2}x_{4}^{-1}+2x_{2}^{-5}x_{3}^{-3}x_{4}^{-2}\end{array}\right)\begin{array}{l}\cdot \frac{1}{120}\left(x_{1} \partial_{1} -x_{2} \partial_{2} \right)^5\left(x_{3} \partial_{3} -x_{4} \partial_{4} \right)\\
		~~\cdot\frac{1}{(1-x_{1} ) (1-x_{1} x_{2} x_{3} x_{4} ) (1-x_{3} ) (1-x_{3} x_{4} ) }\end{array}\\&
	\displaystyle +x_{2}^{-3}x_{3}^{-4}x_{4}^{-3}\cdot \frac{1}{12}\left(x_{1} \partial_{1} \right)^3\left(x_{2} \partial_{2} \right)^2\left(x_{3} \partial_{3} \right)\cdot\frac{1}{(1-x_{1} ) (1-x_{2} ) (1-x_{3} ) (1-x_{4} ) }\end{align*}

	\end{landscape}

\end{document}